\documentclass{amsart}

\newtheorem{theorem}{Theorem}
\newtheorem{definition}{Definition}

\usepackage{amsfonts}
\usepackage[section]{placeins}
\usepackage{comment}
\usepackage{amsmath}
\usepackage{tabularx}
\usepackage{array}
\usepackage{css-colors}
\newtheorem{lemma}{Lemma}
\newtheorem{corollary}{Corollary}

\def\Q{{\mathbb Q}}

\usepackage{pstricks}

\psset{unit=3cm}
\def\ThreeTiling{%
\pspicture(1, 1)
\qline(0,0)(1,0)            
\qline(0,0)(0.5,0.866025404)    
\qline(0.5,0.866025404)(1,0)    
\qline(0,0)(0.5,0.288675135)     
\qline(0.50,0.866025404)(0.50,0.288675135)  
\qline(0.50,0.288675135)(1,0)     
\endpspicture}

\def\SixTiling{%
\pspicture(1, 1)
\qline(0,0)(1,0)            
\qline(0,0)(0.5,0.866025404)    
\qline(0.5,0.866025404)(1,0)    
\qline(0,0)(0.75, 0.433012702)     
\qline(0.50,0.866025404)(0.50,0)  
\qline(0.25, 0.433012702)(1,0)     
\endpspicture}

\def\EightTiling{%
\pspicture(1,1)
\qline(0,0)(1,0)            
\qline(0,0)(0.5,0.866025404)    
\qline(0.5,0.866025404)(1,0)    
\qline(0.25,0.433012602)(0.75,0.433012602)  
\qline(0.25,0.433012602)(0.25,0)  
\qline(0.75,0.433012602)(0.75,0)  
\qline(0.5,0.866025404)(0.5,0)   
\qline(0.25,0.433012602)(0.5,0)  
\qline(0.5,0)(0.75,0.433012602)   
\endpspicture}

\def\FiveTiling{%
\psset{unit=2cm}
\pspicture(2.5,1.2)
\qline(0,0)(2.5,0)              
\qline(0,0)(2,1)    
\qline(2,1)(2.5,0)   
\qline(2,1)(2,0)     
\qline(1,0.5)(2,0.5)  
\qline(2,0.5)(1,0)     
\qline(1,0)(1,0.5)    
\psset{unit=3cm}
\endpspicture}

\def\FiveTilingB{%
\psset{unit=2cm}
\pspicture(2.5,1.2)
\qline(0,0)(2.5,0)              
\qline(0,0)(2,1)    
\qline(2,1)(2.5,0)   
\qline(2,1)(2,0)     
\qline(1,0.5)(2,0.5)  
\qline(1,0.5)(2,0)     
\qline(1,0)(1,0.5)    
\psset{unit=3cm}
\endpspicture}

\def\FourTilingA{%
\psset{unit=2cm}
\pspicture(2,1.2)
\qline(0,0)(1.73205081,0)              
\qline(1.73205081,0)(1.73205081,1)    
\qline(0.866025404,0.5)( 0.866025404,0)   
\qline( 0.866025404,0.5)(1.73205081,0)   
\qline( 0,0)(1.73205081,1) 
\qline( 0.866025404,0.5)(1.73205081,0.5)  
\endpspicture
\psset{unit=3cm}
}

\def\FourTilingB{%
\psset{unit=2cm}
\pspicture(2,1.2)
\qline(0,0)(1.73205081,0)              
\qline(1.73205081,0)(1.73205081,1)    
\qline(0.866025404,0.5)( 0.866025404,0)   
\qline( 0.866025404,0.5)(1.73205081,0)   
\qline( 0,0)(1.73205081,1) 
\qline(1.29903811,0.25)(1.73205081,1)
\endpspicture
\psset{unit=3cm}
}

\def\FourTilingC{%
\psset{unit=2cm}
\pspicture(2,1.2)
\qline(0,0)(1.73205081,0)              
\qline(1.73205081,0)(1.73205081,1)    
\qline(0.866025404,0.5)( 0.866025404,0)   
\qline( 0.866025404,0.5)(1.73205081,0)   
\qline( 0,0)(1.73205081,1) 
\qline(1.73205081,0)(1.29903811,0.75)
\endpspicture
\psset{unit=3cm}
}

\def\FourTiling{%
\pspicture(1,1)
\qline(0,0)(1,0)            
\qline(0,0)(0.5,0.866025404)    
\qline(0.5,0.866025404)(1,0)    
\qline(0.25,0.433012602)(0.75,0.433012602)  
\qline(0.25,0.433012602)(0.5,0)  
\qline(0.5,0)(0.75,0.433012602)   
\endpspicture}

\def\SixteenTiling{%
\pspicture(1,1)
\qline(0,0)(1,0)            
\qline(0,0)(0.5,0.866025404)    
\qline(0.5,0.866025404)(1,0)    
\qline(0.125,0.216506351)(0.875,0.216506351)  
\qline(0.25,0.433012702)(0.75,  0.433012702) 
\qline(0.375,0.649519052)(0.625,0.649519052 ) 
\qline(0.125,0.216506351)(0.25,0)
\qline(0.25,0.433012702)(0.5,0)
\qline(0.375,0.649519052)(0.75,0)
\qline(0.25,0)(0.625,0.649519052 )
\qline(0.5,0)(0.75,  0.433012702)
\qline(0.75,0)(0.875,0.216506351)
\endpspicture}

\def\ThreeTilingA{%
\pspicture(2, 1)
\qline(0,0)(2,0)            
\qline(0,0)(0.5,0.866025404)    
\qline(0.5,0.866025404)(2,0)    
\qline(0,0)(1,0.577350269)     
\qline(1,0.577350269)(1,0)  
\endpspicture}

\def\NineTiling{%
\pspicture(1,1)
\qline(0,0)(1,0)            
\qline(0,0)(0.5,0.866025404)    
\qline(0.5,0.866025404)(1,0)    
\qline(0.167777,0.288675135)(0.833333,0.288675135)  
\qline(0.333333,0.577350269)(0.6777777, 0.577350269) 
\qline(0.1677777,0.288675135)(0.333333,0)
\qline(0.333333,0.577350269)(0.6777777,0)
\qline(0.333333,0)(0.677777,0.577350269)
\qline(0.677777,0)(0.833333,0.288675135)
\endpspicture}

\def\FigureBigQuadratic{%
\hskip0.5in
\pspicture(2.0,0.9)
\psset{unit=0.5cm}
\newrgbcolor{lightblue}{0.8 0.8 1}
\newrgbcolor{pink}{1 0.8 0.8}
\newrgbcolor{lightgreen}{0.8 1 0.8}
\newrgbcolor{lightyellow}{1 1 0.8}
\pspolygon[fillstyle=solid,linewidth=1pt,fillcolor=lightblue](1.00,0.00)(3.00,5.00)(12.00,0.00)
\psline(3.00,5.00)(12.00,0.00)
\psline(2.80,4.50)(10.90,0.00)
\psline(2.60,4.00)(9.80,0.00)
\psline(2.40,3.50)(8.70,0.00)
\psline(2.20,3.00)(7.60,0.00)
\psline(2.00,2.50)(6.50,0.00)
\psline(1.80,2.00)(5.40,0.00)
\psline(1.60,1.50)(4.30,0.00)
\psline(1.40,1.00)(3.20,0.00)
\psline(1.20,0.50)(2.10,0.00)
\psline(1.00,0.00)(1.00,0.00)
\psline(12.00,0.00)(1.00,0.00)
\psline(11.10,0.50)(1.20,0.50)
\psline(10.20,1.00)(1.40,1.00)
\psline(9.30,1.50)(1.60,1.50)
\psline(8.40,2.00)(1.80,2.00)
\psline(7.50,2.50)(2.00,2.50)
\psline(6.60,3.00)(2.20,3.00)
\psline(5.70,3.50)(2.40,3.50)
\psline(4.80,4.00)(2.60,4.00)
\psline(3.90,4.50)(2.80,4.50)
\psline(3.00,5.00)(3.00,5.00)
\psline(1.00,0.00)(3.00,5.00)
\psline(2.10,0.00)(3.90,4.50)
\psline(3.20,0.00)(4.80,4.00)
\psline(4.30,0.00)(5.70,3.50)
\psline(5.40,0.00)(6.60,3.00)
\psline(6.50,0.00)(7.50,2.50)
\psline(7.60,0.00)(8.40,2.00)
\psline(8.70,0.00)(9.30,1.50)
\psline(9.80,0.00)(10.20,1.00)
\psline(10.90,0.00)(11.10,0.50)
\psline(12.00,0.00)(12.00,0.00)
\endpspicture
}

\def\ThirteenTiling{%
\pspicture(1.9,0.8)
\psset{unit=0.4cm}
\newrgbcolor{lightblue}{0.8 0.8 1}
\newrgbcolor{pink}{1 0.8 0.8}
\newrgbcolor{lightgreen}{0.8 1 0.8}
\pspolygon[fillstyle=solid,linewidth=1pt,fillcolor=lightgreen](0.00,0.00)(9.00,6.00)(9.00,0.00)
\psline(9.00,6.00)(9.00,0.00)
\psline(6.00,4.00)(6.00,0.00)
\psline(3.00,2.00)(3.00,0.00)
\psline(0.00,0.00)(0.00,0.00)
\psline(9.00,0.00)(0.00,0.00)
\psline(9.00,2.00)(3.00,2.00)
\psline(9.00,4.00)(6.00,4.00)
\psline(9.00,6.00)(9.00,6.00)
\psline(0.00,0.00)(9.00,6.00)
\psline(3.00,0.00)(9.00,4.00)
\psline(6.00,0.00)(9.00,2.00)
\psline(9.00,0.00)(9.00,0.00)
\pspolygon[fillstyle=solid,linewidth=1pt,fillcolor=lightblue](9.00,6.00)(9.00,0.00)(13.00,0.00)
\psline(9.00,0.00)(13.00,0.00)
\psline(9.00,3.00)(11.00,3.00)
\psline(9.00,6.00)(9.00,6.00)
\psline(13.00,0.00)(9.00,6.00)
\psline(11.00,0.00)(9.00,3.00)
\psline(9.00,0.00)(9.00,0.00)
\psline(9.00,6.00)(9.00,0.00)
\psline(11.00,3.00)(11.00,0.00)
\psline(13.00,0.00)(13.00,0.00)
\endpspicture}

\def\FigureBiquadratic{
\pspicture(2,1)
\psset{unit=0.08cm}
\newrgbcolor{lightblue}{0.8 0.8 1}
\newrgbcolor{pink}{1 0.8 0.8}
\newrgbcolor{lightgreen}{0.8 1 0.8}
\pspolygon[fillstyle=solid,linewidth=1pt,fillcolor=lightgreen](0.00,0.00)(49.00,35.00)(49.00,0.00)
\psline(49.00,35.00)(49.00,0.00)
\psline(42.00,30.00)(42.00,0.00)
\psline(35.00,25.00)(35.00,0.00)
\psline(28.00,20.00)(28.00,0.00)
\psline(21.00,15.00)(21.00,0.00)
\psline(14.00,10.00)(14.00,0.00)
\psline(7.00,5.00)(7.00,0.00)
\psline(0.00,0.00)(0.00,0.00)
\psline(49.00,0.00)(0.00,0.00)
\psline(49.00,5.00)(7.00,5.00)
\psline(49.00,10.00)(14.00,10.00)
\psline(49.00,15.00)(21.00,15.00)
\psline(49.00,20.00)(28.00,20.00)
\psline(49.00,25.00)(35.00,25.00)
\psline(49.00,30.00)(42.00,30.00)
\psline(49.00,35.00)(49.00,35.00)
\psline(0.00,0.00)(49.00,35.00)
\psline(7.00,0.00)(49.00,30.00)
\psline(14.00,0.00)(49.00,25.00)
\psline(21.00,0.00)(49.00,20.00)
\psline(28.00,0.00)(49.00,15.00)
\psline(35.00,0.00)(49.00,10.00)
\psline(42.00,0.00)(49.00,5.00)
\psline(49.00,0.00)(49.00,0.00)
\pspolygon[fillstyle=solid,linewidth=1pt,fillcolor=lightblue](49.00,35.00)(49.00,0.00)(74.00,0.00)
\psline(49.00,0.00)(74.00,0.00)
\psline(49.00,7.00)(69.00,7.00)
\psline(49.00,14.00)(64.00,14.00)
\psline(49.00,21.00)(59.00,21.00)
\psline(49.00,28.00)(54.00,28.00)
\psline(49.00,35.00)(49.00,35.00)
\psline(74.00,0.00)(49.00,35.00)
\psline(69.00,0.00)(49.00,28.00)
\psline(64.00,0.00)(49.00,21.00)
\psline(59.00,0.00)(49.00,14.00)
\psline(54.00,0.00)(49.00,7.00)
\psline(49.00,0.00)(49.00,0.00)
\psline(49.00,35.00)(49.00,0.00)
\psline(54.00,28.00)(54.00,0.00)
\psline(59.00,21.00)(59.00,0.00)
\psline(64.00,14.00)(64.00,0.00)
\psline(69.00,7.00)(69.00,0.00)
\psline(74.00,0.00)(74.00,0.00)
\endpspicture}

\def\TwelveTiling{%
\pspicture(1,1)
\qline(0,0)(1,0)            
\qline(0,0)(0.5,0.866025404)    
\qline(0.5,0.866025404)(1,0)    
\qline(0.25,0.433012602)(0.75,0.433012602)  
\qline(0.25,0.433012602)(0.5,0)  
\qline(0.5,0)(0.75,0.433012602)   
\qline(0,0)(0.25,0.144337567)
\qline(0.25,0.144337567)(0.25,0.433012602)
\qline(0.25,0.144337567)(0.5,0)
\qline(0.5,0)(0.75,0.144337567)
\qline(0.75,0.144337567)(0.75,0.433012602)
\qline(0.75,0.144337567)(1,0)
\qline(0.5,0.288675035)(0.5,0)
\qline(0.5,0.288675035)(0.25,0.433012602)
\qline(0.5,0.288675035)(0.75,0.433012602)
\qline(0.5,0.577350169)(0.5,0.866025404)
\qline(0.5,0.577350169)(0.25,0.433012602)
\qline(0.5,0.577350169)(0.75,0.433012602)
\endpspicture}

\def\HexTilingThree{%
\psset{unit=1cm}
\pspicture(6.1,5.1)(0.8,0.2)
\pspolygon[fillstyle=solid,linewidth=1pt,fillcolor=lightyellow](2.60,1.50)(3.46,1.00)(1.73,1.00)
\pspolygon[fillstyle=solid,linewidth=1pt,fillcolor=lightyellow](2.60,0.50)(3.46,1.00)(1.73,1.00)
\pspolygon[fillstyle=solid,linewidth=1pt,fillcolor=lightyellow](2.60,1.50)(2.60,2.50)(3.46,1.00)
\pspolygon[fillstyle=solid,linewidth=1pt,fillcolor=lightyellow](2.60,1.50)(1.73,1.00)(2.60,2.50)
\pspolygon[fillstyle=solid,linewidth=1pt,fillcolor=lightyellow](3.46,2.00)(2.60,2.50)(3.46,1.00)
\pspolygon[fillstyle=solid,linewidth=1pt,fillcolor=lightyellow](1.73,2.00)(1.73,1.00)(2.60,2.50)
\pspolygon[fillstyle=solid,linewidth=1pt,fillcolor=lightgreen](4.33,1.50)(5.20,1.00)(3.46,1.00)
\pspolygon[fillstyle=solid,linewidth=1pt,fillcolor=lightgreen](4.33,0.50)(5.20,1.00)(3.46,1.00)
\pspolygon[fillstyle=solid,linewidth=1pt,fillcolor=lightgreen](4.33,1.50)(4.33,2.50)(5.20,1.00)
\pspolygon[fillstyle=solid,linewidth=1pt,fillcolor=lightgreen](4.33,1.50)(3.46,1.00)(4.33,2.50)
\pspolygon[fillstyle=solid,linewidth=1pt,fillcolor=lightgreen](5.20,2.00)(4.33,2.50)(5.20,1.00)
\pspolygon[fillstyle=solid,linewidth=1pt,fillcolor=lightgreen](3.46,2.00)(3.46,1.00)(4.33,2.50)
\pspolygon[fillstyle=solid,linewidth=1pt,fillcolor=pink](3.46,3.00)(4.33,2.50)(2.60,2.50)
\pspolygon[fillstyle=solid,linewidth=1pt,fillcolor=pink](3.46,2.00)(4.33,2.50)(2.60,2.50)
\pspolygon[fillstyle=solid,linewidth=1pt,fillcolor=pink](3.46,3.00)(3.46,4.00)(4.33,2.50)
\pspolygon[fillstyle=solid,linewidth=1pt,fillcolor=pink](3.46,3.00)(2.60,2.50)(3.46,4.00)
\pspolygon[fillstyle=solid,linewidth=1pt,fillcolor=pink](4.33,3.50)(3.46,4.00)(4.33,2.50)
\pspolygon[fillstyle=solid,linewidth=1pt,fillcolor=pink](2.60,3.50)(2.60,2.50)(3.46,4.00)
\pspolygon[fillstyle=solid,linewidth=1pt,fillcolor=lightblue](0.87,0.50)(2.60,0.50)(1.73,1.00)
\pspolygon[fillstyle=solid,linewidth=1pt,fillcolor=lightblue](2.60,0.50)(4.33,0.50)(3.46,1.00)
\pspolygon[fillstyle=solid,linewidth=1pt,fillcolor=lightblue](4.33,0.50)(6.06,0.50)(5.20,1.00)
\pspolygon[fillstyle=solid,linewidth=1pt,fillcolor=lightblue](6.06,0.50)(5.20,2.00)(5.20,1.00)
\pspolygon[fillstyle=solid,linewidth=1pt,fillcolor=lightblue](5.20,2.00)(4.33,3.50)(4.33,2.50)
\pspolygon[fillstyle=solid,linewidth=1pt,fillcolor=lightblue](4.33,3.50)(3.46,5.00)(3.46,4.00)
\pspolygon[fillstyle=solid,linewidth=1pt,fillcolor=lightblue](1.73,2.00)(0.87,0.50)(1.73,1.00)
\pspolygon[fillstyle=solid,linewidth=1pt,fillcolor=lightblue](2.60,3.50)(1.73,2.00)(2.60,2.50)
\pspolygon[fillstyle=solid,linewidth=1pt,fillcolor=lightblue](3.46,5.00)(2.60,3.50)(3.46,4.00)
\endpspicture}

\def\HexTilingFour{
\psset{unit=0.8cm}
\pspicture(8,7)(0.8,0)
\pspolygon[fillstyle=solid,linewidth=1pt,fillcolor=lightyellow](2.60,1.50)(3.46,1.00)(1.73,1.00)
\pspolygon[fillstyle=solid,linewidth=1pt,fillcolor=lightyellow](2.60,0.50)(3.46,1.00)(1.73,1.00)
\pspolygon[fillstyle=solid,linewidth=1pt,fillcolor=lightyellow](2.60,1.50)(2.60,2.50)(3.46,1.00)
\pspolygon[fillstyle=solid,linewidth=1pt,fillcolor=lightyellow](2.60,1.50)(1.73,1.00)(2.60,2.50)
\pspolygon[fillstyle=solid,linewidth=1pt,fillcolor=lightyellow](3.46,2.00)(2.60,2.50)(3.46,1.00)
\pspolygon[fillstyle=solid,linewidth=1pt,fillcolor=lightyellow](1.73,2.00)(1.73,1.00)(2.60,2.50)
\pspolygon[fillstyle=solid,linewidth=1pt,fillcolor=lightgreen](4.33,1.50)(5.20,1.00)(3.46,1.00)
\pspolygon[fillstyle=solid,linewidth=1pt,fillcolor=lightgreen](4.33,0.50)(5.20,1.00)(3.46,1.00)
\pspolygon[fillstyle=solid,linewidth=1pt,fillcolor=lightgreen](4.33,1.50)(4.33,2.50)(5.20,1.00)
\pspolygon[fillstyle=solid,linewidth=1pt,fillcolor=lightgreen](4.33,1.50)(3.46,1.00)(4.33,2.50)
\pspolygon[fillstyle=solid,linewidth=1pt,fillcolor=lightgreen](5.20,2.00)(4.33,2.50)(5.20,1.00)
\pspolygon[fillstyle=solid,linewidth=1pt,fillcolor=lightgreen](3.46,2.00)(3.46,1.00)(4.33,2.50)
\pspolygon[fillstyle=solid,linewidth=1pt,fillcolor=pink](6.06,1.50)(6.93,1.00)(5.20,1.00)
\pspolygon[fillstyle=solid,linewidth=1pt,fillcolor=pink](6.06,0.50)(6.93,1.00)(5.20,1.00)
\pspolygon[fillstyle=solid,linewidth=1pt,fillcolor=pink](6.06,1.50)(6.06,2.50)(6.93,1.00)
\pspolygon[fillstyle=solid,linewidth=1pt,fillcolor=pink](6.06,1.50)(5.20,1.00)(6.06,2.50)
\pspolygon[fillstyle=solid,linewidth=1pt,fillcolor=pink](6.93,2.00)(6.06,2.50)(6.93,1.00)
\pspolygon[fillstyle=solid,linewidth=1pt,fillcolor=pink](5.20,2.00)(5.20,1.00)(6.06,2.50)
\pspolygon[fillstyle=solid,linewidth=1pt,fillcolor=red](3.46,3.00)(4.33,2.50)(2.60,2.50)
\pspolygon[fillstyle=solid,linewidth=1pt,fillcolor=red](3.46,2.00)(4.33,2.50)(2.60,2.50)
\pspolygon[fillstyle=solid,linewidth=1pt,fillcolor=red](3.46,3.00)(3.46,4.00)(4.33,2.50)
\pspolygon[fillstyle=solid,linewidth=1pt,fillcolor=red](3.46,3.00)(2.60,2.50)(3.46,4.00)
\pspolygon[fillstyle=solid,linewidth=1pt,fillcolor=red](4.33,3.50)(3.46,4.00)(4.33,2.50)
\pspolygon[fillstyle=solid,linewidth=1pt,fillcolor=red](2.60,3.50)(2.60,2.50)(3.46,4.00)
\pspolygon[fillstyle=solid,linewidth=1pt,fillcolor=blue](5.20,3.00)(6.06,2.50)(4.33,2.50)
\pspolygon[fillstyle=solid,linewidth=1pt,fillcolor=blue](5.20,2.00)(6.06,2.50)(4.33,2.50)
\pspolygon[fillstyle=solid,linewidth=1pt,fillcolor=blue](5.20,3.00)(5.20,4.00)(6.06,2.50)
\pspolygon[fillstyle=solid,linewidth=1pt,fillcolor=blue](5.20,3.00)(4.33,2.50)(5.20,4.00)
\pspolygon[fillstyle=solid,linewidth=1pt,fillcolor=blue](6.06,3.50)(5.20,4.00)(6.06,2.50)
\pspolygon[fillstyle=solid,linewidth=1pt,fillcolor=blue](4.33,3.50)(4.33,2.50)(5.20,4.00)
\pspolygon[fillstyle=solid,linewidth=1pt,fillcolor=Plum](4.33,4.50)(5.20,4.00)(3.46,4.00)
\pspolygon[fillstyle=solid,linewidth=1pt,fillcolor=Plum](4.33,3.50)(5.20,4.00)(3.46,4.00)
\pspolygon[fillstyle=solid,linewidth=1pt,fillcolor=Plum](4.33,4.50)(4.33,5.50)(5.20,4.00)
\pspolygon[fillstyle=solid,linewidth=1pt,fillcolor=Plum](4.33,4.50)(3.46,4.00)(4.33,5.50)
\pspolygon[fillstyle=solid,linewidth=1pt,fillcolor=Plum](5.20,5.00)(4.33,5.50)(5.20,4.00)
\pspolygon[fillstyle=solid,linewidth=1pt,fillcolor=Plum](3.46,5.00)(3.46,4.00)(4.33,5.50)
\pspolygon[fillstyle=solid,linewidth=1pt,fillcolor=lightblue](0.87,0.50)(2.60,0.50)(1.73,1.00)
\pspolygon[fillstyle=solid,linewidth=1pt,fillcolor=lightblue](2.60,0.50)(4.33,0.50)(3.46,1.00)
\pspolygon[fillstyle=solid,linewidth=1pt,fillcolor=lightblue](4.33,0.50)(6.06,0.50)(5.20,1.00)
\pspolygon[fillstyle=solid,linewidth=1pt,fillcolor=lightblue](6.06,0.50)(7.79,0.50)(6.93,1.00)
\pspolygon[fillstyle=solid,linewidth=1pt,fillcolor=lightblue](7.79,0.50)(6.93,2.00)(6.93,1.00)
\pspolygon[fillstyle=solid,linewidth=1pt,fillcolor=lightblue](6.93,2.00)(6.06,3.50)(6.06,2.50)
\pspolygon[fillstyle=solid,linewidth=1pt,fillcolor=lightblue](6.06,3.50)(5.20,5.00)(5.20,4.00)
\pspolygon[fillstyle=solid,linewidth=1pt,fillcolor=lightblue](5.20,5.00)(4.33,6.50)(4.33,5.50)
\pspolygon[fillstyle=solid,linewidth=1pt,fillcolor=lightblue](1.73,2.00)(0.87,0.50)(1.73,1.00)
\pspolygon[fillstyle=solid,linewidth=1pt,fillcolor=lightblue](2.60,3.50)(1.73,2.00)(2.60,2.50)
\pspolygon[fillstyle=solid,linewidth=1pt,fillcolor=lightblue](3.46,5.00)(2.60,3.50)(3.46,4.00)
\pspolygon[fillstyle=solid,linewidth=1pt,fillcolor=lightblue](4.33,6.50)(3.46,5.00)(4.33,5.50)
\endpspicture
}

\def\HexTilingFive{
\psset{unit=0.65cm}
\pspicture(9.5,8)(0,-0.2)
\pspolygon[fillstyle=solid,linewidth=1pt,fillcolor=lightyellow](2.60,1.50)(3.46,1.00)(1.73,1.00)
\pspolygon[fillstyle=solid,linewidth=1pt,fillcolor=lightyellow](2.60,0.50)(3.46,1.00)(1.73,1.00)
\pspolygon[fillstyle=solid,linewidth=1pt,fillcolor=lightyellow](2.60,1.50)(2.60,2.50)(3.46,1.00)
\pspolygon[fillstyle=solid,linewidth=1pt,fillcolor=lightyellow](2.60,1.50)(1.73,1.00)(2.60,2.50)
\pspolygon[fillstyle=solid,linewidth=1pt,fillcolor=lightyellow](3.46,2.00)(2.60,2.50)(3.46,1.00)
\pspolygon[fillstyle=solid,linewidth=1pt,fillcolor=lightyellow](1.73,2.00)(1.73,1.00)(2.60,2.50)
\pspolygon[fillstyle=solid,linewidth=1pt,fillcolor=lightgreen](4.33,1.50)(5.20,1.00)(3.46,1.00)
\pspolygon[fillstyle=solid,linewidth=1pt,fillcolor=lightgreen](4.33,0.50)(5.20,1.00)(3.46,1.00)
\pspolygon[fillstyle=solid,linewidth=1pt,fillcolor=lightgreen](4.33,1.50)(4.33,2.50)(5.20,1.00)
\pspolygon[fillstyle=solid,linewidth=1pt,fillcolor=lightgreen](4.33,1.50)(3.46,1.00)(4.33,2.50)
\pspolygon[fillstyle=solid,linewidth=1pt,fillcolor=lightgreen](5.20,2.00)(4.33,2.50)(5.20,1.00)
\pspolygon[fillstyle=solid,linewidth=1pt,fillcolor=lightgreen](3.46,2.00)(3.46,1.00)(4.33,2.50)
\pspolygon[fillstyle=solid,linewidth=1pt,fillcolor=pink](6.06,1.50)(6.93,1.00)(5.20,1.00)
\pspolygon[fillstyle=solid,linewidth=1pt,fillcolor=pink](6.06,0.50)(6.93,1.00)(5.20,1.00)
\pspolygon[fillstyle=solid,linewidth=1pt,fillcolor=pink](6.06,1.50)(6.06,2.50)(6.93,1.00)
\pspolygon[fillstyle=solid,linewidth=1pt,fillcolor=pink](6.06,1.50)(5.20,1.00)(6.06,2.50)
\pspolygon[fillstyle=solid,linewidth=1pt,fillcolor=pink](6.93,2.00)(6.06,2.50)(6.93,1.00)
\pspolygon[fillstyle=solid,linewidth=1pt,fillcolor=pink](5.20,2.00)(5.20,1.00)(6.06,2.50)
\pspolygon[fillstyle=solid,linewidth=1pt,fillcolor=red](7.79,1.50)(8.66,1.00)(6.93,1.00)
\pspolygon[fillstyle=solid,linewidth=1pt,fillcolor=red](7.79,0.50)(8.66,1.00)(6.93,1.00)
\pspolygon[fillstyle=solid,linewidth=1pt,fillcolor=red](7.79,1.50)(7.79,2.50)(8.66,1.00)
\pspolygon[fillstyle=solid,linewidth=1pt,fillcolor=red](7.79,1.50)(6.93,1.00)(7.79,2.50)
\pspolygon[fillstyle=solid,linewidth=1pt,fillcolor=red](8.66,2.00)(7.79,2.50)(8.66,1.00)
\pspolygon[fillstyle=solid,linewidth=1pt,fillcolor=red](6.93,2.00)(6.93,1.00)(7.79,2.50)
\pspolygon[fillstyle=solid,linewidth=1pt,fillcolor=blue](3.46,3.00)(4.33,2.50)(2.60,2.50)
\pspolygon[fillstyle=solid,linewidth=1pt,fillcolor=blue](3.46,2.00)(4.33,2.50)(2.60,2.50)
\pspolygon[fillstyle=solid,linewidth=1pt,fillcolor=blue](3.46,3.00)(3.46,4.00)(4.33,2.50)
\pspolygon[fillstyle=solid,linewidth=1pt,fillcolor=blue](3.46,3.00)(2.60,2.50)(3.46,4.00)
\pspolygon[fillstyle=solid,linewidth=1pt,fillcolor=blue](4.33,3.50)(3.46,4.00)(4.33,2.50)
\pspolygon[fillstyle=solid,linewidth=1pt,fillcolor=blue](2.60,3.50)(2.60,2.50)(3.46,4.00)
\pspolygon[fillstyle=solid,linewidth=1pt,fillcolor=Plum](5.20,3.00)(6.06,2.50)(4.33,2.50)
\pspolygon[fillstyle=solid,linewidth=1pt,fillcolor=Plum](5.20,2.00)(6.06,2.50)(4.33,2.50)
\pspolygon[fillstyle=solid,linewidth=1pt,fillcolor=Plum](5.20,3.00)(5.20,4.00)(6.06,2.50)
\pspolygon[fillstyle=solid,linewidth=1pt,fillcolor=Plum](5.20,3.00)(4.33,2.50)(5.20,4.00)
\pspolygon[fillstyle=solid,linewidth=1pt,fillcolor=Plum](6.06,3.50)(5.20,4.00)(6.06,2.50)
\pspolygon[fillstyle=solid,linewidth=1pt,fillcolor=Plum](4.33,3.50)(4.33,2.50)(5.20,4.00)
\pspolygon[fillstyle=solid,linewidth=1pt,fillcolor=DarkCyan](6.93,3.00)(7.79,2.50)(6.06,2.50)
\pspolygon[fillstyle=solid,linewidth=1pt,fillcolor=DarkCyan](6.93,2.00)(7.79,2.50)(6.06,2.50)
\pspolygon[fillstyle=solid,linewidth=1pt,fillcolor=DarkCyan](6.93,3.00)(6.93,4.00)(7.79,2.50)
\pspolygon[fillstyle=solid,linewidth=1pt,fillcolor=DarkCyan](6.93,3.00)(6.06,2.50)(6.93,4.00)
\pspolygon[fillstyle=solid,linewidth=1pt,fillcolor=DarkCyan](7.79,3.50)(6.93,4.00)(7.79,2.50)
\pspolygon[fillstyle=solid,linewidth=1pt,fillcolor=DarkCyan](6.06,3.50)(6.06,2.50)(6.93,4.00)
\pspolygon[fillstyle=solid,linewidth=1pt,fillcolor=SlateGray](4.33,4.50)(5.20,4.00)(3.46,4.00)
\pspolygon[fillstyle=solid,linewidth=1pt,fillcolor=SlateGray](4.33,3.50)(5.20,4.00)(3.46,4.00)
\pspolygon[fillstyle=solid,linewidth=1pt,fillcolor=SlateGray](4.33,4.50)(4.33,5.50)(5.20,4.00)
\pspolygon[fillstyle=solid,linewidth=1pt,fillcolor=SlateGray](4.33,4.50)(3.46,4.00)(4.33,5.50)
\pspolygon[fillstyle=solid,linewidth=1pt,fillcolor=SlateGray](5.20,5.00)(4.33,5.50)(5.20,4.00)
\pspolygon[fillstyle=solid,linewidth=1pt,fillcolor=SlateGray](3.46,5.00)(3.46,4.00)(4.33,5.50)
\pspolygon[fillstyle=solid,linewidth=1pt,fillcolor=MediumSpringGreen](6.06,4.50)(6.93,4.00)(5.20,4.00)
\pspolygon[fillstyle=solid,linewidth=1pt,fillcolor=MediumSpringGreen](6.06,3.50)(6.93,4.00)(5.20,4.00)
\pspolygon[fillstyle=solid,linewidth=1pt,fillcolor=MediumSpringGreen](6.06,4.50)(6.06,5.50)(6.93,4.00)
\pspolygon[fillstyle=solid,linewidth=1pt,fillcolor=MediumSpringGreen](6.06,4.50)(5.20,4.00)(6.06,5.50)
\pspolygon[fillstyle=solid,linewidth=1pt,fillcolor=MediumSpringGreen](6.93,5.00)(6.06,5.50)(6.93,4.00)
\pspolygon[fillstyle=solid,linewidth=1pt,fillcolor=MediumSpringGreen](5.20,5.00)(5.20,4.00)(6.06,5.50)
\pspolygon[fillstyle=solid,linewidth=1pt,fillcolor=Azure](5.20,6.00)(6.06,5.50)(4.33,5.50)
\pspolygon[fillstyle=solid,linewidth=1pt,fillcolor=Azure](5.20,5.00)(6.06,5.50)(4.33,5.50)
\pspolygon[fillstyle=solid,linewidth=1pt,fillcolor=Azure](5.20,6.00)(5.20,7.00)(6.06,5.50)
\pspolygon[fillstyle=solid,linewidth=1pt,fillcolor=Azure](5.20,6.00)(4.33,5.50)(5.20,7.00)
\pspolygon[fillstyle=solid,linewidth=1pt,fillcolor=Azure](6.06,6.50)(5.20,7.00)(6.06,5.50)
\pspolygon[fillstyle=solid,linewidth=1pt,fillcolor=Azure](4.33,6.50)(4.33,5.50)(5.20,7.00)
\pspolygon[fillstyle=solid,linewidth=1pt,fillcolor=lightblue](0.87,0.50)(2.60,0.50)(1.73,1.00)
\pspolygon[fillstyle=solid,linewidth=1pt,fillcolor=lightblue](2.60,0.50)(4.33,0.50)(3.46,1.00)
\pspolygon[fillstyle=solid,linewidth=1pt,fillcolor=lightblue](4.33,0.50)(6.06,0.50)(5.20,1.00)
\pspolygon[fillstyle=solid,linewidth=1pt,fillcolor=lightblue](6.06,0.50)(7.79,0.50)(6.93,1.00)
\pspolygon[fillstyle=solid,linewidth=1pt,fillcolor=lightblue](7.79,0.50)(9.53,0.50)(8.66,1.00)
\pspolygon[fillstyle=solid,linewidth=1pt,fillcolor=lightblue](9.53,0.50)(8.66,2.00)(8.66,1.00)
\pspolygon[fillstyle=solid,linewidth=1pt,fillcolor=lightblue](8.66,2.00)(7.79,3.50)(7.79,2.50)
\pspolygon[fillstyle=solid,linewidth=1pt,fillcolor=lightblue](7.79,3.50)(6.93,5.00)(6.93,4.00)
\pspolygon[fillstyle=solid,linewidth=1pt,fillcolor=lightblue](6.93,5.00)(6.06,6.50)(6.06,5.50)
\pspolygon[fillstyle=solid,linewidth=1pt,fillcolor=lightblue](6.06,6.50)(5.20,8.00)(5.20,7.00)
\pspolygon[fillstyle=solid,linewidth=1pt,fillcolor=lightblue](1.73,2.00)(0.87,0.50)(1.73,1.00)
\pspolygon[fillstyle=solid,linewidth=1pt,fillcolor=lightblue](2.60,3.50)(1.73,2.00)(2.60,2.50)
\pspolygon[fillstyle=solid,linewidth=1pt,fillcolor=lightblue](3.46,5.00)(2.60,3.50)(3.46,4.00)
\pspolygon[fillstyle=solid,linewidth=1pt,fillcolor=lightblue](4.33,6.50)(3.46,5.00)(4.33,5.50)
\pspolygon[fillstyle=solid,linewidth=1pt,fillcolor=lightblue](5.20,8.00)(4.33,6.50)(5.20,7.00)
\endpspicture
}

\def\TwentySevenTiling{%
\psset{unit=0.5cm}
\pspicture(10.3923048,9.4)
\qline(0,0)(10.3923048,0)             
\qline(0,0)(5.19615242,9)             
\qline(5.19615242,9)(10.3923048,0)    
\qline(0,0)(1.73205081,1)             
\qline(1.73205081,1)(3.46410162,0)
\qline(3.46410162,0)(5.19615242,1)
\qline(5.19615242,1)(6.92820323,0)
\qline(6.92820323,0)(8.66025404,1)
\qline(8.66025404,1)(10.3923048,0)
\qline(1.73205081,1)(8.66025404,1)  
\qline(1.73205081,1)(1.73205081,3)  
\qline(5.19615242,1)(5.19615242,3)
\qline(8.66025404,1)(8.66025404,3)
\qline(1.73205081,3)(3.46410162,4)  
\qline(3.46410162,4)(5.19615242,3)
\qline(5.19615242,3)(6.92820323,4)
\qline(6.92820323,4)(8.66025404,3)  
\qline(3.46410162,4)(3.46410162,6)  
\qline(3.46410162,6)(5.19615242,7)  
\qline(5.19615242,7)(6.92820323,6)  
\qline(6.92820323,6)(6.92820323,4)  
\qline(5.19615242,7)(5.19615242,9)  
\qline(1.73205081,1)(5.19615242,7)  
\qline(5.19615242,7)(8.66025404,1)
\qline(3.46410162,4)(6.92820323,4)
\qline(3.46410162,4)(5.19615242,1)
\qline(5.19615242,1)(6.92820323,4)
\qline(3.46410162,2)(1.73205081,1)  
\qline(3.46410162,2)(3.46410162,4)
\qline(3.46410162,2)(5.19615242,1)
\qline(6.92820323,2)(5.19615242,1)  
\qline(6.92820323,2)(6.92820323,4)
\qline(6.92820323,2)(8.66025404,1)
\qline(5.19615242,5)(3.46410162,4)  
\qline(5.19615242,5)(6.92820323,4)
\qline(5.19615242,5)(5.19615242,7)
\endpspicture}

\def\FiftyFourTiling{%
\pspicture(3.85,1.2)
\psset{unit=2.2cm}
\pspolygon[fillstyle=solid,linewidth=1pt,fillcolor=lightblue](0.00,0.00)(0.43,0.25)(0.58,0.00)\pspolygon[fillstyle=solid,linewidth=1pt,fillcolor=lightyellow](0.43,0.25)(0.58,0.00)(0.87,0.50)\pspolygon[fillstyle=solid,linewidth=1pt,fillcolor=lightgreen](0.58,0.00)(0.87,0.00)(0.87,0.50)\pspolygon[fillstyle=solid,linewidth=1pt,fillcolor=pink](1.73,0.50)(1.30,0.25)(1.15,0.50)\pspolygon[fillstyle=solid,linewidth=1pt,fillcolor=blue](1.30,0.25)(1.15,0.50)(0.87,0.00)\pspolygon[fillstyle=solid,linewidth=1pt,fillcolor=red](1.15,0.50)(0.87,0.00)(0.87,0.50)\pspolygon[fillstyle=solid,linewidth=1pt,fillcolor=lightblue](0.87,0.00)(1.30,0.25)(1.44,0.00)\pspolygon[fillstyle=solid,linewidth=1pt,fillcolor=lightyellow](1.30,0.25)(1.44,0.00)(1.73,0.50)\pspolygon[fillstyle=solid,linewidth=1pt,fillcolor=lightgreen](1.44,0.00)(1.73,0.00)(1.73,0.50)\pspolygon[fillstyle=solid,linewidth=1pt,fillcolor=pink](2.60,0.50)(2.17,0.25)(2.02,0.50)\pspolygon[fillstyle=solid,linewidth=1pt,fillcolor=blue](2.17,0.25)(2.02,0.50)(1.73,0.00)\pspolygon[fillstyle=solid,linewidth=1pt,fillcolor=red](2.02,0.50)(1.73,0.00)(1.73,0.50)\pspolygon[fillstyle=solid,linewidth=1pt,fillcolor=lightblue](0.87,0.50)(1.30,0.75)(1.44,0.50)\pspolygon[fillstyle=solid,linewidth=1pt,fillcolor=lightyellow](1.30,0.75)(1.44,0.50)(1.73,1.00)\pspolygon[fillstyle=solid,linewidth=1pt,fillcolor=lightgreen](1.44,0.50)(1.73,0.50)(1.73,1.00)\pspolygon[fillstyle=solid,linewidth=1pt,fillcolor=pink](2.60,1.00)(2.17,0.75)(2.02,1.00)\pspolygon[fillstyle=solid,linewidth=1pt,fillcolor=blue](2.17,0.75)(2.02,1.00)(1.73,0.50)\pspolygon[fillstyle=solid,linewidth=1pt,fillcolor=red](2.02,1.00)(1.73,0.50)(1.73,1.00)\pspolygon[fillstyle=solid,linewidth=1pt,fillcolor=lightblue](1.73,0.00)(2.17,0.25)(2.31,0.00)\pspolygon[fillstyle=solid,linewidth=1pt,fillcolor=lightyellow](2.17,0.25)(2.31,0.00)(2.60,0.50)\pspolygon[fillstyle=solid,linewidth=1pt,fillcolor=lightgreen](2.31,0.00)(2.60,0.00)(2.60,0.50)\pspolygon[fillstyle=solid,linewidth=1pt,fillcolor=lightblue](1.73,0.50)(2.17,0.75)(2.31,0.50)\pspolygon[fillstyle=solid,linewidth=1pt,fillcolor=lightyellow](2.17,0.75)(2.31,0.50)(2.60,1.00)\pspolygon[fillstyle=solid,linewidth=1pt,fillcolor=lightgreen](2.31,0.50)(2.60,0.50)(2.60,1.00)\pspolygon[fillstyle=solid,linewidth=1pt,fillcolor=lightblue](1.73,1.00)(2.17,1.25)(2.31,1.00)\pspolygon[fillstyle=solid,linewidth=1pt,fillcolor=lightyellow](2.17,1.25)(2.31,1.00)(2.60,1.50)\pspolygon[fillstyle=solid,linewidth=1pt,fillcolor=lightgreen](2.31,1.00)(2.60,1.00)(2.60,1.50)\pspolygon[fillstyle=solid,linewidth=1pt,fillcolor=lightblue](5.20,0.00)(4.76,0.25)(4.62,0.00)\pspolygon[fillstyle=solid,linewidth=1pt,fillcolor=lightyellow](4.76,0.25)(4.62,0.00)(4.33,0.50)\pspolygon[fillstyle=solid,linewidth=1pt,fillcolor=lightgreen](4.62,0.00)(4.33,0.00)(4.33,0.50)\pspolygon[fillstyle=solid,linewidth=1pt,fillcolor=pink](3.46,0.50)(3.90,0.25)(4.04,0.50)\pspolygon[fillstyle=solid,linewidth=1pt,fillcolor=blue](3.90,0.25)(4.04,0.50)(4.33,0.00)\pspolygon[fillstyle=solid,linewidth=1pt,fillcolor=red](4.04,0.50)(4.33,0.00)(4.33,0.50)\pspolygon[fillstyle=solid,linewidth=1pt,fillcolor=lightblue](4.33,0.00)(3.90,0.25)(3.75,0.00)\pspolygon[fillstyle=solid,linewidth=1pt,fillcolor=lightyellow](3.90,0.25)(3.75,0.00)(3.46,0.50)\pspolygon[fillstyle=solid,linewidth=1pt,fillcolor=lightgreen](3.75,0.00)(3.46,0.00)(3.46,0.50)\pspolygon[fillstyle=solid,linewidth=1pt,fillcolor=pink](2.60,0.50)(3.03,0.25)(3.18,0.50)\pspolygon[fillstyle=solid,linewidth=1pt,fillcolor=blue](3.03,0.25)(3.18,0.50)(3.46,0.00)\pspolygon[fillstyle=solid,linewidth=1pt,fillcolor=red](3.18,0.50)(3.46,0.00)(3.46,0.50)\pspolygon[fillstyle=solid,linewidth=1pt,fillcolor=lightblue](4.33,0.50)(3.90,0.75)(3.75,0.50)\pspolygon[fillstyle=solid,linewidth=1pt,fillcolor=lightyellow](3.90,0.75)(3.75,0.50)(3.46,1.00)\pspolygon[fillstyle=solid,linewidth=1pt,fillcolor=lightgreen](3.75,0.50)(3.46,0.50)(3.46,1.00)\pspolygon[fillstyle=solid,linewidth=1pt,fillcolor=pink](2.60,1.00)(3.03,0.75)(3.18,1.00)\pspolygon[fillstyle=solid,linewidth=1pt,fillcolor=blue](3.03,0.75)(3.18,1.00)(3.46,0.50)\pspolygon[fillstyle=solid,linewidth=1pt,fillcolor=red](3.18,1.00)(3.46,0.50)(3.46,1.00)\pspolygon[fillstyle=solid,linewidth=1pt,fillcolor=lightblue](3.46,0.00)(3.03,0.25)(2.89,0.00)\pspolygon[fillstyle=solid,linewidth=1pt,fillcolor=lightyellow](3.03,0.25)(2.89,0.00)(2.60,0.50)\pspolygon[fillstyle=solid,linewidth=1pt,fillcolor=lightgreen](2.89,0.00)(2.60,0.00)(2.60,0.50)\pspolygon[fillstyle=solid,linewidth=1pt,fillcolor=lightblue](3.46,0.50)(3.03,0.75)(2.89,0.50)\pspolygon[fillstyle=solid,linewidth=1pt,fillcolor=lightyellow](3.03,0.75)(2.89,0.50)(2.60,1.00)\pspolygon[fillstyle=solid,linewidth=1pt,fillcolor=lightgreen](2.89,0.50)(2.60,0.50)(2.60,1.00)\pspolygon[fillstyle=solid,linewidth=1pt,fillcolor=lightblue](3.46,1.00)(3.03,1.25)(2.89,1.00)\pspolygon[fillstyle=solid,linewidth=1pt,fillcolor=lightyellow](3.03,1.25)(2.89,1.00)(2.60,1.50)\pspolygon[fillstyle=solid,linewidth=1pt,fillcolor=lightgreen](2.89,1.00)(2.60,1.00)(2.60,1.50)
\endpspicture
}

\def\TwelveTilingA{
\pspicture(1.4,0.95)
\psset{unit=0.005cm} 
\qline(409.5,10)(609.25,125.3257)
\qline(609.25,125.3257)(675.8334,10)
\qline(675.8334,10)(409.5,10)
\qline(809,240.6514)(809,10)
\qline(809,10)(675.8334,10)
\qline(675.8334,10)(809,240.6514)
\qline(809,240.6514)(609.25,125.3257)
\qline(609.25,125.3257)(675.8334,10)
\qline(675.8334,10)(809,240.6514)
\qline(809,240.6514)(609.25,125.3257)
\qline(609.25,125.3257)(542.6667,240.6514)
\qline(542.6667,240.6514)(809,240.6514)
\qline(409.5,10)(409.5,240.6514)
\qline(409.5,240.6514)(542.6667,240.6514)
\qline(542.6667,240.6514)(409.5,10)
\qline(409.5,10)(609.25,125.3257)
\qline(609.25,125.3257)(542.6667,240.6514)
\qline(542.6667,240.6514)(409.5,10)
\qline(9.999969,10)(209.75,125.3257)
\qline(209.75,125.3257)(276.3333,10)
\qline(276.3333,10)(9.999969,10)
\qline(409.5,240.6514)(409.5,10)
\qline(409.5,10)(276.3333,10)
\qline(276.3333,10)(409.5,240.6514)
\qline(409.5,240.6514)(209.75,125.3257)
\qline(209.75,125.3257)(276.3333,10)
\qline(276.3333,10)(409.5,240.6514)
\qline(409.5,240.6514)(609.25,355.9771)
\qline(609.25,355.9771)(675.8334,240.6514)
\qline(675.8334,240.6514)(409.5,240.6514)
\qline(809,471.3029)(809,240.6514)
\qline(809,240.6514)(675.8334,240.6514)
\qline(675.8334,240.6514)(809,471.3029)
\qline(809,471.3029)(609.25,355.9771)
\qline(609.25,355.9771)(675.8334,240.6514)
\qline(675.8334,240.6514)(809,471.3029)
\endpspicture}

\def\TwelveTilingB{
\pspicture(1.4,0.95)
\psset{unit=0.005cm} 
\qline(276.3333,10)(476.0833,125.3257)
\qline(476.0833,125.3257)(542.6667,10)
\qline(542.6667,10)(276.3333,10)
\qline(675.8333,240.6514)(675.8333,10)
\qline(675.8333,10)(542.6667,10)
\qline(542.6667,10)(675.8333,240.6514)
\qline(675.8333,240.6514)(476.0833,125.3257)
\qline(476.0833,125.3257)(542.6667,10)
\qline(542.6667,10)(675.8333,240.6514)
\qline(9.999969,10)(276.3333,10)
\qline(276.3333,10)(209.75,125.3257)
\qline(209.75,125.3257)(9.999969,10)
\qline(409.5,240.6514)(276.3333,10)
\qline(276.3333,10)(209.75,125.3257)
\qline(209.75,125.3257)(409.5,240.6514)
\qline(409.5,240.6514)(276.3333,10)
\qline(276.3333,10)(476.0833,125.3257)
\qline(476.0833,125.3257)(409.5,240.6514)
\qline(409.5,240.6514)(675.8334,240.6514)
\qline(675.8334,240.6514)(609.2501,355.9771)
\qline(609.2501,355.9771)(409.5,240.6514)
\qline(809,471.3029)(675.8334,240.6514)
\qline(675.8334,240.6514)(609.2501,355.9771)
\qline(609.2501,355.9771)(809,471.3029)
\qline(675.8333,10)(809,10)
\qline(809,10)(675.8333,240.6514)
\qline(675.8333,240.6514)(675.8333,10)
\qline(809,240.6514)(809,10)
\qline(809,10)(675.8333,240.6514)
\qline(675.8333,240.6514)(809,240.6514)
\qline(675.8333,240.6514)(809,240.6514)
\qline(809,240.6514)(809,471.3029)
\qline(809,471.3029)(675.8333,240.6514)
\endpspicture}
\def\TwentySevenTilingA{
\pspicture(1.4,1.0)
\psset{unit=0.005cm} 
\qline(542.6667,10)(675.8334,86.88382)
\qline(675.8334,86.88382)(720.2222,10)
\qline(720.2222,10)(542.6667,10)
\qline(809,163.7676)(809,10)
\qline(809,10)(720.2222,10)
\qline(720.2222,10)(809,163.7676)
\qline(809,163.7676)(675.8334,86.88382)
\qline(675.8334,86.88382)(720.2222,10)
\qline(720.2222,10)(809,163.7676)
\qline(809,163.7676)(675.8334,86.88382)
\qline(675.8334,86.88382)(631.4445,163.7676)
\qline(631.4445,163.7676)(809,163.7676)
\qline(542.6667,10)(542.6667,163.7676)
\qline(542.6667,163.7676)(631.4445,163.7676)
\qline(631.4445,163.7676)(542.6667,10)
\qline(542.6667,10)(675.8334,86.88382)
\qline(675.8334,86.88382)(631.4445,163.7676)
\qline(631.4445,163.7676)(542.6667,10)
\qline(276.3333,10)(409.5,86.88382)
\qline(409.5,86.88382)(453.8889,10)
\qline(453.8889,10)(276.3333,10)
\qline(542.6667,163.7676)(542.6667,10)
\qline(542.6667,10)(453.8889,10)
\qline(453.8889,10)(542.6667,163.7676)
\qline(542.6667,163.7676)(409.5,86.88382)
\qline(409.5,86.88382)(453.8889,10)
\qline(453.8889,10)(542.6667,163.7676)
\qline(542.6667,163.7676)(409.5,86.88382)
\qline(409.5,86.88382)(365.1111,163.7676)
\qline(365.1111,163.7676)(542.6667,163.7676)
\qline(276.3333,10)(276.3333,163.7676)
\qline(276.3333,163.7676)(365.1111,163.7676)
\qline(365.1111,163.7676)(276.3333,10)
\qline(276.3333,10)(409.5,86.88382)
\qline(409.5,86.88382)(365.1111,163.7676)
\qline(365.1111,163.7676)(276.3333,10)
\qline(9.999992,10)(143.1667,86.88382)
\qline(143.1667,86.88382)(187.5555,10)
\qline(187.5555,10)(9.999992,10)
\qline(276.3333,163.7676)(276.3333,10)
\qline(276.3333,10)(187.5555,10)
\qline(187.5555,10)(276.3333,163.7676)
\qline(276.3333,163.7676)(143.1667,86.88382)
\qline(143.1667,86.88382)(187.5555,10)
\qline(187.5555,10)(276.3333,163.7676)
\qline(542.6667,163.7676)(675.8334,240.6514)
\qline(675.8334,240.6514)(720.2222,163.7676)
\qline(720.2222,163.7676)(542.6667,163.7676)
\qline(809,317.5352)(809,163.7676)
\qline(809,163.7676)(720.2222,163.7676)
\qline(720.2222,163.7676)(809,317.5352)
\qline(809,317.5352)(675.8334,240.6514)
\qline(675.8334,240.6514)(720.2222,163.7676)
\qline(720.2222,163.7676)(809,317.5352)
\qline(809,317.5352)(675.8334,240.6514)
\qline(675.8334,240.6514)(631.4445,317.5352)
\qline(631.4445,317.5352)(809,317.5352)
\qline(542.6667,163.7676)(542.6667,317.5353)
\qline(542.6667,317.5353)(631.4445,317.5352)
\qline(631.4445,317.5352)(542.6667,163.7676)
\qline(542.6667,163.7676)(675.8334,240.6514)
\qline(675.8334,240.6514)(631.4445,317.5352)
\qline(631.4445,317.5352)(542.6667,163.7676)
\qline(276.3333,163.7676)(409.5,240.6514)
\qline(409.5,240.6514)(453.8889,163.7676)
\qline(453.8889,163.7676)(276.3333,163.7676)
\qline(542.6667,317.5352)(542.6667,163.7676)
\qline(542.6667,163.7676)(453.8889,163.7676)
\qline(453.8889,163.7676)(542.6667,317.5352)
\qline(542.6667,317.5352)(409.5,240.6514)
\qline(409.5,240.6514)(453.8889,163.7676)
\qline(453.8889,163.7676)(542.6667,317.5352)
\qline(542.6667,317.5352)(675.8334,394.4191)
\qline(675.8334,394.4191)(720.2222,317.5352)
\qline(720.2222,317.5352)(542.6667,317.5352)
\qline(809,471.3029)(809,317.5352)
\qline(809,317.5352)(720.2222,317.5352)
\qline(720.2222,317.5352)(809,471.3029)
\qline(809,471.3029)(675.8334,394.4191)
\qline(675.8334,394.4191)(720.2222,317.5352)
\qline(720.2222,317.5352)(809,471.3029)
\endpspicture}
\def\TwentySevenTilingB{
\pspicture(1.4,1.0)
\psset{unit=0.005cm} 
\qline(453.8889,10)(587.0555,86.88382)
\qline(587.0555,86.88382)(631.4445,10)
\qline(631.4445,10)(453.8889,10)
\qline(720.2222,163.7676)(720.2222,10)
\qline(720.2222,10)(631.4445,10)
\qline(631.4445,10)(720.2222,163.7676)
\qline(720.2222,163.7676)(587.0555,86.88382)
\qline(587.0555,86.88382)(631.4445,10)
\qline(631.4445,10)(720.2222,163.7676)
\qline(720.2222,163.7676)(587.0555,86.88382)
\qline(587.0555,86.88382)(542.6667,163.7676)
\qline(542.6667,163.7676)(720.2222,163.7676)
\qline(453.8889,10)(453.8889,163.7676)
\qline(453.8889,163.7676)(542.6667,163.7676)
\qline(542.6667,163.7676)(453.8889,10)
\qline(453.8889,10)(587.0555,86.88382)
\qline(587.0555,86.88382)(542.6667,163.7676)
\qline(542.6667,163.7676)(453.8889,10)
\qline(187.5555,10)(320.7222,86.88382)
\qline(320.7222,86.88382)(365.1111,10)
\qline(365.1111,10)(187.5555,10)
\qline(453.8889,163.7676)(453.8889,10)
\qline(453.8889,10)(365.1111,10)
\qline(365.1111,10)(453.8889,163.7676)
\qline(453.8889,163.7676)(320.7222,86.88382)
\qline(320.7222,86.88382)(365.1111,10)
\qline(365.1111,10)(453.8889,163.7676)
\qline(453.8889,163.7676)(587.0555,240.6514)
\qline(587.0555,240.6514)(631.4445,163.7676)
\qline(631.4445,163.7676)(453.8889,163.7676)
\qline(720.2222,317.5352)(720.2222,163.7676)
\qline(720.2222,163.7676)(631.4445,163.7676)
\qline(631.4445,163.7676)(720.2222,317.5352)
\qline(720.2222,317.5352)(587.0555,240.6514)
\qline(587.0555,240.6514)(631.4445,163.7676)
\qline(631.4445,163.7676)(720.2222,317.5352)
\qline(9.999992,10)(187.5555,10)
\qline(187.5555,10)(143.1667,86.88382)
\qline(143.1667,86.88382)(9.999992,10)
\qline(276.3333,163.7676)(187.5555,10)
\qline(187.5555,10)(143.1667,86.88382)
\qline(143.1667,86.88382)(276.3333,163.7676)
\qline(276.3333,163.7676)(187.5555,10)
\qline(187.5555,10)(320.7222,86.88382)
\qline(320.7222,86.88382)(276.3333,163.7676)
\qline(276.3333,163.7676)(453.8889,163.7676)
\qline(453.8889,163.7676)(409.5,240.6515)
\qline(409.5,240.6515)(276.3333,163.7676)
\qline(542.6666,317.5353)(453.8889,163.7676)
\qline(453.8889,163.7676)(409.5,240.6515)
\qline(409.5,240.6515)(542.6666,317.5353)
\qline(542.6666,317.5353)(453.8889,163.7676)
\qline(453.8889,163.7676)(587.0555,240.6515)
\qline(587.0555,240.6515)(542.6666,317.5353)
\qline(542.6667,317.5352)(720.2222,317.5352)
\qline(720.2222,317.5352)(675.8333,394.4191)
\qline(675.8333,394.4191)(542.6667,317.5352)
\qline(809,471.3029)(720.2222,317.5352)
\qline(720.2222,317.5352)(675.8333,394.4191)
\qline(675.8333,394.4191)(809,471.3029)
\qline(720.2222,10)(809,10)
\qline(809,10)(809,163.7676)
\qline(809,163.7676)(720.2222,10)
\qline(720.2222,163.7676)(720.2222,10)
\qline(720.2222,10)(809,163.7676)
\qline(809,163.7676)(720.2222,163.7676)
\qline(720.2222,163.7676)(809,163.7676)
\qline(809,163.7676)(720.2222,317.5353)
\qline(720.2222,317.5353)(720.2222,163.7676)
\qline(809,317.5353)(809,163.7676)
\qline(809,163.7676)(720.2222,317.5353)
\qline(720.2222,317.5353)(809,317.5353)
\qline(720.2222,317.5352)(809,317.5352)
\qline(809,317.5352)(809,471.3029)
\qline(809,471.3029)(720.2222,317.5352)
\endpspicture}

\def\NineTilingA{%
\psset{unit=1cm}
\pspicture(6,3.3)
\qline(0,0)(6,0)            
\qline(0,0)(6,3)            
\qline(6,0)(6,3)            
\qline(2,0)(2,1)         
\qline(4,0)(4,2)          
\qline(5,0)(5,2)
\qline(4,2)(6,2)        
\qline(2,1)(4,1)
\qline(2,0)(4,1)        
\qline(4,0)(5,2)
\qline(5,0)(6,2)
\endpspicture
\psset{unit=3cm}
}

\def\FigurePythagorean{
\pspicture(1.7,1.1)
\psset{unit=0.15cm}
\newrgbcolor{lightblue}{0.8 0.8 1}
\newrgbcolor{pink}{1 0.8 0.8}
\newrgbcolor{lightgreen}{0.8 1 0.8}
\newrgbcolor{lightyellow}{1 1 0.8}
\pspolygon[fillstyle=solid,linewidth=1pt,fillcolor=lightyellow](0.00,0.00)(15.00,20.00)(15.00,0.00)
\psline(15.00,20.00)(15.00,0.00)
\psline(12.00,16.00)(12.00,0.00)
\psline(9.00,12.00)(9.00,0.00)
\psline(6.00,8.00)(6.00,0.00)
\psline(3.00,4.00)(3.00,0.00)
\psline(0.00,0.00)(0.00,0.00)
\psline(15.00,0.00)(0.00,0.00)
\psline(15.00,4.00)(3.00,4.00)
\psline(15.00,8.00)(6.00,8.00)
\psline(15.00,12.00)(9.00,12.00)
\psline(15.00,16.00)(12.00,16.00)
\psline(15.00,20.00)(15.00,20.00)
\psline(0.00,0.00)(15.00,20.00)
\psline(3.00,0.00)(15.00,16.00)
\psline(6.00,0.00)(15.00,12.00)
\psline(9.00,0.00)(15.00,8.00)
\psline(12.00,0.00)(15.00,4.00)
\psline(15.00,0.00)(15.00,0.00)
\pspolygon[fillstyle=solid,linewidth=1pt,fillcolor=lightblue](15.00,20.00)(15.00,0.00)(24.60,7.20)
\psline(15.00,0.00)(24.60,7.20)
\psline(15.00,5.00)(22.20,10.40)
\psline(15.00,10.00)(19.80,13.60)
\psline(15.00,15.00)(17.40,16.80)
\psline(15.00,20.00)(15.00,20.00)
\psline(24.60,7.20)(15.00,20.00)
\psline(22.20,5.40)(15.00,15.00)
\psline(19.80,3.60)(15.00,10.00)
\psline(17.40,1.80)(15.00,5.00)
\psline(15.00,0.00)(15.00,0.00)
\psline(15.00,20.00)(15.00,0.00)
\psline(17.40,16.80)(17.40,1.80)
\psline(19.80,13.60)(19.80,3.60)
\psline(22.20,10.40)(22.20,5.40)
\psline(24.60,7.20)(24.60,7.20)
\pspolygon[fillstyle=solid,linewidth=1pt,fillcolor=lightgreen](15.00,0.00)(24.60,7.20)(30.00,0.00)
\psline(24.60,7.20)(30.00,0.00)
\psline(21.40,4.80)(25.00,0.00)
\psline(18.20,2.40)(20.00,0.00)
\psline(15.00,0.00)(15.00,0.00)
\psline(30.00,0.00)(15.00,0.00)
\psline(28.20,2.40)(18.20,2.40)
\psline(26.40,4.80)(21.40,4.80)
\psline(24.60,7.20)(24.60,7.20)
\psline(15.00,0.00)(24.60,7.20)
\psline(20.00,0.00)(26.40,4.80)
\psline(25.00,0.00)(28.20,2.40)
\psline(30.00,0.00)(30.00,0.00)
\endpspicture}

\def\FigureDoubleBiquadratic{
\pspicture(1.7,1.1)
\psset{unit=0.15cm}
\newrgbcolor{lightblue}{0.8 0.8 1}
\newrgbcolor{pink}{1 0.8 0.8}
\newrgbcolor{lightgreen}{0.8 1 0.8}
\newrgbcolor{lightyellow}{1 1 0.8}
\pspolygon[fillstyle=solid,linewidth=1pt,fillcolor=lightblue](15.00,20.00)(15.00,0.00)(24.60,7.20)
\psline(15.00,0.00)(24.60,7.20)
\psline(15.00,5.00)(22.20,10.40)
\psline(15.00,10.00)(19.80,13.60)
\psline(15.00,15.00)(17.40,16.80)
\psline(15.00,20.00)(15.00,20.00)
\psline(24.60,7.20)(15.00,20.00)
\psline(22.20,5.40)(15.00,15.00)
\psline(19.80,3.60)(15.00,10.00)
\psline(17.40,1.80)(15.00,5.00)
\psline(15.00,0.00)(15.00,0.00)
\psline(15.00,20.00)(15.00,0.00)
\psline(17.40,16.80)(17.40,1.80)
\psline(19.80,13.60)(19.80,3.60)
\psline(22.20,10.40)(22.20,5.40)
\psline(24.60,7.20)(24.60,7.20)
\pspolygon[fillstyle=solid,linewidth=1pt,fillcolor=lightgreen](15.00,0.00)(24.60,7.20)(30.00,0.00)
\psline(24.60,7.20)(30.00,0.00)
\psline(21.40,4.80)(25.00,0.00)
\psline(18.20,2.40)(20.00,0.00)
\psline(15.00,0.00)(15.00,0.00)
\psline(30.00,0.00)(15.00,0.00)
\psline(28.20,2.40)(18.20,2.40)
\psline(26.40,4.80)(21.40,4.80)
\psline(24.60,7.20)(24.60,7.20)
\psline(15.00,0.00)(24.60,7.20)
\psline(20.00,0.00)(26.40,4.80)
\psline(25.00,0.00)(28.20,2.40)
\psline(30.00,0.00)(30.00,0.00)
\pspolygon[fillstyle=solid,linewidth=1pt,fillcolor=lightyellow](15.00,20.00)(15.00,0.00)(5.40,7.20)
\psline(15.00,0.00)(5.40,7.20)
\psline(15.00,5.00)(7.80,10.40)
\psline(15.00,10.00)(10.20,13.60)
\psline(15.00,15.00)(12.60,16.80)
\psline(15.00,20.00)(15.00,20.00)
\psline(5.40,7.20)(15.00,20.00)
\psline(7.80,5.40)(15.00,15.00)
\psline(10.20,3.60)(15.00,10.00)
\psline(12.60,1.80)(15.00,5.00)
\psline(15.00,0.00)(15.00,0.00)
\psline(15.00,20.00)(15.00,0.00)
\psline(12.60,16.80)(12.60,1.80)
\psline(10.20,13.60)(10.20,3.60)
\psline(7.80,10.40)(7.80,5.40)
\psline(5.40,7.20)(5.40,7.20)
\pspolygon[fillstyle=solid,linewidth=1pt,fillcolor=pink](15.00,0.00)(5.40,7.20)(0.00,0.00)
\psline(5.40,7.20)(0.00,0.00)
\psline(8.60,4.80)(5.00,0.00)
\psline(11.80,2.40)(10.00,0.00)
\psline(15.00,0.00)(15.00,0.00)
\psline(0.00,0.00)(15.00,0.00)
\psline(1.80,2.40)(11.80,2.40)
\psline(3.60,4.80)(8.60,4.80)
\psline(5.40,7.20)(5.40,7.20)
\psline(15.00,0.00)(5.40,7.20)
\psline(10.00,0.00)(3.60,4.80)
\psline(5.00,0.00)(1.80,2.40)
\psline(0.00,0.00)(0.00,0.00)
\endpspicture}

\def\FigureDoubleQuadratic{
\pspicture(1.7,1.1)
\psset{unit=0.15cm}
\newrgbcolor{lightblue}{0.8 0.8 1}
\newrgbcolor{pink}{1 0.8 0.8}
\newrgbcolor{lightgreen}{0.8 1 0.8}
\newrgbcolor{lightyellow}{1 1 0.8}
\pspolygon[fillstyle=solid,linewidth=1pt,fillcolor=lightblue](15.00,20.00)(15.00,0.00)(30.00,0.00)
\psline(15.00,0.00)(30.00,0.00)
\psline(15.00,4.00)(27.00,4.00)
\psline(15.00,8.00)(24.00,8.00)
\psline(15.00,12.00)(21.00,12.00)
\psline(15.00,16.00)(18.00,16.00)
\psline(15.00,20.00)(15.00,20.00)
\psline(30.00,0.00)(15.00,20.00)
\psline(27.00,0.00)(15.00,16.00)
\psline(24.00,0.00)(15.00,12.00)
\psline(21.00,0.00)(15.00,8.00)
\psline(18.00,0.00)(15.00,4.00)
\psline(15.00,0.00)(15.00,0.00)
\psline(15.00,20.00)(15.00,0.00)
\psline(18.00,16.00)(18.00,0.00)
\psline(21.00,12.00)(21.00,0.00)
\psline(24.00,8.00)(24.00,0.00)
\psline(27.00,4.00)(27.00,0.00)
\psline(30.00,0.00)(30.00,0.00)
\pspolygon[fillstyle=solid,linewidth=1pt,fillcolor=lightyellow](15.00,20.00)(15.00,0.00)(0.00,0.00)
\psline(15.00,0.00)(0.00,0.00)
\psline(15.00,4.00)(3.00,4.00)
\psline(15.00,8.00)(6.00,8.00)
\psline(15.00,12.00)(9.00,12.00)
\psline(15.00,16.00)(12.00,16.00)
\psline(15.00,20.00)(15.00,20.00)
\psline(0.00,0.00)(15.00,20.00)
\psline(3.00,0.00)(15.00,16.00)
\psline(6.00,0.00)(15.00,12.00)
\psline(9.00,0.00)(15.00,8.00)
\psline(12.00,0.00)(15.00,4.00)
\psline(15.00,0.00)(15.00,0.00)
\psline(15.00,20.00)(15.00,0.00)
\psline(12.00,16.00)(12.00,0.00)
\psline(9.00,12.00)(9.00,0.00)
\psline(6.00,8.00)(6.00,0.00)
\psline(3.00,4.00)(3.00,0.00)
\psline(0.00,0.00)(0.00,0.00)
\endpspicture}

\def\FigureBigIsosceles{
\pspicture(4,1.7)
\psset{unit=0.0007cm}
\newrgbcolor{lightblue}{0.8 0.8 1}
\newrgbcolor{pink}{1 0.8 0.8}
\newrgbcolor{lightgreen}{0.8 1 0.8}
\newrgbcolor{lightyellow}{1 1 0.8}
\newrgbcolor{orange}{1 0.5 0}
\newrgbcolor{yellow}{1 1 0}
\psline(0,0)(100,100)
\pspolygon[fillstyle=solid,linewidth=1pt,fillcolor=lightyellow](0.00,0.00)(8235.00,7262.59)(16470.00,0.00)
\pspolygon[fillstyle=solid,linewidth=1pt,fillcolor=pink](2160.00,1904.94)(8235.00,7262.59)(8910.00,1904.94)
\pspolygon[fillstyle=solid,linewidth=1pt,fillcolor=lightblue](8235.00,7262.59)(8910.00,1904.94)(13095.00,2976.47)\pspolygon[fillstyle=solid,linewidth=1pt,fillcolor=red](5214.00,0.00)(7374.00,1904.94)(8670.00,0.00)\pspolygon[fillstyle=solid,linewidth=1pt,fillcolor=orange](7374.00,1904.94)(8910.00,1904.94)(8670.00,0.00)\pspolygon[fillstyle=solid,linewidth=1pt,fillcolor=lightgreen](8910.00,1904.94)(8670.00,0.00)(11070.00,0.00)\pspolygon[fillstyle=solid,linewidth=1pt,fillcolor=blue](8910.00,1904.94)(13095.00,2976.47)(11070.00,0.00)\pspolygon[fillstyle=solid,linewidth=1pt,fillcolor=green](13095.00,2976.47)(16470.00,0.00)(11070.00,0.00)
\pspolygon[fillstyle=solid,linewidth=1pt,fillcolor=violet](0.00,0.00)(2160.00,1904.94)(5214.00,0.00)
\endpspicture}

\def\ColoringTheorem{
\newrgbcolor{lightblue}{0.8 0.8 1}
\psset{unit=0.5cm}
\pspicture(12,14)
\pspolygon[fillstyle=solid,linewidth=1pt,fillcolor=black](7.88,4.36)(0.00,0.00)(12.00,0.00)\pspolygon[fillstyle=solid,linewidth=1pt,fillcolor=white](9.25,2.90)(5.25,2.90)(7.88,4.36)\pspolygon[fillstyle=solid,linewidth=1pt,fillcolor=white](10.62,1.45)(6.62,1.45)(9.25,2.90)\pspolygon[fillstyle=solid,linewidth=1pt,fillcolor=white](6.62,1.45)(2.63,1.45)(5.25,2.90)\pspolygon[fillstyle=solid,linewidth=1pt,fillcolor=white](12.00,0.00)(8.00,0.00)(10.62,1.45)\pspolygon[fillstyle=solid,linewidth=1pt,fillcolor=white](8.00,0.00)(4.00,0.00)(6.62,1.45)\pspolygon[fillstyle=solid,linewidth=1pt,fillcolor=white](4.00,0.00)(0.00,0.00)(2.63,1.45)\pspolygon[fillstyle=solid,linewidth=1pt,fillcolor=white](8.50,13.56)(10.50,5.81)(0.00,0.00)\pspolygon[fillstyle=solid,linewidth=1pt,fillcolor=black](6.38,10.17)(9.00,11.62)(8.50,13.56)\pspolygon[fillstyle=solid,linewidth=1pt,fillcolor=black](4.25,6.78)(6.88,8.23)(6.38,10.17)\pspolygon[fillstyle=solid,linewidth=1pt,fillcolor=black](6.88,8.23)(9.50,9.68)(9.00,11.62)\pspolygon[fillstyle=solid,linewidth=1pt,fillcolor=black](2.13,3.39)(4.75,4.84)(4.25,6.78)\pspolygon[fillstyle=solid,linewidth=1pt,fillcolor=black](4.75,4.84)(7.38,6.29)(6.88,8.23)\pspolygon[fillstyle=solid,linewidth=1pt,fillcolor=black](7.38,6.29)(10.00,7.75)(9.50,9.68)\pspolygon[fillstyle=solid,linewidth=1pt,fillcolor=black](0.00,0.00)(2.62,1.45)(2.13,3.39)\pspolygon[fillstyle=solid,linewidth=1pt,fillcolor=black](2.62,1.45)(5.25,2.90)(4.75,4.84)\pspolygon[fillstyle=solid,linewidth=1pt,fillcolor=black](5.25,2.90)(7.88,4.36)(7.38,6.29)\pspolygon[fillstyle=solid,linewidth=1pt,fillcolor=black](7.88,4.36)(10.50,5.81)(10.00,7.75)\pspolygon[fillstyle=solid,linewidth=1pt,fillcolor=white](12.00,0.00)(10.00,7.75)(7.88,4.36)\pspolygon[fillstyle=solid,linewidth=1pt,fillcolor=black](9.94,2.18)(11.00,3.87)(12.00,0.00)\pspolygon[fillstyle=solid,linewidth=1pt,fillcolor=black](7.88,4.36)(8.94,6.05)(9.94,2.18)\pspolygon[fillstyle=solid,linewidth=1pt,fillcolor=black](8.94,6.05)(10.00,7.75)(11.00,3.87)
\endpspicture}

\def\FigureImpossibleSix{
\newrgbcolor{lightblue}{0.8 0.8 1}
\psset{unit=0.3cm}
\pspicture(14,12)
\pspolygon[fillstyle=solid,linewidth=1pt,fillcolor=white](0.00,0.00)(7.00,12.12)(14.00,0.00)\pspolygon[fillstyle=solid,linewidth=1pt,fillcolor=lightblue](7.00,12.12)(6.29,7.18)(3.50,6.06)\pspolygon[fillstyle=solid,linewidth=1pt,fillcolor=lightblue](3.50,6.06)(2.79,1.11)(0.00,0.00)\pspolygon[fillstyle=solid,linewidth=1pt,fillcolor=lightblue](0.00,0.00)(4.64,1.86)(7.00,0.00)\pspolygon[fillstyle=solid,linewidth=1pt,fillcolor=lightblue](7.00,0.00)(11.64,1.86)(14.00,0.00)\pspolygon[fillstyle=solid,linewidth=1pt,fillcolor=lightblue](14.00,0.00)(10.07,3.09)(10.50,6.06)\pspolygon[fillstyle=solid,linewidth=1pt,fillcolor=lightblue](10.50,6.06)(6.57,9.16)(7.00,12.12)
\endpspicture}

\def\IsoscelesBetaFortyFourTiling{
\pspicture(2.5,1.2)
\psset{unit=0.3cm}
\newrgbcolor{lightblue}{0.8 0.8 1}
\newrgbcolor{pink}{1 0.8 0.8}
\newrgbcolor{lightgreen}{0.8 1 0.8}
\newrgbcolor{lightyellow}{1 1 0.8} 
\newrgbcolor{orange}{1 0.5 0}
\newrgbcolor{violet}{1 0 1}
\pspolygon[fillstyle=solid,linewidth=1pt,fillcolor=lightblue](13.25,2.90)(11.00,11.62)(8.00,0.00)
\psline(11.00,11.62)(8.00,0.00)
\psline(11.75,8.71)(9.75,0.97)
\psline(12.50,5.81)(11.50,1.94)
\psline(13.25,2.90)(13.25,2.90)
\psline(8.00,0.00)(13.25,2.90)
\psline(9.00,3.87)(12.50,5.81)
\psline(10.00,7.75)(11.75,8.71)
\psline(11.00,11.62)(11.00,11.62)
\psline(13.25,2.90)(11.00,11.62)
\psline(11.50,1.94)(10.00,7.75)
\psline(9.75,0.97)(9.00,3.87)
\psline(8.00,0.00)(8.00,0.00)
\pspolygon[fillstyle=solid,linewidth=1pt,fillcolor=lightgreen](14.00,0.00)(11.00,11.62)(22.00,0.00)
\psline(11.00,11.62)(22.00,0.00)
\psline(11.75,8.71)(20.00,0.00)
\psline(12.50,5.81)(18.00,0.00)
\psline(13.25,2.90)(16.00,0.00)
\psline(14.00,0.00)(14.00,0.00)
\psline(22.00,0.00)(14.00,0.00)
\psline(19.25,2.90)(13.25,2.90)
\psline(16.50,5.81)(12.50,5.81)
\psline(13.75,8.71)(11.75,8.71)
\psline(11.00,11.62)(11.00,11.62)
\psline(14.00,0.00)(11.00,11.62)
\psline(16.00,0.00)(13.75,8.71)
\psline(18.00,0.00)(16.50,5.81)
\psline(20.00,0.00)(19.25,2.90)
\psline(22.00,0.00)(22.00,0.00)
\pspolygon[fillstyle=solid,linewidth=1pt,fillcolor=pink](13.25,2.90)(8.00,0.00)(16.00,0.00)
\psline(8.00,0.00)(16.00,0.00)
\psline(10.62,1.45)(14.62,1.45)
\psline(13.25,2.90)(13.25,2.90)
\psline(16.00,0.00)(13.25,2.90)
\psline(12.00,0.00)(10.62,1.45)
\psline(8.00,0.00)(8.00,0.00)
\psline(13.25,2.90)(8.00,0.00)
\psline(14.62,1.45)(12.00,0.00)
\psline(16.00,0.00)(16.00,0.00)
\pspolygon[fillstyle=solid,linewidth=1pt,fillcolor=lightyellow](0.00,0.00)(8.00,0.00)(11.00,11.62)
\psline(8.00,0.00)(11.00,11.62)
\psline(6.00,0.00)(8.25,8.71)
\psline(4.00,0.00)(5.50,5.81)
\psline(2.00,0.00)(2.75,2.90)
\psline(0.00,0.00)(0.00,0.00)
\psline(11.00,11.62)(0.00,0.00)
\psline(10.25,8.71)(2.00,0.00)
\psline(9.50,5.81)(4.00,0.00)
\psline(8.75,2.90)(6.00,0.00)
\psline(8.00,0.00)(8.00,0.00)
\psline(0.00,0.00)(8.00,0.00)
\psline(2.75,2.90)(8.75,2.90)
\psline(5.50,5.81)(9.50,5.81)
\psline(8.25,8.71)(10.25,8.71)
\psline(11.00,11.62)(11.00,11.62)
\endpspicture}

\def\IsoscelesFortyEightTiling{
\pspicture(1,2)
\psset{unit=0.25cm}
\newrgbcolor{lightblue}{0.8 0.8 1}
\newrgbcolor{pink}{1 0.8 0.8}
\newrgbcolor{lightgreen}{0.8 1 0.8}
\newrgbcolor{lightyellow}{1 1 0.8} 
\newrgbcolor{orange}{1 0.5 0}
\newrgbcolor{violet}{1 0 1}
\pspolygon[fillstyle=solid,linewidth=1pt,fillcolor=green](0.00,0.00)(5.25,2.90)(8.00,0.00)
\psline(5.25,2.90)(8.00,0.00)
\psline(2.62,1.45)(4.00,0.00)
\psline(0.00,0.00)(0.00,0.00)
\psline(8.00,0.00)(0.00,0.00)
\psline(6.62,1.45)(2.62,1.45)
\psline(5.25,2.90)(5.25,2.90)
\psline(0.00,0.00)(5.25,2.90)
\psline(4.00,0.00)(6.62,1.45)
\psline(8.00,0.00)(8.00,0.00)
\pspolygon[fillstyle=solid,linewidth=1pt,fillcolor=lightblue](0.00,0.00)(5.25,2.90)(3.00,11.62)
\psline(5.25,2.90)(3.00,11.62)
\psline(3.50,1.94)(2.00,7.75)
\psline(1.75,0.97)(1.00,3.87)
\psline(0.00,0.00)(0.00,0.00)
\psline(3.00,11.62)(0.00,0.00)
\psline(3.75,8.71)(1.75,0.97)
\psline(4.50,5.81)(3.50,1.94)
\psline(5.25,2.90)(5.25,2.90)
\psline(0.00,0.00)(5.25,2.90)
\psline(1.00,3.87)(4.50,5.81)
\psline(2.00,7.75)(3.75,8.71)
\psline(3.00,11.62)(3.00,11.62)
\pspolygon[fillstyle=solid,linewidth=1pt,fillcolor=lightyellow](6.00,23.24)(3.00,11.62)(10.00,7.75)
\psline(3.00,11.62)(10.00,7.75)
\psline(3.75,14.52)(9.00,11.62)
\psline(4.50,17.43)(8.00,15.49)
\psline(5.25,20.33)(7.00,19.36)
\psline(6.00,23.24)(6.00,23.24)
\psline(10.00,7.75)(6.00,23.24)
\psline(8.25,8.71)(5.25,20.33)
\psline(6.50,9.68)(4.50,17.43)
\psline(4.75,10.65)(3.75,14.52)
\psline(3.00,11.62)(3.00,11.62)
\psline(6.00,23.24)(3.00,11.62)
\psline(7.00,19.36)(4.75,10.65)
\psline(8.00,15.49)(6.50,9.68)
\psline(9.00,11.62)(8.25,8.71)
\psline(10.00,7.75)(10.00,7.75)
\pspolygon[fillstyle=solid,linewidth=1pt,fillcolor=red](5.25,2.90)(7.25,2.90)(8.00,0.00)
\psline(7.25,2.90)(8.00,0.00)
\psline(5.25,2.90)(5.25,2.90)
\psline(8.00,0.00)(5.25,2.90)
\psline(7.25,2.90)(7.25,2.90)
\psline(5.25,2.90)(7.25,2.90)
\psline(8.00,0.00)(8.00,0.00)
\pspolygon[fillstyle=solid,linewidth=1pt,fillcolor=red](7.25,2.90)(11.25,2.90)(12.00,0.00)(8.00,0.00)
\psline(11.25,2.90)(12.00,0.00)
\psline(9.25,2.90)(10.00,0.00)
\psline(7.25,2.90)(8.00,0.00)
\psline(8.00,0.00)(12.00,0.00)
\psline(7.25,2.90)(11.25,2.90)
\psline(12.00,0.00)(9.25,2.90)
\psline(10.00,0.00)(7.25,2.90)
\pspolygon[fillstyle=solid,linewidth=1pt,fillcolor=pink](3.00,11.62)(10.00,7.75)(4.00,7.75)
\psline(10.00,7.75)(4.00,7.75)
\psline(6.50,9.68)(3.50,9.68)
\psline(3.00,11.62)(3.00,11.62)
\psline(4.00,7.75)(3.00,11.62)
\psline(7.00,7.75)(6.50,9.68)
\psline(10.00,7.75)(10.00,7.75)
\psline(3.00,11.62)(10.00,7.75)
\psline(3.50,9.68)(7.00,7.75)
\psline(4.00,7.75)(4.00,7.75)
\pspolygon[fillstyle=solid,linewidth=1pt,fillcolor=orange](4.00,7.75)(10.00,7.75)(10.50,5.81)(4.50,5.81)
\psline(10.00,7.75)(10.50,5.81)
\psline(7.00,7.75)(7.50,5.81)
\psline(4.00,7.75)(4.50,5.81)
\psline(4.50,5.81)(10.50,5.81)
\psline(4.00,7.75)(10.00,7.75)
\psline(10.50,5.81)(7.00,7.75)
\psline(7.50,5.81)(4.00,7.75)
\pspolygon[fillstyle=solid,linewidth=1pt,fillcolor=violet](11.25,2.90)(5.25,2.90)(4.50,5.81)(10.50,5.81)
\psline(5.25,2.90)(4.50,5.81)
\psline(7.25,2.90)(6.50,5.81)
\psline(9.25,2.90)(8.50,5.81)
\psline(11.25,2.90)(10.50,5.81)
\psline(10.50,5.81)(4.50,5.81)
\psline(11.25,2.90)(5.25,2.90)
\psline(4.50,5.81)(7.25,2.90)
\psline(6.50,5.81)(9.25,2.90)
\psline(8.50,5.81)(11.25,2.90)
\endpspicture}

\def\TriquadraticTwentyEight{
\pspicture(2,2.4)
\psset{unit=0.5cm}
\newrgbcolor{lightblue}{0.8 0.8 1}
\newrgbcolor{pink}{1 0.8 0.8}
\newrgbcolor{lightgreen}{0.8 1 0.8}
\pspolygon[fillstyle=solid,linewidth=1pt,fillcolor=lightblue](7.88,4.36)(0.00,0.00)(12.00,0.00)
\psline(0.00,0.00)(12.00,0.00)
\psline(2.62,1.45)(10.62,1.45)
\psline(5.25,2.90)(9.25,2.90)
\psline(7.88,4.36)(7.88,4.36)
\psline(12.00,0.00)(7.88,4.36)
\psline(8.00,0.00)(5.25,2.90)
\psline(4.00,0.00)(2.63,1.45)
\psline(0.00,0.00)(0.00,0.00)
\psline(7.88,4.36)(0.00,0.00)
\psline(9.25,2.90)(4.00,0.00)
\psline(10.62,1.45)(8.00,0.00)
\psline(12.00,0.00)(12.00,0.00)
\pspolygon[fillstyle=solid,linewidth=1pt,fillcolor=lightgreen](10.50,5.81)(0.00,0.00)(8.50,13.56)
\psline(0.00,0.00)(8.50,13.56)
\psline(2.62,1.45)(9.00,11.62)
\psline(5.25,2.90)(9.50,9.68)
\psline(7.88,4.36)(10.00,7.75)
\psline(10.50,5.81)(10.50,5.81)
\psline(8.50,13.56)(10.50,5.81)
\psline(6.38,10.17)(7.88,4.36)
\psline(4.25,6.78)(5.25,2.90)
\psline(2.13,3.39)(2.62,1.45)
\psline(0.00,0.00)(0.00,0.00)
\psline(10.50,5.81)(0.00,0.00)
\psline(10.00,7.75)(2.13,3.39)
\psline(9.50,9.68)(4.25,6.78)
\psline(9.00,11.62)(6.38,10.17)
\psline(8.50,13.56)(8.50,13.56)
\pspolygon[fillstyle=solid,linewidth=1pt,fillcolor=pink](7.88,4.36)(12.00,0.00)(10.00,7.75)
\psline(12.00,0.00)(10.00,7.75)
\psline(9.94,2.18)(8.94,6.05)
\psline(7.88,4.36)(7.88,4.36)
\psline(10.00,7.75)(7.88,4.36)
\psline(11.00,3.87)(9.94,2.18)
\psline(12.00,0.00)(12.00,0.00)
\psline(7.88,4.36)(12.00,0.00)
\psline(8.94,6.05)(11.00,3.87)
\psline(10.00,7.75)(10.00,7.75)
\endpspicture
}

\def\EquilateralFigure{
\psset{unit=1cm}
\pspicture(6,7)(-6,-3.5)
\psset{unit=0.03cm}
\psset{linewidth=0.1pt}
\newrgbcolor{lightblue}{0.8 0.8 1}
\newrgbcolor{pink}{1 0.8 0.8}
\newrgbcolor{lightgreen}{0.8 1 0.8}
\newrgbcolor{lightyellow}{1 1 0.8}
\pspolygon[fillstyle=solid,linewidth=0.1pt,fillcolor=lightblue](-22.50,-38.97)(75.00,0.00)(124.50,-38.97)
\psline(75.00,0.00)(124.50,-38.97)
\psline(70.36,-1.86)(117.50,-38.97)
\psline(65.71,-3.71)(110.50,-38.97)
\psline(61.07,-5.57)(103.50,-38.97)
\psline(56.43,-7.42)(96.50,-38.97)
\psline(51.79,-9.28)(89.50,-38.97)
\psline(47.14,-11.13)(82.50,-38.97)
\psline(42.50,-12.99)(75.50,-38.97)
\psline(37.86,-14.85)(68.50,-38.97)
\psline(33.21,-16.70)(61.50,-38.97)
\psline(28.57,-18.56)(54.50,-38.97)
\psline(23.93,-20.41)(47.50,-38.97)
\psline(19.29,-22.27)(40.50,-38.97)
\psline(14.64,-24.12)(33.50,-38.97)
\psline(10.00,-25.98)(26.50,-38.97)
\psline(5.36,-27.84)(19.50,-38.97)
\psline(0.71,-29.69)(12.50,-38.97)
\psline(-3.93,-31.55)(5.50,-38.97)
\psline(-8.57,-33.40)(-1.50,-38.97)
\psline(-13.21,-35.26)(-8.50,-38.97)
\psline(-17.86,-37.12)(-15.50,-38.97)
\psline(-22.50,-38.97)(-22.50,-38.97)
\psline(124.50,-38.97)(-22.50,-38.97)
\psline(122.14,-37.12)(-17.86,-37.12)
\psline(119.79,-35.26)(-13.21,-35.26)
\psline(117.43,-33.40)(-8.57,-33.40)
\psline(115.07,-31.55)(-3.93,-31.55)
\psline(112.71,-29.69)(0.71,-29.69)
\psline(110.36,-27.84)(5.36,-27.84)
\psline(108.00,-25.98)(10.00,-25.98)
\psline(105.64,-24.12)(14.64,-24.12)
\psline(103.29,-22.27)(19.29,-22.27)
\psline(100.93,-20.41)(23.93,-20.41)
\psline(98.57,-18.56)(28.57,-18.56)
\psline(96.21,-16.70)(33.21,-16.70)
\psline(93.86,-14.85)(37.86,-14.85)
\psline(91.50,-12.99)(42.50,-12.99)
\psline(89.14,-11.13)(47.14,-11.13)
\psline(86.79,-9.28)(51.79,-9.28)
\psline(84.43,-7.42)(56.43,-7.42)
\psline(82.07,-5.57)(61.07,-5.57)
\psline(79.71,-3.71)(65.71,-3.71)
\psline(77.36,-1.86)(70.36,-1.86)
\psline(75.00,0.00)(75.00,0.00)
\psline(-22.50,-38.97)(75.00,0.00)
\psline(-15.50,-38.97)(77.36,-1.86)
\psline(-8.50,-38.97)(79.71,-3.71)
\psline(-1.50,-38.97)(82.07,-5.57)
\psline(5.50,-38.97)(84.43,-7.42)
\psline(12.50,-38.97)(86.79,-9.28)
\psline(19.50,-38.97)(89.14,-11.13)
\psline(26.50,-38.97)(91.50,-12.99)
\psline(33.50,-38.97)(93.86,-14.85)
\psline(40.50,-38.97)(96.21,-16.70)
\psline(47.50,-38.97)(98.57,-18.56)
\psline(54.50,-38.97)(100.93,-20.41)
\psline(61.50,-38.97)(103.29,-22.27)
\psline(68.50,-38.97)(105.64,-24.12)
\psline(75.50,-38.97)(108.00,-25.98)
\psline(82.50,-38.97)(110.36,-27.84)
\psline(89.50,-38.97)(112.71,-29.69)
\psline(96.50,-38.97)(115.07,-31.55)
\psline(103.50,-38.97)(117.43,-33.40)
\psline(110.50,-38.97)(119.79,-35.26)
\psline(117.50,-38.97)(122.14,-37.12)
\psline(124.50,-38.97)(124.50,-38.97)
\pspolygon[fillstyle=solid,linewidth=0.1pt,fillcolor=lightblue](-45.00,-77.94)(52.50,-38.97)(102.00,-77.94)
\psline(52.50,-38.97)(102.00,-77.94)
\psline(47.86,-40.83)(95.00,-77.94)
\psline(43.21,-42.68)(88.00,-77.94)
\psline(38.57,-44.54)(81.00,-77.94)
\psline(33.93,-46.39)(74.00,-77.94)
\psline(29.29,-48.25)(67.00,-77.94)
\psline(24.64,-50.11)(60.00,-77.94)
\psline(20.00,-51.96)(53.00,-77.94)
\psline(15.36,-53.82)(46.00,-77.94)
\psline(10.71,-55.67)(39.00,-77.94)
\psline(6.07,-57.53)(32.00,-77.94)
\psline(1.43,-59.38)(25.00,-77.94)
\psline(-3.21,-61.24)(18.00,-77.94)
\psline(-7.86,-63.10)(11.00,-77.94)
\psline(-12.50,-64.95)(4.00,-77.94)
\psline(-17.14,-66.81)(-3.00,-77.94)
\psline(-21.79,-68.66)(-10.00,-77.94)
\psline(-26.43,-70.52)(-17.00,-77.94)
\psline(-31.07,-72.37)(-24.00,-77.94)
\psline(-35.71,-74.23)(-31.00,-77.94)
\psline(-40.36,-76.09)(-38.00,-77.94)
\psline(-45.00,-77.94)(-45.00,-77.94)
\psline(102.00,-77.94)(-45.00,-77.94)
\psline(99.64,-76.09)(-40.36,-76.09)
\psline(97.29,-74.23)(-35.71,-74.23)
\psline(94.93,-72.37)(-31.07,-72.37)
\psline(92.57,-70.52)(-26.43,-70.52)
\psline(90.21,-68.66)(-21.79,-68.66)
\psline(87.86,-66.81)(-17.14,-66.81)
\psline(85.50,-64.95)(-12.50,-64.95)
\psline(83.14,-63.10)(-7.86,-63.10)
\psline(80.79,-61.24)(-3.21,-61.24)
\psline(78.43,-59.38)(1.43,-59.38)
\psline(76.07,-57.53)(6.07,-57.53)
\psline(73.71,-55.67)(10.71,-55.67)
\psline(71.36,-53.82)(15.36,-53.82)
\psline(69.00,-51.96)(20.00,-51.96)
\psline(66.64,-50.11)(24.64,-50.11)
\psline(64.29,-48.25)(29.29,-48.25)
\psline(61.93,-46.39)(33.93,-46.39)
\psline(59.57,-44.54)(38.57,-44.54)
\psline(57.21,-42.68)(43.21,-42.68)
\psline(54.86,-40.83)(47.86,-40.83)
\psline(52.50,-38.97)(52.50,-38.97)
\psline(-45.00,-77.94)(52.50,-38.97)
\psline(-38.00,-77.94)(54.86,-40.83)
\psline(-31.00,-77.94)(57.21,-42.68)
\psline(-24.00,-77.94)(59.57,-44.54)
\psline(-17.00,-77.94)(61.93,-46.39)
\psline(-10.00,-77.94)(64.29,-48.25)
\psline(-3.00,-77.94)(66.64,-50.11)
\psline(4.00,-77.94)(69.00,-51.96)
\psline(11.00,-77.94)(71.36,-53.82)
\psline(18.00,-77.94)(73.71,-55.67)
\psline(25.00,-77.94)(76.07,-57.53)
\psline(32.00,-77.94)(78.43,-59.38)
\psline(39.00,-77.94)(80.79,-61.24)
\psline(46.00,-77.94)(83.14,-63.10)
\psline(53.00,-77.94)(85.50,-64.95)
\psline(60.00,-77.94)(87.86,-66.81)
\psline(67.00,-77.94)(90.21,-68.66)
\psline(74.00,-77.94)(92.57,-70.52)
\psline(81.00,-77.94)(94.93,-72.37)
\psline(88.00,-77.94)(97.29,-74.23)
\psline(95.00,-77.94)(99.64,-76.09)
\psline(102.00,-77.94)(102.00,-77.94)
\pspolygon[fillstyle=solid,linewidth=0.1pt,fillcolor=lightblue](-67.50,-116.91)(79.50,-116.91)(30.00,-77.94)
\psline(79.50,-116.91)(30.00,-77.94)
\psline(72.50,-116.91)(25.36,-79.80)
\psline(65.50,-116.91)(20.71,-81.65)
\psline(58.50,-116.91)(16.07,-83.51)
\psline(51.50,-116.91)(11.43,-85.37)
\psline(44.50,-116.91)(6.79,-87.22)
\psline(37.50,-116.91)(2.14,-89.08)
\psline(30.50,-116.91)(-2.50,-90.93)
\psline(23.50,-116.91)(-7.14,-92.79)
\psline(16.50,-116.91)(-11.79,-94.64)
\psline(9.50,-116.91)(-16.43,-96.50)
\psline(2.50,-116.91)(-21.07,-98.36)
\psline(-4.50,-116.91)(-25.71,-100.21)
\psline(-11.50,-116.91)(-30.36,-102.07)
\psline(-18.50,-116.91)(-35.00,-103.92)
\psline(-25.50,-116.91)(-39.64,-105.78)
\psline(-32.50,-116.91)(-44.29,-107.63)
\psline(-39.50,-116.91)(-48.93,-109.49)
\psline(-46.50,-116.91)(-53.57,-111.35)
\psline(-53.50,-116.91)(-58.21,-113.20)
\psline(-60.50,-116.91)(-62.86,-115.06)
\psline(-67.50,-116.91)(-67.50,-116.91)
\psline(30.00,-77.94)(-67.50,-116.91)
\psline(32.36,-79.80)(-60.50,-116.91)
\psline(34.71,-81.65)(-53.50,-116.91)
\psline(37.07,-83.51)(-46.50,-116.91)
\psline(39.43,-85.37)(-39.50,-116.91)
\psline(41.79,-87.22)(-32.50,-116.91)
\psline(44.14,-89.08)(-25.50,-116.91)
\psline(46.50,-90.93)(-18.50,-116.91)
\psline(48.86,-92.79)(-11.50,-116.91)
\psline(51.21,-94.64)(-4.50,-116.91)
\psline(53.57,-96.50)(2.50,-116.91)
\psline(55.93,-98.36)(9.50,-116.91)
\psline(58.29,-100.21)(16.50,-116.91)
\psline(60.64,-102.07)(23.50,-116.91)
\psline(63.00,-103.92)(30.50,-116.91)
\psline(65.36,-105.78)(37.50,-116.91)
\psline(67.71,-107.63)(44.50,-116.91)
\psline(70.07,-109.49)(51.50,-116.91)
\psline(72.43,-111.35)(58.50,-116.91)
\psline(74.79,-113.20)(65.50,-116.91)
\psline(77.14,-115.06)(72.50,-116.91)
\psline(79.50,-116.91)(79.50,-116.91)
\psline(-67.50,-116.91)(79.50,-116.91)
\psline(-62.86,-115.06)(77.14,-115.06)
\psline(-58.21,-113.20)(74.79,-113.20)
\psline(-53.57,-111.35)(72.43,-111.35)
\psline(-48.93,-109.49)(70.07,-109.49)
\psline(-44.29,-107.63)(67.71,-107.63)
\psline(-39.64,-105.78)(65.36,-105.78)
\psline(-35.00,-103.92)(63.00,-103.92)
\psline(-30.36,-102.07)(60.64,-102.07)
\psline(-25.71,-100.21)(58.29,-100.21)
\psline(-21.07,-98.36)(55.93,-98.36)
\psline(-16.43,-96.50)(53.57,-96.50)
\psline(-11.79,-94.64)(51.21,-94.64)
\psline(-7.14,-92.79)(48.86,-92.79)
\psline(-2.50,-90.93)(46.50,-90.93)
\psline(2.14,-89.08)(44.14,-89.08)
\psline(6.79,-87.22)(41.79,-87.22)
\psline(11.43,-85.37)(39.43,-85.37)
\psline(16.07,-83.51)(37.07,-83.51)
\psline(20.71,-81.65)(34.71,-81.65)
\psline(25.36,-79.80)(32.36,-79.80)
\psline(30.00,-77.94)(30.00,-77.94)
\pspolygon[fillstyle=solid,linewidth=0.1pt,fillcolor=lightyellow](0.00,0.00)(75.00,0.00)(-22.50,-38.97)
\psline(75.00,0.00)(-22.50,-38.97)
\psline(70.00,0.00)(-21.00,-36.37)
\psline(65.00,0.00)(-19.50,-33.77)
\psline(60.00,0.00)(-18.00,-31.18)
\psline(55.00,0.00)(-16.50,-28.58)
\psline(50.00,0.00)(-15.00,-25.98)
\psline(45.00,0.00)(-13.50,-23.38)
\psline(40.00,0.00)(-12.00,-20.78)
\psline(35.00,0.00)(-10.50,-18.19)
\psline(30.00,0.00)(-9.00,-15.59)
\psline(25.00,0.00)(-7.50,-12.99)
\psline(20.00,0.00)(-6.00,-10.39)
\psline(15.00,0.00)(-4.50,-7.79)
\psline(10.00,0.00)(-3.00,-5.20)
\psline(5.00,0.00)(-1.50,-2.60)
\psline(0.00,0.00)(0.00,0.00)
\psline(-22.50,-38.97)(0.00,0.00)
\psline(-16.00,-36.37)(5.00,0.00)
\psline(-9.50,-33.77)(10.00,0.00)
\psline(-3.00,-31.18)(15.00,0.00)
\psline(3.50,-28.58)(20.00,0.00)
\psline(10.00,-25.98)(25.00,0.00)
\psline(16.50,-23.38)(30.00,0.00)
\psline(23.00,-20.78)(35.00,0.00)
\psline(29.50,-18.19)(40.00,0.00)
\psline(36.00,-15.59)(45.00,0.00)
\psline(42.50,-12.99)(50.00,0.00)
\psline(49.00,-10.39)(55.00,0.00)
\psline(55.50,-7.79)(60.00,0.00)
\psline(62.00,-5.20)(65.00,0.00)
\psline(68.50,-2.60)(70.00,0.00)
\psline(75.00,0.00)(75.00,0.00)
\psline(0.00,0.00)(75.00,0.00)
\psline(-1.50,-2.60)(68.50,-2.60)
\psline(-3.00,-5.20)(62.00,-5.20)
\psline(-4.50,-7.79)(55.50,-7.79)
\psline(-6.00,-10.39)(49.00,-10.39)
\psline(-7.50,-12.99)(42.50,-12.99)
\psline(-9.00,-15.59)(36.00,-15.59)
\psline(-10.50,-18.19)(29.50,-18.19)
\psline(-12.00,-20.78)(23.00,-20.78)
\psline(-13.50,-23.38)(16.50,-23.38)
\psline(-15.00,-25.98)(10.00,-25.98)
\psline(-16.50,-28.58)(3.50,-28.58)
\psline(-18.00,-31.18)(-3.00,-31.18)
\psline(-19.50,-33.77)(-9.50,-33.77)
\psline(-21.00,-36.37)(-16.00,-36.37)
\psline(-22.50,-38.97)(-22.50,-38.97)
\pspolygon[fillstyle=solid,linewidth=0.1pt,fillcolor=lightyellow](-45.00,-77.94)(-22.50,-38.97)(52.50,-38.97)
\psline(-22.50,-38.97)(52.50,-38.97)
\psline(-24.00,-41.57)(46.00,-41.57)
\psline(-25.50,-44.17)(39.50,-44.17)
\psline(-27.00,-46.77)(33.00,-46.77)
\psline(-28.50,-49.36)(26.50,-49.36)
\psline(-30.00,-51.96)(20.00,-51.96)
\psline(-31.50,-54.56)(13.50,-54.56)
\psline(-33.00,-57.16)(7.00,-57.16)
\psline(-34.50,-59.76)(0.50,-59.76)
\psline(-36.00,-62.35)(-6.00,-62.35)
\psline(-37.50,-64.95)(-12.50,-64.95)
\psline(-39.00,-67.55)(-19.00,-67.55)
\psline(-40.50,-70.15)(-25.50,-70.15)
\psline(-42.00,-72.75)(-32.00,-72.75)
\psline(-43.50,-75.34)(-38.50,-75.34)
\psline(-45.00,-77.94)(-45.00,-77.94)
\psline(52.50,-38.97)(-45.00,-77.94)
\psline(47.50,-38.97)(-43.50,-75.34)
\psline(42.50,-38.97)(-42.00,-72.75)
\psline(37.50,-38.97)(-40.50,-70.15)
\psline(32.50,-38.97)(-39.00,-67.55)
\psline(27.50,-38.97)(-37.50,-64.95)
\psline(22.50,-38.97)(-36.00,-62.35)
\psline(17.50,-38.97)(-34.50,-59.76)
\psline(12.50,-38.97)(-33.00,-57.16)
\psline(7.50,-38.97)(-31.50,-54.56)
\psline(2.50,-38.97)(-30.00,-51.96)
\psline(-2.50,-38.97)(-28.50,-49.36)
\psline(-7.50,-38.97)(-27.00,-46.77)
\psline(-12.50,-38.97)(-25.50,-44.17)
\psline(-17.50,-38.97)(-24.00,-41.57)
\psline(-22.50,-38.97)(-22.50,-38.97)
\psline(-45.00,-77.94)(-22.50,-38.97)
\psline(-38.50,-75.34)(-17.50,-38.97)
\psline(-32.00,-72.75)(-12.50,-38.97)
\psline(-25.50,-70.15)(-7.50,-38.97)
\psline(-19.00,-67.55)(-2.50,-38.97)
\psline(-12.50,-64.95)(2.50,-38.97)
\psline(-6.00,-62.35)(7.50,-38.97)
\psline(0.50,-59.76)(12.50,-38.97)
\psline(7.00,-57.16)(17.50,-38.97)
\psline(13.50,-54.56)(22.50,-38.97)
\psline(20.00,-51.96)(27.50,-38.97)
\psline(26.50,-49.36)(32.50,-38.97)
\psline(33.00,-46.77)(37.50,-38.97)
\psline(39.50,-44.17)(42.50,-38.97)
\psline(46.00,-41.57)(47.50,-38.97)
\psline(52.50,-38.97)(52.50,-38.97)
\pspolygon[fillstyle=solid,linewidth=0.1pt,fillcolor=lightyellow](-67.50,-116.91)(-45.00,-77.94)(30.00,-77.94)
\psline(-45.00,-77.94)(30.00,-77.94)
\psline(-46.50,-80.54)(23.50,-80.54)
\psline(-48.00,-83.14)(17.00,-83.14)
\psline(-49.50,-85.74)(10.50,-85.74)
\psline(-51.00,-88.33)(4.00,-88.33)
\psline(-52.50,-90.93)(-2.50,-90.93)
\psline(-54.00,-93.53)(-9.00,-93.53)
\psline(-55.50,-96.13)(-15.50,-96.13)
\psline(-57.00,-98.73)(-22.00,-98.73)
\psline(-58.50,-101.32)(-28.50,-101.32)
\psline(-60.00,-103.92)(-35.00,-103.92)
\psline(-61.50,-106.52)(-41.50,-106.52)
\psline(-63.00,-109.12)(-48.00,-109.12)
\psline(-64.50,-111.72)(-54.50,-111.72)
\psline(-66.00,-114.32)(-61.00,-114.32)
\psline(-67.50,-116.91)(-67.50,-116.91)
\psline(30.00,-77.94)(-67.50,-116.91)
\psline(25.00,-77.94)(-66.00,-114.32)
\psline(20.00,-77.94)(-64.50,-111.72)
\psline(15.00,-77.94)(-63.00,-109.12)
\psline(10.00,-77.94)(-61.50,-106.52)
\psline(5.00,-77.94)(-60.00,-103.92)
\psline(0.00,-77.94)(-58.50,-101.32)
\psline(-5.00,-77.94)(-57.00,-98.73)
\psline(-10.00,-77.94)(-55.50,-96.13)
\psline(-15.00,-77.94)(-54.00,-93.53)
\psline(-20.00,-77.94)(-52.50,-90.93)
\psline(-25.00,-77.94)(-51.00,-88.33)
\psline(-30.00,-77.94)(-49.50,-85.74)
\psline(-35.00,-77.94)(-48.00,-83.14)
\psline(-40.00,-77.94)(-46.50,-80.54)
\psline(-45.00,-77.94)(-45.00,-77.94)
\psline(-67.50,-116.91)(-45.00,-77.94)
\psline(-61.00,-114.32)(-40.00,-77.94)
\psline(-54.50,-111.72)(-35.00,-77.94)
\psline(-48.00,-109.12)(-30.00,-77.94)
\psline(-41.50,-106.52)(-25.00,-77.94)
\psline(-35.00,-103.92)(-20.00,-77.94)
\psline(-28.50,-101.32)(-15.00,-77.94)
\psline(-22.00,-98.73)(-10.00,-77.94)
\psline(-15.50,-96.13)(-5.00,-77.94)
\psline(-9.00,-93.53)(0.00,-77.94)
\psline(-2.50,-90.93)(5.00,-77.94)
\psline(4.00,-88.33)(10.00,-77.94)
\psline(10.50,-85.74)(15.00,-77.94)
\psline(17.00,-83.14)(20.00,-77.94)
\psline(23.50,-80.54)(25.00,-77.94)
\psline(30.00,-77.94)(30.00,-77.94)
\pspolygon[fillstyle=solid,linewidth=0.1pt,fillcolor=green](75.00,0.00)(102.00,0.00)(124.50,-38.97)
\psline(102.00,0.00)(124.50,-38.97)
\psline(99.00,0.00)(119.00,-34.64)
\psline(96.00,0.00)(113.50,-30.31)
\psline(93.00,0.00)(108.00,-25.98)
\psline(90.00,0.00)(102.50,-21.65)
\psline(87.00,0.00)(97.00,-17.32)
\psline(84.00,0.00)(91.50,-12.99)
\psline(81.00,0.00)(86.00,-8.66)
\psline(78.00,0.00)(80.50,-4.33)
\psline(75.00,0.00)(75.00,0.00)
\psline(124.50,-38.97)(75.00,0.00)
\psline(122.00,-34.64)(78.00,0.00)
\psline(119.50,-30.31)(81.00,0.00)
\psline(117.00,-25.98)(84.00,0.00)
\psline(114.50,-21.65)(87.00,0.00)
\psline(112.00,-17.32)(90.00,0.00)
\psline(109.50,-12.99)(93.00,0.00)
\psline(107.00,-8.66)(96.00,0.00)
\psline(104.50,-4.33)(99.00,0.00)
\psline(102.00,0.00)(102.00,0.00)
\psline(75.00,0.00)(102.00,0.00)
\psline(80.50,-4.33)(104.50,-4.33)
\psline(86.00,-8.66)(107.00,-8.66)
\psline(91.50,-12.99)(109.50,-12.99)
\psline(97.00,-17.32)(112.00,-17.32)
\psline(102.50,-21.65)(114.50,-21.65)
\psline(108.00,-25.98)(117.00,-25.98)
\psline(113.50,-30.31)(119.50,-30.31)
\psline(119.00,-34.64)(122.00,-34.64)
\psline(124.50,-38.97)(124.50,-38.97)
\pspolygon[fillstyle=solid,linewidth=0.1pt,fillcolor=green](52.50,-38.97)(79.50,-38.97)(102.00,-77.94)
\psline(79.50,-38.97)(102.00,-77.94)
\psline(76.50,-38.97)(96.50,-73.61)
\psline(73.50,-38.97)(91.00,-69.28)
\psline(70.50,-38.97)(85.50,-64.95)
\psline(67.50,-38.97)(80.00,-60.62)
\psline(64.50,-38.97)(74.50,-56.29)
\psline(61.50,-38.97)(69.00,-51.96)
\psline(58.50,-38.97)(63.50,-47.63)
\psline(55.50,-38.97)(58.00,-43.30)
\psline(52.50,-38.97)(52.50,-38.97)
\psline(102.00,-77.94)(52.50,-38.97)
\psline(99.50,-73.61)(55.50,-38.97)
\psline(97.00,-69.28)(58.50,-38.97)
\psline(94.50,-64.95)(61.50,-38.97)
\psline(92.00,-60.62)(64.50,-38.97)
\psline(89.50,-56.29)(67.50,-38.97)
\psline(87.00,-51.96)(70.50,-38.97)
\psline(84.50,-47.63)(73.50,-38.97)
\psline(82.00,-43.30)(76.50,-38.97)
\psline(79.50,-38.97)(79.50,-38.97)
\psline(52.50,-38.97)(79.50,-38.97)
\psline(58.00,-43.30)(82.00,-43.30)
\psline(63.50,-47.63)(84.50,-47.63)
\psline(69.00,-51.96)(87.00,-51.96)
\psline(74.50,-56.29)(89.50,-56.29)
\psline(80.00,-60.62)(92.00,-60.62)
\psline(85.50,-64.95)(94.50,-64.95)
\psline(91.00,-69.28)(97.00,-69.28)
\psline(96.50,-73.61)(99.50,-73.61)
\psline(102.00,-77.94)(102.00,-77.94)
\pspolygon[fillstyle=solid,linewidth=0.1pt,fillcolor=green](30.00,-77.94)(57.00,-77.94)(79.50,-116.91)
\psline(57.00,-77.94)(79.50,-116.91)
\psline(54.00,-77.94)(74.00,-112.58)
\psline(51.00,-77.94)(68.50,-108.25)
\psline(48.00,-77.94)(63.00,-103.92)
\psline(45.00,-77.94)(57.50,-99.59)
\psline(42.00,-77.94)(52.00,-95.26)
\psline(39.00,-77.94)(46.50,-90.93)
\psline(36.00,-77.94)(41.00,-86.60)
\psline(33.00,-77.94)(35.50,-82.27)
\psline(30.00,-77.94)(30.00,-77.94)
\psline(79.50,-116.91)(30.00,-77.94)
\psline(77.00,-112.58)(33.00,-77.94)
\psline(74.50,-108.25)(36.00,-77.94)
\psline(72.00,-103.92)(39.00,-77.94)
\psline(69.50,-99.59)(42.00,-77.94)
\psline(67.00,-95.26)(45.00,-77.94)
\psline(64.50,-90.93)(48.00,-77.94)
\psline(62.00,-86.60)(51.00,-77.94)
\psline(59.50,-82.27)(54.00,-77.94)
\psline(57.00,-77.94)(57.00,-77.94)
\psline(30.00,-77.94)(57.00,-77.94)
\psline(35.50,-82.27)(59.50,-82.27)
\psline(41.00,-86.60)(62.00,-86.60)
\psline(46.50,-90.93)(64.50,-90.93)
\psline(52.00,-95.26)(67.00,-95.26)
\psline(57.50,-99.59)(69.50,-99.59)
\psline(63.00,-103.92)(72.00,-103.92)
\psline(68.50,-108.25)(74.50,-108.25)
\psline(74.00,-112.58)(77.00,-112.58)
\psline(79.50,-116.91)(79.50,-116.91)
\pspolygon[fillstyle=solid,linewidth=0.1pt,fillcolor=red](79.50,-38.97)(157.50,-38.97)(180.00,-77.94)(102.00,-77.94)
\psline(157.50,-38.97)(180.00,-77.94)
\psline(154.50,-38.97)(177.00,-77.94)
\psline(151.50,-38.97)(174.00,-77.94)
\psline(148.50,-38.97)(171.00,-77.94)
\psline(145.50,-38.97)(168.00,-77.94)
\psline(142.50,-38.97)(165.00,-77.94)
\psline(139.50,-38.97)(162.00,-77.94)
\psline(136.50,-38.97)(159.00,-77.94)
\psline(133.50,-38.97)(156.00,-77.94)
\psline(130.50,-38.97)(153.00,-77.94)
\psline(127.50,-38.97)(150.00,-77.94)
\psline(124.50,-38.97)(147.00,-77.94)
\psline(121.50,-38.97)(144.00,-77.94)
\psline(118.50,-38.97)(141.00,-77.94)
\psline(115.50,-38.97)(138.00,-77.94)
\psline(112.50,-38.97)(135.00,-77.94)
\psline(109.50,-38.97)(132.00,-77.94)
\psline(106.50,-38.97)(129.00,-77.94)
\psline(103.50,-38.97)(126.00,-77.94)
\psline(100.50,-38.97)(123.00,-77.94)
\psline(97.50,-38.97)(120.00,-77.94)
\psline(94.50,-38.97)(117.00,-77.94)
\psline(91.50,-38.97)(114.00,-77.94)
\psline(88.50,-38.97)(111.00,-77.94)
\psline(85.50,-38.97)(108.00,-77.94)
\psline(82.50,-38.97)(105.00,-77.94)
\psline(79.50,-38.97)(102.00,-77.94)
\psline(102.00,-77.94)(180.00,-77.94)
\psline(99.50,-73.61)(177.50,-73.61)
\psline(97.00,-69.28)(175.00,-69.28)
\psline(94.50,-64.95)(172.50,-64.95)
\psline(92.00,-60.62)(170.00,-60.62)
\psline(89.50,-56.29)(167.50,-56.29)
\psline(87.00,-51.96)(165.00,-51.96)
\psline(84.50,-47.63)(162.50,-47.63)
\psline(82.00,-43.30)(160.00,-43.30)
\psline(79.50,-38.97)(157.50,-38.97)
\psline(180.00,-77.94)(174.50,-73.61)
\psline(177.50,-73.61)(172.00,-69.28)
\psline(175.00,-69.28)(169.50,-64.95)
\psline(172.50,-64.95)(167.00,-60.62)
\psline(170.00,-60.62)(164.50,-56.29)
\psline(167.50,-56.29)(162.00,-51.96)
\psline(165.00,-51.96)(159.50,-47.63)
\psline(162.50,-47.63)(157.00,-43.30)
\psline(160.00,-43.30)(154.50,-38.97)
\psline(177.00,-77.94)(171.50,-73.61)
\psline(174.50,-73.61)(169.00,-69.28)
\psline(172.00,-69.28)(166.50,-64.95)
\psline(169.50,-64.95)(164.00,-60.62)
\psline(167.00,-60.62)(161.50,-56.29)
\psline(164.50,-56.29)(159.00,-51.96)
\psline(162.00,-51.96)(156.50,-47.63)
\psline(159.50,-47.63)(154.00,-43.30)
\psline(157.00,-43.30)(151.50,-38.97)
\psline(174.00,-77.94)(168.50,-73.61)
\psline(171.50,-73.61)(166.00,-69.28)
\psline(169.00,-69.28)(163.50,-64.95)
\psline(166.50,-64.95)(161.00,-60.62)
\psline(164.00,-60.62)(158.50,-56.29)
\psline(161.50,-56.29)(156.00,-51.96)
\psline(159.00,-51.96)(153.50,-47.63)
\psline(156.50,-47.63)(151.00,-43.30)
\psline(154.00,-43.30)(148.50,-38.97)
\psline(171.00,-77.94)(165.50,-73.61)
\psline(168.50,-73.61)(163.00,-69.28)
\psline(166.00,-69.28)(160.50,-64.95)
\psline(163.50,-64.95)(158.00,-60.62)
\psline(161.00,-60.62)(155.50,-56.29)
\psline(158.50,-56.29)(153.00,-51.96)
\psline(156.00,-51.96)(150.50,-47.63)
\psline(153.50,-47.63)(148.00,-43.30)
\psline(151.00,-43.30)(145.50,-38.97)
\psline(168.00,-77.94)(162.50,-73.61)
\psline(165.50,-73.61)(160.00,-69.28)
\psline(163.00,-69.28)(157.50,-64.95)
\psline(160.50,-64.95)(155.00,-60.62)
\psline(158.00,-60.62)(152.50,-56.29)
\psline(155.50,-56.29)(150.00,-51.96)
\psline(153.00,-51.96)(147.50,-47.63)
\psline(150.50,-47.63)(145.00,-43.30)
\psline(148.00,-43.30)(142.50,-38.97)
\psline(165.00,-77.94)(159.50,-73.61)
\psline(162.50,-73.61)(157.00,-69.28)
\psline(160.00,-69.28)(154.50,-64.95)
\psline(157.50,-64.95)(152.00,-60.62)
\psline(155.00,-60.62)(149.50,-56.29)
\psline(152.50,-56.29)(147.00,-51.96)
\psline(150.00,-51.96)(144.50,-47.63)
\psline(147.50,-47.63)(142.00,-43.30)
\psline(145.00,-43.30)(139.50,-38.97)
\psline(162.00,-77.94)(156.50,-73.61)
\psline(159.50,-73.61)(154.00,-69.28)
\psline(157.00,-69.28)(151.50,-64.95)
\psline(154.50,-64.95)(149.00,-60.62)
\psline(152.00,-60.62)(146.50,-56.29)
\psline(149.50,-56.29)(144.00,-51.96)
\psline(147.00,-51.96)(141.50,-47.63)
\psline(144.50,-47.63)(139.00,-43.30)
\psline(142.00,-43.30)(136.50,-38.97)
\psline(159.00,-77.94)(153.50,-73.61)
\psline(156.50,-73.61)(151.00,-69.28)
\psline(154.00,-69.28)(148.50,-64.95)
\psline(151.50,-64.95)(146.00,-60.62)
\psline(149.00,-60.62)(143.50,-56.29)
\psline(146.50,-56.29)(141.00,-51.96)
\psline(144.00,-51.96)(138.50,-47.63)
\psline(141.50,-47.63)(136.00,-43.30)
\psline(139.00,-43.30)(133.50,-38.97)
\psline(156.00,-77.94)(150.50,-73.61)
\psline(153.50,-73.61)(148.00,-69.28)
\psline(151.00,-69.28)(145.50,-64.95)
\psline(148.50,-64.95)(143.00,-60.62)
\psline(146.00,-60.62)(140.50,-56.29)
\psline(143.50,-56.29)(138.00,-51.96)
\psline(141.00,-51.96)(135.50,-47.63)
\psline(138.50,-47.63)(133.00,-43.30)
\psline(136.00,-43.30)(130.50,-38.97)
\psline(153.00,-77.94)(147.50,-73.61)
\psline(150.50,-73.61)(145.00,-69.28)
\psline(148.00,-69.28)(142.50,-64.95)
\psline(145.50,-64.95)(140.00,-60.62)
\psline(143.00,-60.62)(137.50,-56.29)
\psline(140.50,-56.29)(135.00,-51.96)
\psline(138.00,-51.96)(132.50,-47.63)
\psline(135.50,-47.63)(130.00,-43.30)
\psline(133.00,-43.30)(127.50,-38.97)
\psline(150.00,-77.94)(144.50,-73.61)
\psline(147.50,-73.61)(142.00,-69.28)
\psline(145.00,-69.28)(139.50,-64.95)
\psline(142.50,-64.95)(137.00,-60.62)
\psline(140.00,-60.62)(134.50,-56.29)
\psline(137.50,-56.29)(132.00,-51.96)
\psline(135.00,-51.96)(129.50,-47.63)
\psline(132.50,-47.63)(127.00,-43.30)
\psline(130.00,-43.30)(124.50,-38.97)
\psline(147.00,-77.94)(141.50,-73.61)
\psline(144.50,-73.61)(139.00,-69.28)
\psline(142.00,-69.28)(136.50,-64.95)
\psline(139.50,-64.95)(134.00,-60.62)
\psline(137.00,-60.62)(131.50,-56.29)
\psline(134.50,-56.29)(129.00,-51.96)
\psline(132.00,-51.96)(126.50,-47.63)
\psline(129.50,-47.63)(124.00,-43.30)
\psline(127.00,-43.30)(121.50,-38.97)
\psline(144.00,-77.94)(138.50,-73.61)
\psline(141.50,-73.61)(136.00,-69.28)
\psline(139.00,-69.28)(133.50,-64.95)
\psline(136.50,-64.95)(131.00,-60.62)
\psline(134.00,-60.62)(128.50,-56.29)
\psline(131.50,-56.29)(126.00,-51.96)
\psline(129.00,-51.96)(123.50,-47.63)
\psline(126.50,-47.63)(121.00,-43.30)
\psline(124.00,-43.30)(118.50,-38.97)
\psline(141.00,-77.94)(135.50,-73.61)
\psline(138.50,-73.61)(133.00,-69.28)
\psline(136.00,-69.28)(130.50,-64.95)
\psline(133.50,-64.95)(128.00,-60.62)
\psline(131.00,-60.62)(125.50,-56.29)
\psline(128.50,-56.29)(123.00,-51.96)
\psline(126.00,-51.96)(120.50,-47.63)
\psline(123.50,-47.63)(118.00,-43.30)
\psline(121.00,-43.30)(115.50,-38.97)
\psline(138.00,-77.94)(132.50,-73.61)
\psline(135.50,-73.61)(130.00,-69.28)
\psline(133.00,-69.28)(127.50,-64.95)
\psline(130.50,-64.95)(125.00,-60.62)
\psline(128.00,-60.62)(122.50,-56.29)
\psline(125.50,-56.29)(120.00,-51.96)
\psline(123.00,-51.96)(117.50,-47.63)
\psline(120.50,-47.63)(115.00,-43.30)
\psline(118.00,-43.30)(112.50,-38.97)
\psline(135.00,-77.94)(129.50,-73.61)
\psline(132.50,-73.61)(127.00,-69.28)
\psline(130.00,-69.28)(124.50,-64.95)
\psline(127.50,-64.95)(122.00,-60.62)
\psline(125.00,-60.62)(119.50,-56.29)
\psline(122.50,-56.29)(117.00,-51.96)
\psline(120.00,-51.96)(114.50,-47.63)
\psline(117.50,-47.63)(112.00,-43.30)
\psline(115.00,-43.30)(109.50,-38.97)
\psline(132.00,-77.94)(126.50,-73.61)
\psline(129.50,-73.61)(124.00,-69.28)
\psline(127.00,-69.28)(121.50,-64.95)
\psline(124.50,-64.95)(119.00,-60.62)
\psline(122.00,-60.62)(116.50,-56.29)
\psline(119.50,-56.29)(114.00,-51.96)
\psline(117.00,-51.96)(111.50,-47.63)
\psline(114.50,-47.63)(109.00,-43.30)
\psline(112.00,-43.30)(106.50,-38.97)
\psline(129.00,-77.94)(123.50,-73.61)
\psline(126.50,-73.61)(121.00,-69.28)
\psline(124.00,-69.28)(118.50,-64.95)
\psline(121.50,-64.95)(116.00,-60.62)
\psline(119.00,-60.62)(113.50,-56.29)
\psline(116.50,-56.29)(111.00,-51.96)
\psline(114.00,-51.96)(108.50,-47.63)
\psline(111.50,-47.63)(106.00,-43.30)
\psline(109.00,-43.30)(103.50,-38.97)
\psline(126.00,-77.94)(120.50,-73.61)
\psline(123.50,-73.61)(118.00,-69.28)
\psline(121.00,-69.28)(115.50,-64.95)
\psline(118.50,-64.95)(113.00,-60.62)
\psline(116.00,-60.62)(110.50,-56.29)
\psline(113.50,-56.29)(108.00,-51.96)
\psline(111.00,-51.96)(105.50,-47.63)
\psline(108.50,-47.63)(103.00,-43.30)
\psline(106.00,-43.30)(100.50,-38.97)
\psline(123.00,-77.94)(117.50,-73.61)
\psline(120.50,-73.61)(115.00,-69.28)
\psline(118.00,-69.28)(112.50,-64.95)
\psline(115.50,-64.95)(110.00,-60.62)
\psline(113.00,-60.62)(107.50,-56.29)
\psline(110.50,-56.29)(105.00,-51.96)
\psline(108.00,-51.96)(102.50,-47.63)
\psline(105.50,-47.63)(100.00,-43.30)
\psline(103.00,-43.30)(97.50,-38.97)
\psline(120.00,-77.94)(114.50,-73.61)
\psline(117.50,-73.61)(112.00,-69.28)
\psline(115.00,-69.28)(109.50,-64.95)
\psline(112.50,-64.95)(107.00,-60.62)
\psline(110.00,-60.62)(104.50,-56.29)
\psline(107.50,-56.29)(102.00,-51.96)
\psline(105.00,-51.96)(99.50,-47.63)
\psline(102.50,-47.63)(97.00,-43.30)
\psline(100.00,-43.30)(94.50,-38.97)
\psline(117.00,-77.94)(111.50,-73.61)
\psline(114.50,-73.61)(109.00,-69.28)
\psline(112.00,-69.28)(106.50,-64.95)
\psline(109.50,-64.95)(104.00,-60.62)
\psline(107.00,-60.62)(101.50,-56.29)
\psline(104.50,-56.29)(99.00,-51.96)
\psline(102.00,-51.96)(96.50,-47.63)
\psline(99.50,-47.63)(94.00,-43.30)
\psline(97.00,-43.30)(91.50,-38.97)
\psline(114.00,-77.94)(108.50,-73.61)
\psline(111.50,-73.61)(106.00,-69.28)
\psline(109.00,-69.28)(103.50,-64.95)
\psline(106.50,-64.95)(101.00,-60.62)
\psline(104.00,-60.62)(98.50,-56.29)
\psline(101.50,-56.29)(96.00,-51.96)
\psline(99.00,-51.96)(93.50,-47.63)
\psline(96.50,-47.63)(91.00,-43.30)
\psline(94.00,-43.30)(88.50,-38.97)
\psline(111.00,-77.94)(105.50,-73.61)
\psline(108.50,-73.61)(103.00,-69.28)
\psline(106.00,-69.28)(100.50,-64.95)
\psline(103.50,-64.95)(98.00,-60.62)
\psline(101.00,-60.62)(95.50,-56.29)
\psline(98.50,-56.29)(93.00,-51.96)
\psline(96.00,-51.96)(90.50,-47.63)
\psline(93.50,-47.63)(88.00,-43.30)
\psline(91.00,-43.30)(85.50,-38.97)
\psline(108.00,-77.94)(102.50,-73.61)
\psline(105.50,-73.61)(100.00,-69.28)
\psline(103.00,-69.28)(97.50,-64.95)
\psline(100.50,-64.95)(95.00,-60.62)
\psline(98.00,-60.62)(92.50,-56.29)
\psline(95.50,-56.29)(90.00,-51.96)
\psline(93.00,-51.96)(87.50,-47.63)
\psline(90.50,-47.63)(85.00,-43.30)
\psline(88.00,-43.30)(82.50,-38.97)
\psline(105.00,-77.94)(99.50,-73.61)
\psline(102.50,-73.61)(97.00,-69.28)
\psline(100.00,-69.28)(94.50,-64.95)
\psline(97.50,-64.95)(92.00,-60.62)
\psline(95.00,-60.62)(89.50,-56.29)
\psline(92.50,-56.29)(87.00,-51.96)
\psline(90.00,-51.96)(84.50,-47.63)
\psline(87.50,-47.63)(82.00,-43.30)
\psline(85.00,-43.30)(79.50,-38.97)
\pspolygon[fillstyle=solid,linewidth=0.1pt,fillcolor=red](102.00,0.00)(135.00,0.00)(157.50,-38.97)(124.50,-38.97)
\psline(135.00,0.00)(157.50,-38.97)
\psline(132.00,0.00)(154.50,-38.97)
\psline(129.00,0.00)(151.50,-38.97)
\psline(126.00,0.00)(148.50,-38.97)
\psline(123.00,0.00)(145.50,-38.97)
\psline(120.00,0.00)(142.50,-38.97)
\psline(117.00,0.00)(139.50,-38.97)
\psline(114.00,0.00)(136.50,-38.97)
\psline(111.00,0.00)(133.50,-38.97)
\psline(108.00,0.00)(130.50,-38.97)
\psline(105.00,0.00)(127.50,-38.97)
\psline(102.00,0.00)(124.50,-38.97)
\psline(124.50,-38.97)(157.50,-38.97)
\psline(122.00,-34.64)(155.00,-34.64)
\psline(119.50,-30.31)(152.50,-30.31)
\psline(117.00,-25.98)(150.00,-25.98)
\psline(114.50,-21.65)(147.50,-21.65)
\psline(112.00,-17.32)(145.00,-17.32)
\psline(109.50,-12.99)(142.50,-12.99)
\psline(107.00,-8.66)(140.00,-8.66)
\psline(104.50,-4.33)(137.50,-4.33)
\psline(102.00,0.00)(135.00,0.00)
\psline(157.50,-38.97)(152.00,-34.64)
\psline(155.00,-34.64)(149.50,-30.31)
\psline(152.50,-30.31)(147.00,-25.98)
\psline(150.00,-25.98)(144.50,-21.65)
\psline(147.50,-21.65)(142.00,-17.32)
\psline(145.00,-17.32)(139.50,-12.99)
\psline(142.50,-12.99)(137.00,-8.66)
\psline(140.00,-8.66)(134.50,-4.33)
\psline(137.50,-4.33)(132.00,0.00)
\psline(154.50,-38.97)(149.00,-34.64)
\psline(152.00,-34.64)(146.50,-30.31)
\psline(149.50,-30.31)(144.00,-25.98)
\psline(147.00,-25.98)(141.50,-21.65)
\psline(144.50,-21.65)(139.00,-17.32)
\psline(142.00,-17.32)(136.50,-12.99)
\psline(139.50,-12.99)(134.00,-8.66)
\psline(137.00,-8.66)(131.50,-4.33)
\psline(134.50,-4.33)(129.00,0.00)
\psline(151.50,-38.97)(146.00,-34.64)
\psline(149.00,-34.64)(143.50,-30.31)
\psline(146.50,-30.31)(141.00,-25.98)
\psline(144.00,-25.98)(138.50,-21.65)
\psline(141.50,-21.65)(136.00,-17.32)
\psline(139.00,-17.32)(133.50,-12.99)
\psline(136.50,-12.99)(131.00,-8.66)
\psline(134.00,-8.66)(128.50,-4.33)
\psline(131.50,-4.33)(126.00,0.00)
\psline(148.50,-38.97)(143.00,-34.64)
\psline(146.00,-34.64)(140.50,-30.31)
\psline(143.50,-30.31)(138.00,-25.98)
\psline(141.00,-25.98)(135.50,-21.65)
\psline(138.50,-21.65)(133.00,-17.32)
\psline(136.00,-17.32)(130.50,-12.99)
\psline(133.50,-12.99)(128.00,-8.66)
\psline(131.00,-8.66)(125.50,-4.33)
\psline(128.50,-4.33)(123.00,0.00)
\psline(145.50,-38.97)(140.00,-34.64)
\psline(143.00,-34.64)(137.50,-30.31)
\psline(140.50,-30.31)(135.00,-25.98)
\psline(138.00,-25.98)(132.50,-21.65)
\psline(135.50,-21.65)(130.00,-17.32)
\psline(133.00,-17.32)(127.50,-12.99)
\psline(130.50,-12.99)(125.00,-8.66)
\psline(128.00,-8.66)(122.50,-4.33)
\psline(125.50,-4.33)(120.00,0.00)
\psline(142.50,-38.97)(137.00,-34.64)
\psline(140.00,-34.64)(134.50,-30.31)
\psline(137.50,-30.31)(132.00,-25.98)
\psline(135.00,-25.98)(129.50,-21.65)
\psline(132.50,-21.65)(127.00,-17.32)
\psline(130.00,-17.32)(124.50,-12.99)
\psline(127.50,-12.99)(122.00,-8.66)
\psline(125.00,-8.66)(119.50,-4.33)
\psline(122.50,-4.33)(117.00,0.00)
\psline(139.50,-38.97)(134.00,-34.64)
\psline(137.00,-34.64)(131.50,-30.31)
\psline(134.50,-30.31)(129.00,-25.98)
\psline(132.00,-25.98)(126.50,-21.65)
\psline(129.50,-21.65)(124.00,-17.32)
\psline(127.00,-17.32)(121.50,-12.99)
\psline(124.50,-12.99)(119.00,-8.66)
\psline(122.00,-8.66)(116.50,-4.33)
\psline(119.50,-4.33)(114.00,0.00)
\psline(136.50,-38.97)(131.00,-34.64)
\psline(134.00,-34.64)(128.50,-30.31)
\psline(131.50,-30.31)(126.00,-25.98)
\psline(129.00,-25.98)(123.50,-21.65)
\psline(126.50,-21.65)(121.00,-17.32)
\psline(124.00,-17.32)(118.50,-12.99)
\psline(121.50,-12.99)(116.00,-8.66)
\psline(119.00,-8.66)(113.50,-4.33)
\psline(116.50,-4.33)(111.00,0.00)
\psline(133.50,-38.97)(128.00,-34.64)
\psline(131.00,-34.64)(125.50,-30.31)
\psline(128.50,-30.31)(123.00,-25.98)
\psline(126.00,-25.98)(120.50,-21.65)
\psline(123.50,-21.65)(118.00,-17.32)
\psline(121.00,-17.32)(115.50,-12.99)
\psline(118.50,-12.99)(113.00,-8.66)
\psline(116.00,-8.66)(110.50,-4.33)
\psline(113.50,-4.33)(108.00,0.00)
\psline(130.50,-38.97)(125.00,-34.64)
\psline(128.00,-34.64)(122.50,-30.31)
\psline(125.50,-30.31)(120.00,-25.98)
\psline(123.00,-25.98)(117.50,-21.65)
\psline(120.50,-21.65)(115.00,-17.32)
\psline(118.00,-17.32)(112.50,-12.99)
\psline(115.50,-12.99)(110.00,-8.66)
\psline(113.00,-8.66)(107.50,-4.33)
\psline(110.50,-4.33)(105.00,0.00)
\psline(127.50,-38.97)(122.00,-34.64)
\psline(125.00,-34.64)(119.50,-30.31)
\psline(122.50,-30.31)(117.00,-25.98)
\psline(120.00,-25.98)(114.50,-21.65)
\psline(117.50,-21.65)(112.00,-17.32)
\psline(115.00,-17.32)(109.50,-12.99)
\psline(112.50,-12.99)(107.00,-8.66)
\psline(110.00,-8.66)(104.50,-4.33)
\psline(107.50,-4.33)(102.00,0.00)
\pspolygon[fillstyle=solid,linewidth=0.1pt,fillcolor=red](57.00,-77.94)(180.00,-77.94)(202.50,-116.91)(79.50,-116.91)
\psline(180.00,-77.94)(202.50,-116.91)
\psline(177.00,-77.94)(199.50,-116.91)
\psline(174.00,-77.94)(196.50,-116.91)
\psline(171.00,-77.94)(193.50,-116.91)
\psline(168.00,-77.94)(190.50,-116.91)
\psline(165.00,-77.94)(187.50,-116.91)
\psline(162.00,-77.94)(184.50,-116.91)
\psline(159.00,-77.94)(181.50,-116.91)
\psline(156.00,-77.94)(178.50,-116.91)
\psline(153.00,-77.94)(175.50,-116.91)
\psline(150.00,-77.94)(172.50,-116.91)
\psline(147.00,-77.94)(169.50,-116.91)
\psline(144.00,-77.94)(166.50,-116.91)
\psline(141.00,-77.94)(163.50,-116.91)
\psline(138.00,-77.94)(160.50,-116.91)
\psline(135.00,-77.94)(157.50,-116.91)
\psline(132.00,-77.94)(154.50,-116.91)
\psline(129.00,-77.94)(151.50,-116.91)
\psline(126.00,-77.94)(148.50,-116.91)
\psline(123.00,-77.94)(145.50,-116.91)
\psline(120.00,-77.94)(142.50,-116.91)
\psline(117.00,-77.94)(139.50,-116.91)
\psline(114.00,-77.94)(136.50,-116.91)
\psline(111.00,-77.94)(133.50,-116.91)
\psline(108.00,-77.94)(130.50,-116.91)
\psline(105.00,-77.94)(127.50,-116.91)
\psline(102.00,-77.94)(124.50,-116.91)
\psline(99.00,-77.94)(121.50,-116.91)
\psline(96.00,-77.94)(118.50,-116.91)
\psline(93.00,-77.94)(115.50,-116.91)
\psline(90.00,-77.94)(112.50,-116.91)
\psline(87.00,-77.94)(109.50,-116.91)
\psline(84.00,-77.94)(106.50,-116.91)
\psline(81.00,-77.94)(103.50,-116.91)
\psline(78.00,-77.94)(100.50,-116.91)
\psline(75.00,-77.94)(97.50,-116.91)
\psline(72.00,-77.94)(94.50,-116.91)
\psline(69.00,-77.94)(91.50,-116.91)
\psline(66.00,-77.94)(88.50,-116.91)
\psline(63.00,-77.94)(85.50,-116.91)
\psline(60.00,-77.94)(82.50,-116.91)
\psline(57.00,-77.94)(79.50,-116.91)
\psline(79.50,-116.91)(202.50,-116.91)
\psline(77.00,-112.58)(200.00,-112.58)
\psline(74.50,-108.25)(197.50,-108.25)
\psline(72.00,-103.92)(195.00,-103.92)
\psline(69.50,-99.59)(192.50,-99.59)
\psline(67.00,-95.26)(190.00,-95.26)
\psline(64.50,-90.93)(187.50,-90.93)
\psline(62.00,-86.60)(185.00,-86.60)
\psline(59.50,-82.27)(182.50,-82.27)
\psline(57.00,-77.94)(180.00,-77.94)
\psline(202.50,-116.91)(197.00,-112.58)
\psline(200.00,-112.58)(194.50,-108.25)
\psline(197.50,-108.25)(192.00,-103.92)
\psline(195.00,-103.92)(189.50,-99.59)
\psline(192.50,-99.59)(187.00,-95.26)
\psline(190.00,-95.26)(184.50,-90.93)
\psline(187.50,-90.93)(182.00,-86.60)
\psline(185.00,-86.60)(179.50,-82.27)
\psline(182.50,-82.27)(177.00,-77.94)
\psline(199.50,-116.91)(194.00,-112.58)
\psline(197.00,-112.58)(191.50,-108.25)
\psline(194.50,-108.25)(189.00,-103.92)
\psline(192.00,-103.92)(186.50,-99.59)
\psline(189.50,-99.59)(184.00,-95.26)
\psline(187.00,-95.26)(181.50,-90.93)
\psline(184.50,-90.93)(179.00,-86.60)
\psline(182.00,-86.60)(176.50,-82.27)
\psline(179.50,-82.27)(174.00,-77.94)
\psline(196.50,-116.91)(191.00,-112.58)
\psline(194.00,-112.58)(188.50,-108.25)
\psline(191.50,-108.25)(186.00,-103.92)
\psline(189.00,-103.92)(183.50,-99.59)
\psline(186.50,-99.59)(181.00,-95.26)
\psline(184.00,-95.26)(178.50,-90.93)
\psline(181.50,-90.93)(176.00,-86.60)
\psline(179.00,-86.60)(173.50,-82.27)
\psline(176.50,-82.27)(171.00,-77.94)
\psline(193.50,-116.91)(188.00,-112.58)
\psline(191.00,-112.58)(185.50,-108.25)
\psline(188.50,-108.25)(183.00,-103.92)
\psline(186.00,-103.92)(180.50,-99.59)
\psline(183.50,-99.59)(178.00,-95.26)
\psline(181.00,-95.26)(175.50,-90.93)
\psline(178.50,-90.93)(173.00,-86.60)
\psline(176.00,-86.60)(170.50,-82.27)
\psline(173.50,-82.27)(168.00,-77.94)
\psline(190.50,-116.91)(185.00,-112.58)
\psline(188.00,-112.58)(182.50,-108.25)
\psline(185.50,-108.25)(180.00,-103.92)
\psline(183.00,-103.92)(177.50,-99.59)
\psline(180.50,-99.59)(175.00,-95.26)
\psline(178.00,-95.26)(172.50,-90.93)
\psline(175.50,-90.93)(170.00,-86.60)
\psline(173.00,-86.60)(167.50,-82.27)
\psline(170.50,-82.27)(165.00,-77.94)
\psline(187.50,-116.91)(182.00,-112.58)
\psline(185.00,-112.58)(179.50,-108.25)
\psline(182.50,-108.25)(177.00,-103.92)
\psline(180.00,-103.92)(174.50,-99.59)
\psline(177.50,-99.59)(172.00,-95.26)
\psline(175.00,-95.26)(169.50,-90.93)
\psline(172.50,-90.93)(167.00,-86.60)
\psline(170.00,-86.60)(164.50,-82.27)
\psline(167.50,-82.27)(162.00,-77.94)
\psline(184.50,-116.91)(179.00,-112.58)
\psline(182.00,-112.58)(176.50,-108.25)
\psline(179.50,-108.25)(174.00,-103.92)
\psline(177.00,-103.92)(171.50,-99.59)
\psline(174.50,-99.59)(169.00,-95.26)
\psline(172.00,-95.26)(166.50,-90.93)
\psline(169.50,-90.93)(164.00,-86.60)
\psline(167.00,-86.60)(161.50,-82.27)
\psline(164.50,-82.27)(159.00,-77.94)
\psline(181.50,-116.91)(176.00,-112.58)
\psline(179.00,-112.58)(173.50,-108.25)
\psline(176.50,-108.25)(171.00,-103.92)
\psline(174.00,-103.92)(168.50,-99.59)
\psline(171.50,-99.59)(166.00,-95.26)
\psline(169.00,-95.26)(163.50,-90.93)
\psline(166.50,-90.93)(161.00,-86.60)
\psline(164.00,-86.60)(158.50,-82.27)
\psline(161.50,-82.27)(156.00,-77.94)
\psline(178.50,-116.91)(173.00,-112.58)
\psline(176.00,-112.58)(170.50,-108.25)
\psline(173.50,-108.25)(168.00,-103.92)
\psline(171.00,-103.92)(165.50,-99.59)
\psline(168.50,-99.59)(163.00,-95.26)
\psline(166.00,-95.26)(160.50,-90.93)
\psline(163.50,-90.93)(158.00,-86.60)
\psline(161.00,-86.60)(155.50,-82.27)
\psline(158.50,-82.27)(153.00,-77.94)
\psline(175.50,-116.91)(170.00,-112.58)
\psline(173.00,-112.58)(167.50,-108.25)
\psline(170.50,-108.25)(165.00,-103.92)
\psline(168.00,-103.92)(162.50,-99.59)
\psline(165.50,-99.59)(160.00,-95.26)
\psline(163.00,-95.26)(157.50,-90.93)
\psline(160.50,-90.93)(155.00,-86.60)
\psline(158.00,-86.60)(152.50,-82.27)
\psline(155.50,-82.27)(150.00,-77.94)
\psline(172.50,-116.91)(167.00,-112.58)
\psline(170.00,-112.58)(164.50,-108.25)
\psline(167.50,-108.25)(162.00,-103.92)
\psline(165.00,-103.92)(159.50,-99.59)
\psline(162.50,-99.59)(157.00,-95.26)
\psline(160.00,-95.26)(154.50,-90.93)
\psline(157.50,-90.93)(152.00,-86.60)
\psline(155.00,-86.60)(149.50,-82.27)
\psline(152.50,-82.27)(147.00,-77.94)
\psline(169.50,-116.91)(164.00,-112.58)
\psline(167.00,-112.58)(161.50,-108.25)
\psline(164.50,-108.25)(159.00,-103.92)
\psline(162.00,-103.92)(156.50,-99.59)
\psline(159.50,-99.59)(154.00,-95.26)
\psline(157.00,-95.26)(151.50,-90.93)
\psline(154.50,-90.93)(149.00,-86.60)
\psline(152.00,-86.60)(146.50,-82.27)
\psline(149.50,-82.27)(144.00,-77.94)
\psline(166.50,-116.91)(161.00,-112.58)
\psline(164.00,-112.58)(158.50,-108.25)
\psline(161.50,-108.25)(156.00,-103.92)
\psline(159.00,-103.92)(153.50,-99.59)
\psline(156.50,-99.59)(151.00,-95.26)
\psline(154.00,-95.26)(148.50,-90.93)
\psline(151.50,-90.93)(146.00,-86.60)
\psline(149.00,-86.60)(143.50,-82.27)
\psline(146.50,-82.27)(141.00,-77.94)
\psline(163.50,-116.91)(158.00,-112.58)
\psline(161.00,-112.58)(155.50,-108.25)
\psline(158.50,-108.25)(153.00,-103.92)
\psline(156.00,-103.92)(150.50,-99.59)
\psline(153.50,-99.59)(148.00,-95.26)
\psline(151.00,-95.26)(145.50,-90.93)
\psline(148.50,-90.93)(143.00,-86.60)
\psline(146.00,-86.60)(140.50,-82.27)
\psline(143.50,-82.27)(138.00,-77.94)
\psline(160.50,-116.91)(155.00,-112.58)
\psline(158.00,-112.58)(152.50,-108.25)
\psline(155.50,-108.25)(150.00,-103.92)
\psline(153.00,-103.92)(147.50,-99.59)
\psline(150.50,-99.59)(145.00,-95.26)
\psline(148.00,-95.26)(142.50,-90.93)
\psline(145.50,-90.93)(140.00,-86.60)
\psline(143.00,-86.60)(137.50,-82.27)
\psline(140.50,-82.27)(135.00,-77.94)
\psline(157.50,-116.91)(152.00,-112.58)
\psline(155.00,-112.58)(149.50,-108.25)
\psline(152.50,-108.25)(147.00,-103.92)
\psline(150.00,-103.92)(144.50,-99.59)
\psline(147.50,-99.59)(142.00,-95.26)
\psline(145.00,-95.26)(139.50,-90.93)
\psline(142.50,-90.93)(137.00,-86.60)
\psline(140.00,-86.60)(134.50,-82.27)
\psline(137.50,-82.27)(132.00,-77.94)
\psline(154.50,-116.91)(149.00,-112.58)
\psline(152.00,-112.58)(146.50,-108.25)
\psline(149.50,-108.25)(144.00,-103.92)
\psline(147.00,-103.92)(141.50,-99.59)
\psline(144.50,-99.59)(139.00,-95.26)
\psline(142.00,-95.26)(136.50,-90.93)
\psline(139.50,-90.93)(134.00,-86.60)
\psline(137.00,-86.60)(131.50,-82.27)
\psline(134.50,-82.27)(129.00,-77.94)
\psline(151.50,-116.91)(146.00,-112.58)
\psline(149.00,-112.58)(143.50,-108.25)
\psline(146.50,-108.25)(141.00,-103.92)
\psline(144.00,-103.92)(138.50,-99.59)
\psline(141.50,-99.59)(136.00,-95.26)
\psline(139.00,-95.26)(133.50,-90.93)
\psline(136.50,-90.93)(131.00,-86.60)
\psline(134.00,-86.60)(128.50,-82.27)
\psline(131.50,-82.27)(126.00,-77.94)
\psline(148.50,-116.91)(143.00,-112.58)
\psline(146.00,-112.58)(140.50,-108.25)
\psline(143.50,-108.25)(138.00,-103.92)
\psline(141.00,-103.92)(135.50,-99.59)
\psline(138.50,-99.59)(133.00,-95.26)
\psline(136.00,-95.26)(130.50,-90.93)
\psline(133.50,-90.93)(128.00,-86.60)
\psline(131.00,-86.60)(125.50,-82.27)
\psline(128.50,-82.27)(123.00,-77.94)
\psline(145.50,-116.91)(140.00,-112.58)
\psline(143.00,-112.58)(137.50,-108.25)
\psline(140.50,-108.25)(135.00,-103.92)
\psline(138.00,-103.92)(132.50,-99.59)
\psline(135.50,-99.59)(130.00,-95.26)
\psline(133.00,-95.26)(127.50,-90.93)
\psline(130.50,-90.93)(125.00,-86.60)
\psline(128.00,-86.60)(122.50,-82.27)
\psline(125.50,-82.27)(120.00,-77.94)
\psline(142.50,-116.91)(137.00,-112.58)
\psline(140.00,-112.58)(134.50,-108.25)
\psline(137.50,-108.25)(132.00,-103.92)
\psline(135.00,-103.92)(129.50,-99.59)
\psline(132.50,-99.59)(127.00,-95.26)
\psline(130.00,-95.26)(124.50,-90.93)
\psline(127.50,-90.93)(122.00,-86.60)
\psline(125.00,-86.60)(119.50,-82.27)
\psline(122.50,-82.27)(117.00,-77.94)
\psline(139.50,-116.91)(134.00,-112.58)
\psline(137.00,-112.58)(131.50,-108.25)
\psline(134.50,-108.25)(129.00,-103.92)
\psline(132.00,-103.92)(126.50,-99.59)
\psline(129.50,-99.59)(124.00,-95.26)
\psline(127.00,-95.26)(121.50,-90.93)
\psline(124.50,-90.93)(119.00,-86.60)
\psline(122.00,-86.60)(116.50,-82.27)
\psline(119.50,-82.27)(114.00,-77.94)
\psline(136.50,-116.91)(131.00,-112.58)
\psline(134.00,-112.58)(128.50,-108.25)
\psline(131.50,-108.25)(126.00,-103.92)
\psline(129.00,-103.92)(123.50,-99.59)
\psline(126.50,-99.59)(121.00,-95.26)
\psline(124.00,-95.26)(118.50,-90.93)
\psline(121.50,-90.93)(116.00,-86.60)
\psline(119.00,-86.60)(113.50,-82.27)
\psline(116.50,-82.27)(111.00,-77.94)
\psline(133.50,-116.91)(128.00,-112.58)
\psline(131.00,-112.58)(125.50,-108.25)
\psline(128.50,-108.25)(123.00,-103.92)
\psline(126.00,-103.92)(120.50,-99.59)
\psline(123.50,-99.59)(118.00,-95.26)
\psline(121.00,-95.26)(115.50,-90.93)
\psline(118.50,-90.93)(113.00,-86.60)
\psline(116.00,-86.60)(110.50,-82.27)
\psline(113.50,-82.27)(108.00,-77.94)
\psline(130.50,-116.91)(125.00,-112.58)
\psline(128.00,-112.58)(122.50,-108.25)
\psline(125.50,-108.25)(120.00,-103.92)
\psline(123.00,-103.92)(117.50,-99.59)
\psline(120.50,-99.59)(115.00,-95.26)
\psline(118.00,-95.26)(112.50,-90.93)
\psline(115.50,-90.93)(110.00,-86.60)
\psline(113.00,-86.60)(107.50,-82.27)
\psline(110.50,-82.27)(105.00,-77.94)
\psline(127.50,-116.91)(122.00,-112.58)
\psline(125.00,-112.58)(119.50,-108.25)
\psline(122.50,-108.25)(117.00,-103.92)
\psline(120.00,-103.92)(114.50,-99.59)
\psline(117.50,-99.59)(112.00,-95.26)
\psline(115.00,-95.26)(109.50,-90.93)
\psline(112.50,-90.93)(107.00,-86.60)
\psline(110.00,-86.60)(104.50,-82.27)
\psline(107.50,-82.27)(102.00,-77.94)
\psline(124.50,-116.91)(119.00,-112.58)
\psline(122.00,-112.58)(116.50,-108.25)
\psline(119.50,-108.25)(114.00,-103.92)
\psline(117.00,-103.92)(111.50,-99.59)
\psline(114.50,-99.59)(109.00,-95.26)
\psline(112.00,-95.26)(106.50,-90.93)
\psline(109.50,-90.93)(104.00,-86.60)
\psline(107.00,-86.60)(101.50,-82.27)
\psline(104.50,-82.27)(99.00,-77.94)
\psline(121.50,-116.91)(116.00,-112.58)
\psline(119.00,-112.58)(113.50,-108.25)
\psline(116.50,-108.25)(111.00,-103.92)
\psline(114.00,-103.92)(108.50,-99.59)
\psline(111.50,-99.59)(106.00,-95.26)
\psline(109.00,-95.26)(103.50,-90.93)
\psline(106.50,-90.93)(101.00,-86.60)
\psline(104.00,-86.60)(98.50,-82.27)
\psline(101.50,-82.27)(96.00,-77.94)
\psline(118.50,-116.91)(113.00,-112.58)
\psline(116.00,-112.58)(110.50,-108.25)
\psline(113.50,-108.25)(108.00,-103.92)
\psline(111.00,-103.92)(105.50,-99.59)
\psline(108.50,-99.59)(103.00,-95.26)
\psline(106.00,-95.26)(100.50,-90.93)
\psline(103.50,-90.93)(98.00,-86.60)
\psline(101.00,-86.60)(95.50,-82.27)
\psline(98.50,-82.27)(93.00,-77.94)
\psline(115.50,-116.91)(110.00,-112.58)
\psline(113.00,-112.58)(107.50,-108.25)
\psline(110.50,-108.25)(105.00,-103.92)
\psline(108.00,-103.92)(102.50,-99.59)
\psline(105.50,-99.59)(100.00,-95.26)
\psline(103.00,-95.26)(97.50,-90.93)
\psline(100.50,-90.93)(95.00,-86.60)
\psline(98.00,-86.60)(92.50,-82.27)
\psline(95.50,-82.27)(90.00,-77.94)
\psline(112.50,-116.91)(107.00,-112.58)
\psline(110.00,-112.58)(104.50,-108.25)
\psline(107.50,-108.25)(102.00,-103.92)
\psline(105.00,-103.92)(99.50,-99.59)
\psline(102.50,-99.59)(97.00,-95.26)
\psline(100.00,-95.26)(94.50,-90.93)
\psline(97.50,-90.93)(92.00,-86.60)
\psline(95.00,-86.60)(89.50,-82.27)
\psline(92.50,-82.27)(87.00,-77.94)
\psline(109.50,-116.91)(104.00,-112.58)
\psline(107.00,-112.58)(101.50,-108.25)
\psline(104.50,-108.25)(99.00,-103.92)
\psline(102.00,-103.92)(96.50,-99.59)
\psline(99.50,-99.59)(94.00,-95.26)
\psline(97.00,-95.26)(91.50,-90.93)
\psline(94.50,-90.93)(89.00,-86.60)
\psline(92.00,-86.60)(86.50,-82.27)
\psline(89.50,-82.27)(84.00,-77.94)
\psline(106.50,-116.91)(101.00,-112.58)
\psline(104.00,-112.58)(98.50,-108.25)
\psline(101.50,-108.25)(96.00,-103.92)
\psline(99.00,-103.92)(93.50,-99.59)
\psline(96.50,-99.59)(91.00,-95.26)
\psline(94.00,-95.26)(88.50,-90.93)
\psline(91.50,-90.93)(86.00,-86.60)
\psline(89.00,-86.60)(83.50,-82.27)
\psline(86.50,-82.27)(81.00,-77.94)
\psline(103.50,-116.91)(98.00,-112.58)
\psline(101.00,-112.58)(95.50,-108.25)
\psline(98.50,-108.25)(93.00,-103.92)
\psline(96.00,-103.92)(90.50,-99.59)
\psline(93.50,-99.59)(88.00,-95.26)
\psline(91.00,-95.26)(85.50,-90.93)
\psline(88.50,-90.93)(83.00,-86.60)
\psline(86.00,-86.60)(80.50,-82.27)
\psline(83.50,-82.27)(78.00,-77.94)
\psline(100.50,-116.91)(95.00,-112.58)
\psline(98.00,-112.58)(92.50,-108.25)
\psline(95.50,-108.25)(90.00,-103.92)
\psline(93.00,-103.92)(87.50,-99.59)
\psline(90.50,-99.59)(85.00,-95.26)
\psline(88.00,-95.26)(82.50,-90.93)
\psline(85.50,-90.93)(80.00,-86.60)
\psline(83.00,-86.60)(77.50,-82.27)
\psline(80.50,-82.27)(75.00,-77.94)
\psline(97.50,-116.91)(92.00,-112.58)
\psline(95.00,-112.58)(89.50,-108.25)
\psline(92.50,-108.25)(87.00,-103.92)
\psline(90.00,-103.92)(84.50,-99.59)
\psline(87.50,-99.59)(82.00,-95.26)
\psline(85.00,-95.26)(79.50,-90.93)
\psline(82.50,-90.93)(77.00,-86.60)
\psline(80.00,-86.60)(74.50,-82.27)
\psline(77.50,-82.27)(72.00,-77.94)
\psline(94.50,-116.91)(89.00,-112.58)
\psline(92.00,-112.58)(86.50,-108.25)
\psline(89.50,-108.25)(84.00,-103.92)
\psline(87.00,-103.92)(81.50,-99.59)
\psline(84.50,-99.59)(79.00,-95.26)
\psline(82.00,-95.26)(76.50,-90.93)
\psline(79.50,-90.93)(74.00,-86.60)
\psline(77.00,-86.60)(71.50,-82.27)
\psline(74.50,-82.27)(69.00,-77.94)
\psline(91.50,-116.91)(86.00,-112.58)
\psline(89.00,-112.58)(83.50,-108.25)
\psline(86.50,-108.25)(81.00,-103.92)
\psline(84.00,-103.92)(78.50,-99.59)
\psline(81.50,-99.59)(76.00,-95.26)
\psline(79.00,-95.26)(73.50,-90.93)
\psline(76.50,-90.93)(71.00,-86.60)
\psline(74.00,-86.60)(68.50,-82.27)
\psline(71.50,-82.27)(66.00,-77.94)
\psline(88.50,-116.91)(83.00,-112.58)
\psline(86.00,-112.58)(80.50,-108.25)
\psline(83.50,-108.25)(78.00,-103.92)
\psline(81.00,-103.92)(75.50,-99.59)
\psline(78.50,-99.59)(73.00,-95.26)
\psline(76.00,-95.26)(70.50,-90.93)
\psline(73.50,-90.93)(68.00,-86.60)
\psline(71.00,-86.60)(65.50,-82.27)
\psline(68.50,-82.27)(63.00,-77.94)
\psline(85.50,-116.91)(80.00,-112.58)
\psline(83.00,-112.58)(77.50,-108.25)
\psline(80.50,-108.25)(75.00,-103.92)
\psline(78.00,-103.92)(72.50,-99.59)
\psline(75.50,-99.59)(70.00,-95.26)
\psline(73.00,-95.26)(67.50,-90.93)
\psline(70.50,-90.93)(65.00,-86.60)
\psline(68.00,-86.60)(62.50,-82.27)
\psline(65.50,-82.27)(60.00,-77.94)
\psline(82.50,-116.91)(77.00,-112.58)
\psline(80.00,-112.58)(74.50,-108.25)
\psline(77.50,-108.25)(72.00,-103.92)
\psline(75.00,-103.92)(69.50,-99.59)
\psline(72.50,-99.59)(67.00,-95.26)
\psline(70.00,-95.26)(64.50,-90.93)
\psline(67.50,-90.93)(62.00,-86.60)
\psline(65.00,-86.60)(59.50,-82.27)
\psline(62.50,-82.27)(57.00,-77.94)
\pspolygon[fillstyle=solid,linewidth=0.1pt,fillcolor=lightblue](45.00,0.00)(-37.50,64.95)(-28.50,127.31)
\psline(-37.50,64.95)(-28.50,127.31)
\psline(-33.57,61.86)(-25.00,121.24)
\psline(-29.64,58.77)(-21.50,115.18)
\psline(-25.71,55.67)(-18.00,109.12)
\psline(-21.79,52.58)(-14.50,103.06)
\psline(-17.86,49.49)(-11.00,96.99)
\psline(-13.93,46.39)(-7.50,90.93)
\psline(-10.00,43.30)(-4.00,84.87)
\psline(-6.07,40.21)(-0.50,78.81)
\psline(-2.14,37.12)(3.00,72.75)
\psline(1.79,34.02)(6.50,66.68)
\psline(5.71,30.93)(10.00,60.62)
\psline(9.64,27.84)(13.50,54.56)
\psline(13.57,24.74)(17.00,48.50)
\psline(17.50,21.65)(20.50,42.44)
\psline(21.43,18.56)(24.00,36.37)
\psline(25.36,15.46)(27.50,30.31)
\psline(29.29,12.37)(31.00,24.25)
\psline(33.21,9.28)(34.50,18.19)
\psline(37.14,6.19)(38.00,12.12)
\psline(41.07,3.09)(41.50,6.06)
\psline(45.00,0.00)(45.00,0.00)
\psline(-28.50,127.31)(45.00,0.00)
\psline(-28.93,124.34)(41.07,3.09)
\psline(-29.36,121.37)(37.14,6.19)
\psline(-29.79,118.40)(33.21,9.28)
\psline(-30.21,115.43)(29.29,12.37)
\psline(-30.64,112.46)(25.36,15.46)
\psline(-31.07,109.49)(21.43,18.56)
\psline(-31.50,106.52)(17.50,21.65)
\psline(-31.93,103.55)(13.57,24.74)
\psline(-32.36,100.58)(9.64,27.84)
\psline(-32.79,97.61)(5.71,30.93)
\psline(-33.21,94.64)(1.79,34.02)
\psline(-33.64,91.67)(-2.14,37.12)
\psline(-34.07,88.71)(-6.07,40.21)
\psline(-34.50,85.74)(-10.00,43.30)
\psline(-34.93,82.77)(-13.93,46.39)
\psline(-35.36,79.80)(-17.86,49.49)
\psline(-35.79,76.83)(-21.79,52.58)
\psline(-36.21,73.86)(-25.71,55.67)
\psline(-36.64,70.89)(-29.64,58.77)
\psline(-37.07,67.92)(-33.57,61.86)
\psline(-37.50,64.95)(-37.50,64.95)
\psline(45.00,0.00)(-37.50,64.95)
\psline(41.50,6.06)(-37.07,67.92)
\psline(38.00,12.12)(-36.64,70.89)
\psline(34.50,18.19)(-36.21,73.86)
\psline(31.00,24.25)(-35.79,76.83)
\psline(27.50,30.31)(-35.36,79.80)
\psline(24.00,36.37)(-34.93,82.77)
\psline(20.50,42.44)(-34.50,85.74)
\psline(17.00,48.50)(-34.07,88.71)
\psline(13.50,54.56)(-33.64,91.67)
\psline(10.00,60.62)(-33.21,94.64)
\psline(6.50,66.68)(-32.79,97.61)
\psline(3.00,72.75)(-32.36,100.58)
\psline(-0.50,78.81)(-31.93,103.55)
\psline(-4.00,84.87)(-31.50,106.52)
\psline(-7.50,90.93)(-31.07,109.49)
\psline(-11.00,96.99)(-30.64,112.46)
\psline(-14.50,103.06)(-30.21,115.43)
\psline(-18.00,109.12)(-29.79,118.40)
\psline(-21.50,115.18)(-29.36,121.37)
\psline(-25.00,121.24)(-28.93,124.34)
\psline(-28.50,127.31)(-28.50,127.31)
\pspolygon[fillstyle=solid,linewidth=0.1pt,fillcolor=lightblue](90.00,0.00)(7.50,64.95)(16.50,127.31)
\psline(7.50,64.95)(16.50,127.31)
\psline(11.43,61.86)(20.00,121.24)
\psline(15.36,58.77)(23.50,115.18)
\psline(19.29,55.67)(27.00,109.12)
\psline(23.21,52.58)(30.50,103.06)
\psline(27.14,49.49)(34.00,96.99)
\psline(31.07,46.39)(37.50,90.93)
\psline(35.00,43.30)(41.00,84.87)
\psline(38.93,40.21)(44.50,78.81)
\psline(42.86,37.12)(48.00,72.75)
\psline(46.79,34.02)(51.50,66.68)
\psline(50.71,30.93)(55.00,60.62)
\psline(54.64,27.84)(58.50,54.56)
\psline(58.57,24.74)(62.00,48.50)
\psline(62.50,21.65)(65.50,42.44)
\psline(66.43,18.56)(69.00,36.37)
\psline(70.36,15.46)(72.50,30.31)
\psline(74.29,12.37)(76.00,24.25)
\psline(78.21,9.28)(79.50,18.19)
\psline(82.14,6.19)(83.00,12.12)
\psline(86.07,3.09)(86.50,6.06)
\psline(90.00,0.00)(90.00,0.00)
\psline(16.50,127.31)(90.00,0.00)
\psline(16.07,124.34)(86.07,3.09)
\psline(15.64,121.37)(82.14,6.19)
\psline(15.21,118.40)(78.21,9.28)
\psline(14.79,115.43)(74.29,12.37)
\psline(14.36,112.46)(70.36,15.46)
\psline(13.93,109.49)(66.43,18.56)
\psline(13.50,106.52)(62.50,21.65)
\psline(13.07,103.55)(58.57,24.74)
\psline(12.64,100.58)(54.64,27.84)
\psline(12.21,97.61)(50.71,30.93)
\psline(11.79,94.64)(46.79,34.02)
\psline(11.36,91.67)(42.86,37.12)
\psline(10.93,88.71)(38.93,40.21)
\psline(10.50,85.74)(35.00,43.30)
\psline(10.07,82.77)(31.07,46.39)
\psline(9.64,79.80)(27.14,49.49)
\psline(9.21,76.83)(23.21,52.58)
\psline(8.79,73.86)(19.29,55.67)
\psline(8.36,70.89)(15.36,58.77)
\psline(7.93,67.92)(11.43,61.86)
\psline(7.50,64.95)(7.50,64.95)
\psline(90.00,0.00)(7.50,64.95)
\psline(86.50,6.06)(7.93,67.92)
\psline(83.00,12.12)(8.36,70.89)
\psline(79.50,18.19)(8.79,73.86)
\psline(76.00,24.25)(9.21,76.83)
\psline(72.50,30.31)(9.64,79.80)
\psline(69.00,36.37)(10.07,82.77)
\psline(65.50,42.44)(10.50,85.74)
\psline(62.00,48.50)(10.93,88.71)
\psline(58.50,54.56)(11.36,91.67)
\psline(55.00,60.62)(11.79,94.64)
\psline(51.50,66.68)(12.21,97.61)
\psline(48.00,72.75)(12.64,100.58)
\psline(44.50,78.81)(13.07,103.55)
\psline(41.00,84.87)(13.50,106.52)
\psline(37.50,90.93)(13.93,109.49)
\psline(34.00,96.99)(14.36,112.46)
\psline(30.50,103.06)(14.79,115.43)
\psline(27.00,109.12)(15.21,118.40)
\psline(23.50,115.18)(15.64,121.37)
\psline(20.00,121.24)(16.07,124.34)
\psline(16.50,127.31)(16.50,127.31)
\pspolygon[fillstyle=solid,linewidth=0.1pt,fillcolor=lightblue](135.00,0.00)(61.50,127.31)(52.50,64.95)
\psline(61.50,127.31)(52.50,64.95)
\psline(65.00,121.24)(56.43,61.86)
\psline(68.50,115.18)(60.36,58.77)
\psline(72.00,109.12)(64.29,55.67)
\psline(75.50,103.06)(68.21,52.58)
\psline(79.00,96.99)(72.14,49.49)
\psline(82.50,90.93)(76.07,46.39)
\psline(86.00,84.87)(80.00,43.30)
\psline(89.50,78.81)(83.93,40.21)
\psline(93.00,72.75)(87.86,37.12)
\psline(96.50,66.68)(91.79,34.02)
\psline(100.00,60.62)(95.71,30.93)
\psline(103.50,54.56)(99.64,27.84)
\psline(107.00,48.50)(103.57,24.74)
\psline(110.50,42.44)(107.50,21.65)
\psline(114.00,36.37)(111.43,18.56)
\psline(117.50,30.31)(115.36,15.46)
\psline(121.00,24.25)(119.29,12.37)
\psline(124.50,18.19)(123.21,9.28)
\psline(128.00,12.12)(127.14,6.19)
\psline(131.50,6.06)(131.07,3.09)
\psline(135.00,0.00)(135.00,0.00)
\psline(52.50,64.95)(135.00,0.00)
\psline(52.93,67.92)(131.50,6.06)
\psline(53.36,70.89)(128.00,12.12)
\psline(53.79,73.86)(124.50,18.19)
\psline(54.21,76.83)(121.00,24.25)
\psline(54.64,79.80)(117.50,30.31)
\psline(55.07,82.77)(114.00,36.37)
\psline(55.50,85.74)(110.50,42.44)
\psline(55.93,88.71)(107.00,48.50)
\psline(56.36,91.67)(103.50,54.56)
\psline(56.79,94.64)(100.00,60.62)
\psline(57.21,97.61)(96.50,66.68)
\psline(57.64,100.58)(93.00,72.75)
\psline(58.07,103.55)(89.50,78.81)
\psline(58.50,106.52)(86.00,84.87)
\psline(58.93,109.49)(82.50,90.93)
\psline(59.36,112.46)(79.00,96.99)
\psline(59.79,115.43)(75.50,103.06)
\psline(60.21,118.40)(72.00,109.12)
\psline(60.64,121.37)(68.50,115.18)
\psline(61.07,124.34)(65.00,121.24)
\psline(61.50,127.31)(61.50,127.31)
\psline(135.00,0.00)(61.50,127.31)
\psline(131.07,3.09)(61.07,124.34)
\psline(127.14,6.19)(60.64,121.37)
\psline(123.21,9.28)(60.21,118.40)
\psline(119.29,12.37)(59.79,115.43)
\psline(115.36,15.46)(59.36,112.46)
\psline(111.43,18.56)(58.93,109.49)
\psline(107.50,21.65)(58.50,106.52)
\psline(103.57,24.74)(58.07,103.55)
\psline(99.64,27.84)(57.64,100.58)
\psline(95.71,30.93)(57.21,97.61)
\psline(91.79,34.02)(56.79,94.64)
\psline(87.86,37.12)(56.36,91.67)
\psline(83.93,40.21)(55.93,88.71)
\psline(80.00,43.30)(55.50,85.74)
\psline(76.07,46.39)(55.07,82.77)
\psline(72.14,49.49)(54.64,79.80)
\psline(68.21,52.58)(54.21,76.83)
\psline(64.29,55.67)(53.79,73.86)
\psline(60.36,58.77)(53.36,70.89)
\psline(56.43,61.86)(52.93,67.92)
\psline(52.50,64.95)(52.50,64.95)
\pspolygon[fillstyle=solid,linewidth=0.1pt,fillcolor=lightyellow](-0.00,0.00)(-37.50,64.95)(45.00,0.00)
\psline(-37.50,64.95)(45.00,0.00)
\psline(-35.00,60.62)(42.00,0.00)
\psline(-32.50,56.29)(39.00,0.00)
\psline(-30.00,51.96)(36.00,0.00)
\psline(-27.50,47.63)(33.00,0.00)
\psline(-25.00,43.30)(30.00,0.00)
\psline(-22.50,38.97)(27.00,0.00)
\psline(-20.00,34.64)(24.00,0.00)
\psline(-17.50,30.31)(21.00,0.00)
\psline(-15.00,25.98)(18.00,0.00)
\psline(-12.50,21.65)(15.00,0.00)
\psline(-10.00,17.32)(12.00,0.00)
\psline(-7.50,12.99)(9.00,0.00)
\psline(-5.00,8.66)(6.00,0.00)
\psline(-2.50,4.33)(3.00,0.00)
\psline(-0.00,0.00)(0.00,0.00)
\psline(45.00,0.00)(-0.00,0.00)
\psline(39.50,4.33)(-2.50,4.33)
\psline(34.00,8.66)(-5.00,8.66)
\psline(28.50,12.99)(-7.50,12.99)
\psline(23.00,17.32)(-10.00,17.32)
\psline(17.50,21.65)(-12.50,21.65)
\psline(12.00,25.98)(-15.00,25.98)
\psline(6.50,30.31)(-17.50,30.31)
\psline(1.00,34.64)(-20.00,34.64)
\psline(-4.50,38.97)(-22.50,38.97)
\psline(-10.00,43.30)(-25.00,43.30)
\psline(-15.50,47.63)(-27.50,47.63)
\psline(-21.00,51.96)(-30.00,51.96)
\psline(-26.50,56.29)(-32.50,56.29)
\psline(-32.00,60.62)(-35.00,60.62)
\psline(-37.50,64.95)(-37.50,64.95)
\psline(0.00,0.00)(-37.50,64.95)
\psline(3.00,0.00)(-32.00,60.62)
\psline(6.00,0.00)(-26.50,56.29)
\psline(9.00,0.00)(-21.00,51.96)
\psline(12.00,0.00)(-15.50,47.63)
\psline(15.00,0.00)(-10.00,43.30)
\psline(18.00,0.00)(-4.50,38.97)
\psline(21.00,0.00)(1.00,34.64)
\psline(24.00,0.00)(6.50,30.31)
\psline(27.00,0.00)(12.00,25.98)
\psline(30.00,0.00)(17.50,21.65)
\psline(33.00,0.00)(23.00,17.32)
\psline(36.00,0.00)(28.50,12.99)
\psline(39.00,0.00)(34.00,8.66)
\psline(42.00,0.00)(39.50,4.33)
\psline(45.00,0.00)(45.00,0.00)
\pspolygon[fillstyle=solid,linewidth=0.1pt,fillcolor=lightyellow](90.00,0.00)(45.00,0.00)(7.50,64.95)
\psline(45.00,0.00)(7.50,64.95)
\psline(48.00,0.00)(13.00,60.62)
\psline(51.00,0.00)(18.50,56.29)
\psline(54.00,0.00)(24.00,51.96)
\psline(57.00,0.00)(29.50,47.63)
\psline(60.00,0.00)(35.00,43.30)
\psline(63.00,0.00)(40.50,38.97)
\psline(66.00,0.00)(46.00,34.64)
\psline(69.00,0.00)(51.50,30.31)
\psline(72.00,0.00)(57.00,25.98)
\psline(75.00,0.00)(62.50,21.65)
\psline(78.00,0.00)(68.00,17.32)
\psline(81.00,0.00)(73.50,12.99)
\psline(84.00,0.00)(79.00,8.66)
\psline(87.00,0.00)(84.50,4.33)
\psline(90.00,0.00)(90.00,0.00)
\psline(7.50,64.95)(90.00,0.00)
\psline(10.00,60.62)(87.00,0.00)
\psline(12.50,56.29)(84.00,0.00)
\psline(15.00,51.96)(81.00,0.00)
\psline(17.50,47.63)(78.00,0.00)
\psline(20.00,43.30)(75.00,0.00)
\psline(22.50,38.97)(72.00,0.00)
\psline(25.00,34.64)(69.00,0.00)
\psline(27.50,30.31)(66.00,0.00)
\psline(30.00,25.98)(63.00,0.00)
\psline(32.50,21.65)(60.00,0.00)
\psline(35.00,17.32)(57.00,0.00)
\psline(37.50,12.99)(54.00,0.00)
\psline(40.00,8.66)(51.00,0.00)
\psline(42.50,4.33)(48.00,0.00)
\psline(45.00,0.00)(45.00,0.00)
\psline(90.00,0.00)(45.00,0.00)
\psline(84.50,4.33)(42.50,4.33)
\psline(79.00,8.66)(40.00,8.66)
\psline(73.50,12.99)(37.50,12.99)
\psline(68.00,17.32)(35.00,17.32)
\psline(62.50,21.65)(32.50,21.65)
\psline(57.00,25.98)(30.00,25.98)
\psline(51.50,30.31)(27.50,30.31)
\psline(46.00,34.64)(25.00,34.64)
\psline(40.50,38.97)(22.50,38.97)
\psline(35.00,43.30)(20.00,43.30)
\psline(29.50,47.63)(17.50,47.63)
\psline(24.00,51.96)(15.00,51.96)
\psline(18.50,56.29)(12.50,56.29)
\psline(13.00,60.62)(10.00,60.62)
\psline(7.50,64.95)(7.50,64.95)
\pspolygon[fillstyle=solid,linewidth=0.1pt,fillcolor=lightyellow](135.00,0.00)(90.00,0.00)(52.50,64.95)
\psline(90.00,0.00)(52.50,64.95)
\psline(93.00,0.00)(58.00,60.62)
\psline(96.00,0.00)(63.50,56.29)
\psline(99.00,0.00)(69.00,51.96)
\psline(102.00,0.00)(74.50,47.63)
\psline(105.00,0.00)(80.00,43.30)
\psline(108.00,0.00)(85.50,38.97)
\psline(111.00,0.00)(91.00,34.64)
\psline(114.00,0.00)(96.50,30.31)
\psline(117.00,0.00)(102.00,25.98)
\psline(120.00,0.00)(107.50,21.65)
\psline(123.00,0.00)(113.00,17.32)
\psline(126.00,0.00)(118.50,12.99)
\psline(129.00,0.00)(124.00,8.66)
\psline(132.00,0.00)(129.50,4.33)
\psline(135.00,0.00)(135.00,0.00)
\psline(52.50,64.95)(135.00,0.00)
\psline(55.00,60.62)(132.00,0.00)
\psline(57.50,56.29)(129.00,0.00)
\psline(60.00,51.96)(126.00,0.00)
\psline(62.50,47.63)(123.00,0.00)
\psline(65.00,43.30)(120.00,0.00)
\psline(67.50,38.97)(117.00,0.00)
\psline(70.00,34.64)(114.00,0.00)
\psline(72.50,30.31)(111.00,0.00)
\psline(75.00,25.98)(108.00,0.00)
\psline(77.50,21.65)(105.00,0.00)
\psline(80.00,17.32)(102.00,0.00)
\psline(82.50,12.99)(99.00,0.00)
\psline(85.00,8.66)(96.00,0.00)
\psline(87.50,4.33)(93.00,0.00)
\psline(90.00,0.00)(90.00,0.00)
\psline(135.00,0.00)(90.00,0.00)
\psline(129.50,4.33)(87.50,4.33)
\psline(124.00,8.66)(85.00,8.66)
\psline(118.50,12.99)(82.50,12.99)
\psline(113.00,17.32)(80.00,17.32)
\psline(107.50,21.65)(77.50,21.65)
\psline(102.00,25.98)(75.00,25.98)
\psline(96.50,30.31)(72.50,30.31)
\psline(91.00,34.64)(70.00,34.64)
\psline(85.50,38.97)(67.50,38.97)
\psline(80.00,43.30)(65.00,43.30)
\psline(74.50,47.63)(62.50,47.63)
\psline(69.00,51.96)(60.00,51.96)
\psline(63.50,56.29)(57.50,56.29)
\psline(58.00,60.62)(55.00,60.62)
\psline(52.50,64.95)(52.50,64.95)
\pspolygon[fillstyle=solid,linewidth=0.1pt,fillcolor=green](-37.50,64.95)(-51.00,88.33)(-28.50,127.31)
\psline(-51.00,88.33)(-28.50,127.31)
\psline(-49.50,85.74)(-29.50,120.38)
\psline(-48.00,83.14)(-30.50,113.45)
\psline(-46.50,80.54)(-31.50,106.52)
\psline(-45.00,77.94)(-32.50,99.59)
\psline(-43.50,75.34)(-33.50,92.66)
\psline(-42.00,72.75)(-34.50,85.74)
\psline(-40.50,70.15)(-35.50,78.81)
\psline(-39.00,67.55)(-36.50,71.88)
\psline(-37.50,64.95)(-37.50,64.95)
\psline(-28.50,127.31)(-37.50,64.95)
\psline(-31.00,122.98)(-39.00,67.55)
\psline(-33.50,118.65)(-40.50,70.15)
\psline(-36.00,114.32)(-42.00,72.75)
\psline(-38.50,109.99)(-43.50,75.34)
\psline(-41.00,105.66)(-45.00,77.94)
\psline(-43.50,101.32)(-46.50,80.54)
\psline(-46.00,96.99)(-48.00,83.14)
\psline(-48.50,92.66)(-49.50,85.74)
\psline(-51.00,88.33)(-51.00,88.33)
\psline(-37.50,64.95)(-51.00,88.33)
\psline(-36.50,71.88)(-48.50,92.66)
\psline(-35.50,78.81)(-46.00,96.99)
\psline(-34.50,85.74)(-43.50,101.32)
\psline(-33.50,92.66)(-41.00,105.66)
\psline(-32.50,99.59)(-38.50,109.99)
\psline(-31.50,106.52)(-36.00,114.32)
\psline(-30.50,113.45)(-33.50,118.65)
\psline(-29.50,120.38)(-31.00,122.98)
\psline(-28.50,127.31)(-28.50,127.31)
\pspolygon[fillstyle=solid,linewidth=0.1pt,fillcolor=green](7.50,64.95)(-6.00,88.33)(16.50,127.31)
\psline(-6.00,88.33)(16.50,127.31)
\psline(-4.50,85.74)(15.50,120.38)
\psline(-3.00,83.14)(14.50,113.45)
\psline(-1.50,80.54)(13.50,106.52)
\psline(0.00,77.94)(12.50,99.59)
\psline(1.50,75.34)(11.50,92.66)
\psline(3.00,72.75)(10.50,85.74)
\psline(4.50,70.15)(9.50,78.81)
\psline(6.00,67.55)(8.50,71.88)
\psline(7.50,64.95)(7.50,64.95)
\psline(16.50,127.31)(7.50,64.95)
\psline(14.00,122.98)(6.00,67.55)
\psline(11.50,118.65)(4.50,70.15)
\psline(9.00,114.32)(3.00,72.75)
\psline(6.50,109.99)(1.50,75.34)
\psline(4.00,105.66)(0.00,77.94)
\psline(1.50,101.32)(-1.50,80.54)
\psline(-1.00,96.99)(-3.00,83.14)
\psline(-3.50,92.66)(-4.50,85.74)
\psline(-6.00,88.33)(-6.00,88.33)
\psline(7.50,64.95)(-6.00,88.33)
\psline(8.50,71.88)(-3.50,92.66)
\psline(9.50,78.81)(-1.00,96.99)
\psline(10.50,85.74)(1.50,101.32)
\psline(11.50,92.66)(4.00,105.66)
\psline(12.50,99.59)(6.50,109.99)
\psline(13.50,106.52)(9.00,114.32)
\psline(14.50,113.45)(11.50,118.65)
\psline(15.50,120.38)(14.00,122.98)
\psline(16.50,127.31)(16.50,127.31)
\pspolygon[fillstyle=solid,linewidth=0.1pt,fillcolor=green](52.50,64.95)(39.00,88.33)(61.50,127.31)
\psline(39.00,88.33)(61.50,127.31)
\psline(40.50,85.74)(60.50,120.38)
\psline(42.00,83.14)(59.50,113.45)
\psline(43.50,80.54)(58.50,106.52)
\psline(45.00,77.94)(57.50,99.59)
\psline(46.50,75.34)(56.50,92.66)
\psline(48.00,72.75)(55.50,85.74)
\psline(49.50,70.15)(54.50,78.81)
\psline(51.00,67.55)(53.50,71.88)
\psline(52.50,64.95)(52.50,64.95)
\psline(61.50,127.31)(52.50,64.95)
\psline(59.00,122.98)(51.00,67.55)
\psline(56.50,118.65)(49.50,70.15)
\psline(54.00,114.32)(48.00,72.75)
\psline(51.50,109.99)(46.50,75.34)
\psline(49.00,105.66)(45.00,77.94)
\psline(46.50,101.32)(43.50,80.54)
\psline(44.00,96.99)(42.00,83.14)
\psline(41.50,92.66)(40.50,85.74)
\psline(39.00,88.33)(39.00,88.33)
\psline(52.50,64.95)(39.00,88.33)
\psline(53.50,71.88)(41.50,92.66)
\psline(54.50,78.81)(44.00,96.99)
\psline(55.50,85.74)(46.50,101.32)
\psline(56.50,92.66)(49.00,105.66)
\psline(57.50,99.59)(51.50,109.99)
\psline(58.50,106.52)(54.00,114.32)
\psline(59.50,113.45)(56.50,118.65)
\psline(60.50,120.38)(59.00,122.98)
\psline(61.50,127.31)(61.50,127.31)
\pspolygon[fillstyle=solid,linewidth=0.1pt,fillcolor=red](-6.00,88.33)(-45.00,155.88)(-22.50,194.86)(16.50,127.31)
\psline(-45.00,155.88)(-22.50,194.86)
\psline(-43.50,153.29)(-21.00,192.26)
\psline(-42.00,150.69)(-19.50,189.66)
\psline(-40.50,148.09)(-18.00,187.06)
\psline(-39.00,145.49)(-16.50,184.46)
\psline(-37.50,142.89)(-15.00,181.87)
\psline(-36.00,140.30)(-13.50,179.27)
\psline(-34.50,137.70)(-12.00,176.67)
\psline(-33.00,135.10)(-10.50,174.07)
\psline(-31.50,132.50)(-9.00,171.47)
\psline(-30.00,129.90)(-7.50,168.87)
\psline(-28.50,127.31)(-6.00,166.28)
\psline(-27.00,124.71)(-4.50,163.68)
\psline(-25.50,122.11)(-3.00,161.08)
\psline(-24.00,119.51)(-1.50,158.48)
\psline(-22.50,116.91)(0.00,155.88)
\psline(-21.00,114.32)(1.50,153.29)
\psline(-19.50,111.72)(3.00,150.69)
\psline(-18.00,109.12)(4.50,148.09)
\psline(-16.50,106.52)(6.00,145.49)
\psline(-15.00,103.92)(7.50,142.89)
\psline(-13.50,101.32)(9.00,140.30)
\psline(-12.00,98.73)(10.50,137.70)
\psline(-10.50,96.13)(12.00,135.10)
\psline(-9.00,93.53)(13.50,132.50)
\psline(-7.50,90.93)(15.00,129.90)
\psline(-6.00,88.33)(16.50,127.31)
\psline(16.50,127.31)(-22.50,194.86)
\psline(14.00,122.98)(-25.00,190.53)
\psline(11.50,118.65)(-27.50,186.20)
\psline(9.00,114.32)(-30.00,181.87)
\psline(6.50,109.99)(-32.50,177.54)
\psline(4.00,105.66)(-35.00,173.21)
\psline(1.50,101.32)(-37.50,168.87)
\psline(-1.00,96.99)(-40.00,164.54)
\psline(-3.50,92.66)(-42.50,160.21)
\psline(-6.00,88.33)(-45.00,155.88)
\psline(-22.50,194.86)(-23.50,187.93)
\psline(-25.00,190.53)(-26.00,183.60)
\psline(-27.50,186.20)(-28.50,179.27)
\psline(-30.00,181.87)(-31.00,174.94)
\psline(-32.50,177.54)(-33.50,170.61)
\psline(-35.00,173.21)(-36.00,166.28)
\psline(-37.50,168.87)(-38.50,161.95)
\psline(-40.00,164.54)(-41.00,157.62)
\psline(-42.50,160.21)(-43.50,153.29)
\psline(-21.00,192.26)(-22.00,185.33)
\psline(-23.50,187.93)(-24.50,181.00)
\psline(-26.00,183.60)(-27.00,176.67)
\psline(-28.50,179.27)(-29.50,172.34)
\psline(-31.00,174.94)(-32.00,168.01)
\psline(-33.50,170.61)(-34.50,163.68)
\psline(-36.00,166.28)(-37.00,159.35)
\psline(-38.50,161.95)(-39.50,155.02)
\psline(-41.00,157.62)(-42.00,150.69)
\psline(-19.50,189.66)(-20.50,182.73)
\psline(-22.00,185.33)(-23.00,178.40)
\psline(-24.50,181.00)(-25.50,174.07)
\psline(-27.00,176.67)(-28.00,169.74)
\psline(-29.50,172.34)(-30.50,165.41)
\psline(-32.00,168.01)(-33.00,161.08)
\psline(-34.50,163.68)(-35.50,156.75)
\psline(-37.00,159.35)(-38.00,152.42)
\psline(-39.50,155.02)(-40.50,148.09)
\psline(-18.00,187.06)(-19.00,180.13)
\psline(-20.50,182.73)(-21.50,175.80)
\psline(-23.00,178.40)(-24.00,171.47)
\psline(-25.50,174.07)(-26.50,167.14)
\psline(-28.00,169.74)(-29.00,162.81)
\psline(-30.50,165.41)(-31.50,158.48)
\psline(-33.00,161.08)(-34.00,154.15)
\psline(-35.50,156.75)(-36.50,149.82)
\psline(-38.00,152.42)(-39.00,145.49)
\psline(-16.50,184.46)(-17.50,177.54)
\psline(-19.00,180.13)(-20.00,173.21)
\psline(-21.50,175.80)(-22.50,168.87)
\psline(-24.00,171.47)(-25.00,164.54)
\psline(-26.50,167.14)(-27.50,160.21)
\psline(-29.00,162.81)(-30.00,155.88)
\psline(-31.50,158.48)(-32.50,151.55)
\psline(-34.00,154.15)(-35.00,147.22)
\psline(-36.50,149.82)(-37.50,142.89)
\psline(-15.00,181.87)(-16.00,174.94)
\psline(-17.50,177.54)(-18.50,170.61)
\psline(-20.00,173.21)(-21.00,166.28)
\psline(-22.50,168.87)(-23.50,161.95)
\psline(-25.00,164.54)(-26.00,157.62)
\psline(-27.50,160.21)(-28.50,153.29)
\psline(-30.00,155.88)(-31.00,148.96)
\psline(-32.50,151.55)(-33.50,144.63)
\psline(-35.00,147.22)(-36.00,140.30)
\psline(-13.50,179.27)(-14.50,172.34)
\psline(-16.00,174.94)(-17.00,168.01)
\psline(-18.50,170.61)(-19.50,163.68)
\psline(-21.00,166.28)(-22.00,159.35)
\psline(-23.50,161.95)(-24.50,155.02)
\psline(-26.00,157.62)(-27.00,150.69)
\psline(-28.50,153.29)(-29.50,146.36)
\psline(-31.00,148.96)(-32.00,142.03)
\psline(-33.50,144.63)(-34.50,137.70)
\psline(-12.00,176.67)(-13.00,169.74)
\psline(-14.50,172.34)(-15.50,165.41)
\psline(-17.00,168.01)(-18.00,161.08)
\psline(-19.50,163.68)(-20.50,156.75)
\psline(-22.00,159.35)(-23.00,152.42)
\psline(-24.50,155.02)(-25.50,148.09)
\psline(-27.00,150.69)(-28.00,143.76)
\psline(-29.50,146.36)(-30.50,139.43)
\psline(-32.00,142.03)(-33.00,135.10)
\psline(-10.50,174.07)(-11.50,167.14)
\psline(-13.00,169.74)(-14.00,162.81)
\psline(-15.50,165.41)(-16.50,158.48)
\psline(-18.00,161.08)(-19.00,154.15)
\psline(-20.50,156.75)(-21.50,149.82)
\psline(-23.00,152.42)(-24.00,145.49)
\psline(-25.50,148.09)(-26.50,141.16)
\psline(-28.00,143.76)(-29.00,136.83)
\psline(-30.50,139.43)(-31.50,132.50)
\psline(-9.00,171.47)(-10.00,164.54)
\psline(-11.50,167.14)(-12.50,160.21)
\psline(-14.00,162.81)(-15.00,155.88)
\psline(-16.50,158.48)(-17.50,151.55)
\psline(-19.00,154.15)(-20.00,147.22)
\psline(-21.50,149.82)(-22.50,142.89)
\psline(-24.00,145.49)(-25.00,138.56)
\psline(-26.50,141.16)(-27.50,134.23)
\psline(-29.00,136.83)(-30.00,129.90)
\psline(-7.50,168.87)(-8.50,161.95)
\psline(-10.00,164.54)(-11.00,157.62)
\psline(-12.50,160.21)(-13.50,153.29)
\psline(-15.00,155.88)(-16.00,148.96)
\psline(-17.50,151.55)(-18.50,144.63)
\psline(-20.00,147.22)(-21.00,140.30)
\psline(-22.50,142.89)(-23.50,135.97)
\psline(-25.00,138.56)(-26.00,131.64)
\psline(-27.50,134.23)(-28.50,127.31)
\psline(-6.00,166.28)(-7.00,159.35)
\psline(-8.50,161.95)(-9.50,155.02)
\psline(-11.00,157.62)(-12.00,150.69)
\psline(-13.50,153.29)(-14.50,146.36)
\psline(-16.00,148.96)(-17.00,142.03)
\psline(-18.50,144.63)(-19.50,137.70)
\psline(-21.00,140.30)(-22.00,133.37)
\psline(-23.50,135.97)(-24.50,129.04)
\psline(-26.00,131.64)(-27.00,124.71)
\psline(-4.50,163.68)(-5.50,156.75)
\psline(-7.00,159.35)(-8.00,152.42)
\psline(-9.50,155.02)(-10.50,148.09)
\psline(-12.00,150.69)(-13.00,143.76)
\psline(-14.50,146.36)(-15.50,139.43)
\psline(-17.00,142.03)(-18.00,135.10)
\psline(-19.50,137.70)(-20.50,130.77)
\psline(-22.00,133.37)(-23.00,126.44)
\psline(-24.50,129.04)(-25.50,122.11)
\psline(-3.00,161.08)(-4.00,154.15)
\psline(-5.50,156.75)(-6.50,149.82)
\psline(-8.00,152.42)(-9.00,145.49)
\psline(-10.50,148.09)(-11.50,141.16)
\psline(-13.00,143.76)(-14.00,136.83)
\psline(-15.50,139.43)(-16.50,132.50)
\psline(-18.00,135.10)(-19.00,128.17)
\psline(-20.50,130.77)(-21.50,123.84)
\psline(-23.00,126.44)(-24.00,119.51)
\psline(-1.50,158.48)(-2.50,151.55)
\psline(-4.00,154.15)(-5.00,147.22)
\psline(-6.50,149.82)(-7.50,142.89)
\psline(-9.00,145.49)(-10.00,138.56)
\psline(-11.50,141.16)(-12.50,134.23)
\psline(-14.00,136.83)(-15.00,129.90)
\psline(-16.50,132.50)(-17.50,125.57)
\psline(-19.00,128.17)(-20.00,121.24)
\psline(-21.50,123.84)(-22.50,116.91)
\psline(0.00,155.88)(-1.00,148.96)
\psline(-2.50,151.55)(-3.50,144.63)
\psline(-5.00,147.22)(-6.00,140.30)
\psline(-7.50,142.89)(-8.50,135.97)
\psline(-10.00,138.56)(-11.00,131.64)
\psline(-12.50,134.23)(-13.50,127.31)
\psline(-15.00,129.90)(-16.00,122.98)
\psline(-17.50,125.57)(-18.50,118.65)
\psline(-20.00,121.24)(-21.00,114.32)
\psline(1.50,153.29)(0.50,146.36)
\psline(-1.00,148.96)(-2.00,142.03)
\psline(-3.50,144.63)(-4.50,137.70)
\psline(-6.00,140.30)(-7.00,133.37)
\psline(-8.50,135.97)(-9.50,129.04)
\psline(-11.00,131.64)(-12.00,124.71)
\psline(-13.50,127.31)(-14.50,120.38)
\psline(-16.00,122.98)(-17.00,116.05)
\psline(-18.50,118.65)(-19.50,111.72)
\psline(3.00,150.69)(2.00,143.76)
\psline(0.50,146.36)(-0.50,139.43)
\psline(-2.00,142.03)(-3.00,135.10)
\psline(-4.50,137.70)(-5.50,130.77)
\psline(-7.00,133.37)(-8.00,126.44)
\psline(-9.50,129.04)(-10.50,122.11)
\psline(-12.00,124.71)(-13.00,117.78)
\psline(-14.50,120.38)(-15.50,113.45)
\psline(-17.00,116.05)(-18.00,109.12)
\psline(4.50,148.09)(3.50,141.16)
\psline(2.00,143.76)(1.00,136.83)
\psline(-0.50,139.43)(-1.50,132.50)
\psline(-3.00,135.10)(-4.00,128.17)
\psline(-5.50,130.77)(-6.50,123.84)
\psline(-8.00,126.44)(-9.00,119.51)
\psline(-10.50,122.11)(-11.50,115.18)
\psline(-13.00,117.78)(-14.00,110.85)
\psline(-15.50,113.45)(-16.50,106.52)
\psline(6.00,145.49)(5.00,138.56)
\psline(3.50,141.16)(2.50,134.23)
\psline(1.00,136.83)(0.00,129.90)
\psline(-1.50,132.50)(-2.50,125.57)
\psline(-4.00,128.17)(-5.00,121.24)
\psline(-6.50,123.84)(-7.50,116.91)
\psline(-9.00,119.51)(-10.00,112.58)
\psline(-11.50,115.18)(-12.50,108.25)
\psline(-14.00,110.85)(-15.00,103.92)
\psline(7.50,142.89)(6.50,135.97)
\psline(5.00,138.56)(4.00,131.64)
\psline(2.50,134.23)(1.50,127.31)
\psline(0.00,129.90)(-1.00,122.98)
\psline(-2.50,125.57)(-3.50,118.65)
\psline(-5.00,121.24)(-6.00,114.32)
\psline(-7.50,116.91)(-8.50,109.99)
\psline(-10.00,112.58)(-11.00,105.66)
\psline(-12.50,108.25)(-13.50,101.32)
\psline(9.00,140.30)(8.00,133.37)
\psline(6.50,135.97)(5.50,129.04)
\psline(4.00,131.64)(3.00,124.71)
\psline(1.50,127.31)(0.50,120.38)
\psline(-1.00,122.98)(-2.00,116.05)
\psline(-3.50,118.65)(-4.50,111.72)
\psline(-6.00,114.32)(-7.00,107.39)
\psline(-8.50,109.99)(-9.50,103.06)
\psline(-11.00,105.66)(-12.00,98.73)
\psline(10.50,137.70)(9.50,130.77)
\psline(8.00,133.37)(7.00,126.44)
\psline(5.50,129.04)(4.50,122.11)
\psline(3.00,124.71)(2.00,117.78)
\psline(0.50,120.38)(-0.50,113.45)
\psline(-2.00,116.05)(-3.00,109.12)
\psline(-4.50,111.72)(-5.50,104.79)
\psline(-7.00,107.39)(-8.00,100.46)
\psline(-9.50,103.06)(-10.50,96.13)
\psline(12.00,135.10)(11.00,128.17)
\psline(9.50,130.77)(8.50,123.84)
\psline(7.00,126.44)(6.00,119.51)
\psline(4.50,122.11)(3.50,115.18)
\psline(2.00,117.78)(1.00,110.85)
\psline(-0.50,113.45)(-1.50,106.52)
\psline(-3.00,109.12)(-4.00,102.19)
\psline(-5.50,104.79)(-6.50,97.86)
\psline(-8.00,100.46)(-9.00,93.53)
\psline(13.50,132.50)(12.50,125.57)
\psline(11.00,128.17)(10.00,121.24)
\psline(8.50,123.84)(7.50,116.91)
\psline(6.00,119.51)(5.00,112.58)
\psline(3.50,115.18)(2.50,108.25)
\psline(1.00,110.85)(0.00,103.92)
\psline(-1.50,106.52)(-2.50,99.59)
\psline(-4.00,102.19)(-5.00,95.26)
\psline(-6.50,97.86)(-7.50,90.93)
\psline(15.00,129.90)(14.00,122.98)
\psline(12.50,125.57)(11.50,118.65)
\psline(10.00,121.24)(9.00,114.32)
\psline(7.50,116.91)(6.50,109.99)
\psline(5.00,112.58)(4.00,105.66)
\psline(2.50,108.25)(1.50,101.32)
\psline(0.00,103.92)(-1.00,96.99)
\psline(-2.50,99.59)(-3.50,92.66)
\psline(-5.00,95.26)(-6.00,88.33)
\pspolygon[fillstyle=solid,linewidth=0.1pt,fillcolor=red](-51.00,88.33)(-67.50,116.91)(-45.00,155.88)(-28.50,127.31)
\psline(-67.50,116.91)(-45.00,155.88)
\psline(-66.00,114.32)(-43.50,153.29)
\psline(-64.50,111.72)(-42.00,150.69)
\psline(-63.00,109.12)(-40.50,148.09)
\psline(-61.50,106.52)(-39.00,145.49)
\psline(-60.00,103.92)(-37.50,142.89)
\psline(-58.50,101.32)(-36.00,140.30)
\psline(-57.00,98.73)(-34.50,137.70)
\psline(-55.50,96.13)(-33.00,135.10)
\psline(-54.00,93.53)(-31.50,132.50)
\psline(-52.50,90.93)(-30.00,129.90)
\psline(-51.00,88.33)(-28.50,127.31)
\psline(-28.50,127.31)(-45.00,155.88)
\psline(-31.00,122.98)(-47.50,151.55)
\psline(-33.50,118.65)(-50.00,147.22)
\psline(-36.00,114.32)(-52.50,142.89)
\psline(-38.50,109.99)(-55.00,138.56)
\psline(-41.00,105.66)(-57.50,134.23)
\psline(-43.50,101.32)(-60.00,129.90)
\psline(-46.00,96.99)(-62.50,125.57)
\psline(-48.50,92.66)(-65.00,121.24)
\psline(-51.00,88.33)(-67.50,116.91)
\psline(-45.00,155.88)(-46.00,148.96)
\psline(-47.50,151.55)(-48.50,144.63)
\psline(-50.00,147.22)(-51.00,140.30)
\psline(-52.50,142.89)(-53.50,135.97)
\psline(-55.00,138.56)(-56.00,131.64)
\psline(-57.50,134.23)(-58.50,127.31)
\psline(-60.00,129.90)(-61.00,122.98)
\psline(-62.50,125.57)(-63.50,118.65)
\psline(-65.00,121.24)(-66.00,114.32)
\psline(-43.50,153.29)(-44.50,146.36)
\psline(-46.00,148.96)(-47.00,142.03)
\psline(-48.50,144.63)(-49.50,137.70)
\psline(-51.00,140.30)(-52.00,133.37)
\psline(-53.50,135.97)(-54.50,129.04)
\psline(-56.00,131.64)(-57.00,124.71)
\psline(-58.50,127.31)(-59.50,120.38)
\psline(-61.00,122.98)(-62.00,116.05)
\psline(-63.50,118.65)(-64.50,111.72)
\psline(-42.00,150.69)(-43.00,143.76)
\psline(-44.50,146.36)(-45.50,139.43)
\psline(-47.00,142.03)(-48.00,135.10)
\psline(-49.50,137.70)(-50.50,130.77)
\psline(-52.00,133.37)(-53.00,126.44)
\psline(-54.50,129.04)(-55.50,122.11)
\psline(-57.00,124.71)(-58.00,117.78)
\psline(-59.50,120.38)(-60.50,113.45)
\psline(-62.00,116.05)(-63.00,109.12)
\psline(-40.50,148.09)(-41.50,141.16)
\psline(-43.00,143.76)(-44.00,136.83)
\psline(-45.50,139.43)(-46.50,132.50)
\psline(-48.00,135.10)(-49.00,128.17)
\psline(-50.50,130.77)(-51.50,123.84)
\psline(-53.00,126.44)(-54.00,119.51)
\psline(-55.50,122.11)(-56.50,115.18)
\psline(-58.00,117.78)(-59.00,110.85)
\psline(-60.50,113.45)(-61.50,106.52)
\psline(-39.00,145.49)(-40.00,138.56)
\psline(-41.50,141.16)(-42.50,134.23)
\psline(-44.00,136.83)(-45.00,129.90)
\psline(-46.50,132.50)(-47.50,125.57)
\psline(-49.00,128.17)(-50.00,121.24)
\psline(-51.50,123.84)(-52.50,116.91)
\psline(-54.00,119.51)(-55.00,112.58)
\psline(-56.50,115.18)(-57.50,108.25)
\psline(-59.00,110.85)(-60.00,103.92)
\psline(-37.50,142.89)(-38.50,135.97)
\psline(-40.00,138.56)(-41.00,131.64)
\psline(-42.50,134.23)(-43.50,127.31)
\psline(-45.00,129.90)(-46.00,122.98)
\psline(-47.50,125.57)(-48.50,118.65)
\psline(-50.00,121.24)(-51.00,114.32)
\psline(-52.50,116.91)(-53.50,109.99)
\psline(-55.00,112.58)(-56.00,105.66)
\psline(-57.50,108.25)(-58.50,101.32)
\psline(-36.00,140.30)(-37.00,133.37)
\psline(-38.50,135.97)(-39.50,129.04)
\psline(-41.00,131.64)(-42.00,124.71)
\psline(-43.50,127.31)(-44.50,120.38)
\psline(-46.00,122.98)(-47.00,116.05)
\psline(-48.50,118.65)(-49.50,111.72)
\psline(-51.00,114.32)(-52.00,107.39)
\psline(-53.50,109.99)(-54.50,103.06)
\psline(-56.00,105.66)(-57.00,98.73)
\psline(-34.50,137.70)(-35.50,130.77)
\psline(-37.00,133.37)(-38.00,126.44)
\psline(-39.50,129.04)(-40.50,122.11)
\psline(-42.00,124.71)(-43.00,117.78)
\psline(-44.50,120.38)(-45.50,113.45)
\psline(-47.00,116.05)(-48.00,109.12)
\psline(-49.50,111.72)(-50.50,104.79)
\psline(-52.00,107.39)(-53.00,100.46)
\psline(-54.50,103.06)(-55.50,96.13)
\psline(-33.00,135.10)(-34.00,128.17)
\psline(-35.50,130.77)(-36.50,123.84)
\psline(-38.00,126.44)(-39.00,119.51)
\psline(-40.50,122.11)(-41.50,115.18)
\psline(-43.00,117.78)(-44.00,110.85)
\psline(-45.50,113.45)(-46.50,106.52)
\psline(-48.00,109.12)(-49.00,102.19)
\psline(-50.50,104.79)(-51.50,97.86)
\psline(-53.00,100.46)(-54.00,93.53)
\psline(-31.50,132.50)(-32.50,125.57)
\psline(-34.00,128.17)(-35.00,121.24)
\psline(-36.50,123.84)(-37.50,116.91)
\psline(-39.00,119.51)(-40.00,112.58)
\psline(-41.50,115.18)(-42.50,108.25)
\psline(-44.00,110.85)(-45.00,103.92)
\psline(-46.50,106.52)(-47.50,99.59)
\psline(-49.00,102.19)(-50.00,95.26)
\psline(-51.50,97.86)(-52.50,90.93)
\psline(-30.00,129.90)(-31.00,122.98)
\psline(-32.50,125.57)(-33.50,118.65)
\psline(-35.00,121.24)(-36.00,114.32)
\psline(-37.50,116.91)(-38.50,109.99)
\psline(-40.00,112.58)(-41.00,105.66)
\psline(-42.50,108.25)(-43.50,101.32)
\psline(-45.00,103.92)(-46.00,96.99)
\psline(-47.50,99.59)(-48.50,92.66)
\psline(-50.00,95.26)(-51.00,88.33)
\pspolygon[fillstyle=solid,linewidth=0.1pt,fillcolor=red](39.00,88.33)(-22.50,194.86)(0.00,233.83)(61.50,127.31)
\psline(-22.50,194.86)(0.00,233.83)
\psline(-21.00,192.26)(1.50,231.23)
\psline(-19.50,189.66)(3.00,228.63)
\psline(-18.00,187.06)(4.50,226.03)
\psline(-16.50,184.46)(6.00,223.43)
\psline(-15.00,181.87)(7.50,220.84)
\psline(-13.50,179.27)(9.00,218.24)
\psline(-12.00,176.67)(10.50,215.64)
\psline(-10.50,174.07)(12.00,213.04)
\psline(-9.00,171.47)(13.50,210.44)
\psline(-7.50,168.87)(15.00,207.85)
\psline(-6.00,166.28)(16.50,205.25)
\psline(-4.50,163.68)(18.00,202.65)
\psline(-3.00,161.08)(19.50,200.05)
\psline(-1.50,158.48)(21.00,197.45)
\psline(0.00,155.88)(22.50,194.86)
\psline(1.50,153.29)(24.00,192.26)
\psline(3.00,150.69)(25.50,189.66)
\psline(4.50,148.09)(27.00,187.06)
\psline(6.00,145.49)(28.50,184.46)
\psline(7.50,142.89)(30.00,181.87)
\psline(9.00,140.30)(31.50,179.27)
\psline(10.50,137.70)(33.00,176.67)
\psline(12.00,135.10)(34.50,174.07)
\psline(13.50,132.50)(36.00,171.47)
\psline(15.00,129.90)(37.50,168.87)
\psline(16.50,127.31)(39.00,166.28)
\psline(18.00,124.71)(40.50,163.68)
\psline(19.50,122.11)(42.00,161.08)
\psline(21.00,119.51)(43.50,158.48)
\psline(22.50,116.91)(45.00,155.88)
\psline(24.00,114.32)(46.50,153.29)
\psline(25.50,111.72)(48.00,150.69)
\psline(27.00,109.12)(49.50,148.09)
\psline(28.50,106.52)(51.00,145.49)
\psline(30.00,103.92)(52.50,142.89)
\psline(31.50,101.32)(54.00,140.30)
\psline(33.00,98.73)(55.50,137.70)
\psline(34.50,96.13)(57.00,135.10)
\psline(36.00,93.53)(58.50,132.50)
\psline(37.50,90.93)(60.00,129.90)
\psline(39.00,88.33)(61.50,127.31)
\psline(61.50,127.31)(0.00,233.83)
\psline(59.00,122.98)(-2.50,229.50)
\psline(56.50,118.65)(-5.00,225.17)
\psline(54.00,114.32)(-7.50,220.84)
\psline(51.50,109.99)(-10.00,216.51)
\psline(49.00,105.66)(-12.50,212.18)
\psline(46.50,101.32)(-15.00,207.85)
\psline(44.00,96.99)(-17.50,203.52)
\psline(41.50,92.66)(-20.00,199.19)
\psline(39.00,88.33)(-22.50,194.86)
\psline(0.00,233.83)(-1.00,226.90)
\psline(-2.50,229.50)(-3.50,222.57)
\psline(-5.00,225.17)(-6.00,218.24)
\psline(-7.50,220.84)(-8.50,213.91)
\psline(-10.00,216.51)(-11.00,209.58)
\psline(-12.50,212.18)(-13.50,205.25)
\psline(-15.00,207.85)(-16.00,200.92)
\psline(-17.50,203.52)(-18.50,196.59)
\psline(-20.00,199.19)(-21.00,192.26)
\psline(1.50,231.23)(0.50,224.30)
\psline(-1.00,226.90)(-2.00,219.97)
\psline(-3.50,222.57)(-4.50,215.64)
\psline(-6.00,218.24)(-7.00,211.31)
\psline(-8.50,213.91)(-9.50,206.98)
\psline(-11.00,209.58)(-12.00,202.65)
\psline(-13.50,205.25)(-14.50,198.32)
\psline(-16.00,200.92)(-17.00,193.99)
\psline(-18.50,196.59)(-19.50,189.66)
\psline(3.00,228.63)(2.00,221.70)
\psline(0.50,224.30)(-0.50,217.37)
\psline(-2.00,219.97)(-3.00,213.04)
\psline(-4.50,215.64)(-5.50,208.71)
\psline(-7.00,211.31)(-8.00,204.38)
\psline(-9.50,206.98)(-10.50,200.05)
\psline(-12.00,202.65)(-13.00,195.72)
\psline(-14.50,198.32)(-15.50,191.39)
\psline(-17.00,193.99)(-18.00,187.06)
\psline(4.50,226.03)(3.50,219.10)
\psline(2.00,221.70)(1.00,214.77)
\psline(-0.50,217.37)(-1.50,210.44)
\psline(-3.00,213.04)(-4.00,206.11)
\psline(-5.50,208.71)(-6.50,201.78)
\psline(-8.00,204.38)(-9.00,197.45)
\psline(-10.50,200.05)(-11.50,193.12)
\psline(-13.00,195.72)(-14.00,188.79)
\psline(-15.50,191.39)(-16.50,184.46)
\psline(6.00,223.43)(5.00,216.51)
\psline(3.50,219.10)(2.50,212.18)
\psline(1.00,214.77)(0.00,207.85)
\psline(-1.50,210.44)(-2.50,203.52)
\psline(-4.00,206.11)(-5.00,199.19)
\psline(-6.50,201.78)(-7.50,194.86)
\psline(-9.00,197.45)(-10.00,190.53)
\psline(-11.50,193.12)(-12.50,186.20)
\psline(-14.00,188.79)(-15.00,181.87)
\psline(7.50,220.84)(6.50,213.91)
\psline(5.00,216.51)(4.00,209.58)
\psline(2.50,212.18)(1.50,205.25)
\psline(0.00,207.85)(-1.00,200.92)
\psline(-2.50,203.52)(-3.50,196.59)
\psline(-5.00,199.19)(-6.00,192.26)
\psline(-7.50,194.86)(-8.50,187.93)
\psline(-10.00,190.53)(-11.00,183.60)
\psline(-12.50,186.20)(-13.50,179.27)
\psline(9.00,218.24)(8.00,211.31)
\psline(6.50,213.91)(5.50,206.98)
\psline(4.00,209.58)(3.00,202.65)
\psline(1.50,205.25)(0.50,198.32)
\psline(-1.00,200.92)(-2.00,193.99)
\psline(-3.50,196.59)(-4.50,189.66)
\psline(-6.00,192.26)(-7.00,185.33)
\psline(-8.50,187.93)(-9.50,181.00)
\psline(-11.00,183.60)(-12.00,176.67)
\psline(10.50,215.64)(9.50,208.71)
\psline(8.00,211.31)(7.00,204.38)
\psline(5.50,206.98)(4.50,200.05)
\psline(3.00,202.65)(2.00,195.72)
\psline(0.50,198.32)(-0.50,191.39)
\psline(-2.00,193.99)(-3.00,187.06)
\psline(-4.50,189.66)(-5.50,182.73)
\psline(-7.00,185.33)(-8.00,178.40)
\psline(-9.50,181.00)(-10.50,174.07)
\psline(12.00,213.04)(11.00,206.11)
\psline(9.50,208.71)(8.50,201.78)
\psline(7.00,204.38)(6.00,197.45)
\psline(4.50,200.05)(3.50,193.12)
\psline(2.00,195.72)(1.00,188.79)
\psline(-0.50,191.39)(-1.50,184.46)
\psline(-3.00,187.06)(-4.00,180.13)
\psline(-5.50,182.73)(-6.50,175.80)
\psline(-8.00,178.40)(-9.00,171.47)
\psline(13.50,210.44)(12.50,203.52)
\psline(11.00,206.11)(10.00,199.19)
\psline(8.50,201.78)(7.50,194.86)
\psline(6.00,197.45)(5.00,190.53)
\psline(3.50,193.12)(2.50,186.20)
\psline(1.00,188.79)(0.00,181.87)
\psline(-1.50,184.46)(-2.50,177.54)
\psline(-4.00,180.13)(-5.00,173.21)
\psline(-6.50,175.80)(-7.50,168.87)
\psline(15.00,207.85)(14.00,200.92)
\psline(12.50,203.52)(11.50,196.59)
\psline(10.00,199.19)(9.00,192.26)
\psline(7.50,194.86)(6.50,187.93)
\psline(5.00,190.53)(4.00,183.60)
\psline(2.50,186.20)(1.50,179.27)
\psline(0.00,181.87)(-1.00,174.94)
\psline(-2.50,177.54)(-3.50,170.61)
\psline(-5.00,173.21)(-6.00,166.28)
\psline(16.50,205.25)(15.50,198.32)
\psline(14.00,200.92)(13.00,193.99)
\psline(11.50,196.59)(10.50,189.66)
\psline(9.00,192.26)(8.00,185.33)
\psline(6.50,187.93)(5.50,181.00)
\psline(4.00,183.60)(3.00,176.67)
\psline(1.50,179.27)(0.50,172.34)
\psline(-1.00,174.94)(-2.00,168.01)
\psline(-3.50,170.61)(-4.50,163.68)
\psline(18.00,202.65)(17.00,195.72)
\psline(15.50,198.32)(14.50,191.39)
\psline(13.00,193.99)(12.00,187.06)
\psline(10.50,189.66)(9.50,182.73)
\psline(8.00,185.33)(7.00,178.40)
\psline(5.50,181.00)(4.50,174.07)
\psline(3.00,176.67)(2.00,169.74)
\psline(0.50,172.34)(-0.50,165.41)
\psline(-2.00,168.01)(-3.00,161.08)
\psline(19.50,200.05)(18.50,193.12)
\psline(17.00,195.72)(16.00,188.79)
\psline(14.50,191.39)(13.50,184.46)
\psline(12.00,187.06)(11.00,180.13)
\psline(9.50,182.73)(8.50,175.80)
\psline(7.00,178.40)(6.00,171.47)
\psline(4.50,174.07)(3.50,167.14)
\psline(2.00,169.74)(1.00,162.81)
\psline(-0.50,165.41)(-1.50,158.48)
\psline(21.00,197.45)(20.00,190.53)
\psline(18.50,193.12)(17.50,186.20)
\psline(16.00,188.79)(15.00,181.87)
\psline(13.50,184.46)(12.50,177.54)
\psline(11.00,180.13)(10.00,173.21)
\psline(8.50,175.80)(7.50,168.87)
\psline(6.00,171.47)(5.00,164.54)
\psline(3.50,167.14)(2.50,160.21)
\psline(1.00,162.81)(0.00,155.88)
\psline(22.50,194.86)(21.50,187.93)
\psline(20.00,190.53)(19.00,183.60)
\psline(17.50,186.20)(16.50,179.27)
\psline(15.00,181.87)(14.00,174.94)
\psline(12.50,177.54)(11.50,170.61)
\psline(10.00,173.21)(9.00,166.28)
\psline(7.50,168.87)(6.50,161.95)
\psline(5.00,164.54)(4.00,157.62)
\psline(2.50,160.21)(1.50,153.29)
\psline(24.00,192.26)(23.00,185.33)
\psline(21.50,187.93)(20.50,181.00)
\psline(19.00,183.60)(18.00,176.67)
\psline(16.50,179.27)(15.50,172.34)
\psline(14.00,174.94)(13.00,168.01)
\psline(11.50,170.61)(10.50,163.68)
\psline(9.00,166.28)(8.00,159.35)
\psline(6.50,161.95)(5.50,155.02)
\psline(4.00,157.62)(3.00,150.69)
\psline(25.50,189.66)(24.50,182.73)
\psline(23.00,185.33)(22.00,178.40)
\psline(20.50,181.00)(19.50,174.07)
\psline(18.00,176.67)(17.00,169.74)
\psline(15.50,172.34)(14.50,165.41)
\psline(13.00,168.01)(12.00,161.08)
\psline(10.50,163.68)(9.50,156.75)
\psline(8.00,159.35)(7.00,152.42)
\psline(5.50,155.02)(4.50,148.09)
\psline(27.00,187.06)(26.00,180.13)
\psline(24.50,182.73)(23.50,175.80)
\psline(22.00,178.40)(21.00,171.47)
\psline(19.50,174.07)(18.50,167.14)
\psline(17.00,169.74)(16.00,162.81)
\psline(14.50,165.41)(13.50,158.48)
\psline(12.00,161.08)(11.00,154.15)
\psline(9.50,156.75)(8.50,149.82)
\psline(7.00,152.42)(6.00,145.49)
\psline(28.50,184.46)(27.50,177.54)
\psline(26.00,180.13)(25.00,173.21)
\psline(23.50,175.80)(22.50,168.87)
\psline(21.00,171.47)(20.00,164.54)
\psline(18.50,167.14)(17.50,160.21)
\psline(16.00,162.81)(15.00,155.88)
\psline(13.50,158.48)(12.50,151.55)
\psline(11.00,154.15)(10.00,147.22)
\psline(8.50,149.82)(7.50,142.89)
\psline(30.00,181.87)(29.00,174.94)
\psline(27.50,177.54)(26.50,170.61)
\psline(25.00,173.21)(24.00,166.28)
\psline(22.50,168.87)(21.50,161.95)
\psline(20.00,164.54)(19.00,157.62)
\psline(17.50,160.21)(16.50,153.29)
\psline(15.00,155.88)(14.00,148.96)
\psline(12.50,151.55)(11.50,144.63)
\psline(10.00,147.22)(9.00,140.30)
\psline(31.50,179.27)(30.50,172.34)
\psline(29.00,174.94)(28.00,168.01)
\psline(26.50,170.61)(25.50,163.68)
\psline(24.00,166.28)(23.00,159.35)
\psline(21.50,161.95)(20.50,155.02)
\psline(19.00,157.62)(18.00,150.69)
\psline(16.50,153.29)(15.50,146.36)
\psline(14.00,148.96)(13.00,142.03)
\psline(11.50,144.63)(10.50,137.70)
\psline(33.00,176.67)(32.00,169.74)
\psline(30.50,172.34)(29.50,165.41)
\psline(28.00,168.01)(27.00,161.08)
\psline(25.50,163.68)(24.50,156.75)
\psline(23.00,159.35)(22.00,152.42)
\psline(20.50,155.02)(19.50,148.09)
\psline(18.00,150.69)(17.00,143.76)
\psline(15.50,146.36)(14.50,139.43)
\psline(13.00,142.03)(12.00,135.10)
\psline(34.50,174.07)(33.50,167.14)
\psline(32.00,169.74)(31.00,162.81)
\psline(29.50,165.41)(28.50,158.48)
\psline(27.00,161.08)(26.00,154.15)
\psline(24.50,156.75)(23.50,149.82)
\psline(22.00,152.42)(21.00,145.49)
\psline(19.50,148.09)(18.50,141.16)
\psline(17.00,143.76)(16.00,136.83)
\psline(14.50,139.43)(13.50,132.50)
\psline(36.00,171.47)(35.00,164.54)
\psline(33.50,167.14)(32.50,160.21)
\psline(31.00,162.81)(30.00,155.88)
\psline(28.50,158.48)(27.50,151.55)
\psline(26.00,154.15)(25.00,147.22)
\psline(23.50,149.82)(22.50,142.89)
\psline(21.00,145.49)(20.00,138.56)
\psline(18.50,141.16)(17.50,134.23)
\psline(16.00,136.83)(15.00,129.90)
\psline(37.50,168.87)(36.50,161.95)
\psline(35.00,164.54)(34.00,157.62)
\psline(32.50,160.21)(31.50,153.29)
\psline(30.00,155.88)(29.00,148.96)
\psline(27.50,151.55)(26.50,144.63)
\psline(25.00,147.22)(24.00,140.30)
\psline(22.50,142.89)(21.50,135.97)
\psline(20.00,138.56)(19.00,131.64)
\psline(17.50,134.23)(16.50,127.31)
\psline(39.00,166.28)(38.00,159.35)
\psline(36.50,161.95)(35.50,155.02)
\psline(34.00,157.62)(33.00,150.69)
\psline(31.50,153.29)(30.50,146.36)
\psline(29.00,148.96)(28.00,142.03)
\psline(26.50,144.63)(25.50,137.70)
\psline(24.00,140.30)(23.00,133.37)
\psline(21.50,135.97)(20.50,129.04)
\psline(19.00,131.64)(18.00,124.71)
\psline(40.50,163.68)(39.50,156.75)
\psline(38.00,159.35)(37.00,152.42)
\psline(35.50,155.02)(34.50,148.09)
\psline(33.00,150.69)(32.00,143.76)
\psline(30.50,146.36)(29.50,139.43)
\psline(28.00,142.03)(27.00,135.10)
\psline(25.50,137.70)(24.50,130.77)
\psline(23.00,133.37)(22.00,126.44)
\psline(20.50,129.04)(19.50,122.11)
\psline(42.00,161.08)(41.00,154.15)
\psline(39.50,156.75)(38.50,149.82)
\psline(37.00,152.42)(36.00,145.49)
\psline(34.50,148.09)(33.50,141.16)
\psline(32.00,143.76)(31.00,136.83)
\psline(29.50,139.43)(28.50,132.50)
\psline(27.00,135.10)(26.00,128.17)
\psline(24.50,130.77)(23.50,123.84)
\psline(22.00,126.44)(21.00,119.51)
\psline(43.50,158.48)(42.50,151.55)
\psline(41.00,154.15)(40.00,147.22)
\psline(38.50,149.82)(37.50,142.89)
\psline(36.00,145.49)(35.00,138.56)
\psline(33.50,141.16)(32.50,134.23)
\psline(31.00,136.83)(30.00,129.90)
\psline(28.50,132.50)(27.50,125.57)
\psline(26.00,128.17)(25.00,121.24)
\psline(23.50,123.84)(22.50,116.91)
\psline(45.00,155.88)(44.00,148.96)
\psline(42.50,151.55)(41.50,144.63)
\psline(40.00,147.22)(39.00,140.30)
\psline(37.50,142.89)(36.50,135.97)
\psline(35.00,138.56)(34.00,131.64)
\psline(32.50,134.23)(31.50,127.31)
\psline(30.00,129.90)(29.00,122.98)
\psline(27.50,125.57)(26.50,118.65)
\psline(25.00,121.24)(24.00,114.32)
\psline(46.50,153.29)(45.50,146.36)
\psline(44.00,148.96)(43.00,142.03)
\psline(41.50,144.63)(40.50,137.70)
\psline(39.00,140.30)(38.00,133.37)
\psline(36.50,135.97)(35.50,129.04)
\psline(34.00,131.64)(33.00,124.71)
\psline(31.50,127.31)(30.50,120.38)
\psline(29.00,122.98)(28.00,116.05)
\psline(26.50,118.65)(25.50,111.72)
\psline(48.00,150.69)(47.00,143.76)
\psline(45.50,146.36)(44.50,139.43)
\psline(43.00,142.03)(42.00,135.10)
\psline(40.50,137.70)(39.50,130.77)
\psline(38.00,133.37)(37.00,126.44)
\psline(35.50,129.04)(34.50,122.11)
\psline(33.00,124.71)(32.00,117.78)
\psline(30.50,120.38)(29.50,113.45)
\psline(28.00,116.05)(27.00,109.12)
\psline(49.50,148.09)(48.50,141.16)
\psline(47.00,143.76)(46.00,136.83)
\psline(44.50,139.43)(43.50,132.50)
\psline(42.00,135.10)(41.00,128.17)
\psline(39.50,130.77)(38.50,123.84)
\psline(37.00,126.44)(36.00,119.51)
\psline(34.50,122.11)(33.50,115.18)
\psline(32.00,117.78)(31.00,110.85)
\psline(29.50,113.45)(28.50,106.52)
\psline(51.00,145.49)(50.00,138.56)
\psline(48.50,141.16)(47.50,134.23)
\psline(46.00,136.83)(45.00,129.90)
\psline(43.50,132.50)(42.50,125.57)
\psline(41.00,128.17)(40.00,121.24)
\psline(38.50,123.84)(37.50,116.91)
\psline(36.00,119.51)(35.00,112.58)
\psline(33.50,115.18)(32.50,108.25)
\psline(31.00,110.85)(30.00,103.92)
\psline(52.50,142.89)(51.50,135.97)
\psline(50.00,138.56)(49.00,131.64)
\psline(47.50,134.23)(46.50,127.31)
\psline(45.00,129.90)(44.00,122.98)
\psline(42.50,125.57)(41.50,118.65)
\psline(40.00,121.24)(39.00,114.32)
\psline(37.50,116.91)(36.50,109.99)
\psline(35.00,112.58)(34.00,105.66)
\psline(32.50,108.25)(31.50,101.32)
\psline(54.00,140.30)(53.00,133.37)
\psline(51.50,135.97)(50.50,129.04)
\psline(49.00,131.64)(48.00,124.71)
\psline(46.50,127.31)(45.50,120.38)
\psline(44.00,122.98)(43.00,116.05)
\psline(41.50,118.65)(40.50,111.72)
\psline(39.00,114.32)(38.00,107.39)
\psline(36.50,109.99)(35.50,103.06)
\psline(34.00,105.66)(33.00,98.73)
\psline(55.50,137.70)(54.50,130.77)
\psline(53.00,133.37)(52.00,126.44)
\psline(50.50,129.04)(49.50,122.11)
\psline(48.00,124.71)(47.00,117.78)
\psline(45.50,120.38)(44.50,113.45)
\psline(43.00,116.05)(42.00,109.12)
\psline(40.50,111.72)(39.50,104.79)
\psline(38.00,107.39)(37.00,100.46)
\psline(35.50,103.06)(34.50,96.13)
\psline(57.00,135.10)(56.00,128.17)
\psline(54.50,130.77)(53.50,123.84)
\psline(52.00,126.44)(51.00,119.51)
\psline(49.50,122.11)(48.50,115.18)
\psline(47.00,117.78)(46.00,110.85)
\psline(44.50,113.45)(43.50,106.52)
\psline(42.00,109.12)(41.00,102.19)
\psline(39.50,104.79)(38.50,97.86)
\psline(37.00,100.46)(36.00,93.53)
\psline(58.50,132.50)(57.50,125.57)
\psline(56.00,128.17)(55.00,121.24)
\psline(53.50,123.84)(52.50,116.91)
\psline(51.00,119.51)(50.00,112.58)
\psline(48.50,115.18)(47.50,108.25)
\psline(46.00,110.85)(45.00,103.92)
\psline(43.50,106.52)(42.50,99.59)
\psline(41.00,102.19)(40.00,95.26)
\psline(38.50,97.86)(37.50,90.93)
\psline(60.00,129.90)(59.00,122.98)
\psline(57.50,125.57)(56.50,118.65)
\psline(55.00,121.24)(54.00,114.32)
\psline(52.50,116.91)(51.50,109.99)
\psline(50.00,112.58)(49.00,105.66)
\psline(47.50,108.25)(46.50,101.32)
\psline(45.00,103.92)(44.00,96.99)
\psline(42.50,99.59)(41.50,92.66)
\psline(40.00,95.26)(39.00,88.33)
\pspolygon[fillstyle=solid,linewidth=0.1pt,fillcolor=lightblue](-22.50,38.97)(-37.50,-64.95)(-96.00,-88.33)
\psline(-37.50,-64.95)(-96.00,-88.33)
\psline(-36.79,-60.00)(-92.50,-82.27)
\psline(-36.07,-55.05)(-89.00,-76.21)
\psline(-35.36,-50.11)(-85.50,-70.15)
\psline(-34.64,-45.16)(-82.00,-64.09)
\psline(-33.93,-40.21)(-78.50,-58.02)
\psline(-33.21,-35.26)(-75.00,-51.96)
\psline(-32.50,-30.31)(-71.50,-45.90)
\psline(-31.79,-25.36)(-68.00,-39.84)
\psline(-31.07,-20.41)(-64.50,-33.77)
\psline(-30.36,-15.46)(-61.00,-27.71)
\psline(-29.64,-10.52)(-57.50,-21.65)
\psline(-28.93,-5.57)(-54.00,-15.59)
\psline(-28.21,-0.62)(-50.50,-9.53)
\psline(-27.50,4.33)(-47.00,-3.46)
\psline(-26.79,9.28)(-43.50,2.60)
\psline(-26.07,14.23)(-40.00,8.66)
\psline(-25.36,19.18)(-36.50,14.72)
\psline(-24.64,24.12)(-33.00,20.78)
\psline(-23.93,29.07)(-29.50,26.85)
\psline(-23.21,34.02)(-26.00,32.91)
\psline(-22.50,38.97)(-22.50,38.97)
\psline(-96.00,-88.33)(-22.50,38.97)
\psline(-93.21,-87.22)(-23.21,34.02)
\psline(-90.43,-86.11)(-23.93,29.07)
\psline(-87.64,-84.99)(-24.64,24.12)
\psline(-84.86,-83.88)(-25.36,19.18)
\psline(-82.07,-82.77)(-26.07,14.23)
\psline(-79.29,-81.65)(-26.79,9.28)
\psline(-76.50,-80.54)(-27.50,4.33)
\psline(-73.71,-79.43)(-28.21,-0.62)
\psline(-70.93,-78.31)(-28.93,-5.57)
\psline(-68.14,-77.20)(-29.64,-10.52)
\psline(-65.36,-76.09)(-30.36,-15.46)
\psline(-62.57,-74.97)(-31.07,-20.41)
\psline(-59.79,-73.86)(-31.79,-25.36)
\psline(-57.00,-72.75)(-32.50,-30.31)
\psline(-54.21,-71.63)(-33.21,-35.26)
\psline(-51.43,-70.52)(-33.93,-40.21)
\psline(-48.64,-69.41)(-34.64,-45.16)
\psline(-45.86,-68.29)(-35.36,-50.11)
\psline(-43.07,-67.18)(-36.07,-55.05)
\psline(-40.29,-66.07)(-36.79,-60.00)
\psline(-37.50,-64.95)(-37.50,-64.95)
\psline(-22.50,38.97)(-37.50,-64.95)
\psline(-26.00,32.91)(-40.29,-66.07)
\psline(-29.50,26.85)(-43.07,-67.18)
\psline(-33.00,20.78)(-45.86,-68.29)
\psline(-36.50,14.72)(-48.64,-69.41)
\psline(-40.00,8.66)(-51.43,-70.52)
\psline(-43.50,2.60)(-54.21,-71.63)
\psline(-47.00,-3.46)(-57.00,-72.75)
\psline(-50.50,-9.53)(-59.79,-73.86)
\psline(-54.00,-15.59)(-62.57,-74.97)
\psline(-57.50,-21.65)(-65.36,-76.09)
\psline(-61.00,-27.71)(-68.14,-77.20)
\psline(-64.50,-33.77)(-70.93,-78.31)
\psline(-68.00,-39.84)(-73.71,-79.43)
\psline(-71.50,-45.90)(-76.50,-80.54)
\psline(-75.00,-51.96)(-79.29,-81.65)
\psline(-78.50,-58.02)(-82.07,-82.77)
\psline(-82.00,-64.09)(-84.86,-83.88)
\psline(-85.50,-70.15)(-87.64,-84.99)
\psline(-89.00,-76.21)(-90.43,-86.11)
\psline(-92.50,-82.27)(-93.21,-87.22)
\psline(-96.00,-88.33)(-96.00,-88.33)
\pspolygon[fillstyle=solid,linewidth=0.1pt,fillcolor=lightblue](-45.00,77.94)(-60.00,-25.98)(-118.50,-49.36)
\psline(-60.00,-25.98)(-118.50,-49.36)
\psline(-59.29,-21.03)(-115.00,-43.30)
\psline(-58.57,-16.08)(-111.50,-37.24)
\psline(-57.86,-11.13)(-108.00,-31.18)
\psline(-57.14,-6.19)(-104.50,-25.11)
\psline(-56.43,-1.24)(-101.00,-19.05)
\psline(-55.71,3.71)(-97.50,-12.99)
\psline(-55.00,8.66)(-94.00,-6.93)
\psline(-54.29,13.61)(-90.50,-0.87)
\psline(-53.57,18.56)(-87.00,5.20)
\psline(-52.86,23.51)(-83.50,11.26)
\psline(-52.14,28.46)(-80.00,17.32)
\psline(-51.43,33.40)(-76.50,23.38)
\psline(-50.71,38.35)(-73.00,29.44)
\psline(-50.00,43.30)(-69.50,35.51)
\psline(-49.29,48.25)(-66.00,41.57)
\psline(-48.57,53.20)(-62.50,47.63)
\psline(-47.86,58.15)(-59.00,53.69)
\psline(-47.14,63.10)(-55.50,59.76)
\psline(-46.43,68.04)(-52.00,65.82)
\psline(-45.71,72.99)(-48.50,71.88)
\psline(-45.00,77.94)(-45.00,77.94)
\psline(-118.50,-49.36)(-45.00,77.94)
\psline(-115.71,-48.25)(-45.71,72.99)
\psline(-112.93,-47.14)(-46.43,68.04)
\psline(-110.14,-46.02)(-47.14,63.10)
\psline(-107.36,-44.91)(-47.86,58.15)
\psline(-104.57,-43.80)(-48.57,53.20)
\psline(-101.79,-42.68)(-49.29,48.25)
\psline(-99.00,-41.57)(-50.00,43.30)
\psline(-96.21,-40.46)(-50.71,38.35)
\psline(-93.43,-39.34)(-51.43,33.40)
\psline(-90.64,-38.23)(-52.14,28.46)
\psline(-87.86,-37.12)(-52.86,23.51)
\psline(-85.07,-36.00)(-53.57,18.56)
\psline(-82.29,-34.89)(-54.29,13.61)
\psline(-79.50,-33.77)(-55.00,8.66)
\psline(-76.71,-32.66)(-55.71,3.71)
\psline(-73.93,-31.55)(-56.43,-1.24)
\psline(-71.14,-30.43)(-57.14,-6.19)
\psline(-68.36,-29.32)(-57.86,-11.13)
\psline(-65.57,-28.21)(-58.57,-16.08)
\psline(-62.79,-27.09)(-59.29,-21.03)
\psline(-60.00,-25.98)(-60.00,-25.98)
\psline(-45.00,77.94)(-60.00,-25.98)
\psline(-48.50,71.88)(-62.79,-27.09)
\psline(-52.00,65.82)(-65.57,-28.21)
\psline(-55.50,59.76)(-68.36,-29.32)
\psline(-59.00,53.69)(-71.14,-30.43)
\psline(-62.50,47.63)(-73.93,-31.55)
\psline(-66.00,41.57)(-76.71,-32.66)
\psline(-69.50,35.51)(-79.50,-33.77)
\psline(-73.00,29.44)(-82.29,-34.89)
\psline(-76.50,23.38)(-85.07,-36.00)
\psline(-80.00,17.32)(-87.86,-37.12)
\psline(-83.50,11.26)(-90.64,-38.23)
\psline(-87.00,5.20)(-93.43,-39.34)
\psline(-90.50,-0.87)(-96.21,-40.46)
\psline(-94.00,-6.93)(-99.00,-41.57)
\psline(-97.50,-12.99)(-101.79,-42.68)
\psline(-101.00,-19.05)(-104.57,-43.80)
\psline(-104.50,-25.11)(-107.36,-44.91)
\psline(-108.00,-31.18)(-110.14,-46.02)
\psline(-111.50,-37.24)(-112.93,-47.14)
\psline(-115.00,-43.30)(-115.71,-48.25)
\psline(-118.50,-49.36)(-118.50,-49.36)
\pspolygon[fillstyle=solid,linewidth=0.1pt,fillcolor=lightblue](-67.50,116.91)(-141.00,-10.39)(-82.50,12.99)
\psline(-141.00,-10.39)(-82.50,12.99)
\psline(-137.50,-4.33)(-81.79,17.94)
\psline(-134.00,1.73)(-81.07,22.89)
\psline(-130.50,7.79)(-80.36,27.84)
\psline(-127.00,13.86)(-79.64,32.79)
\psline(-123.50,19.92)(-78.93,37.73)
\psline(-120.00,25.98)(-78.21,42.68)
\psline(-116.50,32.04)(-77.50,47.63)
\psline(-113.00,38.11)(-76.79,52.58)
\psline(-109.50,44.17)(-76.07,57.53)
\psline(-106.00,50.23)(-75.36,62.48)
\psline(-102.50,56.29)(-74.64,67.43)
\psline(-99.00,62.35)(-73.93,72.37)
\psline(-95.50,68.42)(-73.21,77.32)
\psline(-92.00,74.48)(-72.50,82.27)
\psline(-88.50,80.54)(-71.79,87.22)
\psline(-85.00,86.60)(-71.07,92.17)
\psline(-81.50,92.66)(-70.36,97.12)
\psline(-78.00,98.73)(-69.64,102.07)
\psline(-74.50,104.79)(-68.93,107.02)
\psline(-71.00,110.85)(-68.21,111.96)
\psline(-67.50,116.91)(-67.50,116.91)
\psline(-82.50,12.99)(-67.50,116.91)
\psline(-85.29,11.88)(-71.00,110.85)
\psline(-88.07,10.76)(-74.50,104.79)
\psline(-90.86,9.65)(-78.00,98.73)
\psline(-93.64,8.54)(-81.50,92.66)
\psline(-96.43,7.42)(-85.00,86.60)
\psline(-99.21,6.31)(-88.50,80.54)
\psline(-102.00,5.20)(-92.00,74.48)
\psline(-104.79,4.08)(-95.50,68.42)
\psline(-107.57,2.97)(-99.00,62.35)
\psline(-110.36,1.86)(-102.50,56.29)
\psline(-113.14,0.74)(-106.00,50.23)
\psline(-115.93,-0.37)(-109.50,44.17)
\psline(-118.71,-1.48)(-113.00,38.11)
\psline(-121.50,-2.60)(-116.50,32.04)
\psline(-124.29,-3.71)(-120.00,25.98)
\psline(-127.07,-4.82)(-123.50,19.92)
\psline(-129.86,-5.94)(-127.00,13.86)
\psline(-132.64,-7.05)(-130.50,7.79)
\psline(-135.43,-8.17)(-134.00,1.73)
\psline(-138.21,-9.28)(-137.50,-4.33)
\psline(-141.00,-10.39)(-141.00,-10.39)
\psline(-67.50,116.91)(-141.00,-10.39)
\psline(-68.21,111.96)(-138.21,-9.28)
\psline(-68.93,107.02)(-135.43,-8.17)
\psline(-69.64,102.07)(-132.64,-7.05)
\psline(-70.36,97.12)(-129.86,-5.94)
\psline(-71.07,92.17)(-127.07,-4.82)
\psline(-71.79,87.22)(-124.29,-3.71)
\psline(-72.50,82.27)(-121.50,-2.60)
\psline(-73.21,77.32)(-118.71,-1.48)
\psline(-73.93,72.37)(-115.93,-0.37)
\psline(-74.64,67.43)(-113.14,0.74)
\psline(-75.36,62.48)(-110.36,1.86)
\psline(-76.07,57.53)(-107.57,2.97)
\psline(-76.79,52.58)(-104.79,4.08)
\psline(-77.50,47.63)(-102.00,5.20)
\psline(-78.21,42.68)(-99.21,6.31)
\psline(-78.93,37.73)(-96.43,7.42)
\psline(-79.64,32.79)(-93.64,8.54)
\psline(-80.36,27.84)(-90.86,9.65)
\psline(-81.07,22.89)(-88.07,10.76)
\psline(-81.79,17.94)(-85.29,11.88)
\psline(-82.50,12.99)(-82.50,12.99)
\pspolygon[fillstyle=solid,linewidth=0.1pt,fillcolor=lightyellow](0.00,-0.00)(-37.50,-64.95)(-22.50,38.97)
\psline(-37.50,-64.95)(-22.50,38.97)
\psline(-35.00,-60.62)(-21.00,36.37)
\psline(-32.50,-56.29)(-19.50,33.77)
\psline(-30.00,-51.96)(-18.00,31.18)
\psline(-27.50,-47.63)(-16.50,28.58)
\psline(-25.00,-43.30)(-15.00,25.98)
\psline(-22.50,-38.97)(-13.50,23.38)
\psline(-20.00,-34.64)(-12.00,20.78)
\psline(-17.50,-30.31)(-10.50,18.19)
\psline(-15.00,-25.98)(-9.00,15.59)
\psline(-12.50,-21.65)(-7.50,12.99)
\psline(-10.00,-17.32)(-6.00,10.39)
\psline(-7.50,-12.99)(-4.50,7.79)
\psline(-5.00,-8.66)(-3.00,5.20)
\psline(-2.50,-4.33)(-1.50,2.60)
\psline(0.00,-0.00)(0.00,0.00)
\psline(-22.50,38.97)(0.00,-0.00)
\psline(-23.50,32.04)(-2.50,-4.33)
\psline(-24.50,25.11)(-5.00,-8.66)
\psline(-25.50,18.19)(-7.50,-12.99)
\psline(-26.50,11.26)(-10.00,-17.32)
\psline(-27.50,4.33)(-12.50,-21.65)
\psline(-28.50,-2.60)(-15.00,-25.98)
\psline(-29.50,-9.53)(-17.50,-30.31)
\psline(-30.50,-16.45)(-20.00,-34.64)
\psline(-31.50,-23.38)(-22.50,-38.97)
\psline(-32.50,-30.31)(-25.00,-43.30)
\psline(-33.50,-37.24)(-27.50,-47.63)
\psline(-34.50,-44.17)(-30.00,-51.96)
\psline(-35.50,-51.10)(-32.50,-56.29)
\psline(-36.50,-58.02)(-35.00,-60.62)
\psline(-37.50,-64.95)(-37.50,-64.95)
\psline(0.00,0.00)(-37.50,-64.95)
\psline(-1.50,2.60)(-36.50,-58.02)
\psline(-3.00,5.20)(-35.50,-51.10)
\psline(-4.50,7.79)(-34.50,-44.17)
\psline(-6.00,10.39)(-33.50,-37.24)
\psline(-7.50,12.99)(-32.50,-30.31)
\psline(-9.00,15.59)(-31.50,-23.38)
\psline(-10.50,18.19)(-30.50,-16.45)
\psline(-12.00,20.78)(-29.50,-9.53)
\psline(-13.50,23.38)(-28.50,-2.60)
\psline(-15.00,25.98)(-27.50,4.33)
\psline(-16.50,28.58)(-26.50,11.26)
\psline(-18.00,31.18)(-25.50,18.19)
\psline(-19.50,33.77)(-24.50,25.11)
\psline(-21.00,36.37)(-23.50,32.04)
\psline(-22.50,38.97)(-22.50,38.97)
\pspolygon[fillstyle=solid,linewidth=0.1pt,fillcolor=lightyellow](-45.00,77.94)(-22.50,38.97)(-60.00,-25.98)
\psline(-22.50,38.97)(-60.00,-25.98)
\psline(-24.00,41.57)(-59.00,-19.05)
\psline(-25.50,44.17)(-58.00,-12.12)
\psline(-27.00,46.77)(-57.00,-5.20)
\psline(-28.50,49.36)(-56.00,1.73)
\psline(-30.00,51.96)(-55.00,8.66)
\psline(-31.50,54.56)(-54.00,15.59)
\psline(-33.00,57.16)(-53.00,22.52)
\psline(-34.50,59.76)(-52.00,29.44)
\psline(-36.00,62.35)(-51.00,36.37)
\psline(-37.50,64.95)(-50.00,43.30)
\psline(-39.00,67.55)(-49.00,50.23)
\psline(-40.50,70.15)(-48.00,57.16)
\psline(-42.00,72.75)(-47.00,64.09)
\psline(-43.50,75.34)(-46.00,71.01)
\psline(-45.00,77.94)(-45.00,77.94)
\psline(-60.00,-25.98)(-45.00,77.94)
\psline(-57.50,-21.65)(-43.50,75.34)
\psline(-55.00,-17.32)(-42.00,72.75)
\psline(-52.50,-12.99)(-40.50,70.15)
\psline(-50.00,-8.66)(-39.00,67.55)
\psline(-47.50,-4.33)(-37.50,64.95)
\psline(-45.00,-0.00)(-36.00,62.35)
\psline(-42.50,4.33)(-34.50,59.76)
\psline(-40.00,8.66)(-33.00,57.16)
\psline(-37.50,12.99)(-31.50,54.56)
\psline(-35.00,17.32)(-30.00,51.96)
\psline(-32.50,21.65)(-28.50,49.36)
\psline(-30.00,25.98)(-27.00,46.77)
\psline(-27.50,30.31)(-25.50,44.17)
\psline(-25.00,34.64)(-24.00,41.57)
\psline(-22.50,38.97)(-22.50,38.97)
\psline(-45.00,77.94)(-22.50,38.97)
\psline(-46.00,71.01)(-25.00,34.64)
\psline(-47.00,64.09)(-27.50,30.31)
\psline(-48.00,57.16)(-30.00,25.98)
\psline(-49.00,50.23)(-32.50,21.65)
\psline(-50.00,43.30)(-35.00,17.32)
\psline(-51.00,36.37)(-37.50,12.99)
\psline(-52.00,29.44)(-40.00,8.66)
\psline(-53.00,22.52)(-42.50,4.33)
\psline(-54.00,15.59)(-45.00,-0.00)
\psline(-55.00,8.66)(-47.50,-4.33)
\psline(-56.00,1.73)(-50.00,-8.66)
\psline(-57.00,-5.20)(-52.50,-12.99)
\psline(-58.00,-12.12)(-55.00,-17.32)
\psline(-59.00,-19.05)(-57.50,-21.65)
\psline(-60.00,-25.98)(-60.00,-25.98)
\pspolygon[fillstyle=solid,linewidth=0.1pt,fillcolor=lightyellow](-67.50,116.91)(-45.00,77.94)(-82.50,12.99)
\psline(-45.00,77.94)(-82.50,12.99)
\psline(-46.50,80.54)(-81.50,19.92)
\psline(-48.00,83.14)(-80.50,26.85)
\psline(-49.50,85.74)(-79.50,33.77)
\psline(-51.00,88.33)(-78.50,40.70)
\psline(-52.50,90.93)(-77.50,47.63)
\psline(-54.00,93.53)(-76.50,54.56)
\psline(-55.50,96.13)(-75.50,61.49)
\psline(-57.00,98.73)(-74.50,68.42)
\psline(-58.50,101.32)(-73.50,75.34)
\psline(-60.00,103.92)(-72.50,82.27)
\psline(-61.50,106.52)(-71.50,89.20)
\psline(-63.00,109.12)(-70.50,96.13)
\psline(-64.50,111.72)(-69.50,103.06)
\psline(-66.00,114.32)(-68.50,109.99)
\psline(-67.50,116.91)(-67.50,116.91)
\psline(-82.50,12.99)(-67.50,116.91)
\psline(-80.00,17.32)(-66.00,114.32)
\psline(-77.50,21.65)(-64.50,111.72)
\psline(-75.00,25.98)(-63.00,109.12)
\psline(-72.50,30.31)(-61.50,106.52)
\psline(-70.00,34.64)(-60.00,103.92)
\psline(-67.50,38.97)(-58.50,101.32)
\psline(-65.00,43.30)(-57.00,98.73)
\psline(-62.50,47.63)(-55.50,96.13)
\psline(-60.00,51.96)(-54.00,93.53)
\psline(-57.50,56.29)(-52.50,90.93)
\psline(-55.00,60.62)(-51.00,88.33)
\psline(-52.50,64.95)(-49.50,85.74)
\psline(-50.00,69.28)(-48.00,83.14)
\psline(-47.50,73.61)(-46.50,80.54)
\psline(-45.00,77.94)(-45.00,77.94)
\psline(-67.50,116.91)(-45.00,77.94)
\psline(-68.50,109.99)(-47.50,73.61)
\psline(-69.50,103.06)(-50.00,69.28)
\psline(-70.50,96.13)(-52.50,64.95)
\psline(-71.50,89.20)(-55.00,60.62)
\psline(-72.50,82.27)(-57.50,56.29)
\psline(-73.50,75.34)(-60.00,51.96)
\psline(-74.50,68.42)(-62.50,47.63)
\psline(-75.50,61.49)(-65.00,43.30)
\psline(-76.50,54.56)(-67.50,38.97)
\psline(-77.50,47.63)(-70.00,34.64)
\psline(-78.50,40.70)(-72.50,30.31)
\psline(-79.50,33.77)(-75.00,25.98)
\psline(-80.50,26.85)(-77.50,21.65)
\psline(-81.50,19.92)(-80.00,17.32)
\psline(-82.50,12.99)(-82.50,12.99)
\pspolygon[fillstyle=solid,linewidth=0.1pt,fillcolor=green](-37.50,-64.95)(-51.00,-88.33)(-96.00,-88.33)
\psline(-51.00,-88.33)(-96.00,-88.33)
\psline(-49.50,-85.74)(-89.50,-85.74)
\psline(-48.00,-83.14)(-83.00,-83.14)
\psline(-46.50,-80.54)(-76.50,-80.54)
\psline(-45.00,-77.94)(-70.00,-77.94)
\psline(-43.50,-75.34)(-63.50,-75.34)
\psline(-42.00,-72.75)(-57.00,-72.75)
\psline(-40.50,-70.15)(-50.50,-70.15)
\psline(-39.00,-67.55)(-44.00,-67.55)
\psline(-37.50,-64.95)(-37.50,-64.95)
\psline(-96.00,-88.33)(-37.50,-64.95)
\psline(-91.00,-88.33)(-39.00,-67.55)
\psline(-86.00,-88.33)(-40.50,-70.15)
\psline(-81.00,-88.33)(-42.00,-72.75)
\psline(-76.00,-88.33)(-43.50,-75.34)
\psline(-71.00,-88.33)(-45.00,-77.94)
\psline(-66.00,-88.33)(-46.50,-80.54)
\psline(-61.00,-88.33)(-48.00,-83.14)
\psline(-56.00,-88.33)(-49.50,-85.74)
\psline(-51.00,-88.33)(-51.00,-88.33)
\psline(-37.50,-64.95)(-51.00,-88.33)
\psline(-44.00,-67.55)(-56.00,-88.33)
\psline(-50.50,-70.15)(-61.00,-88.33)
\psline(-57.00,-72.75)(-66.00,-88.33)
\psline(-63.50,-75.34)(-71.00,-88.33)
\psline(-70.00,-77.94)(-76.00,-88.33)
\psline(-76.50,-80.54)(-81.00,-88.33)
\psline(-83.00,-83.14)(-86.00,-88.33)
\psline(-89.50,-85.74)(-91.00,-88.33)
\psline(-96.00,-88.33)(-96.00,-88.33)
\pspolygon[fillstyle=solid,linewidth=0.1pt,fillcolor=green](-60.00,-25.98)(-73.50,-49.36)(-118.50,-49.36)
\psline(-73.50,-49.36)(-118.50,-49.36)
\psline(-72.00,-46.77)(-112.00,-46.77)
\psline(-70.50,-44.17)(-105.50,-44.17)
\psline(-69.00,-41.57)(-99.00,-41.57)
\psline(-67.50,-38.97)(-92.50,-38.97)
\psline(-66.00,-36.37)(-86.00,-36.37)
\psline(-64.50,-33.77)(-79.50,-33.77)
\psline(-63.00,-31.18)(-73.00,-31.18)
\psline(-61.50,-28.58)(-66.50,-28.58)
\psline(-60.00,-25.98)(-60.00,-25.98)
\psline(-118.50,-49.36)(-60.00,-25.98)
\psline(-113.50,-49.36)(-61.50,-28.58)
\psline(-108.50,-49.36)(-63.00,-31.18)
\psline(-103.50,-49.36)(-64.50,-33.77)
\psline(-98.50,-49.36)(-66.00,-36.37)
\psline(-93.50,-49.36)(-67.50,-38.97)
\psline(-88.50,-49.36)(-69.00,-41.57)
\psline(-83.50,-49.36)(-70.50,-44.17)
\psline(-78.50,-49.36)(-72.00,-46.77)
\psline(-73.50,-49.36)(-73.50,-49.36)
\psline(-60.00,-25.98)(-73.50,-49.36)
\psline(-66.50,-28.58)(-78.50,-49.36)
\psline(-73.00,-31.18)(-83.50,-49.36)
\psline(-79.50,-33.77)(-88.50,-49.36)
\psline(-86.00,-36.37)(-93.50,-49.36)
\psline(-92.50,-38.97)(-98.50,-49.36)
\psline(-99.00,-41.57)(-103.50,-49.36)
\psline(-105.50,-44.17)(-108.50,-49.36)
\psline(-112.00,-46.77)(-113.50,-49.36)
\psline(-118.50,-49.36)(-118.50,-49.36)
\pspolygon[fillstyle=solid,linewidth=0.1pt,fillcolor=green](-82.50,12.99)(-96.00,-10.39)(-141.00,-10.39)
\psline(-96.00,-10.39)(-141.00,-10.39)
\psline(-94.50,-7.79)(-134.50,-7.79)
\psline(-93.00,-5.20)(-128.00,-5.20)
\psline(-91.50,-2.60)(-121.50,-2.60)
\psline(-90.00,-0.00)(-115.00,-0.00)
\psline(-88.50,2.60)(-108.50,2.60)
\psline(-87.00,5.20)(-102.00,5.20)
\psline(-85.50,7.79)(-95.50,7.79)
\psline(-84.00,10.39)(-89.00,10.39)
\psline(-82.50,12.99)(-82.50,12.99)
\psline(-141.00,-10.39)(-82.50,12.99)
\psline(-136.00,-10.39)(-84.00,10.39)
\psline(-131.00,-10.39)(-85.50,7.79)
\psline(-126.00,-10.39)(-87.00,5.20)
\psline(-121.00,-10.39)(-88.50,2.60)
\psline(-116.00,-10.39)(-90.00,-0.00)
\psline(-111.00,-10.39)(-91.50,-2.60)
\psline(-106.00,-10.39)(-93.00,-5.20)
\psline(-101.00,-10.39)(-94.50,-7.79)
\psline(-96.00,-10.39)(-96.00,-10.39)
\psline(-82.50,12.99)(-96.00,-10.39)
\psline(-89.00,10.39)(-101.00,-10.39)
\psline(-95.50,7.79)(-106.00,-10.39)
\psline(-102.00,5.20)(-111.00,-10.39)
\psline(-108.50,2.60)(-116.00,-10.39)
\psline(-115.00,-0.00)(-121.00,-10.39)
\psline(-121.50,-2.60)(-126.00,-10.39)
\psline(-128.00,-5.20)(-131.00,-10.39)
\psline(-134.50,-7.79)(-136.00,-10.39)
\psline(-141.00,-10.39)(-141.00,-10.39)
\pspolygon[fillstyle=solid,linewidth=0.1pt,fillcolor=red](-73.50,-49.36)(-112.50,-116.91)(-157.50,-116.91)(-118.50,-49.36)
\psline(-112.50,-116.91)(-157.50,-116.91)
\psline(-111.00,-114.32)(-156.00,-114.32)
\psline(-109.50,-111.72)(-154.50,-111.72)
\psline(-108.00,-109.12)(-153.00,-109.12)
\psline(-106.50,-106.52)(-151.50,-106.52)
\psline(-105.00,-103.92)(-150.00,-103.92)
\psline(-103.50,-101.32)(-148.50,-101.32)
\psline(-102.00,-98.73)(-147.00,-98.73)
\psline(-100.50,-96.13)(-145.50,-96.13)
\psline(-99.00,-93.53)(-144.00,-93.53)
\psline(-97.50,-90.93)(-142.50,-90.93)
\psline(-96.00,-88.33)(-141.00,-88.33)
\psline(-94.50,-85.74)(-139.50,-85.74)
\psline(-93.00,-83.14)(-138.00,-83.14)
\psline(-91.50,-80.54)(-136.50,-80.54)
\psline(-90.00,-77.94)(-135.00,-77.94)
\psline(-88.50,-75.34)(-133.50,-75.34)
\psline(-87.00,-72.75)(-132.00,-72.75)
\psline(-85.50,-70.15)(-130.50,-70.15)
\psline(-84.00,-67.55)(-129.00,-67.55)
\psline(-82.50,-64.95)(-127.50,-64.95)
\psline(-81.00,-62.35)(-126.00,-62.35)
\psline(-79.50,-59.76)(-124.50,-59.76)
\psline(-78.00,-57.16)(-123.00,-57.16)
\psline(-76.50,-54.56)(-121.50,-54.56)
\psline(-75.00,-51.96)(-120.00,-51.96)
\psline(-73.50,-49.36)(-118.50,-49.36)
\psline(-118.50,-49.36)(-157.50,-116.91)
\psline(-113.50,-49.36)(-152.50,-116.91)
\psline(-108.50,-49.36)(-147.50,-116.91)
\psline(-103.50,-49.36)(-142.50,-116.91)
\psline(-98.50,-49.36)(-137.50,-116.91)
\psline(-93.50,-49.36)(-132.50,-116.91)
\psline(-88.50,-49.36)(-127.50,-116.91)
\psline(-83.50,-49.36)(-122.50,-116.91)
\psline(-78.50,-49.36)(-117.50,-116.91)
\psline(-73.50,-49.36)(-112.50,-116.91)
\psline(-157.50,-116.91)(-151.00,-114.32)
\psline(-152.50,-116.91)(-146.00,-114.32)
\psline(-147.50,-116.91)(-141.00,-114.32)
\psline(-142.50,-116.91)(-136.00,-114.32)
\psline(-137.50,-116.91)(-131.00,-114.32)
\psline(-132.50,-116.91)(-126.00,-114.32)
\psline(-127.50,-116.91)(-121.00,-114.32)
\psline(-122.50,-116.91)(-116.00,-114.32)
\psline(-117.50,-116.91)(-111.00,-114.32)
\psline(-156.00,-114.32)(-149.50,-111.72)
\psline(-151.00,-114.32)(-144.50,-111.72)
\psline(-146.00,-114.32)(-139.50,-111.72)
\psline(-141.00,-114.32)(-134.50,-111.72)
\psline(-136.00,-114.32)(-129.50,-111.72)
\psline(-131.00,-114.32)(-124.50,-111.72)
\psline(-126.00,-114.32)(-119.50,-111.72)
\psline(-121.00,-114.32)(-114.50,-111.72)
\psline(-116.00,-114.32)(-109.50,-111.72)
\psline(-154.50,-111.72)(-148.00,-109.12)
\psline(-149.50,-111.72)(-143.00,-109.12)
\psline(-144.50,-111.72)(-138.00,-109.12)
\psline(-139.50,-111.72)(-133.00,-109.12)
\psline(-134.50,-111.72)(-128.00,-109.12)
\psline(-129.50,-111.72)(-123.00,-109.12)
\psline(-124.50,-111.72)(-118.00,-109.12)
\psline(-119.50,-111.72)(-113.00,-109.12)
\psline(-114.50,-111.72)(-108.00,-109.12)
\psline(-153.00,-109.12)(-146.50,-106.52)
\psline(-148.00,-109.12)(-141.50,-106.52)
\psline(-143.00,-109.12)(-136.50,-106.52)
\psline(-138.00,-109.12)(-131.50,-106.52)
\psline(-133.00,-109.12)(-126.50,-106.52)
\psline(-128.00,-109.12)(-121.50,-106.52)
\psline(-123.00,-109.12)(-116.50,-106.52)
\psline(-118.00,-109.12)(-111.50,-106.52)
\psline(-113.00,-109.12)(-106.50,-106.52)
\psline(-151.50,-106.52)(-145.00,-103.92)
\psline(-146.50,-106.52)(-140.00,-103.92)
\psline(-141.50,-106.52)(-135.00,-103.92)
\psline(-136.50,-106.52)(-130.00,-103.92)
\psline(-131.50,-106.52)(-125.00,-103.92)
\psline(-126.50,-106.52)(-120.00,-103.92)
\psline(-121.50,-106.52)(-115.00,-103.92)
\psline(-116.50,-106.52)(-110.00,-103.92)
\psline(-111.50,-106.52)(-105.00,-103.92)
\psline(-150.00,-103.92)(-143.50,-101.32)
\psline(-145.00,-103.92)(-138.50,-101.32)
\psline(-140.00,-103.92)(-133.50,-101.32)
\psline(-135.00,-103.92)(-128.50,-101.32)
\psline(-130.00,-103.92)(-123.50,-101.32)
\psline(-125.00,-103.92)(-118.50,-101.32)
\psline(-120.00,-103.92)(-113.50,-101.32)
\psline(-115.00,-103.92)(-108.50,-101.32)
\psline(-110.00,-103.92)(-103.50,-101.32)
\psline(-148.50,-101.32)(-142.00,-98.73)
\psline(-143.50,-101.32)(-137.00,-98.73)
\psline(-138.50,-101.32)(-132.00,-98.73)
\psline(-133.50,-101.32)(-127.00,-98.73)
\psline(-128.50,-101.32)(-122.00,-98.73)
\psline(-123.50,-101.32)(-117.00,-98.73)
\psline(-118.50,-101.32)(-112.00,-98.73)
\psline(-113.50,-101.32)(-107.00,-98.73)
\psline(-108.50,-101.32)(-102.00,-98.73)
\psline(-147.00,-98.73)(-140.50,-96.13)
\psline(-142.00,-98.73)(-135.50,-96.13)
\psline(-137.00,-98.73)(-130.50,-96.13)
\psline(-132.00,-98.73)(-125.50,-96.13)
\psline(-127.00,-98.73)(-120.50,-96.13)
\psline(-122.00,-98.73)(-115.50,-96.13)
\psline(-117.00,-98.73)(-110.50,-96.13)
\psline(-112.00,-98.73)(-105.50,-96.13)
\psline(-107.00,-98.73)(-100.50,-96.13)
\psline(-145.50,-96.13)(-139.00,-93.53)
\psline(-140.50,-96.13)(-134.00,-93.53)
\psline(-135.50,-96.13)(-129.00,-93.53)
\psline(-130.50,-96.13)(-124.00,-93.53)
\psline(-125.50,-96.13)(-119.00,-93.53)
\psline(-120.50,-96.13)(-114.00,-93.53)
\psline(-115.50,-96.13)(-109.00,-93.53)
\psline(-110.50,-96.13)(-104.00,-93.53)
\psline(-105.50,-96.13)(-99.00,-93.53)
\psline(-144.00,-93.53)(-137.50,-90.93)
\psline(-139.00,-93.53)(-132.50,-90.93)
\psline(-134.00,-93.53)(-127.50,-90.93)
\psline(-129.00,-93.53)(-122.50,-90.93)
\psline(-124.00,-93.53)(-117.50,-90.93)
\psline(-119.00,-93.53)(-112.50,-90.93)
\psline(-114.00,-93.53)(-107.50,-90.93)
\psline(-109.00,-93.53)(-102.50,-90.93)
\psline(-104.00,-93.53)(-97.50,-90.93)
\psline(-142.50,-90.93)(-136.00,-88.33)
\psline(-137.50,-90.93)(-131.00,-88.33)
\psline(-132.50,-90.93)(-126.00,-88.33)
\psline(-127.50,-90.93)(-121.00,-88.33)
\psline(-122.50,-90.93)(-116.00,-88.33)
\psline(-117.50,-90.93)(-111.00,-88.33)
\psline(-112.50,-90.93)(-106.00,-88.33)
\psline(-107.50,-90.93)(-101.00,-88.33)
\psline(-102.50,-90.93)(-96.00,-88.33)
\psline(-141.00,-88.33)(-134.50,-85.74)
\psline(-136.00,-88.33)(-129.50,-85.74)
\psline(-131.00,-88.33)(-124.50,-85.74)
\psline(-126.00,-88.33)(-119.50,-85.74)
\psline(-121.00,-88.33)(-114.50,-85.74)
\psline(-116.00,-88.33)(-109.50,-85.74)
\psline(-111.00,-88.33)(-104.50,-85.74)
\psline(-106.00,-88.33)(-99.50,-85.74)
\psline(-101.00,-88.33)(-94.50,-85.74)
\psline(-139.50,-85.74)(-133.00,-83.14)
\psline(-134.50,-85.74)(-128.00,-83.14)
\psline(-129.50,-85.74)(-123.00,-83.14)
\psline(-124.50,-85.74)(-118.00,-83.14)
\psline(-119.50,-85.74)(-113.00,-83.14)
\psline(-114.50,-85.74)(-108.00,-83.14)
\psline(-109.50,-85.74)(-103.00,-83.14)
\psline(-104.50,-85.74)(-98.00,-83.14)
\psline(-99.50,-85.74)(-93.00,-83.14)
\psline(-138.00,-83.14)(-131.50,-80.54)
\psline(-133.00,-83.14)(-126.50,-80.54)
\psline(-128.00,-83.14)(-121.50,-80.54)
\psline(-123.00,-83.14)(-116.50,-80.54)
\psline(-118.00,-83.14)(-111.50,-80.54)
\psline(-113.00,-83.14)(-106.50,-80.54)
\psline(-108.00,-83.14)(-101.50,-80.54)
\psline(-103.00,-83.14)(-96.50,-80.54)
\psline(-98.00,-83.14)(-91.50,-80.54)
\psline(-136.50,-80.54)(-130.00,-77.94)
\psline(-131.50,-80.54)(-125.00,-77.94)
\psline(-126.50,-80.54)(-120.00,-77.94)
\psline(-121.50,-80.54)(-115.00,-77.94)
\psline(-116.50,-80.54)(-110.00,-77.94)
\psline(-111.50,-80.54)(-105.00,-77.94)
\psline(-106.50,-80.54)(-100.00,-77.94)
\psline(-101.50,-80.54)(-95.00,-77.94)
\psline(-96.50,-80.54)(-90.00,-77.94)
\psline(-135.00,-77.94)(-128.50,-75.34)
\psline(-130.00,-77.94)(-123.50,-75.34)
\psline(-125.00,-77.94)(-118.50,-75.34)
\psline(-120.00,-77.94)(-113.50,-75.34)
\psline(-115.00,-77.94)(-108.50,-75.34)
\psline(-110.00,-77.94)(-103.50,-75.34)
\psline(-105.00,-77.94)(-98.50,-75.34)
\psline(-100.00,-77.94)(-93.50,-75.34)
\psline(-95.00,-77.94)(-88.50,-75.34)
\psline(-133.50,-75.34)(-127.00,-72.75)
\psline(-128.50,-75.34)(-122.00,-72.75)
\psline(-123.50,-75.34)(-117.00,-72.75)
\psline(-118.50,-75.34)(-112.00,-72.75)
\psline(-113.50,-75.34)(-107.00,-72.75)
\psline(-108.50,-75.34)(-102.00,-72.75)
\psline(-103.50,-75.34)(-97.00,-72.75)
\psline(-98.50,-75.34)(-92.00,-72.75)
\psline(-93.50,-75.34)(-87.00,-72.75)
\psline(-132.00,-72.75)(-125.50,-70.15)
\psline(-127.00,-72.75)(-120.50,-70.15)
\psline(-122.00,-72.75)(-115.50,-70.15)
\psline(-117.00,-72.75)(-110.50,-70.15)
\psline(-112.00,-72.75)(-105.50,-70.15)
\psline(-107.00,-72.75)(-100.50,-70.15)
\psline(-102.00,-72.75)(-95.50,-70.15)
\psline(-97.00,-72.75)(-90.50,-70.15)
\psline(-92.00,-72.75)(-85.50,-70.15)
\psline(-130.50,-70.15)(-124.00,-67.55)
\psline(-125.50,-70.15)(-119.00,-67.55)
\psline(-120.50,-70.15)(-114.00,-67.55)
\psline(-115.50,-70.15)(-109.00,-67.55)
\psline(-110.50,-70.15)(-104.00,-67.55)
\psline(-105.50,-70.15)(-99.00,-67.55)
\psline(-100.50,-70.15)(-94.00,-67.55)
\psline(-95.50,-70.15)(-89.00,-67.55)
\psline(-90.50,-70.15)(-84.00,-67.55)
\psline(-129.00,-67.55)(-122.50,-64.95)
\psline(-124.00,-67.55)(-117.50,-64.95)
\psline(-119.00,-67.55)(-112.50,-64.95)
\psline(-114.00,-67.55)(-107.50,-64.95)
\psline(-109.00,-67.55)(-102.50,-64.95)
\psline(-104.00,-67.55)(-97.50,-64.95)
\psline(-99.00,-67.55)(-92.50,-64.95)
\psline(-94.00,-67.55)(-87.50,-64.95)
\psline(-89.00,-67.55)(-82.50,-64.95)
\psline(-127.50,-64.95)(-121.00,-62.35)
\psline(-122.50,-64.95)(-116.00,-62.35)
\psline(-117.50,-64.95)(-111.00,-62.35)
\psline(-112.50,-64.95)(-106.00,-62.35)
\psline(-107.50,-64.95)(-101.00,-62.35)
\psline(-102.50,-64.95)(-96.00,-62.35)
\psline(-97.50,-64.95)(-91.00,-62.35)
\psline(-92.50,-64.95)(-86.00,-62.35)
\psline(-87.50,-64.95)(-81.00,-62.35)
\psline(-126.00,-62.35)(-119.50,-59.76)
\psline(-121.00,-62.35)(-114.50,-59.76)
\psline(-116.00,-62.35)(-109.50,-59.76)
\psline(-111.00,-62.35)(-104.50,-59.76)
\psline(-106.00,-62.35)(-99.50,-59.76)
\psline(-101.00,-62.35)(-94.50,-59.76)
\psline(-96.00,-62.35)(-89.50,-59.76)
\psline(-91.00,-62.35)(-84.50,-59.76)
\psline(-86.00,-62.35)(-79.50,-59.76)
\psline(-124.50,-59.76)(-118.00,-57.16)
\psline(-119.50,-59.76)(-113.00,-57.16)
\psline(-114.50,-59.76)(-108.00,-57.16)
\psline(-109.50,-59.76)(-103.00,-57.16)
\psline(-104.50,-59.76)(-98.00,-57.16)
\psline(-99.50,-59.76)(-93.00,-57.16)
\psline(-94.50,-59.76)(-88.00,-57.16)
\psline(-89.50,-59.76)(-83.00,-57.16)
\psline(-84.50,-59.76)(-78.00,-57.16)
\psline(-123.00,-57.16)(-116.50,-54.56)
\psline(-118.00,-57.16)(-111.50,-54.56)
\psline(-113.00,-57.16)(-106.50,-54.56)
\psline(-108.00,-57.16)(-101.50,-54.56)
\psline(-103.00,-57.16)(-96.50,-54.56)
\psline(-98.00,-57.16)(-91.50,-54.56)
\psline(-93.00,-57.16)(-86.50,-54.56)
\psline(-88.00,-57.16)(-81.50,-54.56)
\psline(-83.00,-57.16)(-76.50,-54.56)
\psline(-121.50,-54.56)(-115.00,-51.96)
\psline(-116.50,-54.56)(-110.00,-51.96)
\psline(-111.50,-54.56)(-105.00,-51.96)
\psline(-106.50,-54.56)(-100.00,-51.96)
\psline(-101.50,-54.56)(-95.00,-51.96)
\psline(-96.50,-54.56)(-90.00,-51.96)
\psline(-91.50,-54.56)(-85.00,-51.96)
\psline(-86.50,-54.56)(-80.00,-51.96)
\psline(-81.50,-54.56)(-75.00,-51.96)
\psline(-120.00,-51.96)(-113.50,-49.36)
\psline(-115.00,-51.96)(-108.50,-49.36)
\psline(-110.00,-51.96)(-103.50,-49.36)
\psline(-105.00,-51.96)(-98.50,-49.36)
\psline(-100.00,-51.96)(-93.50,-49.36)
\psline(-95.00,-51.96)(-88.50,-49.36)
\psline(-90.00,-51.96)(-83.50,-49.36)
\psline(-85.00,-51.96)(-78.50,-49.36)
\psline(-80.00,-51.96)(-73.50,-49.36)
\pspolygon[fillstyle=solid,linewidth=0.1pt,fillcolor=red](-51.00,-88.33)(-67.50,-116.91)(-112.50,-116.91)(-96.00,-88.33)
\psline(-67.50,-116.91)(-112.50,-116.91)
\psline(-66.00,-114.32)(-111.00,-114.32)
\psline(-64.50,-111.72)(-109.50,-111.72)
\psline(-63.00,-109.12)(-108.00,-109.12)
\psline(-61.50,-106.52)(-106.50,-106.52)
\psline(-60.00,-103.92)(-105.00,-103.92)
\psline(-58.50,-101.32)(-103.50,-101.32)
\psline(-57.00,-98.73)(-102.00,-98.73)
\psline(-55.50,-96.13)(-100.50,-96.13)
\psline(-54.00,-93.53)(-99.00,-93.53)
\psline(-52.50,-90.93)(-97.50,-90.93)
\psline(-51.00,-88.33)(-96.00,-88.33)
\psline(-96.00,-88.33)(-112.50,-116.91)
\psline(-91.00,-88.33)(-107.50,-116.91)
\psline(-86.00,-88.33)(-102.50,-116.91)
\psline(-81.00,-88.33)(-97.50,-116.91)
\psline(-76.00,-88.33)(-92.50,-116.91)
\psline(-71.00,-88.33)(-87.50,-116.91)
\psline(-66.00,-88.33)(-82.50,-116.91)
\psline(-61.00,-88.33)(-77.50,-116.91)
\psline(-56.00,-88.33)(-72.50,-116.91)
\psline(-51.00,-88.33)(-67.50,-116.91)
\psline(-112.50,-116.91)(-106.00,-114.32)
\psline(-107.50,-116.91)(-101.00,-114.32)
\psline(-102.50,-116.91)(-96.00,-114.32)
\psline(-97.50,-116.91)(-91.00,-114.32)
\psline(-92.50,-116.91)(-86.00,-114.32)
\psline(-87.50,-116.91)(-81.00,-114.32)
\psline(-82.50,-116.91)(-76.00,-114.32)
\psline(-77.50,-116.91)(-71.00,-114.32)
\psline(-72.50,-116.91)(-66.00,-114.32)
\psline(-111.00,-114.32)(-104.50,-111.72)
\psline(-106.00,-114.32)(-99.50,-111.72)
\psline(-101.00,-114.32)(-94.50,-111.72)
\psline(-96.00,-114.32)(-89.50,-111.72)
\psline(-91.00,-114.32)(-84.50,-111.72)
\psline(-86.00,-114.32)(-79.50,-111.72)
\psline(-81.00,-114.32)(-74.50,-111.72)
\psline(-76.00,-114.32)(-69.50,-111.72)
\psline(-71.00,-114.32)(-64.50,-111.72)
\psline(-109.50,-111.72)(-103.00,-109.12)
\psline(-104.50,-111.72)(-98.00,-109.12)
\psline(-99.50,-111.72)(-93.00,-109.12)
\psline(-94.50,-111.72)(-88.00,-109.12)
\psline(-89.50,-111.72)(-83.00,-109.12)
\psline(-84.50,-111.72)(-78.00,-109.12)
\psline(-79.50,-111.72)(-73.00,-109.12)
\psline(-74.50,-111.72)(-68.00,-109.12)
\psline(-69.50,-111.72)(-63.00,-109.12)
\psline(-108.00,-109.12)(-101.50,-106.52)
\psline(-103.00,-109.12)(-96.50,-106.52)
\psline(-98.00,-109.12)(-91.50,-106.52)
\psline(-93.00,-109.12)(-86.50,-106.52)
\psline(-88.00,-109.12)(-81.50,-106.52)
\psline(-83.00,-109.12)(-76.50,-106.52)
\psline(-78.00,-109.12)(-71.50,-106.52)
\psline(-73.00,-109.12)(-66.50,-106.52)
\psline(-68.00,-109.12)(-61.50,-106.52)
\psline(-106.50,-106.52)(-100.00,-103.92)
\psline(-101.50,-106.52)(-95.00,-103.92)
\psline(-96.50,-106.52)(-90.00,-103.92)
\psline(-91.50,-106.52)(-85.00,-103.92)
\psline(-86.50,-106.52)(-80.00,-103.92)
\psline(-81.50,-106.52)(-75.00,-103.92)
\psline(-76.50,-106.52)(-70.00,-103.92)
\psline(-71.50,-106.52)(-65.00,-103.92)
\psline(-66.50,-106.52)(-60.00,-103.92)
\psline(-105.00,-103.92)(-98.50,-101.32)
\psline(-100.00,-103.92)(-93.50,-101.32)
\psline(-95.00,-103.92)(-88.50,-101.32)
\psline(-90.00,-103.92)(-83.50,-101.32)
\psline(-85.00,-103.92)(-78.50,-101.32)
\psline(-80.00,-103.92)(-73.50,-101.32)
\psline(-75.00,-103.92)(-68.50,-101.32)
\psline(-70.00,-103.92)(-63.50,-101.32)
\psline(-65.00,-103.92)(-58.50,-101.32)
\psline(-103.50,-101.32)(-97.00,-98.73)
\psline(-98.50,-101.32)(-92.00,-98.73)
\psline(-93.50,-101.32)(-87.00,-98.73)
\psline(-88.50,-101.32)(-82.00,-98.73)
\psline(-83.50,-101.32)(-77.00,-98.73)
\psline(-78.50,-101.32)(-72.00,-98.73)
\psline(-73.50,-101.32)(-67.00,-98.73)
\psline(-68.50,-101.32)(-62.00,-98.73)
\psline(-63.50,-101.32)(-57.00,-98.73)
\psline(-102.00,-98.73)(-95.50,-96.13)
\psline(-97.00,-98.73)(-90.50,-96.13)
\psline(-92.00,-98.73)(-85.50,-96.13)
\psline(-87.00,-98.73)(-80.50,-96.13)
\psline(-82.00,-98.73)(-75.50,-96.13)
\psline(-77.00,-98.73)(-70.50,-96.13)
\psline(-72.00,-98.73)(-65.50,-96.13)
\psline(-67.00,-98.73)(-60.50,-96.13)
\psline(-62.00,-98.73)(-55.50,-96.13)
\psline(-100.50,-96.13)(-94.00,-93.53)
\psline(-95.50,-96.13)(-89.00,-93.53)
\psline(-90.50,-96.13)(-84.00,-93.53)
\psline(-85.50,-96.13)(-79.00,-93.53)
\psline(-80.50,-96.13)(-74.00,-93.53)
\psline(-75.50,-96.13)(-69.00,-93.53)
\psline(-70.50,-96.13)(-64.00,-93.53)
\psline(-65.50,-96.13)(-59.00,-93.53)
\psline(-60.50,-96.13)(-54.00,-93.53)
\psline(-99.00,-93.53)(-92.50,-90.93)
\psline(-94.00,-93.53)(-87.50,-90.93)
\psline(-89.00,-93.53)(-82.50,-90.93)
\psline(-84.00,-93.53)(-77.50,-90.93)
\psline(-79.00,-93.53)(-72.50,-90.93)
\psline(-74.00,-93.53)(-67.50,-90.93)
\psline(-69.00,-93.53)(-62.50,-90.93)
\psline(-64.00,-93.53)(-57.50,-90.93)
\psline(-59.00,-93.53)(-52.50,-90.93)
\psline(-97.50,-90.93)(-91.00,-88.33)
\psline(-92.50,-90.93)(-86.00,-88.33)
\psline(-87.50,-90.93)(-81.00,-88.33)
\psline(-82.50,-90.93)(-76.00,-88.33)
\psline(-77.50,-90.93)(-71.00,-88.33)
\psline(-72.50,-90.93)(-66.00,-88.33)
\psline(-67.50,-90.93)(-61.00,-88.33)
\psline(-62.50,-90.93)(-56.00,-88.33)
\psline(-57.50,-90.93)(-51.00,-88.33)
\pspolygon[fillstyle=solid,linewidth=0.1pt,fillcolor=red](-96.00,-10.39)(-157.50,-116.91)(-202.50,-116.91)(-141.00,-10.39)
\psline(-157.50,-116.91)(-202.50,-116.91)
\psline(-156.00,-114.32)(-201.00,-114.32)
\psline(-154.50,-111.72)(-199.50,-111.72)
\psline(-153.00,-109.12)(-198.00,-109.12)
\psline(-151.50,-106.52)(-196.50,-106.52)
\psline(-150.00,-103.92)(-195.00,-103.92)
\psline(-148.50,-101.32)(-193.50,-101.32)
\psline(-147.00,-98.73)(-192.00,-98.73)
\psline(-145.50,-96.13)(-190.50,-96.13)
\psline(-144.00,-93.53)(-189.00,-93.53)
\psline(-142.50,-90.93)(-187.50,-90.93)
\psline(-141.00,-88.33)(-186.00,-88.33)
\psline(-139.50,-85.74)(-184.50,-85.74)
\psline(-138.00,-83.14)(-183.00,-83.14)
\psline(-136.50,-80.54)(-181.50,-80.54)
\psline(-135.00,-77.94)(-180.00,-77.94)
\psline(-133.50,-75.34)(-178.50,-75.34)
\psline(-132.00,-72.75)(-177.00,-72.75)
\psline(-130.50,-70.15)(-175.50,-70.15)
\psline(-129.00,-67.55)(-174.00,-67.55)
\psline(-127.50,-64.95)(-172.50,-64.95)
\psline(-126.00,-62.35)(-171.00,-62.35)
\psline(-124.50,-59.76)(-169.50,-59.76)
\psline(-123.00,-57.16)(-168.00,-57.16)
\psline(-121.50,-54.56)(-166.50,-54.56)
\psline(-120.00,-51.96)(-165.00,-51.96)
\psline(-118.50,-49.36)(-163.50,-49.36)
\psline(-117.00,-46.77)(-162.00,-46.77)
\psline(-115.50,-44.17)(-160.50,-44.17)
\psline(-114.00,-41.57)(-159.00,-41.57)
\psline(-112.50,-38.97)(-157.50,-38.97)
\psline(-111.00,-36.37)(-156.00,-36.37)
\psline(-109.50,-33.77)(-154.50,-33.77)
\psline(-108.00,-31.18)(-153.00,-31.18)
\psline(-106.50,-28.58)(-151.50,-28.58)
\psline(-105.00,-25.98)(-150.00,-25.98)
\psline(-103.50,-23.38)(-148.50,-23.38)
\psline(-102.00,-20.78)(-147.00,-20.78)
\psline(-100.50,-18.19)(-145.50,-18.19)
\psline(-99.00,-15.59)(-144.00,-15.59)
\psline(-97.50,-12.99)(-142.50,-12.99)
\psline(-96.00,-10.39)(-141.00,-10.39)
\psline(-141.00,-10.39)(-202.50,-116.91)
\psline(-136.00,-10.39)(-197.50,-116.91)
\psline(-131.00,-10.39)(-192.50,-116.91)
\psline(-126.00,-10.39)(-187.50,-116.91)
\psline(-121.00,-10.39)(-182.50,-116.91)
\psline(-116.00,-10.39)(-177.50,-116.91)
\psline(-111.00,-10.39)(-172.50,-116.91)
\psline(-106.00,-10.39)(-167.50,-116.91)
\psline(-101.00,-10.39)(-162.50,-116.91)
\psline(-96.00,-10.39)(-157.50,-116.91)
\psline(-202.50,-116.91)(-196.00,-114.32)
\psline(-197.50,-116.91)(-191.00,-114.32)
\psline(-192.50,-116.91)(-186.00,-114.32)
\psline(-187.50,-116.91)(-181.00,-114.32)
\psline(-182.50,-116.91)(-176.00,-114.32)
\psline(-177.50,-116.91)(-171.00,-114.32)
\psline(-172.50,-116.91)(-166.00,-114.32)
\psline(-167.50,-116.91)(-161.00,-114.32)
\psline(-162.50,-116.91)(-156.00,-114.32)
\psline(-201.00,-114.32)(-194.50,-111.72)
\psline(-196.00,-114.32)(-189.50,-111.72)
\psline(-191.00,-114.32)(-184.50,-111.72)
\psline(-186.00,-114.32)(-179.50,-111.72)
\psline(-181.00,-114.32)(-174.50,-111.72)
\psline(-176.00,-114.32)(-169.50,-111.72)
\psline(-171.00,-114.32)(-164.50,-111.72)
\psline(-166.00,-114.32)(-159.50,-111.72)
\psline(-161.00,-114.32)(-154.50,-111.72)
\psline(-199.50,-111.72)(-193.00,-109.12)
\psline(-194.50,-111.72)(-188.00,-109.12)
\psline(-189.50,-111.72)(-183.00,-109.12)
\psline(-184.50,-111.72)(-178.00,-109.12)
\psline(-179.50,-111.72)(-173.00,-109.12)
\psline(-174.50,-111.72)(-168.00,-109.12)
\psline(-169.50,-111.72)(-163.00,-109.12)
\psline(-164.50,-111.72)(-158.00,-109.12)
\psline(-159.50,-111.72)(-153.00,-109.12)
\psline(-198.00,-109.12)(-191.50,-106.52)
\psline(-193.00,-109.12)(-186.50,-106.52)
\psline(-188.00,-109.12)(-181.50,-106.52)
\psline(-183.00,-109.12)(-176.50,-106.52)
\psline(-178.00,-109.12)(-171.50,-106.52)
\psline(-173.00,-109.12)(-166.50,-106.52)
\psline(-168.00,-109.12)(-161.50,-106.52)
\psline(-163.00,-109.12)(-156.50,-106.52)
\psline(-158.00,-109.12)(-151.50,-106.52)
\psline(-196.50,-106.52)(-190.00,-103.92)
\psline(-191.50,-106.52)(-185.00,-103.92)
\psline(-186.50,-106.52)(-180.00,-103.92)
\psline(-181.50,-106.52)(-175.00,-103.92)
\psline(-176.50,-106.52)(-170.00,-103.92)
\psline(-171.50,-106.52)(-165.00,-103.92)
\psline(-166.50,-106.52)(-160.00,-103.92)
\psline(-161.50,-106.52)(-155.00,-103.92)
\psline(-156.50,-106.52)(-150.00,-103.92)
\psline(-195.00,-103.92)(-188.50,-101.32)
\psline(-190.00,-103.92)(-183.50,-101.32)
\psline(-185.00,-103.92)(-178.50,-101.32)
\psline(-180.00,-103.92)(-173.50,-101.32)
\psline(-175.00,-103.92)(-168.50,-101.32)
\psline(-170.00,-103.92)(-163.50,-101.32)
\psline(-165.00,-103.92)(-158.50,-101.32)
\psline(-160.00,-103.92)(-153.50,-101.32)
\psline(-155.00,-103.92)(-148.50,-101.32)
\psline(-193.50,-101.32)(-187.00,-98.73)
\psline(-188.50,-101.32)(-182.00,-98.73)
\psline(-183.50,-101.32)(-177.00,-98.73)
\psline(-178.50,-101.32)(-172.00,-98.73)
\psline(-173.50,-101.32)(-167.00,-98.73)
\psline(-168.50,-101.32)(-162.00,-98.73)
\psline(-163.50,-101.32)(-157.00,-98.73)
\psline(-158.50,-101.32)(-152.00,-98.73)
\psline(-153.50,-101.32)(-147.00,-98.73)
\psline(-192.00,-98.73)(-185.50,-96.13)
\psline(-187.00,-98.73)(-180.50,-96.13)
\psline(-182.00,-98.73)(-175.50,-96.13)
\psline(-177.00,-98.73)(-170.50,-96.13)
\psline(-172.00,-98.73)(-165.50,-96.13)
\psline(-167.00,-98.73)(-160.50,-96.13)
\psline(-162.00,-98.73)(-155.50,-96.13)
\psline(-157.00,-98.73)(-150.50,-96.13)
\psline(-152.00,-98.73)(-145.50,-96.13)
\psline(-190.50,-96.13)(-184.00,-93.53)
\psline(-185.50,-96.13)(-179.00,-93.53)
\psline(-180.50,-96.13)(-174.00,-93.53)
\psline(-175.50,-96.13)(-169.00,-93.53)
\psline(-170.50,-96.13)(-164.00,-93.53)
\psline(-165.50,-96.13)(-159.00,-93.53)
\psline(-160.50,-96.13)(-154.00,-93.53)
\psline(-155.50,-96.13)(-149.00,-93.53)
\psline(-150.50,-96.13)(-144.00,-93.53)
\psline(-189.00,-93.53)(-182.50,-90.93)
\psline(-184.00,-93.53)(-177.50,-90.93)
\psline(-179.00,-93.53)(-172.50,-90.93)
\psline(-174.00,-93.53)(-167.50,-90.93)
\psline(-169.00,-93.53)(-162.50,-90.93)
\psline(-164.00,-93.53)(-157.50,-90.93)
\psline(-159.00,-93.53)(-152.50,-90.93)
\psline(-154.00,-93.53)(-147.50,-90.93)
\psline(-149.00,-93.53)(-142.50,-90.93)
\psline(-187.50,-90.93)(-181.00,-88.33)
\psline(-182.50,-90.93)(-176.00,-88.33)
\psline(-177.50,-90.93)(-171.00,-88.33)
\psline(-172.50,-90.93)(-166.00,-88.33)
\psline(-167.50,-90.93)(-161.00,-88.33)
\psline(-162.50,-90.93)(-156.00,-88.33)
\psline(-157.50,-90.93)(-151.00,-88.33)
\psline(-152.50,-90.93)(-146.00,-88.33)
\psline(-147.50,-90.93)(-141.00,-88.33)
\psline(-186.00,-88.33)(-179.50,-85.74)
\psline(-181.00,-88.33)(-174.50,-85.74)
\psline(-176.00,-88.33)(-169.50,-85.74)
\psline(-171.00,-88.33)(-164.50,-85.74)
\psline(-166.00,-88.33)(-159.50,-85.74)
\psline(-161.00,-88.33)(-154.50,-85.74)
\psline(-156.00,-88.33)(-149.50,-85.74)
\psline(-151.00,-88.33)(-144.50,-85.74)
\psline(-146.00,-88.33)(-139.50,-85.74)
\psline(-184.50,-85.74)(-178.00,-83.14)
\psline(-179.50,-85.74)(-173.00,-83.14)
\psline(-174.50,-85.74)(-168.00,-83.14)
\psline(-169.50,-85.74)(-163.00,-83.14)
\psline(-164.50,-85.74)(-158.00,-83.14)
\psline(-159.50,-85.74)(-153.00,-83.14)
\psline(-154.50,-85.74)(-148.00,-83.14)
\psline(-149.50,-85.74)(-143.00,-83.14)
\psline(-144.50,-85.74)(-138.00,-83.14)
\psline(-183.00,-83.14)(-176.50,-80.54)
\psline(-178.00,-83.14)(-171.50,-80.54)
\psline(-173.00,-83.14)(-166.50,-80.54)
\psline(-168.00,-83.14)(-161.50,-80.54)
\psline(-163.00,-83.14)(-156.50,-80.54)
\psline(-158.00,-83.14)(-151.50,-80.54)
\psline(-153.00,-83.14)(-146.50,-80.54)
\psline(-148.00,-83.14)(-141.50,-80.54)
\psline(-143.00,-83.14)(-136.50,-80.54)
\psline(-181.50,-80.54)(-175.00,-77.94)
\psline(-176.50,-80.54)(-170.00,-77.94)
\psline(-171.50,-80.54)(-165.00,-77.94)
\psline(-166.50,-80.54)(-160.00,-77.94)
\psline(-161.50,-80.54)(-155.00,-77.94)
\psline(-156.50,-80.54)(-150.00,-77.94)
\psline(-151.50,-80.54)(-145.00,-77.94)
\psline(-146.50,-80.54)(-140.00,-77.94)
\psline(-141.50,-80.54)(-135.00,-77.94)
\psline(-180.00,-77.94)(-173.50,-75.34)
\psline(-175.00,-77.94)(-168.50,-75.34)
\psline(-170.00,-77.94)(-163.50,-75.34)
\psline(-165.00,-77.94)(-158.50,-75.34)
\psline(-160.00,-77.94)(-153.50,-75.34)
\psline(-155.00,-77.94)(-148.50,-75.34)
\psline(-150.00,-77.94)(-143.50,-75.34)
\psline(-145.00,-77.94)(-138.50,-75.34)
\psline(-140.00,-77.94)(-133.50,-75.34)
\psline(-178.50,-75.34)(-172.00,-72.75)
\psline(-173.50,-75.34)(-167.00,-72.75)
\psline(-168.50,-75.34)(-162.00,-72.75)
\psline(-163.50,-75.34)(-157.00,-72.75)
\psline(-158.50,-75.34)(-152.00,-72.75)
\psline(-153.50,-75.34)(-147.00,-72.75)
\psline(-148.50,-75.34)(-142.00,-72.75)
\psline(-143.50,-75.34)(-137.00,-72.75)
\psline(-138.50,-75.34)(-132.00,-72.75)
\psline(-177.00,-72.75)(-170.50,-70.15)
\psline(-172.00,-72.75)(-165.50,-70.15)
\psline(-167.00,-72.75)(-160.50,-70.15)
\psline(-162.00,-72.75)(-155.50,-70.15)
\psline(-157.00,-72.75)(-150.50,-70.15)
\psline(-152.00,-72.75)(-145.50,-70.15)
\psline(-147.00,-72.75)(-140.50,-70.15)
\psline(-142.00,-72.75)(-135.50,-70.15)
\psline(-137.00,-72.75)(-130.50,-70.15)
\psline(-175.50,-70.15)(-169.00,-67.55)
\psline(-170.50,-70.15)(-164.00,-67.55)
\psline(-165.50,-70.15)(-159.00,-67.55)
\psline(-160.50,-70.15)(-154.00,-67.55)
\psline(-155.50,-70.15)(-149.00,-67.55)
\psline(-150.50,-70.15)(-144.00,-67.55)
\psline(-145.50,-70.15)(-139.00,-67.55)
\psline(-140.50,-70.15)(-134.00,-67.55)
\psline(-135.50,-70.15)(-129.00,-67.55)
\psline(-174.00,-67.55)(-167.50,-64.95)
\psline(-169.00,-67.55)(-162.50,-64.95)
\psline(-164.00,-67.55)(-157.50,-64.95)
\psline(-159.00,-67.55)(-152.50,-64.95)
\psline(-154.00,-67.55)(-147.50,-64.95)
\psline(-149.00,-67.55)(-142.50,-64.95)
\psline(-144.00,-67.55)(-137.50,-64.95)
\psline(-139.00,-67.55)(-132.50,-64.95)
\psline(-134.00,-67.55)(-127.50,-64.95)
\psline(-172.50,-64.95)(-166.00,-62.35)
\psline(-167.50,-64.95)(-161.00,-62.35)
\psline(-162.50,-64.95)(-156.00,-62.35)
\psline(-157.50,-64.95)(-151.00,-62.35)
\psline(-152.50,-64.95)(-146.00,-62.35)
\psline(-147.50,-64.95)(-141.00,-62.35)
\psline(-142.50,-64.95)(-136.00,-62.35)
\psline(-137.50,-64.95)(-131.00,-62.35)
\psline(-132.50,-64.95)(-126.00,-62.35)
\psline(-171.00,-62.35)(-164.50,-59.76)
\psline(-166.00,-62.35)(-159.50,-59.76)
\psline(-161.00,-62.35)(-154.50,-59.76)
\psline(-156.00,-62.35)(-149.50,-59.76)
\psline(-151.00,-62.35)(-144.50,-59.76)
\psline(-146.00,-62.35)(-139.50,-59.76)
\psline(-141.00,-62.35)(-134.50,-59.76)
\psline(-136.00,-62.35)(-129.50,-59.76)
\psline(-131.00,-62.35)(-124.50,-59.76)
\psline(-169.50,-59.76)(-163.00,-57.16)
\psline(-164.50,-59.76)(-158.00,-57.16)
\psline(-159.50,-59.76)(-153.00,-57.16)
\psline(-154.50,-59.76)(-148.00,-57.16)
\psline(-149.50,-59.76)(-143.00,-57.16)
\psline(-144.50,-59.76)(-138.00,-57.16)
\psline(-139.50,-59.76)(-133.00,-57.16)
\psline(-134.50,-59.76)(-128.00,-57.16)
\psline(-129.50,-59.76)(-123.00,-57.16)
\psline(-168.00,-57.16)(-161.50,-54.56)
\psline(-163.00,-57.16)(-156.50,-54.56)
\psline(-158.00,-57.16)(-151.50,-54.56)
\psline(-153.00,-57.16)(-146.50,-54.56)
\psline(-148.00,-57.16)(-141.50,-54.56)
\psline(-143.00,-57.16)(-136.50,-54.56)
\psline(-138.00,-57.16)(-131.50,-54.56)
\psline(-133.00,-57.16)(-126.50,-54.56)
\psline(-128.00,-57.16)(-121.50,-54.56)
\psline(-166.50,-54.56)(-160.00,-51.96)
\psline(-161.50,-54.56)(-155.00,-51.96)
\psline(-156.50,-54.56)(-150.00,-51.96)
\psline(-151.50,-54.56)(-145.00,-51.96)
\psline(-146.50,-54.56)(-140.00,-51.96)
\psline(-141.50,-54.56)(-135.00,-51.96)
\psline(-136.50,-54.56)(-130.00,-51.96)
\psline(-131.50,-54.56)(-125.00,-51.96)
\psline(-126.50,-54.56)(-120.00,-51.96)
\psline(-165.00,-51.96)(-158.50,-49.36)
\psline(-160.00,-51.96)(-153.50,-49.36)
\psline(-155.00,-51.96)(-148.50,-49.36)
\psline(-150.00,-51.96)(-143.50,-49.36)
\psline(-145.00,-51.96)(-138.50,-49.36)
\psline(-140.00,-51.96)(-133.50,-49.36)
\psline(-135.00,-51.96)(-128.50,-49.36)
\psline(-130.00,-51.96)(-123.50,-49.36)
\psline(-125.00,-51.96)(-118.50,-49.36)
\psline(-163.50,-49.36)(-157.00,-46.77)
\psline(-158.50,-49.36)(-152.00,-46.77)
\psline(-153.50,-49.36)(-147.00,-46.77)
\psline(-148.50,-49.36)(-142.00,-46.77)
\psline(-143.50,-49.36)(-137.00,-46.77)
\psline(-138.50,-49.36)(-132.00,-46.77)
\psline(-133.50,-49.36)(-127.00,-46.77)
\psline(-128.50,-49.36)(-122.00,-46.77)
\psline(-123.50,-49.36)(-117.00,-46.77)
\psline(-162.00,-46.77)(-155.50,-44.17)
\psline(-157.00,-46.77)(-150.50,-44.17)
\psline(-152.00,-46.77)(-145.50,-44.17)
\psline(-147.00,-46.77)(-140.50,-44.17)
\psline(-142.00,-46.77)(-135.50,-44.17)
\psline(-137.00,-46.77)(-130.50,-44.17)
\psline(-132.00,-46.77)(-125.50,-44.17)
\psline(-127.00,-46.77)(-120.50,-44.17)
\psline(-122.00,-46.77)(-115.50,-44.17)
\psline(-160.50,-44.17)(-154.00,-41.57)
\psline(-155.50,-44.17)(-149.00,-41.57)
\psline(-150.50,-44.17)(-144.00,-41.57)
\psline(-145.50,-44.17)(-139.00,-41.57)
\psline(-140.50,-44.17)(-134.00,-41.57)
\psline(-135.50,-44.17)(-129.00,-41.57)
\psline(-130.50,-44.17)(-124.00,-41.57)
\psline(-125.50,-44.17)(-119.00,-41.57)
\psline(-120.50,-44.17)(-114.00,-41.57)
\psline(-159.00,-41.57)(-152.50,-38.97)
\psline(-154.00,-41.57)(-147.50,-38.97)
\psline(-149.00,-41.57)(-142.50,-38.97)
\psline(-144.00,-41.57)(-137.50,-38.97)
\psline(-139.00,-41.57)(-132.50,-38.97)
\psline(-134.00,-41.57)(-127.50,-38.97)
\psline(-129.00,-41.57)(-122.50,-38.97)
\psline(-124.00,-41.57)(-117.50,-38.97)
\psline(-119.00,-41.57)(-112.50,-38.97)
\psline(-157.50,-38.97)(-151.00,-36.37)
\psline(-152.50,-38.97)(-146.00,-36.37)
\psline(-147.50,-38.97)(-141.00,-36.37)
\psline(-142.50,-38.97)(-136.00,-36.37)
\psline(-137.50,-38.97)(-131.00,-36.37)
\psline(-132.50,-38.97)(-126.00,-36.37)
\psline(-127.50,-38.97)(-121.00,-36.37)
\psline(-122.50,-38.97)(-116.00,-36.37)
\psline(-117.50,-38.97)(-111.00,-36.37)
\psline(-156.00,-36.37)(-149.50,-33.77)
\psline(-151.00,-36.37)(-144.50,-33.77)
\psline(-146.00,-36.37)(-139.50,-33.77)
\psline(-141.00,-36.37)(-134.50,-33.77)
\psline(-136.00,-36.37)(-129.50,-33.77)
\psline(-131.00,-36.37)(-124.50,-33.77)
\psline(-126.00,-36.37)(-119.50,-33.77)
\psline(-121.00,-36.37)(-114.50,-33.77)
\psline(-116.00,-36.37)(-109.50,-33.77)
\psline(-154.50,-33.77)(-148.00,-31.18)
\psline(-149.50,-33.77)(-143.00,-31.18)
\psline(-144.50,-33.77)(-138.00,-31.18)
\psline(-139.50,-33.77)(-133.00,-31.18)
\psline(-134.50,-33.77)(-128.00,-31.18)
\psline(-129.50,-33.77)(-123.00,-31.18)
\psline(-124.50,-33.77)(-118.00,-31.18)
\psline(-119.50,-33.77)(-113.00,-31.18)
\psline(-114.50,-33.77)(-108.00,-31.18)
\psline(-153.00,-31.18)(-146.50,-28.58)
\psline(-148.00,-31.18)(-141.50,-28.58)
\psline(-143.00,-31.18)(-136.50,-28.58)
\psline(-138.00,-31.18)(-131.50,-28.58)
\psline(-133.00,-31.18)(-126.50,-28.58)
\psline(-128.00,-31.18)(-121.50,-28.58)
\psline(-123.00,-31.18)(-116.50,-28.58)
\psline(-118.00,-31.18)(-111.50,-28.58)
\psline(-113.00,-31.18)(-106.50,-28.58)
\psline(-151.50,-28.58)(-145.00,-25.98)
\psline(-146.50,-28.58)(-140.00,-25.98)
\psline(-141.50,-28.58)(-135.00,-25.98)
\psline(-136.50,-28.58)(-130.00,-25.98)
\psline(-131.50,-28.58)(-125.00,-25.98)
\psline(-126.50,-28.58)(-120.00,-25.98)
\psline(-121.50,-28.58)(-115.00,-25.98)
\psline(-116.50,-28.58)(-110.00,-25.98)
\psline(-111.50,-28.58)(-105.00,-25.98)
\psline(-150.00,-25.98)(-143.50,-23.38)
\psline(-145.00,-25.98)(-138.50,-23.38)
\psline(-140.00,-25.98)(-133.50,-23.38)
\psline(-135.00,-25.98)(-128.50,-23.38)
\psline(-130.00,-25.98)(-123.50,-23.38)
\psline(-125.00,-25.98)(-118.50,-23.38)
\psline(-120.00,-25.98)(-113.50,-23.38)
\psline(-115.00,-25.98)(-108.50,-23.38)
\psline(-110.00,-25.98)(-103.50,-23.38)
\psline(-148.50,-23.38)(-142.00,-20.78)
\psline(-143.50,-23.38)(-137.00,-20.78)
\psline(-138.50,-23.38)(-132.00,-20.78)
\psline(-133.50,-23.38)(-127.00,-20.78)
\psline(-128.50,-23.38)(-122.00,-20.78)
\psline(-123.50,-23.38)(-117.00,-20.78)
\psline(-118.50,-23.38)(-112.00,-20.78)
\psline(-113.50,-23.38)(-107.00,-20.78)
\psline(-108.50,-23.38)(-102.00,-20.78)
\psline(-147.00,-20.78)(-140.50,-18.19)
\psline(-142.00,-20.78)(-135.50,-18.19)
\psline(-137.00,-20.78)(-130.50,-18.19)
\psline(-132.00,-20.78)(-125.50,-18.19)
\psline(-127.00,-20.78)(-120.50,-18.19)
\psline(-122.00,-20.78)(-115.50,-18.19)
\psline(-117.00,-20.78)(-110.50,-18.19)
\psline(-112.00,-20.78)(-105.50,-18.19)
\psline(-107.00,-20.78)(-100.50,-18.19)
\psline(-145.50,-18.19)(-139.00,-15.59)
\psline(-140.50,-18.19)(-134.00,-15.59)
\psline(-135.50,-18.19)(-129.00,-15.59)
\psline(-130.50,-18.19)(-124.00,-15.59)
\psline(-125.50,-18.19)(-119.00,-15.59)
\psline(-120.50,-18.19)(-114.00,-15.59)
\psline(-115.50,-18.19)(-109.00,-15.59)
\psline(-110.50,-18.19)(-104.00,-15.59)
\psline(-105.50,-18.19)(-99.00,-15.59)
\psline(-144.00,-15.59)(-137.50,-12.99)
\psline(-139.00,-15.59)(-132.50,-12.99)
\psline(-134.00,-15.59)(-127.50,-12.99)
\psline(-129.00,-15.59)(-122.50,-12.99)
\psline(-124.00,-15.59)(-117.50,-12.99)
\psline(-119.00,-15.59)(-112.50,-12.99)
\psline(-114.00,-15.59)(-107.50,-12.99)
\psline(-109.00,-15.59)(-102.50,-12.99)
\psline(-104.00,-15.59)(-97.50,-12.99)
\psline(-142.50,-12.99)(-136.00,-10.39)
\psline(-137.50,-12.99)(-131.00,-10.39)
\psline(-132.50,-12.99)(-126.00,-10.39)
\psline(-127.50,-12.99)(-121.00,-10.39)
\psline(-122.50,-12.99)(-116.00,-10.39)
\psline(-117.50,-12.99)(-111.00,-10.39)
\psline(-112.50,-12.99)(-106.00,-10.39)
\psline(-107.50,-12.99)(-101.00,-10.39)
\psline(-102.50,-12.99)(-96.00,-10.39)
\endpspicture}

\def\FigureZeroLimitsTwo{
\pspicture(3.4,1.6)
\psset{unit=1cm}
\newrgbcolor{lightblue}{0.8 0.8 1}
\newrgbcolor{pink}{1 0.8 0.8}
\newrgbcolor{lightgreen}{0.8 1 0.8}
\newrgbcolor{lightyellow}{1 1 0.8} 
\newrgbcolor{orange}{1 0.5 0}
\newrgbcolor{violet}{1 0 1}
\pspolygon[fillstyle=solid,linewidth=1pt,fillcolor=lightblue](0.00,0.00)(3.50,1.94)(3.00,0.00)\pspolygon[fillstyle=solid,linewidth=1pt,fillcolor=lightblue](3.00,0.00)(5.62,1.45)(7.00,0.00)\pspolygon[fillstyle=solid,linewidth=1pt,fillcolor=lightblue](7.00,0.00)(6.50,1.94)(10.00,0.00)\pspolygon[fillstyle=solid,linewidth=1pt,fillcolor=pink](4.94,2.18)(7.00,0.00)(6.00,3.87)\pspolygon[fillstyle=solid,linewidth=1pt,fillcolor=lightyellow](4.00,3.87)(3.00,0.00)(4.75,0.97)\put(2.17,0.65){$0$}
\put(3.92,1.61){$2$}
\put(5.21,0.48){$1$}
\put(5.8,1.78){$4$}
\put(7.83,0.65){$3$}
\psline[linestyle=dashed](4.75,0.97)(5.75,4.84)
\psline[linestyle=dashed](5.75,4.84)(6.50,1.94)
\psdot(0.00,0.00)\put(-0.30,-0.35){$E$}
\psdot(3.00,0.00)\put(2.70,-0.35){$P$}
\psdot(7.00,0.00)\put(6.70,-0.35){$Q$}
\psdot(5.62,1.45)\put(5.4,1.07){$R$}
\psdot(10.00,0.00)\put(9.70,-0.35){$F$}
\endpspicture}

\def\FigureTheoremSeven{
\pspicture(2.5,1.2)(0,-0.05)
\psset{unit=2cm}
\newrgbcolor{lightblue}{0.8 0.8 1}
\pspolygon[fillstyle=solid,linewidth=1pt,fillcolor=lightblue](0.00,0.00)(1.80,0.87)(1.30,0.00)
\put(1.04,0.29){$1$}
\pspolygon[fillstyle=solid,linewidth=1pt,fillcolor=lightblue](1.30,0.00)(1.80,0.87)(3.11,0.87)
\put(2.07,0.58){$2$}
\pspolygon[fillstyle=solid,linewidth=1pt,fillcolor=lightblue](1.80,0.87)(3.61,1.73)(3.11,0.87)
\put(2.84,1.15){$3$}
\pspolygon[fillstyle=solid,linewidth=1pt,fillcolor=lightblue](1.30,0.00)(2.61,0.00)(3.11,0.87)
\put(2.34,0.29){$5$}
\pspolygon[fillstyle=solid,linewidth=1pt,fillcolor=lightblue](2.61,0.00)(3.68,0.73)(3.61,1.73)
\put(3.30,0.82){$4$}
\psline(2.61,0.00)(3.91,0.00)
\psdot(0.00,0.00)
\psdot(3.61,1.73)
\psdot(1.30,0.00)
\psdot(2.61,0.00)
\put(-0.130000,-0.170000){$A$}
\put(3.7,1.7){$B$}
\put(1.2,-0.170000){$P$}
\put(2.5,-0.170000){$Q$}
\endpspicture}

\author{Michael Beeson}
\title{No triangle can be cut into \goodbreak seven congruent triangles}         
\address{Professor of Mathematics (emeritus), San Jos\'e State University, San Jos\'e, California}

\begin{document}

\begin{abstract} We prove the 
theorem in the title, and prove the theorem for 11 as well as 7.
By previous work of others, the problem reduces to a number of 
cases.  The cases not solved already are solved here.
\medskip

\end{abstract}

\maketitle


\section{Introduction}
An $N$-tiling of triangle $ABC$ by triangle $T$ is a way of writing $ABC$ as a union of $N$ triangles
congruent to $T$, overlapping only at their boundaries.   The triangle $T$ is the ``tile''.
We consider the problem of cutting a triangle into $N$ congruent triangles.  
Shortly we shall give a number of examples of $N$-tilings, for 
various small values of $N$.  These examples will be tilings that have,
for the most part, been known a very long time.  But it will be 
obvious that $N=7$ is not in this list of examples, nor is $N=11$.

These two values of $N$ are the focus of this   
paper.  We will prove 
here that there are no 7-tilings or 11-tilings. 
Originally it was the question of $7$-tilings that attracted us to 
this subject.  This question could easily have
 been understood by the Greek geometers working in Alexandria with 
Euclid three centuries BCE, and possibly could have been solved by them too.
We were able to give a purely Euclidean proof, but it was very long and 
complicated.  Once sufficient machinery is developed, non-existence of
tilings for
many $N$, including all primes congruent to $3$ mod $4$, 
is a consequence, but also 
that development is long and complicated.  Therefore we were happy, in 
October 2018, to discover a relatively simple proof of the non-existence of 
any 7-tiling, which we present here.  It was also possible to treat $N=11$
with very little extra work--something we could not do with a purely 
Euclidean proof.  One might say that here Descartes is victorious over 
Euclid, as algebra and computation is shorter and more efficient than geometry.
Following Euclid we could do $N=7$, but not 11. 

We checked most of the algebra both by hand and by computer, using SageMath \cite{Sage}, 
and we provide the short snippets of code we used.
In only one place is it too laborious to do by hand.%
\footnote{SageMath code, being written in Python, needs to contain 
tabs for indentation.  When you cut and paste from a pdf file, you will
get spaces, not tabs.  Therefore you must paste into a file and 
supply tabs using the Unix utility {\tt unexpand}.  Also single quote
marks are a different character in pdf.  For very short snippets you 
may find it easier just to retype.}

These results fit into a larger research program, begun by Lazkovich
\cite{laczkovich1995}.  He studied the possible shapes of
tiles and triangles that can possibly be used in tilings, and obtained
results that will be described below.  It is our contribution to 
focus attention on $N$ as well.  One may say that Laczkovich studied
the pair $(ABC, T)$,  and we want to refine his work to study the triple
$(ABC,T,N)$. 
 
\section{Some examples of tilings} 
 Figures 1
 through 4 show the simplest examples of $N$-tilings.
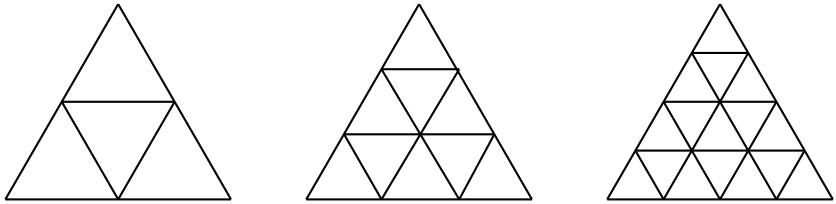
\begin{figure}[ht]    
\caption{Two 3-tilings   }
\label{figure:3-tilings}
\begin{center}
\ThreeTilingA
\hskip1cm
\ThreeTiling
\end{center}
\end{figure}

\begin{figure}[ht]
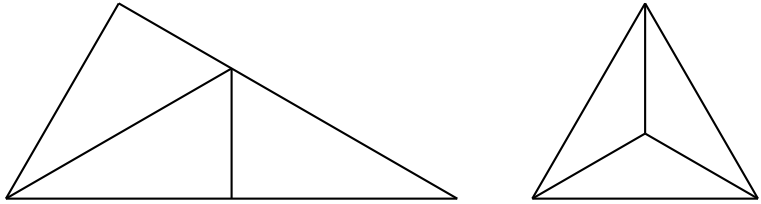
   
\caption{A 4-tiling, a 9-tiling, and a 16-tiling}
\label{figure:nsquared}
\begin{center}
\FourTiling  
\hskip 1cm
\NineTiling
\hskip1cm
\SixteenTiling
\end{center}
\end{figure}

\begin{figure}[ht]
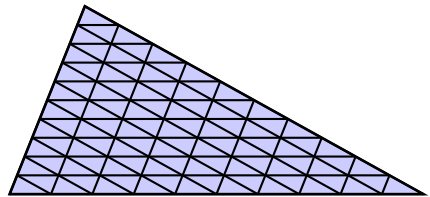

\caption{A quadratic tiling of an arbitrary triangle}
\label{figure:bigquadratic}
\begin{center}
\FigureBigQuadratic
\end{center}
\end{figure}

The method illustrated for $N=4$ ,$9$, and $16$  generalizes to any perfect square $N$. 
While the two exhibited 3-tilings 
clearly depend on the exact angles of the triangle, 
{\em any} triangle can be decomposed into $n^2$ congruent triangles by  drawing $n-1$ equally spaced lines parallel to each of the three sides of 
the triangle, as illustrated in Fig.~\ref{figure:bigquadratic}.
Moreover, the large (tiled) triangle is similar to the 
small triangle (the ``tile'').  We call such a tiling a {\em quadratic tiling.} 
\FloatBarrier

 It follows that if we have a tiling of a triangle $ABC$ into $N$ congruent triangles, and $m$ is any 
integer,  we can tile $ABC$ into $Nm^2$ triangles by subdividing the first tiling, replacing each of the $N$ triangles by $m^2$ smaller ones.
Hence the set of $N$ for which an $N$-tiling of some triangle exists is closed under multiplication by squares.

\begin{figure}[ht]
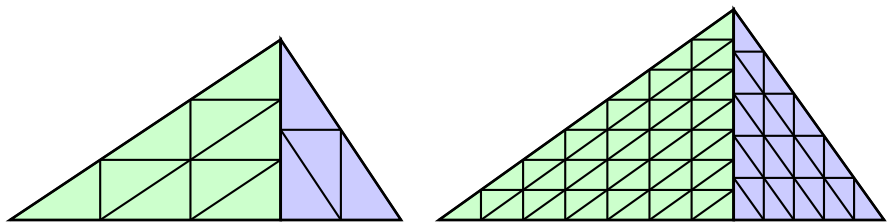

\caption{Biquadratic tilings with $N = 13 = 3^2 + 2^2$ and $N=74 = 5^2 + 7^2$}
\label{figure:biquadratic}
\ThirteenTiling
\FigureBiquadratic
\end{figure}

Sometimes it is possible to combine two quadratic tilings (using the same 
tile) into a single tiling, as shown in Fig.~\ref{figure:biquadratic}.
We will explain how these tilings are constructed.  We start with 
a big right triangle resting on its hypotenuse, and divide it into 
two right triangles by an altitude.  Then we quadratically tile each 
of those triangles.  The trick is to choose the dimensions in such a 
way that the same tile can be used throughout.  If that can be 
done then evidently $N$, the total number of tiles, will be the
sum of two squares, $N = n^2 + m^2$,  one square for each of the two 
quadratic tilings.  On the other hand, if we start with an $N$ of 
that form,  and we choose the 
tile to be an $n$ by $m$ right triangle, then we can construct 
such a tiling.  We call these tilings ``biquadratic.''
More generally, a {\em biquadratic tiling} of triangle $ABC$ is one in which $ABC$ has a right angle at $C$, and can be divided by an altitude 
from $C$ to $AB$ into two triangles, each similar to $ABC$, which can be tiled respectively by $n^2$ and $m^2$ copies of a triangle similar to $ABC$. 
A larger biquadratic tiling, with  $n=5$ and $m=7$ and hence $N=74$, is shown in at the right of  Fig.~\ref{figure:biquadratic}. 

\begin{figure}[ht]
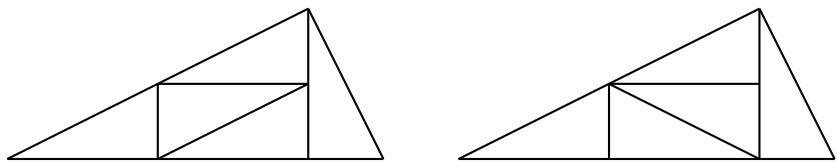

\caption{Two 5-tilings}
\label{figure:5-tilings}
\begin{center}
\FiveTiling  
\hskip 1 cm
\FiveTilingB
\end{center}
\end{figure}

\begin{figure}[ht]
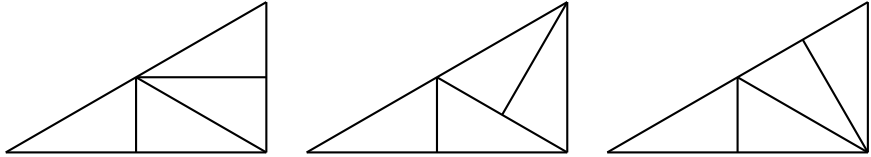
     
\caption{Three 4-tilings}
\label{figure:4-tilings}
\begin{center}
\FourTilingA
\FourTilingB
\FourTilingC
\end{center}
\end{figure}

Since $5 = 2^2 + 1^2$,  the simplest case of a biquadratic tiling is $N=5$. 
The second 5-tiling in Fig.~\ref{figure:5-tilings}  shows that a
biquadratic tiling can sometimes be  more complicated than a combination of two quadratic tilings.  Symmetry can permit rearranging some of the tiles.
The symmetrical tile used in Fig.~\ref{figure:4-tilings} also allows 
for variety.

If the original triangle $ABC$ is chosen to be isosceles, and is 
then quadratically tiled, 
then each of the $n^2$ triangles can be divided in half by an altitude;  hence any isosceles triangle can be decomposed into $2n^2$ congruent 
triangles.   If the original triangle is equilateral,  then it can be first decomposed into $n^2$ equilateral triangles, and then these 
triangles can be decomposed into 3 or 6 triangles each,  showing that any equilateral triangle can be decomposed into $3n^2$ or $6n^2$
congruent triangles. 
For example we can 12-tile an equilateral triangle in two different ways,  starting with a 3-tiling and then subdividing each triangle 
into 4 triangles (``subdividing by 4''),  or starting with a 4-tiling and then subdividing by 3.

\begin{figure}[ht]
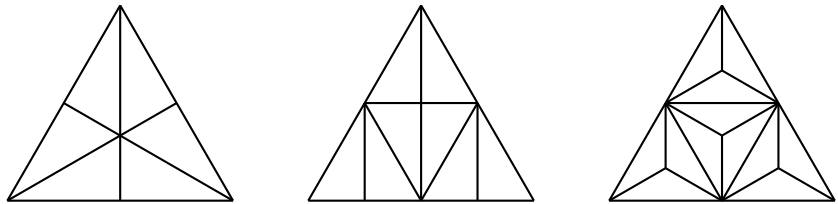
   
\caption{A 6-tiling, an 8-tiling, and a 12-tiling}
\label{figure:6-8-12-tilings}
\begin{center}
\SixTiling
\hskip 1cm
\EightTiling
\hskip 1cm
\TwelveTiling
\end{center}
\end{figure}

\begin{figure}[ht]
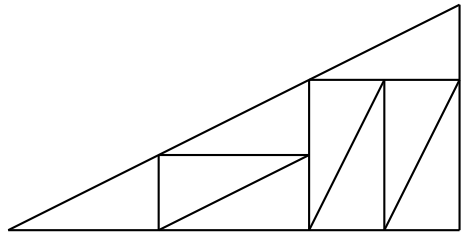
   
\caption{Another 9-tiling}
\label{figure:9-tiling}
\begin{center}
\NineTilingA
\end{center}
\end{figure}

The examples above do not exhaust all possible tilings, even when $N$ is a square.   For example, Fig.~\ref{figure:9-tiling} shows a 9-tiling that is not produced by 
those methods; again this seems attributable to symmetry permitting a 
rearrangement of tiles in a quadratic tiling. 

There is another family of $N$-tilings, in which $N$ is of the form $3m^2$, and both the tile and the tiled triangle are 30-60-90 triangles.  We call these the ``triple-square'' tilings. 
The case $m=1$ is given in Fig.~\ref{figure:3-tilings}; the case $m=2$ makes $N=12$.  There are two ways to 12-tile a 30-60-90 triangle with 30-60-90 triangle.  
One is to first quadratically 4-tile it, and then subtile the four triangles with the 3-tiling of Figure 1.  This produces the first 
12-tiling in Fig.~\ref{figure:12-tilings}.  Somewhat surprisingly, there is another way to tile the same triangle with the same 12 tiles, also shown in Fig.~\ref{figure:12-tilings}.
The next member of this family is $m=3$, which makes $N=27$.  
Two 27-tilings are shown in Fig.~\ref{figure:27-tilings}.
 Similarly, there are two 48-tilings (not shown).

\begin{figure}[ht]    
\caption{Two 12-tilings}
\label{figure:12-tilings}
\begin{center}
\TwelveTilingA
\hskip 2cm
\TwelveTilingB
\end{center}
\end{figure}

\begin{figure}[ht]
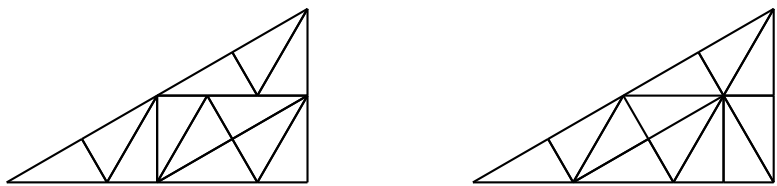
    
\caption{Two 27-tilings   }
\label{figure:27-tilings}
\begin{center}
\TwentySevenTilingA
\hskip 2cm
\TwentySevenTilingB
\end{center}
\end{figure}

\begin{figure}[ht]
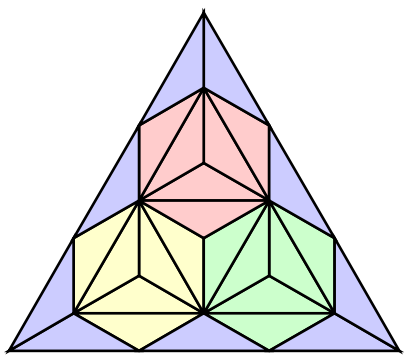
     
\caption{A 27-tiling due to Major MacMahon 1921, rediscovered 2011}
\label{figure:prime27}
\begin{center}
\HexTilingThree
\end{center}
\end{figure}

\begin{figure}[ht]
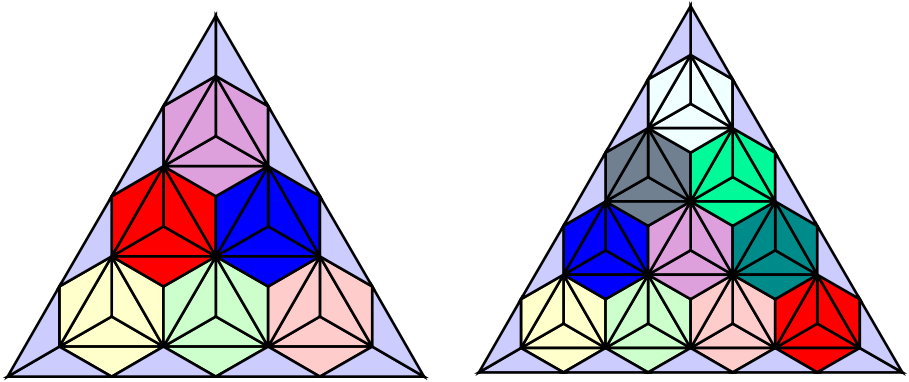
    
\caption{$3m^2$ (hexagonal) tilings for $m=4$ and $m=5$}
\label{figure:hexagonaltilings}
\begin{center}
\HexTilingFour
\HexTilingFive 
\end{center}
\end{figure}

Until October 12, 2008,  we did not know any more complicated tilings than those illustrated above (and
there also none in \cite{soifer}).    Then we found the 
beautiful 27-tiling shown in Fig.~\ref{figure:prime27}.  This tiling is one of a family of $3k^2$ tilings (the case $k=3$).
The next case is a 48-tiling, made from six hexagons (each containing 6 tiles) bordered by 
4 tiles on each of 3 sides.  In general one can arrange $1+2+\ldots+k$ hexagons in bowling-pin fashion, and add $k+1$ tiles 
on each of three sides, for a total number of tiles of $6(1+2+\ldots+k)+3(k+1) = 3k(k+1) + 3(k+1) = 3(k+1)^2$. 
Fig.~\ref{figure:hexagonaltilings} shows more members of this family, which we call the ``hexagonal tilings.''%
\footnote{In January, 2012, I bought a puzzle at the exhibition at the AMS meeting, which contained the 
tiling in Fig.~\ref{figure:prime27} as part of a tiling of a larger hexagon.  The tiling is attributed 
to Major Percy Alexander MacMahon (1854-1929) \cite{macmahon}.}

Whenever there is an $N$-tiling of the right triangle $ABM$, there is a $2N$-tiling of
the isosoceles triangle $ABC$.  Using the biquadratic tilings (see Fig.~\ref{figure:5-tilings}
and Fig.~\ref{figure:biquadratic}) and  triple-square tilings (see Fig~\ref{figure:12-tilings} and Fig.~\ref{figure:27-tilings}),
we can produce $2N$-tilings when $N$ is a sum of squares or three times a sum of squares.  We call these tilings ``double biquadratic''
and ``hexquadratic''.  
  For example,  one has two 10-tilings and two 26-tilings, obtained by reflecting Figs. 4 and 5 about either of the 
sides of the triangles shown in those figures;  and one has 24-tilings and 54-tilings obtained from Figs. 8 and 9.
Note that in the latter two cases, $ABC$ is equilateral.   

In the case when the sides of the tile $T$ form a Pythogorean triple $n^2 + m^2 + k^2 = N/2$, then we can tile one 
half of  $ABC$ with a quadratic tiling and the other half with a biquadratic tiling.  The smallest example is when the tile
has sides 3, 4, and 5, and $N = 50$.  See Fig.~\ref{figure:severalisosceles2}.
 One half is 25-tiled quadratically, and the other half is divided into two smaller
right triangles which are 9-tiled and 16-tiled quadratically.  This shows that the tiling of $ABC$ does not have to be 
symmetric about the altitude.
\FloatBarrier

\section{Definitions and notation}
We give a mathematically precise definition of ``tiling'' and fix some terminology and notation.
Given a triangle $T$ and a larger triangle $ABC$,  a ``tiling'' of triangle $ABC$ by triangle $T$ is 
a set of triangles $T_1,\ldots,T_n$ congruent to $T$, whose interiors are disjoint, and the closure of whose union is triangle $ABC$. 

Let $a$, $b$, and $c$ be the sides of the tile $T$, and angles $\alpha$, $\beta$, and $\gamma$ be the angles 
opposite sides $a$, $b$, and $c$.
The letter ``$N$'' will always be used for the number of triangles used in the tiling.  An $N$-tiling of $ABC$ is a tiling that uses
$N$ copies of some triangle $T$.  
The meanings of $N$, $\alpha$, $\beta$, $\gamma$, $a$, $b$,$c$, $A$, $B$, and $C$ will be fixed throughout this paper, and  
we assume $\alpha \le \beta \le \gamma$, when there is no 
other assumption about $\alpha$ and $\beta$, such as $3\alpha + 2\beta = \pi$.

\section{History}

In our gallery of examples, we saw quadratic and biquadratic tilings
in which the tile is similar to $ABC$, and also hexagonal tilings. These
involve $N$ being square, a sum of two squares, or three times a square.
The biquadratic tilings were known in 1964, when the 
paper \cite{golomb} was published.  This is the earliest paper on 
the subject of which I am aware.%
\footnote{The simplest hexagonal tiling
is attributed to Major MacMahon (1921) in the notes accompanying a 
plastic toy I purchased at an AMS meeting in 2012.} 
 Snover
{\em et. al.} \cite{snover1991} took up the challenge of showing that 
these are the only possible values of $N$.
The following theorem completely answers the question, ``for which $N$ does there exist an $N$-tiling
in which the tile is similar to the tiled triangle?''
\begin{theorem} [Snover {\em et. al.} \cite{snover1991}] 
\label{theorem:snover} 
 Suppose $ABC$ is $N$-tiled by tile $T$ similar to $ABC$.
If $N$ is not a square, then $T$ and $ABC$ are right triangles.  Then either 
\smallskip

(i) $N$ is three times
a square and $T$ is a 30-60-90 triangle, or 

(ii)  $N$ is a sum of squares $e^2 + f^2$,  the right angle of $ABC$ is split 
by the tiling, and the acute angles of $ABC$ have rational tangents $e/f$ and $f/e$,
\smallskip

\noindent
  and these two alternatives are 
mutually exclusive. 
\end{theorem}

Soifer's book \cite{soifer} appeared in 1990, with a second
edition in 2009.  He considered two ``Grand Problems'': for which 
$N$ can {\em every} triangle be $N$-tiled, and for which $N$ can 
{\em every} triangle be dissected into similar, but not necessarily
congruent triangles.  (The latter eventually became a Mathematics
Olympiad problem.)  The 2009 edition has an added chapter in 
which the biquadratic tilings and a theorem of Laczkovich occur.

Miklos Laczkovich published six papers \cite{laczkovich1990,laczkovich1995, laczkovich1998, laczkovich2010, laczkovich-szekeres1995, laczkovich2012}
 on triangle and polygon tilings.  According to Soifer, the 1995 paper
was submitted in 1992.  Laczkovich, like Soifer, 
studied dissecting a triangle into smaller {\em similar} triangles,
not {\em congruent} triangles as we require here.  
If those similar triangles are rational (i.e., the ratios of their sides
are rational) then if we divide each of them into small enough 
quadratic subtilings, we can achieve an $N$-tiling into {\em congruent}
triangles, but of course $N$ may be large.  Laczkovich proved little  
about $N$, focusing instead on the shapes of $ABC$ (or more 
 generally, convex polygons) and 
of the tile.  His theorems, for example, do not address the possibility 
of an $N$-tiling (of some $ABC$ by some tile) for any particular $N$,
but they do give us an exhaustive list of the possible shapes of $ABC$
and the tile, which we will need in our proof that there is no 7-tiling.
This list can be found in \S\ref{section:laczkovich} (of this paper).  
However, his theorem published in the last chapter of \cite{soifer}
does mention
$N$.  It states that given an integer $k$,  there exists an $N$-tiling for 
some $N$ whose square-free part is $k$.   

\section{Laczkovich}\label{section:laczkovich}

A basic fact is that, apart from a small number of cases that can be 
explicitly enumerated, if there is an $N$-tiling of $ABC$ by a tile 
with angles $(\alpha, \beta, \gamma)$, then
 the angles $\alpha$ and $\beta$ are not rational 
multiples of $\pi$.  This theorem is 
Theorem~5.1 of \cite{laczkovich1995}. Laczkovich calls
the angles of the tile commensurable if 
 each of them is a rational multiple of $\pi$.  He states his 
theorem conversely to the way we just described it:  
 if there is a tiling of $ABC$ by a tile $T$ with 
commensurable angles, then the pair $(ABC, T)$ belongs to a specific,
fairly short list.  It is important to note that Laczkovich's list
in Theorem~5.1 is about dissections of $ABC$ into {\em similar},
not necessarily congruent, triangles.  His subsequent Theorem~5.3 
shows that three possibilities for dissecting the right isosceles
triangle $ABC$ into similar triangles are impossible with 
congruent tiles.

Laczkovich's list of 
possibilities from the cited 
1995 paper is given in Table~\ref{table:laczkovich-original}.   In the table, 
the triples giving the angles of the tile are $(\alpha, \beta, \gamma)$
after a suitable permutation, i.e., they are unordered triples. 
The reader who checks with \cite{laczkovich1995} will need to remember
that we have deleted the entries for the right isosceles $ABC$ mentioned
above.   

\begin{table}[ht]
\caption{Laczkovich's 1995 list of tilings by tiles with commensurable angles}
\label{table:laczkovich-original}
\begin{center}
\setlength{\extrarowheight}{.5em}
\begin{tabular}{rr}
$ABC$  &  the tile     \\
\hline
$(\alpha,\beta,\gamma)$  &  similar to $ABC$ \\
$(\alpha,\alpha,2\beta)$ &  $\gamma = \pi/2$ \\
equilateral  & {\large $(\frac \pi 6, \frac \pi 6, \frac {2\pi} 3)$ } \\
equilateral  & {\large$(\frac \pi 3, \frac \pi {12}, \frac {7\pi}{12})$} \\
equilateral & {\large$( \frac \pi 3, \frac \pi {30}, \frac {19\pi}{30})$}\\
equilateral & {\large$( \frac \pi 3, \frac {7\pi} {30}, \frac {13\pi}{30})$}
\end{tabular}
\end{center}
\end{table}

In subsequent work, specifically Theorem~3.3 of \cite{laczkovich2012},
Laczkovich proved that the table can be considerably shortened: the 
tilings of the equilateral triangle mentioned in the last three rows
cannot occur (when the tiles are required, as in this paper, to 
be congruent rather than just similar).  Thus the final version
is as shown in 

\begin{table}[ht]
\caption{Laczkovich's 2012 list of tilings by tiles with commensurable angles}
\label{table:laczkovich}
\begin{center}
\setlength{\extrarowheight}{.5em}
\begin{tabular}{rr}
$ABC$  &  the tile     \\
\hline
$(\alpha,\beta,\gamma)$  &  similar to $ABC$ \\
$(\alpha,\alpha,2\beta)$ &  $\gamma = \pi/2$ \\
equilateral  & {\large $(\frac \pi 6, \frac \pi 6, \frac {2\pi} 3)$ }
\end{tabular}
\end{center}
\end{table} 

It is possible to prove by direct computation that the last 
two rows of Table~\ref{table:laczkovich-original} do not correspond 
to actual tilings.   Namely, the area equation for the equilateral 
triangle with side $X$ tells us $X^2 = Nbc$, if angle $\alpha = \pi/3$.
Then writing $X = pa + qb + rc$ and calculating  
$(a,b,c) = (\sin \alpha, \sin \beta, \sin \gamma)$ for the specific 
angles involved, we get equations in certain algebraic number fields,
that one then has to show impossible.  For example,
\begin{eqnarray}\left( p\left(\xi - \frac {1 + \sqrt 5} 8\right) 
     + q \frac {\sqrt 3} 2 + 
r\left(\xi + \frac {1 + \sqrt 5} 8\right)\right)^2 &=& N\bigg(\frac 3 8  - \frac {\sqrt 5} 8 \bigg)  \label{eq:1797}
\end{eqnarray}
One interesting thing about this approach is that SageMath is fully
capable of performing all the required calculations, including determining
whether certain expressions lie in certain algebraic number fields or not.
We did not succeed, however, in entirely eliminating the row 
mentioning $\pi/12$ by computation; in that case, using the area equation as described
only tells us that $N$ is six times a square.  That would be enough 
for this paper, where we only need that $N$ cannot be 7 or 11; but 
Laczkovich entirely eliminated that possible tiling as well as the other two.

\section{The coloring equation}
\begin{figure} [ht]
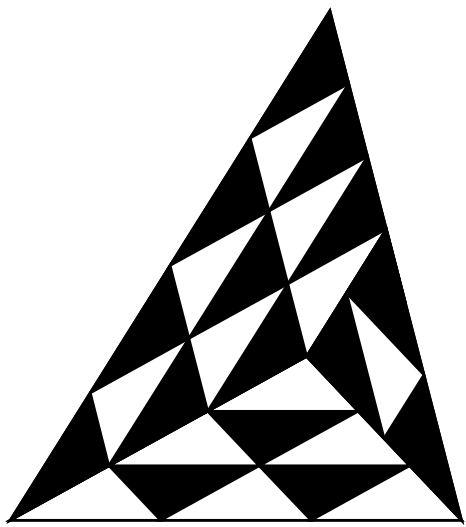

\caption{A tiling colored so that touching tiles have different colors.}
\label{figure:coloring}
\begin{center} \ColoringTheorem
\end{center}
\end{figure}

In this section we introduce a tool that is useful for some, but not all, 
tiling problems.  Suppose that triangle $ABC$ is tiled
by a tile with 
angles $(\alpha,\beta,\gamma)$ and sides $(a,b,c)$, and suppose
there is just one tile at vertex $A$.  We color that tile black, 
and then we color each tile black or white, changing colors as 
we cross tile boundaries.  Under certain conditions this coloring
can be defined unambiguously, and then, we define the ``coloring number''
to be the number of black tiles minus the number of white tiles.
An example
of such a coloring is given in Fig.~\ref{figure:coloring}.

The following theorem spells out the conditions under which this 
can be done.  In the theorem, ``boundary vertex'' refers to a vertex
that lies on the boundary of $ABC$ or on an edge of another tile,
so that the sum of the angles of tiles at that vertex is $\pi$.
At an ``interior vertex'' the sum of the angles is $2\pi$.

\begin{theorem} \label{theorem:coloring}
Suppose that triangle $ABC$ is tiled by the tile $(a,b,c)$ in 
such a way that 
\smallskip

(i) There is just one tile at $A$.
\smallskip

(ii) At every boundary vertex an odd number of tiles meet.
\smallskip

(iii) At every interior vertex an even number of tiles meet.
\smallskip

(iv) The numbers of tiles at $B$ and $C$ are both even, or both odd.
\smallskip

Then every tile can be assigned a color (black or white)
 in such a way that 
colors change across tile boundaries, and the tile at $A$ is black.
  Let $M$ be the 
number of black tiles minus the number of white tiles.
Then the coloring equation 

$$ X \pm Y + Z = M(a+b+c)$$
holds, where $Y$ is the side of $ABC$ opposite $A$, and $X$ and $Z$
are the other two sides.  The sign is $+$ or $-$ according as the 
number of tiles at $B$ and $C$ is odd or even.
\end{theorem}

\noindent{\em Proof}.  Each tile is colored black or white
according as the number of tile boundaries crossed in reaching it 
from $A$ without passing through a vertex is even or odd.  The 
hypotheses of the theorem guarantee that color so defined 
is independent of the path chosen to reach the tile from $A$.   
The total length of black edges, minus the total length of white 
edges,  is $M(a+b+c)$, since $a+b+c$ is the perimeter of each tile.
Each interior edge makes a contribution 
of zero to this sum, since it is black on one side and white on the other.
Therefore only the edges on the boundary of $ABC$ contribute. 
Now sides $X$ and $Y$ contain only edges of black tiles, by 
hypotheses (i) and (ii).  Side $Y$ is also black if the
number of tiles at $B$ and $C$ is odd,  and white if it is even.
Hence the difference in the total length of black and white tiles
is $X \pm Y + Z$, with the sign determined as described.  That completes
the proof.

\section{Possible values of $N$ in tilings with commensurable angles}
We wish to add a third column to Laczkovich's Table~\ref{table:laczkovich},
giving the possible forms of $N$ if there is an $N$-tiling of $ABC$
by the tile in that row.  For example, when $ABC$ is similar to the 
tile, then $N$ must be a square, so we put $n^2$ in the third column.
While we are at it, we add a fourth column with a citation to the 
result,  and delete the rows corresponding to the tilings
of the equilateral triangle that we have proved impossible.
  The revised and extended table is Table~\ref{table:laczkovich-extended}.
All the entries in this table except the last one give necessary and
sufficient conditions on $N$ for the tilings to exist.  The last one
gives necessary conditions for certain tilings that probably do not 
actually exist.

\begin{table}[ht]
\caption{$N$-tilings by tiles with commensurable angles, with form of $N$}
\label{table:laczkovich-extended}
\begin{center}
\setlength{\extrarowheight}{.5em}
\begin{tabular}{rrrr}
$ABC$  &  the tile & form of $N$ & citation    \\
\hline
$(\alpha,\beta,\gamma)$  &  similar to $ABC$ & $n^2$ & \cite{snover1991} \\
$(\alpha,\beta,\gamma)$  & similar to $ABC$,  $\gamma = \pi/2$  & $e^2 + f^2$ & \cite{snover1991} \\
$(\frac \pi 6, \frac \pi 3, \frac \pi 2)$  &  similar to $ABC$  & $3n^2$ & \cite{snover1991} \\
$(\alpha,\alpha,2\beta)$ &  $\gamma = \pi/2$ & $2n^2$ & \cite{beeson-isosceles}\\
$(\alpha,\alpha,2\beta)$ &  $(\frac \pi 4,\frac \pi 4, \frac \pi 2)$ & $n^2$ & \cite{beeson-isosceles}\\
{\large $(\frac \pi 6, \frac \pi 6, \frac {2 \pi} 3)$} &{\large  $(\frac \pi 6,\frac \pi 3, \frac \pi 2)$} & $6n^2$ & \cite{beeson-isosceles}\\
equilateral  & {\large $(\frac \pi 6, \frac \pi 3, \frac \pi 2)$} & $6n^2$ & \cite{beeson-isosceles} \\
equilateral  & {\large $(\frac \pi 6, \frac \pi 6, \frac {2\pi} 3)$ } 
                  & $3n^2$ & \cite{beeson-isosceles} 
\end{tabular}
\end{center}
\end{table}

\begin{theorem} \label{theorem:laczkovich-extended}
Suppose $(\alpha,\beta,\gamma)$ are all rational multiples of 
$2\pi$, and triangle $ABC$ is $N$-tiled by a tile with angles 
$(\alpha,\beta,\gamma)$.  Then $ABC$, $(\alpha,\beta,\gamma)$,
and $N$ correspond to one of the lines in Table~\ref{table:laczkovich-extended}.
\end{theorem}

\noindent{\em Proof}.  As discussed above, Laczkovich characterized
the pairs of tiled triangle and tile, as given in
Table~\ref{table:laczkovich}.%
\footnote{ Again, we remind readers who may check with 
\cite{laczkovich1995} that there are three entries in Laczkovich's 
Theorem~5.1 that are shown in the subsequent Theorem~5.3 not to apply
to tilings by congruent triangles, so they do not appear in our tables.
}
It remains to characterize the 
possible $N$ for each line.  In several cases lines in 
Table~\ref{table:laczkovich} split into two or more lines in 
Table~\ref{table:laczkovich-extended}, which supplies the required
possible forms of values of $N$.  That table
 lists in its last column citations
to the literature or theorems in this paper for each line.  Finally,
we have deleted the rows of Table~\ref{table:laczkovich}
corresponding to the tilings that are impossible 
by Theorem~3.3 of \cite{laczkovich2012}.
 That completes the proof.

\section{Laczkovich's second table}
Laczkovich also studied the case when not all the angles of the tile
are rational multiples of $\pi$. Again a finite number of cases
can arise.  This is Theorem~4.1 of \cite{laczkovich1995}, and the 
list of cases is given in Table~\ref{table:laczkovich2}.

\begin{table}[ht]
\caption{Tilings when not all angles are rational multiples of $\pi$.}
\label{table:laczkovich2}
\begin{center}
\setlength{\extrarowheight}{.5em}
\begin{tabular}{rr}
$ABC$  &  the tile     \\
\hline
$(\alpha,\beta,\gamma)$  &  similar to $ABC$ \\
equilateral              & $\alpha = \pi/3$ \\
$(\alpha,\alpha,2\beta)$ &  $\gamma = \pi/2$ \\
$(\alpha,\alpha,\pi-2\alpha)$  & $\gamma = 2\alpha$\\
$(2\alpha,\beta,\alpha+\beta)$ & $3\alpha + 2\beta = \pi$\\
$(2\alpha,\alpha,2\beta)$ &$3\alpha + 2\beta = \pi$\\
isosceles &$3\alpha + 2\beta = \pi$\\
$(\alpha,\alpha,\pi-2\alpha)$ & $\gamma = 2\pi/3$ \\
$(\alpha,2\alpha,\pi-3\alpha)$ & $\gamma = 2\pi/3$ \\
$(\alpha, 2 \beta, 2\alpha + \beta)$ & $\gamma = 2\pi/3$\\
$(\alpha, \alpha+\beta,\alpha+2\beta)$ & $\gamma = 2\pi/3$\\
$(2\alpha, 2\beta, \alpha+\beta)$ & $\gamma = 2\pi/3$ \\
equilateral  & $\gamma = 2\pi/3$ 
\end{tabular}
\end{center}
\end{table}

Table~\ref{table:laczkovich} and~\ref{table:laczkovich2} together
constitute an exhaustive list of tilings.  If we have some conditions
on the tile, such as for example $3\alpha + 2\beta= \pi$, then 
we look to see what entries in Table~\ref{table:laczkovich} satisfy 
those conditions.  That gives some tilings with commensurable angles.
Then we look in the other table for tilings in which not all the angles 
are rational multiples of $2\pi$.  To fix the ideas we spell
out the details for the case $3\alpha + 2\beta = \pi$.

\begin{lemma} \label{lemma:angles}
Let $3\alpha+2\beta = \pi$. Suppose there is an $N$-tiling of triangle $ABC$
by tile $T$ with angles
$(\alpha,\beta,\gamma)$.  Suppose also that $ABC$ is not similar to $T$.
Then
  $\alpha$ and $\beta$ are not rational 
multiples of $\pi$, and every linear relation between $\pi$, $\alpha$,
and $\beta$ is a multiple of $3\alpha + 2\beta = \pi$.  
\end{lemma}

\noindent{\em Proof}.  Suppose there is an $N$-tiling as in the 
statement of the lemma.  Then if angles of the tile are all 
rational multiples of $\pi$, the pair $ABC$ and the tile must 
occur in Table~\ref{table:laczkovich}. 
 So we have to check if 
any of the triples in that table satisfy $3\alpha+2 \beta = \pi$.
And they do not, so that completes the proof.
\medskip

\noindent{\em Remark}.  The reader of Laczkovich's paper \cite{laczkovich1995}
should beware:  Theorem~5.1 includes the triple
$(\pi/4,\pi/8,5\pi/8)$, which does satisfy
 $3\alpha + 2\beta = \pi$.  But as discussed above, it is included
 since the theorem is about dissections into {\em similar} triangles,
 and Theorem~5.3 of \cite{laczkovich1995} rules it out for tilings into
{\em congruent} triangles.  Hence we have deleted it from Table~\ref{table:laczkovich} and do not need to consider it here.

\section{Adding a column for $N$ to Laczkovich's second table}
The research program that we have been pursuing in this subject 
is to study triples $(ABC, T, N)$ instead of just pairs $(ABC, T)$,
where there is an $N$-tiling of $ABC$ by tile $T$.  Another way to 
say that is that we wish to add a third column to Laczkovich's 
second table, entering the possible forms of $N$ in that column,
as we did to the first table.  This has proved to be a longer
business than I had originally imagined, although also more interesting, 
since several new tilings have been discovered in the process, and 
this research program is not complete.  The point of the present
paper is that we have pursued it far enough to reach the goal 
of showing that 7-tilings are impossible.  Presently, we can supply
entries in the third column down to the cases with $\gamma = 2\pi/3$,
but some of them are only necessary conditions, not necessary and sufficient,
leaving open many questions about particular values of $N$ that are not
ruled out by those necessary conditions.

\section{Some number-theoretic facts}
The facts in this section may not be 
well-known to all our readers, so we collect them here with short proofs or citations. 

\begin{lemma}\label{lemma:sumsofsquares} 
An integer $N$ can be written as a sum of two integer squares if and only if the 
squarefree part of $N$ is not divisible by any prime of the form $4n+3$.
\end{lemma}

\noindent{\em Proof}.  See for example \cite{hardy-wright}, Theorem~366, p.~299.

\begin{lemma}\label{lemma:qsums} A quotient of sums of two rational squares
is a sum of two rational squares.
\end{lemma}

\noindent{\em Proof}. A sum of two rational squares is the square of the
absolute 
value   of some complex number. The quotient of the absolute values is the 
absolute value of the quotient.  Explicitly:
\begin{eqnarray*}
\frac {a^2 + b^2}{c^2 + d^2} &=& \frac {\vert a + bi \vert^2}{\vert c + di \vert^2}\\
&=& \left\vert \frac {a+bi}{c+di} \right\vert^2 \\
&=& \left\vert \frac {(a+bi)(c-di)}{c^2 + d^2} \right\vert^2 \\
&=& \left( \frac {ac + bd}{c^2+d^2} \right)^2 + 
    \left( \frac {bc -ad}{c^2+d^2} \right)^2
\end{eqnarray*}
That completes the proof of the lemma.
\medskip

The following lemma identifies those relatively few rational multiples of $\pi$ that have rational tangents or 
whose sine and cosine satisfy a polynomial of low degree over $\Q$.

\begin{lemma}  \label{lemma:euler} 
Let $\theta = 2m \pi/n$,  where $m$ and $n$ have no common factor.  
Suppose $\cos \theta$ is algebraic of degree 1 or 2 over $\Q$.  
Then $n$ is one of $5,6,8,10,12$.   If both $\cos \theta$ and $\sin \theta$ have
degree 1 or 2 over $\Q$, then $n$ is $6,8$, or $12$.
\end{lemma}

\noindent{\em Proof}.  Let $\varphi$ be the Euler totient function.  
Assume $\cos \theta$  has degree 1 or 2.  By \cite{niven},  Theorem~3.9, p.~37,
$\varphi(n) = 2$ or $4$.  The stated conclusion follows from the well-known formula for $\varphi(n)$.
The second part of Theorem~3.9 of \cite{niven} rules out $n=5$ or $10$ when $\sin \theta$ is 
also of degree 1 or 2.

\section{Isosceles $ABC$ tiled by a right triangle}
\begin{figure}[ht]
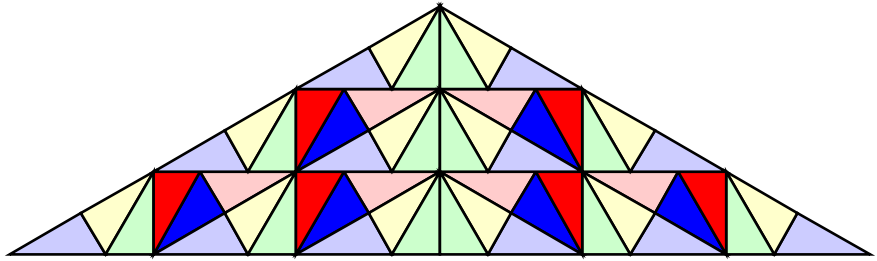

\begin{center}
\caption{A 54-tiling; $N/2$ is three times a square. Tile is 30-60-90.}
\label{figure:54}
\FiftyFourTiling
\end{center}
\end{figure}

The possible ways of tiling an isosceles triangle  
have been studied in \cite{beeson-isosceles}.  
Let the angles of the tile be  $(\alpha,\beta,\gamma)$.
Laczkovich's
results imply that there are only four cases to consider:
$\gamma =  \pi/2$,  $\gamma = 2\alpha$, $2\alpha + 3\beta = \pi$,
 or $\gamma = 2\pi/3$,
with $\alpha$ not a rational multiple of $\pi/2$ in the last two 
cases.   The first three cases have been studied in
\cite{beeson-isosceles}.  Here is the result for the case of 
a right-angled tile:

\begin{theorem} \label{theorem:isosceles2}
  Suppose isosceles triangle $ABC$ with base angles 
$\beta$ is $N$-tiled by a right-angled tile.
Then $N$ is a square, or a sum of two squares, 
or six times a square.   
\end{theorem}

\begin{corollary}  $N$ is not a prime congruent to 3 mod 4;
in particular it is not 7, 11, or 19; nor can $N$ be twice such a prime.
\end{corollary}

\noindent{\em Proof of Corollary}.  Primes congruent to 3 mod 4
cannot be sums of squares, by Lemma~\ref{lemma:sumsofsquares}.
\medskip

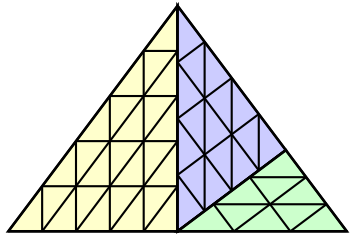
\begin{figure}[ht]
\caption{$N$ is a twice a square or a twice a sum of squares. 50 is both. }
\label{figure:severalisosceles}
\begin{center}
\FigureDoubleQuadratic
\FigureDoubleBiquadratic
\end{center}
\end{figure}

\begin{figure}[ht]
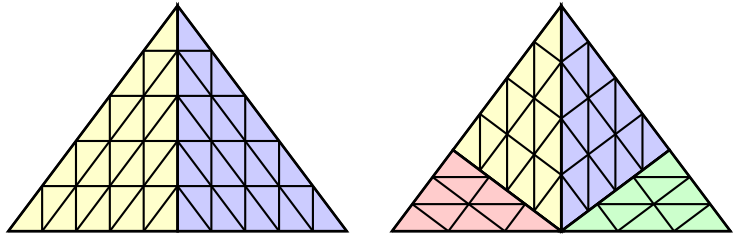

\caption{50 is both twice a square and twice a sum of squares. }
\label{figure:severalisosceles2}
\begin{center}
\FigurePythagorean
\end{center}
\end{figure}

{\em Remark}.  For $N$ even, all the listed possibilities can occur.
See Figs.~\ref{figure:severalisosceles} and \ref{figure:severalisosceles2},
and additional examples in \cite{beeson-isosceles}.
We conjecture that there are no such tilings with $N$ odd.

{\em Remark}. Our original plan was to avoid citing 
\cite{beeson-isosceles} by exhibiting simple SageMath code
that verifies directly that there are no such tilings for $N = 7, 11, 14, 19$,
by showing that the area equation $(pa + qb + r)^2 = N/2$ together
with the right-angle equation $a^2+b^2 = 1$ has no solutions in 
non-negative integers
$(p,q,r)$.  
Although this code is straightforward, to justify it we had to prove
two lemmas: no tile except the one at the vertex $B$ can touch both
equal sides, and we can assume $p+q+r < N/3$.  In total then 
the computational approach is not much shorter than 
the proof in \cite{beeson-isosceles}. 

 \section{Useful lemmas}
In this section, we collect some facts that will be applied when 
we start eliminating the possibilities for 7-tilings and 11-tilings
case by case, according to the cases of Laczkovich's second table.

\subsection{Angles}
\begin{lemma}\label{lemma:gammagreaterpiover2}  Suppose triangle $ABC$ is $N$-tiled by a tile in which $3\alpha + 2\beta = \pi$.   Then 
$\gamma > \pi/2$. 
\end{lemma}

\noindent{\em Proof}.  
\begin{eqnarray*}
\pi &=& 3\alpha + 2\beta \\
    &=& \alpha + 2(\alpha + \beta) \\
    &=& \alpha + 2 (\pi - \gamma) \\
 \gamma &=& \frac \pi 2 + \frac \alpha 2 \ > \ \frac \pi 2
\end{eqnarray*}
That completes the proof.

\begin{lemma} \label{lemma:possibleangles}
Let triangle $ABC$ be $N$-tiled by a tile with angles 
$(\alpha,\beta,\gamma)$.  Suppose that either $3\alpha + 2\beta = \pi$
and $ABC$ is not isosceles with base angles $\alpha$, 
or $\gamma = 2\pi/3$.  Then no tile has its $\gamma$ angle 
at a vertex of $ABC$.
\end{lemma}

\noindent{\em Proof}.  
 By Lemma~\ref{lemma:angles}, $\alpha$ and $\beta$
are not rational multiples of $\pi$.  Hence the 
angles of $ABC$ are linear integral combinations of 
$\alpha$, $\beta$, and $\gamma$.  First assume
$3\alpha + 2\beta = \pi$.  Then 
the angles of $ABC$ are each
equal to $\alpha$, $2\alpha$, $\alpha+\beta$, $\beta$, or 
$2\beta$. Of these angles, all but $2\beta$ are less than $\gamma$,
as we now show.   Then 
$\gamma = \beta + 2\alpha$, and
\begin{eqnarray*}
\alpha &<& \beta + 2\alpha\ =\ \gamma \\ 
\beta &<& \beta + 2\alpha \ < \ \gamma \\
\alpha + \beta &<& \beta + 2\alpha \ = \ \gamma\\
2\alpha &<& \beta + 2\alpha \ = \ \gamma.
\end{eqnarray*} 

 Since $ABC$ is not similar to the tile, 
there cannot be a $\gamma$ angle alone at any vertex, since that 
would leave $\alpha + \beta$ for the other two vertices, making 
$ABC$ similar to the tile, since $\alpha$ is not a rational
multiple of $\beta$.  

   Since all the possible angles but $2\beta$
are less than $\gamma$, it only remains to deal with the case where
angle $C$ is equal to $2\beta$ and 
$\gamma < 2\beta$, and there is a tile with its $\gamma$ angle at $C$.  
We do not have $2\beta = \gamma$, by Lemma~\ref{lemma:angles}.
Then there must be another tile at $C$ as well.  If the angle 
of that tile at $C$ is $\alpha$, then the total angle at $C$
is at least $\gamma + \alpha = 2\alpha + \beta + \alpha =  3\alpha + \beta$,
leaving only $\beta$ for the other two angles of $ABC$. But that 
is impossible, since $\alpha$ is not a rational multiple of $\beta$.
 If the second angle at $C$
is $\beta$, then the total angle at $C$ is at least
 $\gamma + \beta = 2\alpha + 2\beta$,  leaving just $\alpha$
for the other two angles, which is again impossible.
 Hence the second angle at $C$ cannot be $\beta$.
That completes the proof under the assumption $3\alpha + 2\beta = \pi$.

We now take up the case $\gamma = 2\pi/3$. Then 
the possible angles of $ABC$ are
$\alpha$, $\beta$, $\alpha+\beta$, $\alpha + 2\beta$,
$2\alpha+\beta$, $3\alpha$, and $3\beta$.  All but $3\alpha$ and 
$3\beta$ are less than $2\alpha + 2\beta = \gamma$, so a $\gamma$
tile can occur, if at all, only at a vertex angle of $3\alpha$ or 
$3\beta$.  Suppose vertex $C$ has angle $3\alpha$ and there is a 
$\gamma$ angle of a tile at $C$.  Then $\gamma < 3\alpha$ and
angles $A$ and $B$ together
are $\pi-3\alpha < \pi-\gamma$, which is impossible since the three
angles of $ABC$ add up to $\pi$.  Similarly if vertex $C$
has angle $3\beta$ and $\gamma < 3\beta$. 
That completes the proof of the lemma. 

\subsection{Two $c$ edges on each side of $ABC$}

\begin{lemma} \label{lemma:zerolimits1} Suppose triangle $ABC$ is $N$-tiled by a tile with angles $(\alpha,\beta,\gamma)$
and $\gamma > \pi/2$. 
 Suppose all the tiles along one side of $ABC$ do not have their $c$ sides
along that side of $ABC$.  Then there is a tile with a $\gamma$ angle at one of the endpoints of that side of $ABC$.
\end{lemma}

\noindent{\em Proof}.  
 Let $PQ$ be the side of $ABC$ with no $c$ sides of tiles along it.  Then the 
$\gamma$ angle of each of those tiles   occurs at a vertex on $PQ$, since the angle opposite the side 
of the tile on $PQ$ must be $\alpha$ or $\beta$.  Let $n$ be the number of tiles along $PQ$;  then there 
are $n-1$ vertices of these tiles on the interior of $PQ$.  
Since $\gamma > \pi/2$, no
 vertex  on the boundary has more 
than one $\gamma$ angle.  By the pigeonhole principle, there is at least one tile whose $\gamma$ angle is 
not at one of those $n-1$ interior vertices;  that angle must be at $P$ or $Q$.  That completes the proof of the lemma.
\medskip

\begin{lemma} \label{lemma:zerolimits2}
Suppose triangle $ABC$  is 
$N$-tiled by a tile $T$ with angles $(\alpha,\beta,\gamma)$.  Suppose 
\smallskip

(i) $\gamma > \pi/2$, and
\smallskip
 
(ii) $\alpha$ is not a rational multiple of $\pi$, and 
\smallskip

(iii) every angle of triangle $ABC$ is less than $\gamma$.
\smallskip

Then there 
are at least two $c$ edges of tiles on side $AC$.  
\end{lemma}

\noindent{\em Remarks.}  One can prove by the same method that the 
$c$ edges must occur in adjacent blocks of at least two edges, but
we found no use for that result. 
\medskip

\begin{figure}[ht]
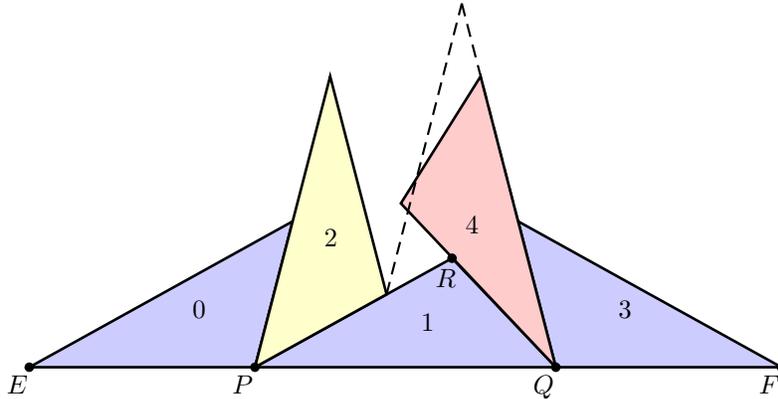

\caption{Proof of Lemma~\ref{lemma:zerolimits2}: another tile won't fit next to Tile~2}
\label{figure:zerolimits2}
\begin{center}
\FigureZeroLimitsTwo
\end{center}
\end{figure}

\noindent{\em Proof}.  By hypothesis (ii), every boundary vertex $P$
(except $A$, $B$, and $C$) that has a $\gamma$ angle (i.e., some 
tile with a vertex at $P$ has its $\gamma$ angle at $P$) 
touches exactly three tiles, which contribute
angles of $\alpha$, $\beta$, and $\gamma$.  By  
Lemma \ref{lemma:zerolimits1}, each side of $ABC$
has at least one $c$ edge. 
The present lemma, however, claims more: there must be 
at least two $c$ edges.  Suppose, to the contrary, that there is 
just one $c$ tile, Tile~1, with an edge on one side  $EF$ of triangle $ABC$.  Then all the 
other tiles with an edge on $EF$ have a $\gamma$ angle on $EF$. 
We visualize $EF$ as horizontal with triangle $ABC$ above, 
and use the word ``north'' and ``northwest'' accordingly. 
See Fig.~\ref{figure:zerolimits2}.

\noindent
Since there cannot be a $\gamma$ angle at the vertices of $ABC$,
it follows that both the tiles on $AC$ adjacent to Tile~1 (if there
are two, or otherwise, only the one) have
their $\gamma$ angles adjacent to Tile~1.  Let $PQ$ be the 
$c$ edge of Tile~1 lying on $AC$.  Let $R$ be the 
northern vertex of Tile~1.  Suppose (without loss 
of generality) that Tile~1 has its $\beta$ angle at $Q$. Then 
the side $PR$ of Tile~1, opposite $Q$,  has length $b$.  Let Tile~2 
be the tile adjacent to $PR$.  

Since the hypotheses of the theorem remain true if (the names of)
$\alpha$ and $\beta$
are interchanged, we may assume without loss of generality that 
$\alpha < \beta$.  Then by the law of sines, $a < b$.  Since $\gamma > \pi/2$
we also have $a < c$ (by the law of cosines).   

Assume, for proof by contradiction, that neither $P$ nor $Q$ is a 
vertex of $ABC$.  Then there exist Tile~0 and Tile~3 on $AC$
sharing vertices $P$ and $Q$ with Tile~1.  
 Tile~2, between Tile~0 and Tile~1, must 
have its $\beta$ angle at $P$, since Tile~1 has its $\alpha$ angle 
there and Tile~0 has its $\gamma$ angle at $P$.
There is an open
$\alpha$ angle between Tile~1 and Tile~3; let Tile~4 be the tile 
that fills that notch.  Then Tile~4 has its $b$ or $c$ edge 
along $QR$.  Since Tile~1 has its $a$ edge along $QR$ and $a < b$
and $a < c$, the edge of Tile~4 on $QR$ extends past $R$.  Then 
the segment $PR$ is of length $b$ and its northwest side is 
composed of a number of tile edges, starting with Tile~2 at $P$.
These must all be $a$ edges, since $a$ is the only edge less than 
$b$.  Since the tiles northwest of $PR$ all have their $a$ edges
on $PR$, they all have a $\gamma$ angle on $PR$.  But Tile~2
does not have its $\gamma$ angle at $P$, since Tile~0 has its 
$\gamma$ angle at $P$.  And the last tile cannot have its $\gamma$
angle at $R$, since Tile~4 extends along $QR$ past $R$, and 
Tile~1 has its $\gamma$ angle at $R$.  So if there are $n$ tiles
northwest of $PR$, there are only $n-1$ possible places for their 
$\gamma$ angles, contradicting the pigeon-hole principle.  
This contradiction proves that one of $P$ or $Q$ is a vertex of $ABC$.
\smallskip

Now we argue by cases. 
\smallskip

Case 1:  $Q$ is a vertex of $ABC$, i.e., $Q=C$. 
If the angle of $ABC$ at $Q$ is strictly between $\beta$ and $2\beta$,
 then Tile~4 must have its 
 $\alpha$ angle at $Q$, and we argue exactly as before.
If  the angle 
of $ABC$ at $Q$ is exactly $\beta$, then we argue as above, except 
that $RQ$ is now extended past $R$ by one side of $ABC$ rather than 
an edge of Tile~4.  The argument about the $\gamma$ angles of the 
tiles northwest of $PR$ is unchanged, if $P$ is not a vertex
of $ABC$.  If $P$ is a vertex of $ABC$, then we still 
can argue that Tile~2 must have its $a$ side on $PR$, because it cannot fit
next to Tile~1 with its $b$ or $c$ side on $PR$.  

Therefore we may assume that the angle of $ABC$ at $Q$ is at least $2\beta$,
and that Tile~4 has its $\beta$ angle at $Q$ and its $a$ edge against 
Tile~1.  Hence there is a double angle at $Q$.  Then by 
hypothesis (iv), $b$ is not a multiple of $a$. 
Tile~4 cannot have its $\gamma$ angle at $Q$, by hypothesis~(iii).
Therefore Tile~4 has its $\gamma$ angle at $R$, and since $\gamma > \pi/2$
by hypothesis~(i), $PR$ does not extend past $R$ as part of the tiling.
The tiles northwest of $PR$ must all have their $a$ edges on $PR$, 
since $a$ is the only edge less than $b$.  Similarly, the tiles 
supported by the west edge of Tile~4 must all have their $a$ edges
against that west edge, which has length $b$.  All 
those tiles northwest of $PR$ have their $\gamma$ angles on $PR$
(since they have their $a$ edges on $PR$), and by the pigeonhole
principle those $\gamma$ angles are all 
at the northwest. Therefore the tile supported by $PR$ at $R$,
call it Tile~5, 
has its $\gamma$ angle there. Since Tiles~1 and 4 already have 
their $\gamma$ angles at $R$, Tile~5 shares
an edge with Tile~4, and as just shown that edge has length $a$.
But now Tile~5 has two $a$ edges, contradiction.
 That completes Case~1.     
\smallskip

Case 2:  $P$ is a vertex of $ABC$, and $Q$ is not.  Then Tile~4
is placed as shown in the figure. Therefore the 
angle of $ABC$ at vertex $P$ must be greater than $\alpha$, 
since if it were equal to $\alpha$, Tile~4 would not lie 
inside $ABC$.    Then Tile~2 exists, 
and Tile~2 must have its $a$ side on $PR$, because it cannot fit
next to Tile~1 with its $b$ or $c$ side on $PR$.  From there
the argument proceeds as before. That completes Case~2. 
\smallskip

 That completes the 
proof of the lemma. 
\medskip

If there are enough tiles on the boundary of $ABC$ then $N$ must 
be at least 12.  How many is ``enough''?  As it turns out we do 
not need a precise answer; the following lemma is helpful enough 
and easy to prove.  No doubt the number 10 can be improved, but this 
is good enough.

\begin{lemma} \label{lemma:enough}
Let $ABC$ be $N$-tiled, and suppose the total number of tiles 
with an edge on the boundary of $ABC$ is at least k, with at 
least two tiles on each side of $ABC$, and only one tile at $B$,
and a total of five tiles at the vertices of $ABC$. 
Suppose $\gamma \neq \pi/2$.  Then $N \ge k+2$.
\end{lemma}

\noindent{\em Proof}.  We must produce at least two non-boundary tiles.
 Case~1, two vertices, say $A$ and $B$, of 
$ABC$ have only one tile each.  Since $\gamma \neq \pi/2$,
at three tiles (at least) meet at
each boundary vertex.  Therefore, 
 the tile that shares an edge with the tile at $A$
is not a boundary tile, and the same for the tile next to the tile at $B$.
That makes at least $k+2$ tiles.

Case~2, only $B$ has a single tile, while vertices $A$ and $C$ have 
two tiles each.  Then the tile adjacent to the tile at $B$ is a
non-boundary tile.  Consider the two tiles at vertex $A$, say Tile~1 and 
Tile~2.  If they do not share a common edge then one of them, say 
Tile~1,  has a shorter edge along their common boundary. Then the 
tile adjacent to that edge is not a boundary tile, and hence it is 
a second non-boundary tile.   If they do share a common edge, then let Tile~3
and Tile~4 be the tiles adjacent to Tile~1 and Tile~2, respectively.
At most one of Tile~3 and Tile~4 can have a boundary extending 
past the common interior vertex $E$ of Tile~1 and Tile~2, and the 
one that does not cannot be a boundary tile.  Hence it is a second
non-boundary tile.  That completes the proof of the lemma.

\section{The case $3\alpha + 2\beta = \pi$}
Three of the rows of Table~\ref{table:laczkovich2} fall under 
the case $3\alpha + 2\beta = \pi$, with $\alpha$ not a rational 
multiple of $\pi$.  For some of those cases we have proved
necessary and sufficient conditions for the existence of an 
$N$-tiling; and for all of them we have strong necessary conditions.
In other words,  we have added a fourth column to Table~\ref{table:laczkovich2},
at least for the rows corresponding to $3\alpha + 2\beta = \pi$. 
From those entries we can simply read off that $N=7$ and $N=11$ are 
impossible.  In fact $N=28$ is the smallest possible $N$.  
But the proofs, which are unpublished, occupy approximately
a hundred pages.
(except for its dependence on \cite{laczkovich1995} and \cite{snover1991}).
Therefore we give a short, self-contained, algebraic and computational proof
that $N$-tilings do not exist when $N < 12$ and $3\alpha+2\beta = \pi$ and
$\alpha$ is not a rational multiple of $\pi$. 

An important tool in the analysis of these tilings is the ``coloring 
equation'' given in Theorem~\ref{theorem:coloring}.  That theorem 
applies here, as we now show.
  If $3\alpha+2\beta = \pi$ and $\alpha$ is not a rational 
multiple of $\pi$, then every boundary vertex is composed of 
three tiles ($\alpha + \beta + \gamma$) or five tiles ($3\alpha + 2\beta$),
and every interior vertex is either a ``center'' with four tiles ($3\gamma +\beta$)
or has six tiles ($2\alpha + 2\beta + 2 \gamma$) or 
eight tiles ($4\alpha + 3 \beta + \gamma$) or ten tiles ($6 \alpha + 4 \beta$).

 Since
there are five tiles at the angles of $ABC$, by renaming the vertices we 
may assume that only one tile is at $B$.  Let $(X,Y,Z)$ be the lengths
of sides $AB$, $BC$, and $AC$.  Then we have the 
 ``coloring equation''
\begin{eqnarray}
M(a+b+c) &=& X + Z \pm Y \label{eq:coloring}
\end{eqnarray} 
where the $+$ sign is taken if the angles at $A$ and $C$ have an odd
number of tiles, and the $-$ sign is taken if they have an even number.

Besides the coloring equation, we have the ``area equation'', which 
says that the area of $ABC$ is equal to $N$ times the area of the tile.
We use the formula for the area of a triangle that says twice the area
is the product of two adjacent sides and the sine of the included angle.
By the law of sines, $a/c = \sin \alpha/ \sin \gamma$.
Then the area equation can be written
\begin{eqnarray}
XZ\sin \alpha &=& Nbc \sin \alpha \\
XZ &=& Nbc  \mbox{\qquad if angle $B = \alpha$} \label{eq:area1} \\
XZ &=& Nac   \mbox{\qquad if angle $B = \beta$} \label{eq:area2}
\end{eqnarray}

\begin{definition} \label{defn:s}
Let a triangle have angles $(\alpha,\beta,\gamma)$.  
We define
$$s=2\sin(\alpha/2).$$
\end{definition}
This definition is useful because 
the ratios $a/c$ and $b/c$ can be expressed simply in terms of $s$, as shown in the following
lemma.

\begin{lemma} \label{lemma:singamma}
  Suppose $3\alpha + 2\beta = \pi$.  Let $s = 2\sin \alpha/2$.  Then we have
\begin{eqnarray*}
\sin \gamma &=& \cos \frac \alpha 2 \\
\frac a c &=& s \\
\frac b c &=& 1-s^2  
\end{eqnarray*}
\end{lemma} 

\noindent{\em Proof}.   
Since $\gamma = \pi-(\alpha + \beta)$, we have 
\begin{eqnarray*}
\sin \gamma &=& \sin(\pi-(\alpha + \beta))  \\
 &=& \sin(\alpha + \beta) \nonumber \\
 &=& \cos(\pi/2- (\alpha + \beta))   \\
 &=& \cos \frac \alpha 2  \mbox{\qquad since $\pi/2 - \beta = 3\alpha/2$} 
 \end{eqnarray*}
Then $c = \sin \gamma = \cos \alpha/2$, and $a = \sin \alpha = 2 \sin(\alpha/2) \cos(\alpha/2)$.
Hence 
$$\frac a c = 2\sin \alpha/2.$$
Since $3 \alpha + 2\beta = \pi$, we have 
\begin{eqnarray*}
\sin \beta &=& \sin (\pi/2 - 3 \alpha/2) \\
&=& \cos(3 \alpha/2) \\
&=& 4 \cos^3 \frac \alpha 2 - 3 \cos \frac \alpha 2 
\end{eqnarray*}
Hence 
\begin{eqnarray*}
b/c &=& 4 \cos^2(\alpha/2) - 3 \\
&=& 4(1-\sin^2 \alpha/2) - 3 \\
&=& 1-4\sin^2 \alpha/2
\end{eqnarray*}
 Then we have 
\begin{eqnarray*}
\frac a c &=& s \\
\frac b c &=& 1-s^2
\end{eqnarray*}
establishing the second equation of the lemma.   
That completes the proof of the lemma.

\begin{theorem} \label{theorem:no7-11a} 
Suppose $3\alpha + 2\beta = \pi$, and triangle $ABC$ is 
$N$-tiled by a tile with angles $(\alpha,\beta,\gamma)$ not 
similar to $ABC$, and $\alpha$ is not a rational multiple of $\pi$.
  Then $N \ge 12$.
\end{theorem}

\noindent{\em Proof}.  We first discuss 
the possibility of applying 
of Lemma~\ref{lemma:zerolimits2}.  Do the hypotheses hold? 
By Lemma~\ref{lemma:possibleangles},
no tile has a $\gamma$ angle at a vertex of $ABC$; and by
 Lemma~\ref{lemma:gammagreaterpiover2}, $\gamma > \pi/2$. 
Since the tile is not similar to $ABC$, and $\alpha$ is not a 
rational multiple of $\pi$, each angle of $ABC$ is 
less than $\gamma$. Therefore 
Lemma~\ref{lemma:zerolimits2} is applicable.

We now explain the idea of the proof.
The tiling provides an expression for each side of $ABC$ as a 
linear combination of $abc$.  Thus 
\begin{eqnarray*}
X &=& pa + qb + rc \\
Z &=& ua + vb + wc \\
Y &=& ka + \ell b + mc
\end{eqnarray*}
Substitute these expressions for $(X,Y,Z)$ in the coloring equation.
With $P = p+u\pm k$, $Q = q+v \pm \ell$, $R = r+w\pm m$
we have $M(a+b+c) = Pa + Qb + Rc$.
 Dividing by $c$ and use $a/c = s$ and 
$b/c = 1-s^2$ we have
\begin{eqnarray*}
M(2+s-s^2) &=& Ps + Q(1-s^2) + R
\end{eqnarray*}
For given $(M,P,Q,R)$ that quadratic can be solved for $s$
(provided its discriminant is nonnegative).  The area equation 
too can be expressed in terms of $s$, and we can check if it 
is satisfied for the $s$ from the coloring equation. 
For a given $N$, we need to consider only values of the integer 
parameters between 0 and $N$, so this search will terminate.
Moreover, as discussed in the first paragraph of this proof, Lemma~\ref{lemma:zerolimits2} tells us that we 
can restrict the search by only examining values of $r$, $w$, and $m$
that are at least 2, provided $ABC$ is isosceles or $s$ is rational. 
Finally, Lemma~\ref{lemma:enough}  allows us to 
not consider cases in which there would be ten or more boundary tiles, 
i.e.,  when $p+q+r+u+v+w+k+\ell+m \ge 10$.   
SageMath code to carry out this 
plan for isosceles $ABC$ with base angles $\alpha$ or $\beta$
is exhibited in Fig.~\ref{figure:oct22}.  Run that code
passing 7 as the function parameter, and then again passing 11.
It runs in about 12 seconds, and produces no output except the reassuring 
progress reports as $M$ changes.   That shows that there is no 
7 or 11 tiling in the case of isosceles $ABC$ with angles $\alpha$
at $A$ and $B$, or $\beta$ at $A$ and $B$, i.e., when the  
coloring equation is 
$M(a+b+c) = X + Y + Z$.   

The other possible shapes of $ABC$ satisfy the 
coloring equation $M(a+b+c) = X-Y+Z$.  That code differs from 
the code in Fig.~\ref{figure:oct22} in two respects.
First, because of the minus sign in the coloring equation,  
negative values of
$(P,Q,R)$ are allowed, and 
the upper limits of $(P,Q,R)$ go up to $N$,
$N-\vert P \vert$, and $N-\vert P \vert - \vert R \vert$, respectively, 
and the values of $(k,\ell,m)$ are preceded by a minus sign, 
with a continue statement inserted to reject negative values.
Second, we are only allowed to assume each side contains at least
two $c$ edges in case $s$ is rational,  so the variable {\tt looplimit}
has to be recalculated each time $s$ is recalculated, and set to 
2 if $s$ is rational, and otherwise to 1.  Although it adds 
a page to the length of the paper, we also include enough of this code so 
that any reader can reproduce our results.  See Fig.~\ref{figure:oct22b}.

To prove the theorem
as stated, we ran both programs for all $N$ between 3 and 11, inclusive.
The second program is slower, requiring 27 seconds for $N=7$,
over three minutes for $N=11$, and about 8 minutes for all values 3 to 11.
But it gets the answer:  no solutions are found.  
That completes the proof. 
\medskip

\noindent{\em Remark}.  This method cannot be used for $N=14$ or $19$, as 
some solutions are found.  As discussed in \S\ref{section:last} below,
we do have a proof that there is no 14-tiling or 
19-tiling with $3\alpha+2\beta = \pi$;
but the short direct computational proof given here will not work.    

\begin{figure}[ht]
\caption{SageMath code used in the proof of Theorem~\ref{theorem:no7-11a}, $ABC$ isosceles}
\label{figure:oct22}
\begin{verbatim}
def oct22(N):  # ABC isosceles 
  var('P,Q,R,M,s,p,q,r,u,v,w,k,ell,m')
  epsilon = 0.0000001
  lowerlimit = 2  # each side has at least 2 c edges
  for M in range(1,N):
    print("M=%d" %M)
    for P in range(0,N):
      for Q in range(0,N-P):
        for R in range(6,N-P-Q):
          eq1 = M*(2+s-s^2) - P*s - Q*(1-s^2) - R
          discriminant = (M-P)^2 - 4*(Q-M)*(2*M-Q-R)
          if discriminant < 0:
            continue
          answers = solve(eq1,s)
          for x in answers:
            if x.rhs() <= 0 or x.rhs() >= 1:
              continue             
            for r in range(lowerlimit,R+1):
              for w in range(lowerlimit,R-r):
                m = R -r-w
                if m < lowerlimit:
                  continue
                for p in range(0,P+1):
                  for u in range(0,P-p):
                    k = P-p-u;
                    for q in range(0,Q+1):
                      for v in range(0,Q-q):
                        ell = Q-q-v
                        boundarytiles = p+q+r+u+v+w+k+ell+m
                        if boundarytiles >= N-2:
                          continue
                        X = r + p*S + q*(1-S^2)
                        Y = w + u*S + v*(1-S^2)
                        area1 = abs( X*Y - N*(1-S^2))
                        area2 = abs(X*Y - N*S)      
                        if n(area1) < epsilon:   #  B = alpha
                          print("alpha",N,M,p,q,r,u,v,w,k,ell,m)
                          print(area1)
                        if n(area2) < epsilon:  #  B = beta
                          print("beta",N,M,p,q,r,u,v,w,k,ell,m)              
\end{verbatim}
\end{figure} 

\begin{figure}[ht]
\caption{SageMath code used in the proof of Theorem~\ref{theorem:no7-11a}, $ABC$ not isosceles}
\label{figure:oct22b}
\begin{verbatim}
def oct22b(N):  # case when ABC is not isosceles
  var('P,Q,R,M,s,p,q,r,u,v,w,k,ell,m')
  epsilon = 0.0000001
  for M in range(1,N):
    print("M=%d" %M)
    for P in range(-N,N+1):
      for Q in range(-(N-abs(P)),N-abs(P)+1):
        for R in range(-(N-abs(P)-abs(Q)),N-abs(P)-abs(Q)+1):
          eq1 = M*(2+s-s^2) - P*s - Q*(1-s^2) - R
          discriminant = (M-P)^2 - 4*(Q-M)*(2*M-Q-R)
          if discriminant < 0:
            continue
          answers = solve(eq1,s)
          for x in answers:
            S = x.rhs()
            if S <= 0 or S >= 1:
              continue
            lowerlimit=1; #will be set to 2 when S is rational
            if S < 1-S^2 and (S in QQ or not (1-S^2)/S in ZZ):
              lowerlimit = 2
            else:
              if 1-S^2 < S and (not (S/(1-S^2) in ZZ)):
                lowerlimit = 2
              else:
                lowerlimit = 1 
            for r in range(lowerlimit,R+1):
              for w in range(lowerlimit,R-r):
                m = -(R -r-w)
                if m < lower3limit:
                  continue
                for p in range(0,P+1):
                  for u in range(0,P-p):
                    k = -(P-p-u);
                    if k < 0:
                      continue
                    for q in range(0,Q+1):
                      for v in range(0,Q-q):                                   
                        ell = -(Q-q-v)
                        if ell < 0:
                          continue
                        # ... the rest as in the previous figure
                                    
\end{verbatim}
\end{figure}
\FloatBarrier

\section{The case $\gamma = 2\pi/3$ and $\alpha$ not a rational multiple of $\pi$}
In this case, $\alpha + \beta = \pi/3$, so a boundary vertex can be composed
of angles contributed by 3 or 6 tiles.  Hence it is not in general possible 
to color the tiles black and white in a way that leads to
a ``coloring equation.'' 

There are several shapes possible for $ABC$, listed in
Table~\ref{table:laczkovich2}, but for our purposes  there are 
just two cases to consider:  either one of the vertices of $ABC$ has  
just one tile (in which case we rename $\alpha$ and $\beta$ so 
that the standalone angle is $\alpha$, and we rename the vertices
so it occurs at $A$),  or there are two tiles at each of the three vertices,
in which case we may assume that the angle at $A$ is $\alpha + \beta$.
We do not need to consider the case $ABC$ similar to the tile, so 
no $\gamma$ angles occur at the vertices of $ABC$.  The shape of 
the tile can be expressed using the law of cosines, since 
$\cos(2\pi/3) = - \frac 1 2$, by the equation 
\begin{eqnarray}
 c^2 &=& a^2 + b^2 + ab. \label{eq:3084}
\end{eqnarray}
For example, $(3,5,7)$ and $(8,7,13)$ are rational tiles satisfying
this equation.  

Although there is no coloring equation, we
 do still have the ``area equation'' that 
says the area of $ABC$ is $N$ times the area of the tile.  That equation
takes different forms depending on the shape of $ABC$.  In case the 
angle at $A$ is $\alpha$, and the sides $AB$ and $AC$ have length
$X$ and $Y$, the area equation is $XY\sin \alpha = Nab \sin \alpha$.
After canceling $\sin \alpha$ we have
\begin{eqnarray} \label{eq:2872}
XY &=& Nab
\end{eqnarray}
We do have some general results about this kind of tiling, but 
the theory is incomplete and we do not go into it here.
For our present purposes it suffices to show that any such 
tiling requires at least 12 tiles; that is Theorem~\ref{theorem:no7-11b} 
below.   

\begin{theorem} \label{theorem:no7-11b}
Let triangle $ABC$ be $N$-tiled by a tile  with angles $(\alpha,\beta,2\pi/3)$,
not similar to $ABC$, and suppose $\alpha$ is not a rational multiple of $\pi$.
Then $N \ge 12$.  In particular, $N$ is not equal to 7 or 11.
\end{theorem}

\begin{figure}[ht]
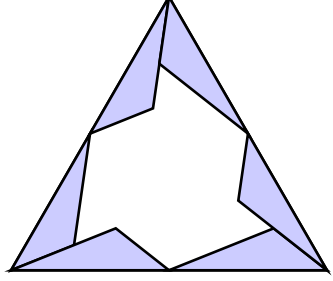

\caption{The case of Theorem~\ref{theorem:no7-11b} when there 
are exactly six boundary tiles} 
\label{figure:impossible6}
\begin{center}
\FigureImpossibleSix
\end{center}
\end{figure}

\noindent{\em Remark}. 
The idea of the proof of this case is that, because
of Lemma~\ref{lemma:zerolimits2}, each side of $ABC$ is at least $2c$
in length, and that makes the area more than the area of 12 tiles.
See Fig.~\ref{figure:impossible6}, which illustrates an equilateral $ABC$
with six tiles placed, and more area remaining than
six tiles can cover.  (This figure is only illustrative.) 
  We first proved this theorem
by a geometrical argument about placing tiles, but algebra is shorter
and simpler. Both ideas can be seen in Fig.~\ref{figure:impossible6}:
It is geometrically impossible to complete the tiling, and also 
the untiled area is more than the area of six tiles.

\FloatBarrier
\medskip
\noindent{\em Proof}.   Since the tile is 
not similar to $ABC$, and $\alpha$ is not a rational multiple of $\pi$,
there can be no $\gamma$ angle at a vertex of $ABC$.  Then there must
occur a total of six tiles at the vertices of $ABC$, 
contributing three $\alpha$
angles and three $\beta$ angles to make up the angles of $ABC$. 
Lemma~\ref{lemma:zerolimits2} is applicable, since no vertex of 
$ABC$ can have a $\gamma = 2\pi/3$ angle, so each side of $ABC$
has at least two $c$ edges.  That is the key idea of this proof. 

We divide the proof into two cases.  Case~1:  One vertex of $ABC$
has an angle $\delta$ with  $\pi/3 \le \delta \le 2\pi/3$.  The 
point of that inequality is that it implies $\sin \pi/3 < \sin \delta$.
Let $X$ and $Y$ be the lengths of the 
sides adjacent to that angle.  Then we have the area equation
\begin{eqnarray*}
XY \sin \delta &=& Nab \sin \frac {2\pi} 3
\end{eqnarray*}
Since $\sin \frac {2\pi} 3 = \sin \frac{\pi} 3 \le \sin \delta$, 
we have
\begin{eqnarray*}
XY &\le & Nab 
\end{eqnarray*}
According to Lemma~\ref{lemma:zerolimits2}, there are at least
two $c$ edges on each side of $ABC$.  Hence $X \ge 2c$ and $Y \ge 2c$.
Therefore $Nab \ge XY \ge 4c^2$.  Therefore
\begin{eqnarray}
3 ab &\ge& \frac {12c^2} {N} \label{eq:3177}
\end{eqnarray}
 Recall (\ref{eq:3084}):
\begin{eqnarray}
c^2 &=& a^2 + b^2 + ab \nonumber \\
c^2 - 3ab &=& a^2 + b^2 - 2ab \ = \ (a-b)^2 \ > \ 0 \label{eq:3188}
\end{eqnarray}
Substituting on the left from (\ref{eq:3177}) we have
\begin{eqnarray*}
c^2\left(1 - \frac {12} N\right) &>& 0
\end{eqnarray*}
We have strict inequality since the hypothesis that
$\alpha$ is not a rational multiple of $\pi$ implies $a \neq b$.
Since $c^2 > 0$  we have
\begin{eqnarray*}
1 - \frac {12}{ N} & > & 0 \\
N &>& 12
\end{eqnarray*}
That completes the proof in Case~1.
\smallskip

Case~2: Every vertex angle of $ABC$ is either more than $2\pi/3$ 
or less than $\pi/3$.  They cannot all be less than $\pi/3$ since 
they add up to $\pi$.  Therefore one angle is more than $2\pi/3$.
Renaming the vertices if necessary, we can assume the angle at $B$
is more than $2\pi/3$.  Renaming $\alpha$ and $\beta$ if necessary,
we can assume $\alpha < \beta$.  Then the angles at $A$ and $C$
are 
 either $(\alpha,\alpha)$ or $(\alpha, 2\alpha)$ or $(2\alpha,\alpha)$,
  since  
otherwise the angle at $B$ is  $\le\pi-(\alpha+\beta)= 2\pi/3$.

I say that no tile has one vertex on $AB$ and another vertex on 
$BC$.  Suppose, to the contrary, that Tile~1 has vertex $P$ on $AB$
and vertex $Q$ on $BC$.  Then Tile~1 does not have a vertex at $B$,
since the tiles at $B$ have at least one vertex interior to $ABC$.
Consider triangle $PBQ$.  Angle $B$ is either $3\beta$ or 
$3\beta+\alpha$, either of 
which is more than $2\pi/3$, so $\vert \cos B \vert > 1/2$.
Hence angles $P$ and $Q$ of triangle $PBQ$ are acute. 
Consider the rays emanating from $B$ along the dividing lines between
tiles with vertices at $B$.  All three or four of the tiles at $B$
have one $c$ edge emanating along such a ray, and the whole tile is
contained in $PBQ$.  If one of the tiles at $B$ supported by $PB$
or $BQ$ has its $c$ edge there, then $PB \ge c$ or $BQ \ge c$.  
Otherwise let $BR$ be the ray containing the other edge of the tile 
supported by $PB$ at $B$, with 
$R$ on $PQ$.  In triangle $BRP$, $BR$ is opposite an acute angle $P$,
and $PB$ is opposite an angle greater than $\pi/2$, since that angle 
is $\pi$ minus angle $P$ minus angle $PBR$, and both the subtracted
angles are acute.  Since the greater side is opposite the greater angle,
$BP > BR \ge c$.   Thus either $BP \ge c$ or $BQ \ge c$.  Relabeling
$P$ and $Q$ if necessary, we can assume $BP \ge c$.   Since $BQ$
is also composed of tile edges, at least $BQ \ge a$.     
Then by the law of cosines
we have 
\begin{eqnarray*}
PQ^2 &=& BP^2 + BQ^2 - 2PB\cdot QB \cos B \\
&=& BP^2 + BQ^2 + 2PB\cdot QB \vert \cos B \vert \\
&\ge& c^2 + a^2 + 2ac \vert \cos B \vert \\
&=&  c^2 + a^2 \\
&>& c^2
\end{eqnarray*}
Hence the length of $PQ$ is greater than $c$, and hence
cannot be just one tile edge.  Hence, as I said, no tile
has one vertex on $AB$ and another vertex on $BC$. 

I say that $c$ is not a linear integral combination of $(a,b)$
unless it is a multiple of $a$.  For suppose  $c = ua + vb$.  Then
\begin{eqnarray*}
c^2 &=& a^2 + b^2 + ab \\
(ua + vb)^2 &=& a^2 + b^2 + ab \\
(u^2-1)a^2 + (v^2-1)b^2 + (2uv-1)ab &=& 0
\end{eqnarray*}
which is a contradiction if both $u$ and $v$ 
are positive.  Therefore $u=0$ or $v=0$.  If $v=0$,
then $c$ is a multiple of $a$ as claimed. Therefore
we may assume $c = vb$ with $v > 1$.
Then by the law of cosines, 
\begin{eqnarray*}
a^2 &=& b^2 + c^2 - 2bc \cos \alpha \\
    &>& b^2 + c^2 - 2bc \\
    &=& b^2(v^2 + 1 -2v) \mbox{\qquad since $c = vb$} \\
    &=& b^2(v-1)^2 \\
    &>& b^2  \mbox{\qquad since $v>1$}
\end{eqnarray*}
Hence $a > b$, contradiction.  Hence, as I said,
$c$ is not a linear integral combination of $(a,b)$ unless
it is a multiple of $a$.

We now divide into further cases.  Case~2A:  each of the two 
edges $AB$ and $BC$ supports at least three tiles.  Then, there
are two ``notches'' between the three tiles, so there are five
tiles touching $AB$ in at least one point, and five tiles 
touching $BC$.  None of these have been double-counted, since 
no tile has a vertex on $AB$ and a vertex on $BC$, so that is 
ten tiles.  If $ABC$ is isosceles then there are four tiles with
vertices at $B$, two of which we have not yet counted, so that 
makes 12.  If $ABC$ is not isosceles then there are three tiles 
with vertices at $B$, and one more with a vertex at $A$ or $C$,
making again two uncounted tiles for 12 total.  That completes
the proof in Case~2A. 
\smallskip

\begin{figure} [ht]
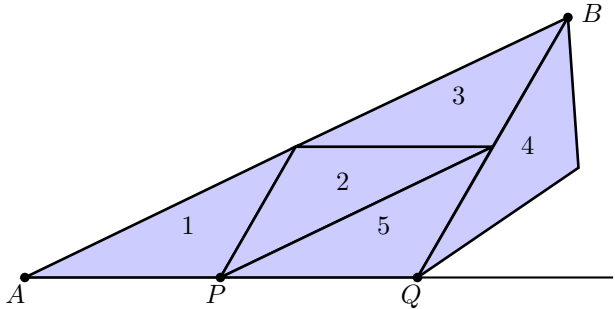

\caption{Case 2B of Theorem 7}
\label{figure:theorem7}
\begin{center}
\FigureTheoremSeven
\end{center}
\end{figure}

Case~2B: $AB$ supports exactly two tiles, and $ABC$ has 
angle $\alpha$ at $A$.   
As already noted, by Lemma~\ref{lemma:zerolimits2}, there are at least
two $c$ edges on each side of $ABC$, so both the tiles supported
by $AB$ have their $c$ edges on $AB$.  See Fig.~\ref{figure:theorem7}
for an illustration of the following argument. 
Let Tile~1 be the tile 
at $A$; then Tile~1 has its $b$ edge on $AC$.  Let Tile~2 be the 
tile east of Tile~1.  Tile~2 shares the $a$ edge of Tile~1, since
that edge terminates at both ends on the boundary of $ABC$.  Tile~2
cannot have its $\gamma$ angle to the south, since that would make
two $\gamma$ angles at a vertex on the boundary.  Therefore it has
its $\gamma$ angle on $AB$ and its $\beta$ angle on $AC$.  Let 
Tile~3 be the tile north of Tile~2.  Then Tile~3 is supported
by $AB$ and hence has its $c$ edge on $AB$, and shares the $b$ edge
of Tile~2 on its southern border.  Let Tile~4 be the tile east of 
Tile~3.  Then Tile~4 has a vertex at $B$.
Either Tile~4  has its $b$ or $c$ edge along
Tile~3 (either of which extends beyond Tile~3, since $a < b < c$),
or it has its $a$ edge shared with Tile~3 and its $\gamma$ angle 
to the southwest.  In all three of those cases, it terminates the 
line forming the southeast border of Tile~2.  Let Tile~5 be the tile 
south of Tile~2, with its southwest vertex $P$ on $AC$ shared with Tile~1.

Suppose, for proof by contradiction, that
Tile~5 does not share the $c$ edge of Tile~2. Since
 its $\alpha$ angle is towards the west,
 it has its $b$ edge along Tile~2. Since the tiles south 
 of Tile~2 must terminate at the eastern vertex of Tile~2, 
 the remaining $c-b$ of the southeast border of Tile~2
must be filled by some tile edges.  But then 
 $c$ would be an integral
 linear combination of $a$ and $b$, including at least one $b$,
 which is impossible, as proved above.
 That contradiction completes
the proof that Tile~5 shares the $c$ edge of Tile~2.  Hence
Tile~5 has its $b$ edge on $AC$, as shown in Fig.~\ref{figure:theorem7}.
Then Tiles~1,2,3,5 are definitely as shown in  Fig.~\ref{figure:theorem7}.

Now consider $BQ$, the eastern border of Tiles~3 and 5.  On the 
west side of $BQ$ are two $a$ edges.  These cannot be matched on the east by 
two $a$ edges, since then the $\gamma$ angles of those two tiles 
would occur on $BQ$, either both to the north, or both to the south,
but both are impossible.  If Tile~4 has its $b$ edge on $BQ$ then 
$2a -b$ is an integer linear combination of $(a,b,c)$, which is 
impossible. Hence Tile~4 has its $c$ edge on $BQ$. Hence $c = 2a$ 
and Tile~4 has a vertex on $AC$ as shown in Fig.~\ref{figure:theorem7}.
Let Tile~6 be south of Tile~4.  Then Tile~6 has its $\beta$ angle 
at $Q$ and hence does not have its $b$ edge along Tile~4.  Since
$2a = c > b$, the tile or tiles south of Tile~4 do not terminate at the 
eastern boundary of Tile~4, but continue to the east.  Let Tile~7
be the tile east of Tile~4.  Then Tile~7 must share its $a$ edge 
with Tile~4, since it cannot extend to the south.  But that is 
impossible, as then it would have a $\gamma$ angle on the shared 
$a$ edge, but that cannot occur at $B$ on the north or at the 
southeast vertex of Tile~4 either.  We have reached a 
contradiction.  That completes the proof in Case~2B.

\smallskip

Case~2C: The angle of $ABC$ at $A$ is $2\alpha$ and $AB$ supports
exactly two tiles.  Then again the two tiles on $AB$ have their $c$
edges on $AB$.  As before let Tile~1 have a vertex at $A$. Let Tile~8
be the other tile with a vertex at $A$, south of Tile~1. Then Tile~8
either extends east of Tile~1, or has its $\gamma$ angle at the 
shared eastern vertex $P$.  Tile~2 cannot be placed with its
$\gamma$ angle also at $P$, making $P$ a ``center'' with three 
$\gamma$ angles, since in that case the $b$ side of Tile~2 would
extend past Tile~8, which is impossible as that would go outside $ABC$.
Therefore, whatever the position of Tile~8, Tile~2 must be placed as before,
with its $a$ edge shared with Tile~1.  Then Tile~3
must be placed as before also, and as before Tile~4 must terminate the 
southern boundary of Tile~2 from extending eastwards.
 Let Tile~5 be the tile south of Tile~2.  Since the southern boundary
 of Tile~2 is terminated at both ends.  Assume,
 for proof by contradiction, that Tile~5 does not share its $c$
 edge with Tile~2.  Then $c$ is a combination of tile edges $a$ and $b$,
 which implies that $c$ is a multiple of $a$. Then Tile~5 
has its $a$ edge against Tile~2.  Then its $\beta$ angle 
or $\gamma$ angle is at $P$, which implies that Tile~8 shares its
$b$ edge with Tile~1 and has its $\gamma$ angle at $P$.  But then,
there is not room for Tile~5 to also have its $\gamma$ angle at $P$.
Hence Tile~5 has its $\beta$ angle at $P$.
Then there exists Tile~9 between Tile~8 and Tile~5, with the $\alpha$
angle of Tile~9 at $P$.  Tile~9  lies next to the $a$ edge of Tile~8,
but since Tile~9 has its $\alpha$ angle at $P$, it must have its
$b$ or $c$ side next to Tile~9, which is impossible as that side would 
extend outside $ABC$.  That contradiction completes the proof 
that Tile~5 does share its $c$ edge with Tile~2. 
Hence Tile~5 must occur in the position shown 
in Fig.~\ref{figure:theorem7}.   That is, Fig.~\ref{figure:theorem7} 
correctly shows Tiles 1,2,3,5 also in Case~2C,  regardless of the position of Tile~8 (which
is not shown in the figure).

Now, in Case~2B, the southern vertex of Tile~4 had to be the southeastern
vertex of Tile~5.  If that is so, then as before $c = 2a$, and the 
proof is completed as in Case~2B. 
Let $P$ and $Q$ be 
the southeastern vertices of Tiles~1 and 5, respectively.  If 
some tile south of $APQ$ blocks line $BQ$ from continuing south of $PQ$,
then the southern vertex of Tile~4 is $Q$, and as in Case~2B,
$c = 2a$ and the proof can be completed. Hence, we can assume
$BQ$ does extend through $APQ$.  Then the tiles south of Tiles~1
and 5 must share the $b$ edges of those tiles.  But that is impossible,
as those two tiles would have their $\gamma$ angles to the east,
crossing $BQ$.   That contradiction completes the proof 
in Case~2C.

If Case~2 holds then either $AB$ or $BC$ supports exactly two tiles. 
Renaming $A$ and $C$ if necessary, we can assume it is $AB$ that 
supports exactly two tiles.  Then either Case~2B or Case~2C applies.
That completes the proof of the theorem.
\medskip

\begin{figure}[ht]
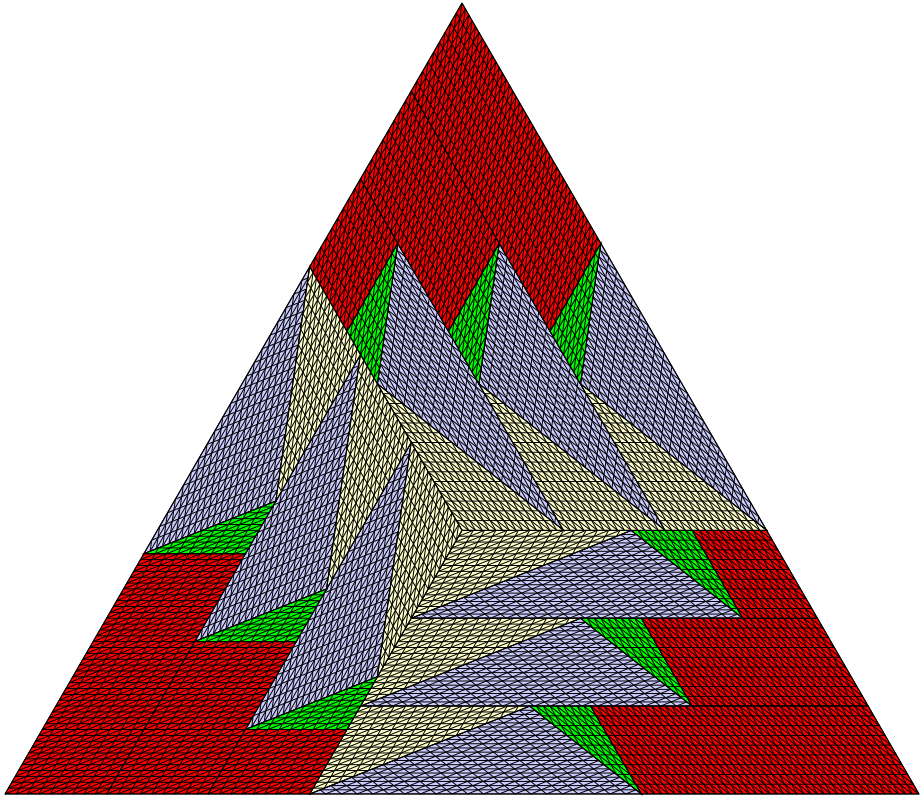

\caption{$N=10935$.  The tile is $(3,5,7)$. $ABC$ is equilateral.}
\label{figure:10935}
\begin{center}
\EquilateralFigure
\end{center}
\end{figure}

Laczkovich proved that $N$-tilings of the kind discussed in this 
section exist, but did not actually exhibit any.  Although in 
\cite{laczkovich1995}, he did not explicitly consider $N$ at all, 
he did consider $N$ 
in a theorem that he allowed Soifer to publish in \cite{soifer}.
In that theorem, he proved
that for tilings of an equilateral triangle by a tile $(\alpha,\beta,2\pi/3)$,
the square-free part of $N$ could be anything desired.  Following
Laczkovich's ideas, we found the tiling of Fig.~\ref{figure:10935}, 
with $N = 10935$.  Other tilings that we found require more than $32,000$
tiles (and so are too big to draw nicely on a normal page).  What the 
smallest possible $N$ is, we have no idea.  For all we know,
the construction method used for this tiling might yield a smaller 
$N$ for some tile with very large sides; or there might be 
a much more efficient tiling construction yet to be discovered.
In 2012 it was not known if there is an $N$-tiling of the equilateral 
triangle for every sufficiently large $N$,  or if instead there
are arbitrarily large $N$ for which, like 7 and 11, there is no 
$N$-tiling at all.  In unpublished work, we have proved $N$ cannot 
be prime. Therefore there are arbitrarily large $N$ for which there is 
no $N$-tiling.

\section{Tilings of an isosceles triangle with $\gamma = 2\alpha$}
In this section we take up the row of Laczkovich's second table in 
which $ABC$ is isosceles with base angles $\alpha$ and is tiled
by a tile with $\gamma = 2\alpha$, and $\alpha$ is not a rational
multiple of $\pi$.    The condition $\gamma = 2\alpha$ can also be 
written as $3\alpha + \beta = \pi$.  Unlike the similar-looking
condition $3\alpha + 2\beta = \pi$, this condition does not imply 
$\gamma > \pi/2$.  The vertex angle of $ABC$ is then $\pi-2\alpha = \alpha +\beta$.

\begin{figure}[ht]
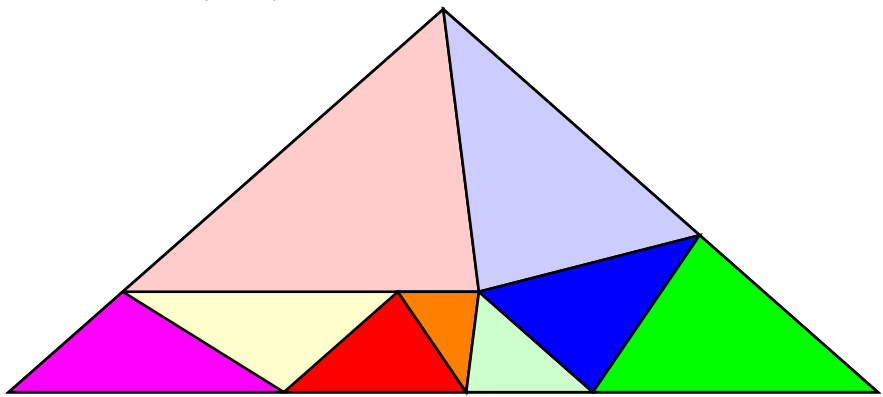

\caption{Laczkovich's dissection of isosceles $ABC$ into triangles similar to $(4,5,6)$ can be used to produce a $5861172$-tiling.}
\label{figure:bigisosceles}
\begin{center}
\FigureBigIsosceles
\end{center}
\end{figure}

Laczkovich \cite{laczkovich1995} proves that, given any
tile with $\gamma = 2\alpha$ and $\alpha$ not a rational multiple 
of $\pi$, an isosceles triangle can 
be dissected into triangles {\em similar} to the tile. 
Following the steps of his proof with the tile $(4,5,6)$, one finds the dissection shown in
 Fig.~\ref{figure:bigisosceles}.
  To make an $N$-tiling, we have 
to tile each of these triangles and the parallelogram with many copies
of the same tile.  Along each edge in the figure there is an arithmetical
condition to satisfy.  Working out those conditions, we find that 
more than five million tiles are required:  5861172 to be precise.
  It is not possible to 
print such a large tiling (unless one could use the side of a large building), and we do not know a smaller one.  But at least,
some such tilings do exist.  Indeed, {\em many} such tilings exist.

The theory of such tilings has progressed far enough in 
\cite{beeson-isosceles}
to prove that that tile is necessary rational,
and from that plus the characterization of the tile given above,
one can show that $N$ cannot be prime.  That is proved
in \cite{beeson-isosceles}; in fact, it is proved that $N$ 
cannot be twice a prime either.  In particular, $N$ cannot be 
7, 11, 14, or 19; so we need not consider the case $\gamma = 2\alpha$
further in this paper.   We merely remark that there is still a big gap
between the smallest $N$ for which a tiling of this kind is known to 
exist (over five million) and the lowest $N$ for which we have not 
ruled out the existence of a tiling, currently $N=20$.  We are 
far from a characterization of the $N$ for which there are $N$-tilings;
but we did prove $N$ cannot be prime.

\section{Tilings of an equilateral triangle with $\alpha/\pi$  irrational}

 Laczkovich's Theorem~3.1
\cite{laczkovich1995} says that, given a 
 rational tile with an angle $\pi/3$, and the other angles
 not rational multiples of $\pi$,  that tile will tile
 {\em some}  an equilateral $ABC$.    There are
infinitely many such rational tiles, as Laczkovich proves. 
The two simplest ones are $(7,5,8)$ and $(7,3,8)$.  Similarly
when the tile has a $2\pi/3$ angle.  Hence there are plenty
of tilings of the kind considered here.

Laczkovich's second table has an entries for the cases when $ABC$
is equilateral and the tile has either a $\pi/3$ or a $2\pi/3$ 
angle.  The second table assumes
not all the angles are rational multiples of $\pi$, so this 
entry also assumes that $\beta$ is not a 
rational multiple of $\pi/3$.  In Theorem~\ref{theorem:no7-11b},
we proved that if the tile has a $2\pi/3$ angle then $N > 12$.
Thus, for the purpose of proving $N$ cannot be $7$ or $11$, we
are already done with that case.   Nevertheless, for very little
extra work, we can push the lower limit on $N$ higher for both
cases of tilings of an equilateral triangle.

The main difference between the two cases is that when the tile has
a $\pi/3$ angle, we can color the tiles black and white in such a 
way that the coloring theorem applies.  That is not possible when 
the tile has a $2\pi/3$ angle, for example because it is possible
(and necessary) for three tiles to meet at some vertex where all 
three have a $2\pi/3$ angle, and three tiles at one vertex cannot 
be colored.

In unpublished work, we have interesting results about the 
existence or non-existence of such tilings, including a proof
that the tile ratios $b/a$ 
and $c/a$ can be computed from $N$ and $M$, and that $N$ cannot 
be prime, which certainly covers the cases $7$ and $11$.  
These proofs will not be presented here; instead we treat
the problem computationally.
Of course that covers only sufficiently small $N$.  But it 
at least deals with $N = 7$ and $11$.

There are two computational approaches to the problem of tiling
an equilateral triangle; one uses the coloring equation and hence
applies only to the case of a tile with a $\pi/3$ angle, but the 
other uses only the area equation, and applies just as well to 
both cases.  The latter method, however, needs to use the fact that
the tile is rational,  which Laczkovich proved in 2012 
\cite{laczkovich2012}, Theorem~3.3.  The approach via the coloring
equation does not use that result, which we may count in its favor,
but the method based on the area equation works better, so 
we present the code for that method.

We will label the angles of the tile
 so that $\gamma$ is the $\pi/3$ angle or the $2\pi/3$ angle.
Then we have 
 the law of cosines:
\begin{eqnarray}
a^2 &=& b^2 + c^2 - 2bc \cos \gamma \nonumber\\
a^2 &=& b^2 + c^2 \pm bc  \label{eq:2605}
\end{eqnarray}
since $\cos \gamma = \pm \pi/3$.  The plus sign 
corresponds to the case $\gamma = 2\pi/3$ and the minus
sign corresponds to $\gamma = \pi/3$.
The area equation is the same in both cases, since
$\sin \gamma = \sin \pi/3$ in either case.  
If $X$ is the length of each side of $ABC$, 
\begin{eqnarray}
X^2 \sin \pi/3 &=& Nab \sin \gamma \nonumber \\
X^2 &=& Nab \label{eq:2511} 
\end{eqnarray}

In case $\gamma = \pi/3$, we also have the coloring
equation
\begin{eqnarray}
M(a+b+c) &=& 3X \label{eq:2500}
\end{eqnarray}
where $M$ is the coloring number of the tiling.

Both computational approaches depend on writing 
$$ X = pa + qb + rc$$
for non-negative integers $(p,q,r)$; this expression
describes how the sides of $ABC$ are composed of tile edges.

The following lemma gives some useful
 restrictions on the possibilities
for $(p,q,r)$.

\begin{lemma} \label{lemma:pqr}
Let the equilateral triangle $ABC$ be tiled by a tile with
$(\alpha,\beta,\gamma)$ and sides $(a,b,c)$, where $\gamma$
is either $(\pi/3)$ or $(2\pi/3)$.  Possibly after a relabeling of 
the vertices of $ABC$, let $X = pa + qb + rc$
where $X$ is the length of $AB$. Then
\smallskip

(i) If $\gamma = \pi/3$, then $p \ge 1$ and $q \ge 1$ 
\smallskip

(ii) $r \ge 2$
\end{lemma}

\noindent{\em Proof}.  Ad (i). If $\gamma = \pi/3$, then
at the vertices of $ABC$ there are 
tiles with altogether three $a$ edges and three $b$ edges.
Therefore at least one side of $ABC$ has one $a$ edge and one 
$b$ edge at its endpoints.  Choosing that side for the decomposition
of $X$, we have $p \ge 1$ and $q \ge 1$.  We relabel the 
vertices so that side is $AB$.   

If $\gamma = 2\pi/3$,
 we do not assert $p \ge 1$ and $q \ge 1$.  But we do have, 
by Lemma~\ref{lemma:zerolimits2}, that there are at least two
$c$ edges on each side of $ABC$, so $r \ge 2$.  

Now suppose $\gamma = \pi/3$. 
We have
$r \ge 1$, since if there are no $c$ edges on $AB$,
then every tile supported by $AB$ has a $\gamma$ angle
on $AB$, which contradicts the pigeonhole principle since the
edges at the endpoints are $a$ or $b$ edges. 
I say that $r \ge 2$ holds also in case $\gamma = \pi/3$.
Indeed if there
is only one $c$ edge on $AB$, and a total of $n$ tiles supported
on $AB$, then there are $n-1$ tiles with a $\gamma$ angle on $AB$,
and $n-1$ possible vertices for them, so the one tile with a 
$c$ edge is bordered by two tiles with their $\gamma$ angles
adjacent.  Let $PQ$ be the $c$ edge of Tile~3 supported by $AB$,
and let Tiles~1,2,3,4,5 be adjacent tiles in order, so that Tiles~1,3,5
are supported by $AB$. Tiles~1 and 5 have their $\gamma$ angles at $P$
and $Q$ respectively.  Since $\gamma = \pi/3$ we have $\alpha + \beta = 2\pi/3$,
so $\alpha < \gamma < \beta$.  Hence $a < c < b$.  Renaming $P$ and $Q$
if necessary, we may assume that Tile~3 has its $\alpha$ angle at $P$
and its $\beta$ angle at $Q$.  Then Tile~2 has its $\beta$ angle at $P$
and Tile~4 has its $\alpha$ angle at $Q$. Let $R$ be the third vertex
of Tile~3.  Then $RQ$ has length $a$ since it is opposite the $\alpha$
angle of Tile~3.  The adjacent edge of Tile~4 is not the $a$ edge, 
since Tile~4 has its $\alpha$ angle at $Q$. Hence Tile~2 has its 
$b$ edge matching the $b$ edge of Tile~3 along $PR$, since the $c$ 
edge is too long to fit, and if the $a$ edge were there, the remaining
part $b-a$ cannot be tiled unless $b$ is an integer multiple of $a$. 
But $b$ cannot be an integer multiple of $a$, for then by the 
law of cosines,
\begin{eqnarray*}
c^2 &=& a^2 + b^2 - ab \\
&=& a^2 + (ma)^2 - a(ma) \\
&=& a^2
\end{eqnarray*}
Hence $c=a$.  Hence by the law of sines $\alpha = \gamma$,
contradicting the assumption that $\alpha$ is not a rational multiple of $\pi$.
Hence Tile~2 has its $b$ edge along $PR$.  But that is impossible, since
it has its $\beta$ angle at $P$.  
That completes the proof of the lemma.

Our algorithm is going to check all possible values of $(p,q,r)$ and 
try to solve the area equation.  The time that takes will clearly 
depend on how large $(p,q,r)$ can be in terms of $N$. 
Of course each of $(p,q,r)$ is at most $N$ since there
are at most $N$ tiles altogether.  But for efficiency of computation,
we want a better bound.  We improve the bound in the next 
two lemmas by a factor of 6; 
still crude, but enough for our purposes.

\begin{lemma} \label{lemma:pqrbound_helper} Let the equilateral triangle
$ABC$ be $N$-tiled by a tile with sides $(a,b,c)$ and $\alpha/\pi$
irrational.  If $\gamma = 2\pi/3$
then no tile touches two different sides of $ABC$.  If $\gamma = \pi/3$
then exactly three tiles touch two different sides of $ABC$.
\end{lemma}

\noindent{\em Proof}. First assume $\gamma= 2\pi/3$.  Since we can
relabel the vertices of $ABC$, it suffices to show that no tile
touches both $AB$ and $BC$.    We have $a < c$
and $b < c$ since $\alpha + \beta = \pi/3 < \gamma$.  Renaming
$\alpha$ and $\beta$ if necessary, we may suppose $\alpha < \beta$.
Suppose 
$PQ$ is a tile edge with $P$ on $AB$ and $Q$ on $AC$.  
  Then angle $APQ$ plus angle $AQP$ is
$2\pi/3$, since the sum of the angles in triangle $APQ$ is $\pi$.
Then one of those two angles is $\le \pi/3$.  Without loss of generality
we can assume it is angle $APQ$.  Then angle $APQ$ is either 
equal to $\pi/3$ or to $\alpha$.  
If $AQ$ does not support a tile with its $c$ edge on $AQ$, then each 
tile supported by $AQ$ has its $\gamma$ angle on $AQ$, and by 
the pigeonhole principle, the tile at $Q$ has its $\gamma$ angle
at $Q$, contradiction.  Therefore $\overline{AQ} \ge c$.
In triangle $APQ$, the angle 
opposite $PQ$ (namely $\pi/3$) is greater than or equal to 
the angle opposite $AQ$
(namely angle $APQ$).  Therefore the length $x$ of $PQ$ is greater
than or equal to the length of $AQ$.  Hence $x \ge c$.  But $x$ is the length 
of a tile edge, and $c$ is the longest tile edge, so $x = c$ and
we have equality throughout, i.e., $AQ = c$ and triangle $APQ$ is 
equilateral. Now we have an equilateral triangle $APQ$ tiled by 
some number $n$ of tiles, with the side of $APQ$ equal to $c$. 
The area equation tells us $c^2 = nab$.  Let $g = gcd(a,b)$. 
Then $g^2$ divides $c$; but $(a,b,c)$ have no common factor, so $g = 1$.
By the law of cosines, 
\begin{eqnarray*}
c^2 &=& a^2 + b^2 + ab \\
(n-1)ab &=& a^2 + b^2 \mbox{\qquad since $c^2 = nab$}
\end{eqnarray*}
Taking this equation mod $a$ we find $b^2 \equiv 0$ mod $a$.
Since $b$ and $a$ are relatively prime, that is a contradiction.
That completes the proof in case $\gamma = 2\pi/3$.

Now we assume $\gamma = \pi/3$.  Then $\alpha < \gamma < \beta$,
so $a < c < b$.   At each vertex of $ABC$,
there is a single tile with its $\gamma$ angle at the vertex. 
Its $c$ side therefore does touch two sides of $ABC$.  We have
to prove that no other tile touches both $AB$ and $AC$. Suppose
to the contrary that $P$ lies on $AB$ and $Q$ on $AC$ and $PQ$
is an edge of a tile in the tiling.  
Let $x$ be the length of $PQ$, so $x$ is one of $(a,b,c)$.
 Assume, for proof by contradiction,
that triangle $APQ$ is 
equilateral.  What is the length $x$ of the sides of $APQ$? 
Since the tiles at $A$ have their $\gamma$ angles at $A$,
they have their $a$ or $b$ edges on $AP$ and $AQ$. If either
has its $b$ edge there, then $x \ge b$, so $x=b$. If neither 
has its $b$ edge there, then they both have their $b$ edges
along a line from $A$ to an interior point on $PQ$.  Such a 
line is shorter than the side of an equilateral triangle, so 
again $x \ge b$.  In this case $x > b$, which is impossible
since $x$ is one of $(a,b,c)$ and $b$ is the largest of these.
Therefore $x=b$.  Now triangle $APQ$ is 
$n$-tiled for some $n$, and the area equation tells us 
\begin{eqnarray*}
b^2 &=& nab \\
b &=& na
\end{eqnarray*}
By the law of cosines
\begin{eqnarray*}
c^2 &=& a^2 + b^2 -ab \\
   &=& a^2 + (na)^2 -a(na) \\
   &=& a^2
\end{eqnarray*}
Hence $a=c$, contradiction.  Therefore triangle $APQ$ is 
not equilateral.

Now consider angles $APQ$ and $AQP$; one of these angle is less than
or equal to $\pi/3$, since the sum of the angles of triangle $APQ$
is $\pi$ and angle $A$ is $\pi/3$.  Since triangle $APQ$ is not 
equilateral, one of those two angles is strictly less than $\pi/3$.
Without loss of generality we may assume it is angle $APQ$.
Since $PQ$ is part of the tiling, angle $APQ$ must be $\alpha$.
Then angle $AQP$ is $\beta$, and triangle $APQ$ is similar to the tile.
Then for some $\lambda > 0$ we have $x = \lambda c$ and $AQ = \lambda b$
and $AP = \lambda a$.  Since $x$ is one of $(a,b,c)$ we consider 
the possibilities one by one. If $x=c$ then $\lambda = 1$ and 
$AQ = b$ and $AP = a$.  Then the two tiles at $A$ form a parallelogram
with sides $a$ and $b$, and $PQ$ is the diagonal, hence not a part 
of the tiling.  That rules out the case $x=c$.  If $x=b$ then 
$\lambda = b/c$ and $AQ = b^2/c$ and $AP = ab/c$.  Since $a < b < c$
that would make $AP < a$, which is impossible since $AP$ must support
at least one tile.  That rules out $x=b$.  Therefore $x=a$. 
Then $\lambda = a/c$ and $AQ = ab/c < a$, which is also impossible.
That completes the proof of the lemma.

\begin{lemma} \label{lemma:pqrbound} Let the equilateral triangle
$ABC$ be $N$-tiled by a tile with sides $(a,b,c)$,
 and let $X$ be the length of $AB$.  Suppose
$X = pa + qb + rc$.  Then $p+q+r \le N/6 + 1$.
\end{lemma}

\noindent{\em Proof}.  There are three tiles at each 
boundary vertex of $ABC$.  Suppose there are $k$ tiles 
supported by $AB$.  Then there are $2k-1$ tiles with an 
edge or vertex on $AB$. 

First we assume $\gamma= 2\pi/3$. Then there are two 
tiles at each vertex, and by Lemma~\ref{lemma:pqrbound_helper},
no tile has a vertex on two 
different sides of $ABC$.  Then there will be no double-counting
of tiles when we triple that number:  there are at least $3(2k-1)$ 
different tiles with an edge or vertex on the boundary.
Hence $6k-3 \le N$, so $k = p+q+r le (N+3)/6 \le N/6 + 1$. 

Now assume $\gamma = \pi/3$.  Then when we triple the number
$3(2k-1)$, we have double-counted the single tiles at the 
vertices of $ABC$, and also there are three tiles with a vertex
on two sides of $ABC$ that will be double-counted.  But by 
Lemma~\ref{lemma:pqrbound_helper}, only those three tiles touch 
two different sides of $ABC$.  Hence 
\begin{eqnarray*}
3(2k-1)-6 &\le& N \\
6k &\le& N+6 \\
p+q+r &\le& N/6 + 1
\end{eqnarray*}

That
completes the proof of the lemma.

\begin{lemma} \label{lemma:equilateral7} 
Let the equilateral triangle $ABC$ be $N$-tiled by a tile with
$\gamma = \pi/3$ or $2\pi/3$, and $\alpha$
not a rational multiple of 
$\pi$.  Then $N \ge 40$.  If $N \le 75$ then the possible
values of $N$ and the associated tiles are given in Table~\ref{table:equilateral}.
\end{lemma}

\begin{table}[ht]
\caption{Possible $N$ and $(a,b,c)$ for equilateral tilings}
\label{table:equilateral}
\begin{center}
\setlength{\extrarowheight}{.2em}
\begin{tabular}{rrr}
$N$ & $\gamma$ & the tile \\
\hline
40 & $\frac {\pi} 3$ & $(5,8,7)$ \\
54 & $\frac \pi 3$ & $(3,8,7)$ \\
56 & $\frac{2\pi} 3$ & $(7,8,13)$ \\
60 & $\frac{2\pi} 3$ & $(3,5,7)$ \\
65 & $\frac \pi 3 $ &  $(9,65,61)$ \\
66 & $\frac {2\pi} 3$ & $(11,24,31)$ \\
70 & $\frac {\pi} 3$ & $(7,40,37)$ \\
80 & $\frac{2\pi} 3$ & $(5,16,19)$ \\
84 & $\frac{\pi} 3 $ & $(16,20,19)$\\
85 & $\frac \pi 3 $ & $(17,80,73)$ 
\end{tabular}
\end{center}
\end{table}
\medskip

\noindent{\em Proof}.  After a suitable relabeling of the 
vertices of $ABC$ so that Lemma~\ref{lemma:pqr} will apply,
let $X$ be the length of $AB$ and let the tiling determine
the integers $(p,q,r)$ such that
 $X = pa+qb+rc$.  According to the area equation (\ref{eq:2511})
 we have 
\begin{eqnarray*}
X^2 &=& Nab \\
(pa+qb+rc)^2 &=& Nab 
\end{eqnarray*}
By the law of cosines (\ref{eq:2605}), we have 
\begin{eqnarray*}
c &=& \sqrt{a^2 + b^2 \pm ab}
\end{eqnarray*}
Putting that into the area equation, we have
\begin{eqnarray*}
(pa + qb + r \sqrt{a^2 + b^2 \pm ab})^2 &=& Nab
\end{eqnarray*}

Define $s:= a/b$ and divide the equation by $b^2$:
\begin{eqnarray*}
(ps + q + r \sqrt{s^2 + 1 \pm s})^2 &=& Ns
\end{eqnarray*}
Expanding the left side we have 
\begin{eqnarray*}
(ps+q)^2 + r^2(s^2 + 1 \pm s) + 2(ps+q)r\sqrt{s^2+1\pm s} &=& Ns
\end{eqnarray*}
\begin{eqnarray*}
2(ps+q)r\sqrt{s^2+1\pm s} &=& Ns - (ps+q)^2 - r^2(s^2 + 1 \pm s)
\end{eqnarray*}
Squaring both sides we have
\begin{eqnarray}
4(ps+q)^2 r^2(s^2 + 1\pm s) &=& (Ns - (ps+q)^2 - r^2(s^2+1\pm s))^2
\label{eq:8696}
\end{eqnarray}
That is a fourth-degree polynomial equation in $s$.  Since the tile
is known to be rational,  and since we can assume $a < b$, we are
looking for rational solutions $s$ with $0 < s < 1$.  If there is 
a tiling,  there will be such a solution.  Since SageMath can 
solve quartic equations, and test whether the solutions are rational,
we are almost finished:  we just need to bound the possible values
of $(p,q,r)$. 

 We make use of Lemma~\ref{lemma:pqrbound} to
use the (still crude) bound $N/6+1$.  Now the algorithm is simple:
given $N$, check all the possibilities $(p,q,r)$ that satisfy
the conditions of Lemma~\ref{lemma:pqr} and satisfy $p+q+r \le N/6+1$.
Solve (\ref{eq:8696}). Reject solutions that are not real and 
solutions that are not rational and solutions that are not 
between 0 and 1.   If there are no solutions remaining, then there
is no $N$-tiling.   If there is a solution, the method is inconclusive.
The code is shown in Fig.~\ref{figure:equilateralpqr}.  We put 
the function defined there in an outer loop and ran it over $N$ from 
$3$ to $47$.  It finds solutions only for the values of $N$
mentioned in the statement of the lemma.  That completes the proof,
at least, if you believe the code is correctly written and correctly
executed. 

We took the following steps to check this code for correctness.
We wrote this code independently in SageMath and in C, and,
commenting out the line to reject rational solutions,  printed out 
the solutions that are found. The C code only works with floating-point
numbers, so to compare the results, we had SageMath numerically 
evaluate its exact solutions.  The C and SageMath code printed out 
the same decimal values,  which we found reassuring.  Since the C 
code works with floating-point numbers, it is not so easy to have 
it reject irrational solutions, so we could not replace the SageMath
code with more efficient C code.

\begin{figure}
\caption{SageMath code for ruling out equilateral $N$-tilings}
\label{figure:equilateralpqr}
\begin{verbatim}
def feb24(N,an):  # an = -1 for pi/3, 1 for 2pi/3
# search for equilateral gamma = pi/3 or 2pi/3 solutions 
# of area equation and boundary conditions.
   var('s')
   if an == -1:    # gamma at A,B,C
      Abound = 1   # at least one a edge  
      Bbound = 1    # at least one b edge 
   else:           # alpha + beta at A,B,C
      Abound = 0
      Bbound = 0
   upperbound = N/6 + 2  
   for p in range(Abound,upperbound):   
      for q in range(Bbound,upperbound):     
         for r in range(2,upperbound):  # at least 2 c edges in either case
            if p+q+r >= floor(upperbound):
               continue
            eq = 4*(p*s + q)^2*r^2*(s^2+1+an*s)
                - (N*s-(p*s+q)^2 - r^2*(s^2+1+ an*s))^2
            answers = solve(eq,s)
            for t in answers:
               S = t.rhs();
               if not S in RR:
                  continue
               if not S in QQ:
                  continue      # tile is known to be rational
               if S <= 0 or 1 <= S:  # we can assume a < b so s<1
                  continue
               A = S.numerator()
               B = S.denominator()
               C = B * sqrt(S^2+1+ an*S)
               if not C in QQ:
                  continue;
               if not C in ZZ:
                  print("oops")   # it better be in ZZ!
                  print(A,B,C);
               g = gcd(A,gcd(B,C));
               a = A/g
               b = B/g
               C = C/g
               print("found (%d, %d, %d)" %(a,b,c))
   return true
\end{verbatim}
\end{figure}

\section{No $7$-tilings}

We break the proof that there are no 7-tilings or 11-tilings into
two cases, according as the angles are commensurable or not.
All the required cases have already been dealt with: it only 
remains to put the pieces together.

\begin{theorem} \label{theorem:no7-commensurable}
Suppose $(\alpha,\beta,\gamma)$ are all rational multiples of 
$2\pi$.  Then there is no 7-tiling of any triangle $ABC$ by a
tile with angles $(\alpha,\beta,\gamma)$.  Moreover, there is 
no $N$-tiling by such a tile for $N = 11,14,19,31$
or any number which is neither a square, sum of squares, or 
$2$, $3$, or $6$ times a square.
\end{theorem}

\noindent{\em Remark.} Any odd $N$ which is not divisible by 3
but whose squarefree part is divisible by some prime congruent to 3 mod 4 meets the 
conditions of the theorem.
\medskip

\noindent{\em Proof}. Assume, for proof by contradiction,
that there is such a tiling.  By Theorem~\ref{theorem:laczkovich-extended},
 the pair $ABC$ and $(\alpha,\beta,\gamma)$  (after a suitable renaming
of the angles) occurs in Table~\ref{table:laczkovich-extended}.  
But 7 does not match any of the forms of $N$ listed in that table,
which are the forms listed in the final sentence of the theorem.
That completes the proof.

Finally we have arrived at the main theorem.

\begin{theorem} There are no 7-tilings or 11-tilings.  
\end{theorem}

\noindent{\em Proof.}  Suppose, for proof by contradiction,
that triangle $ABC$ is N-tiled by a tile with angles 
$(\alpha, \beta,\gamma)$, with $N = 7$ or $11$.  We also
note the cases where the proof works for $N=14$ and $19$. 
Since $N$ is not a square or a sum of two squares, 
then by \cite{snover1991}, $ABC$ is not similar to the tile. 
Then according to Theorem~\ref{theorem:no7-commensurable},
not all the angles $(\alpha,\beta,\gamma)$ are rational multiples of $\pi$.
Then according to \cite{laczkovich1995}, the tiling must correspond
to one of the rows in Table~\ref{table:laczkovich2} in this paper.

For our purpose these rows will be combined into five cases:
Either $3\alpha + 2\beta = \pi$, or $\gamma = 2\pi/3$, 
or $ABC$ is isosceles with base angles $\alpha$ and $\gamma = \pi/2$,
or $ABC$ is isosceles with base angles $\alpha$ and $\gamma = 2\alpha$,
or $ABC$ is equilateral and $\alpha = \pi/3$.

In case $ABC$ is equilateral and $\alpha = \pi/3$,  Lemma~\ref{lemma:equilateral7} tells us there is no 
7-tiling or 11-tiling (but the proof does not work for $N=14$ or $N=19$).

In case $ABC$ is isosceles with base angles $\alpha$ and $\gamma = \pi/2$,
then by \cite{beeson-isosceles}, $N$ is a square, or a sum of squares,
or six times a square.  Hence 
$N$ cannot be $7, 11, 14$, or $19$. 

In case $ABC$ is isosceles with base angles $\alpha$ and 
$\gamma = 2\alpha$,  by \cite{beeson-isosceles},
$N$ is not prime or twice a prime. In particular $N$ is not
$7, 11, 14$, or $19$. 

In case  $3\alpha + 2\beta = \pi$,  Theorem~\ref{theorem:no7-11a}
tells us there is no 7-tiling or 11-tiling. 
In case $\gamma = 2\pi/3$, by Theorem~\ref{theorem:no7-11b} 
there is no 7-tiling or 11-tiling.  That completes
the proof of the theorem.
\medskip

\section{Concluding Remarks} \label{section:last}

This paper has successfully avoided the need to appeal to 
the hundred pages of theoretical work on the case $3\alpha + 2\beta = \pi$,
as well as more than thirty pages on the equilateral
case, although we did appeal to our work on the isosceles case,
since the computational approach was not much shorter.
 Instead we have used
algebraic and computation shortcuts that work only for 
small values of $N$. 

In this section we nevertheless mention
some results of those lengthier investigations.  First,
in each of the three cases ($3\alpha + 2\beta = \pi$, isosceles,
equilateral), we used techniques pioneered by Laczkovich to prove
that the tile has to be rational.  Then we used the area equation
and (for the $3\alpha + 2\beta= \pi$ case and one equilateral case)
the ``coloring equation'' to derive necessary conditions. 
We used these equations to prove that $N$ cannot be prime in 
some cases.  The exceptions are as follow: 
There are the biquadratic tilings of a right isosceles triangle, in which 
case if $N$ is prime it must be congruent to 1 mod 4.  Also, for an
isosceles triangle tiled by a right triangle,  we could only prove
$N$ cannot be a prime congruent to 3 mod 4, although we conjecture
$N$ has to be even (and hence not prime at all).  And
in case $ABC$ is isosceles with base angles $\alpha$, and $\gamma = 2\pi/3$,
and $\alpha/\pi$ is not rational, we were unable to settle the question,
although if $N$ is prime in that case, $N = 2b+a$.  

Generally we want to we expand the lines with in
Table~\ref{table:laczkovich2} 
 by adding a third column with restrictions on the possible 
form of $N$.   In some cases this is a necessary and sufficient condition;
in others it is only a necessary condition.  Where the necessary and 
sufficient conditions do not match,  there are open questions.  Whether
there are yet-undiscovered tilings, or our necessary 
conditions are too weak, we do not know.    
Table~\ref{table:finalsummary} gives a summary of what we know about 
$N \le 100$.  

\begin{figure}
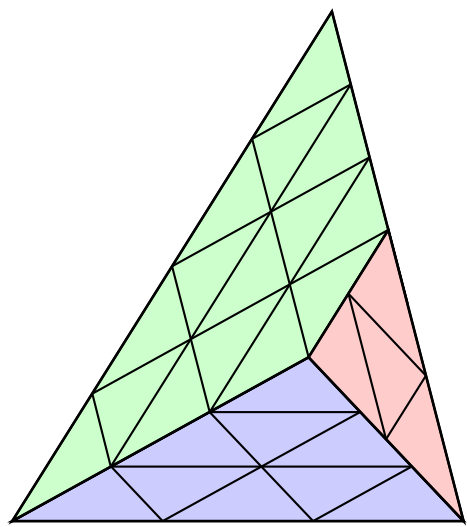

\caption{A tiling with $N=28$ and $3\alpha + \beta = \pi$, and tile $(2,3,4)$}
\label{figure:triquadratics}
\begin{center}
\TriquadraticTwentyEight
\end{center}
\end{figure}

\begin{figure}
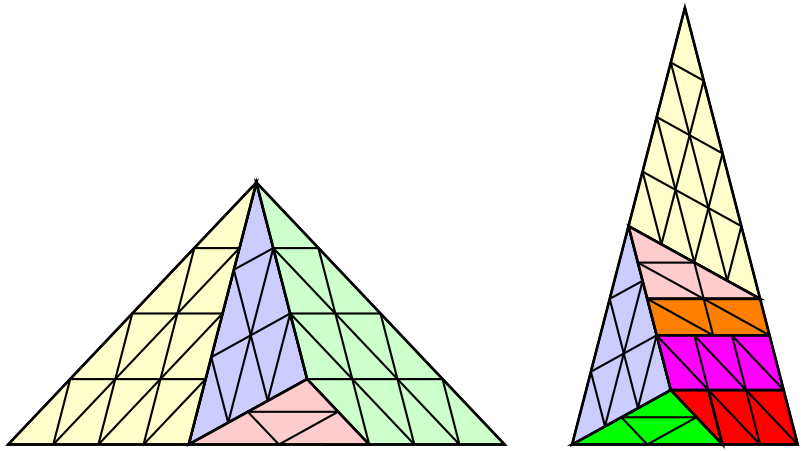

\caption{Tilings with $N=44$ and $48$, with $3\alpha + \beta = \pi$ and tile $(2,3,4)$}
\label{figure:endfigure}
\begin{center}
\IsoscelesBetaFortyFourTiling
\IsoscelesFortyEightTiling
\end{center}
\end{figure}

 Figures~\ref{figure:triquadratics} and \ref{figure:endfigure} show some examples of tilings
unknown before 2011 and 2108, respectively.  We also have many examples
of larger tilings. 

Very little is known about the possible values of $N$ for tilings
of isosceles and equilateral triangles;  proving that $N$ cannot 
be prime is an advance, since for example in 2012 it was not known whether 
there are arbitrarily large $N$ such that no equilateral triangle 
can be $N$-tiled by a tile whose angles are not all rational multiples of $\pi$.
Similar for the case of tiling an isosceles $ABC$ by a tile with 
$\gamma = 2\alpha$.  Now we know $N$ can't be prime, but we still 
don't know if $N$ can be even or not, and we don't know if $N$
can be less than five million, although Laczkovich proved that 
any such tile must tile {\em some} isosceles $ABC$, so there do 
exist a lot of such tilings--we just don't know how big (or small)
$N$ can be.

\begin{table}[ht]
\caption{Knowledge about tilings with $N \le 100$ as of March, 2019}
\label{table:finalsummary}
\begin{center}
\setlength{\extrarowheight}{.2em}
\begin{tabular}{lrrll}
 $ABC$ shape & {\bf the tile} & {\bf known $N$}  & {\bf values $<= 100$} & {\bf least unknown $N$}\\
\hline
{\bf equilateral }& $(\frac \pi 6, \frac \pi 6, \frac {2\pi} 3)$ & $6n^2$ & $24,72,96$ &  \\
&$(\frac \pi 3, \frac \pi 2, \frac {2\pi} 3)$ & $3n^2$ & $12,27,48,75$ &  \\
&$(\alpha, \beta, \frac {\pi} 3)$ &  $5861172$ & ? & 40? \\
&$(\alpha, \beta, \frac {2\pi} 3)$ &  $10395$ & ? & 40? \\
\hline
{\bf right}$(\frac \pi 6, \frac \pi 3, \frac \pi 2)$
         &$(\frac \pi 6, \frac \pi 3, \frac \pi 2)$ & $n^2, 3n^2$ & 
                   $4,9\ldots 81,100$ \\
{\bf right}$(\alpha,\beta,\frac \pi 2)$ & $(\alpha,\beta,\frac \pi 2)$& $N = e^2 + f^2$ &  & \\
\hline 
{\bf isosceles}-$\alpha$  &$(\frac \pi 6, \frac \pi 3, \frac \pi 2)$ & $6n^2$ & $24,72,96$ &  \\
&$(\alpha,\beta,\frac \pi 2)$ & $2n^2$ & $2,8,18\ldots$\\
&$(\frac \pi 4, \frac \pi 4, \frac \pi 2)$ & $6n^2$ & $24,72,96$ &  \\ 
&$(\alpha,\beta, 2\alpha)$ &   & not $p$ or $2p$ & 20? \\ 
&$(\alpha,\beta,\frac {2\pi}3) $ &  1878500 & &33? \\
&$3\alpha + 2\beta = \pi$   && 84 &70? \\
{\bf isosceles}-$\beta$&$3\alpha + 2\beta = \pi$ && 44  & 59? 66? 71? 74?\\
&&&&83? 92? 99? \\ 
{\bf isosceles}-$\alpha +\beta$&$3\alpha + 2\beta = \pi$ & & 48 & 45? 72? 75?  99? \\
\hline
$(\alpha, 2\alpha,2\beta)$ & $3\alpha + 2\beta = \pi$ & & 77&  \\
$(2\alpha,\beta,\alpha+\beta)$&  $3\alpha + 2\beta = \pi$& & 28 &  \\
$(\alpha,\alpha+\beta,\alpha+2\beta)$ & $(\alpha,\beta,\frac {2\pi} 3)$ & 
& & 13? \\
$(\alpha,2\alpha,3\beta)$ &  $(\alpha,\beta,\frac {2\pi} 3)$ &
    & & 13? \\
$(2\alpha,2\beta,\alpha +\beta)$ &  $(\alpha,\beta,\frac {2\pi} 3)$ &
    & & 13? \\
\hline
any $ABC$ & similar to $ABC$ & 	$n^2$ & 4,9,16,\ldots 100&
\end{tabular}
\end{center}
\end{table}

\FloatBarrier

\bibliographystyle{plain-annote}

\end{document}